\documentclass[preprint,review,3p,times]{elsarticle}

\usepackage{amsmath,amssymb}
\usepackage{pgfplots}
\usepackage{marginnote}
\usepackage{textpos}
\usetikzlibrary{positioning}
\usepackage{commath}
\usepackage{mathtools}
\usepackage{subfig}
\usepackage{tabularx}
\usepackage{graphicx}
\usepackage{array} 
\usepackage{multirow}
\usepackage{multicol}
\usepackage{ulem}
\usepackage{hyperref}
\newcolumntype{M}[1]{>{\centering\arraybackslash}m{#1}}
 
\usepackage{xcolor}
\definecolor{Tblue}{RGB}{55,126,184}
\definecolor{Tred}{RGB}{228,26,28}
\definecolor{Tgreen}{RGB}{77,175,74}
\definecolor{Torange}{RGB}{255,127,0}
\definecolor{Tpurple}{RGB}{152,78,163}
\definecolor{Tbrown}{RGB}{166,86,40}
 
\usepackage{tikz}
\usetikzlibrary{patterns,calc,shapes,arrows}

\newcommand{\R}{\mathbb{R}}

\newcommand{\logLogSlopeTriangle}[5]
{
	\pgfplotsextra
	{
		\pgfkeysgetvalue{/pgfplots/xmin}{\xmin}
		\pgfkeysgetvalue{/pgfplots/xmax}{\xmax}
		\pgfkeysgetvalue{/pgfplots/ymin}{\ymin}
		\pgfkeysgetvalue{/pgfplots/ymax}{\ymax}
		
		\pgfmathsetmacro{\xArel}{#1}
		\pgfmathsetmacro{\yArel}{#3}
		\pgfmathsetmacro{\xBrel}{#1-#2}
		\pgfmathsetmacro{\yBrel}{\yArel}
		\pgfmathsetmacro{\xCrel}{\xBrel}
		
		\pgfmathsetmacro{\lnxB}{\xmin*(1-(#1-#2))+\xmax*(#1-#2)} 
		\pgfmathsetmacro{\lnxA}{\xmin*(1-#1)+\xmax*#1} 
		\pgfmathsetmacro{\lnyA}{\ymin*(1-#3)+\ymax*#3} 
		\pgfmathsetmacro{\lnyC}{\lnyA+#4*(\lnxA-\lnxB)}
		\pgfmathsetmacro{\yCrel}{\lnyC-\ymin)/(\ymax-\ymin)} 
		
		\coordinate (A) at (rel axis cs:\xArel,\yArel);
		\coordinate (B) at (rel axis cs:\xBrel,\yBrel);
		\coordinate (C) at (rel axis cs:\xCrel,\yCrel);
		
		\draw[#5]   (B)-- node[pos=0.5,anchor=north] {1}
		(A)-- 
		(C)-- node[pos=0.5,anchor=east] {#4}
		cycle;
	}
}

\journal{Computer-Aided Design}

\begin{document}
\begin{textblock*}{\textwidth}(0cm,-3cm)
	\noindent\small{This article a preprint of an article published in the journal \textit{Computer-Aided Design}. The final authenticated version is available online at: \url{https://doi.org/10.1016/j.cad.2022.103409}.}
\end{textblock*}\vspace{-1cm}
 


\begin{frontmatter}

	\author[1]{Thibault Jacquemin}
	\author[1]{Pratik Suchde}
	\author[1,2]{St\'ephane P.A. Bordas}

	\affiliation[1]{organization={University of Luxembourg, Institute of Computational Engineering},
		addressline={Maison du Nombre, 6 Avenue de la Fonte},
		city={Esch-sur-Alzette},
		postcode={L-4364},
		country={Luxembourg}}

	\affiliation[2]{organization={Department of Medical Research, China Medical University Hospital},
		addressline={China Medical University},
		city={Taichung},
		country={Taiwan}}

	\title{Smart cloud collocation: geometry-aware adaptivity directly from CAD}
	\author{}

	\begin{abstract}
		Computer Aided Design (CAD) is widely used in the creation and optimization of various industrial systems and processes. Transforming a CAD geometry into a computational discretization that be used to solve PDEs requires care and a deep knowledge of the selected computational method. In this article, we present a novel integrated collocation scheme based on smart clouds. It allows us to transform a CAD geometry into a complete point collocation model, aware of the base geometry, with minimum effort. For this process, only the geometry of the domain, in the form of a STEP file, and the boundary conditions are needed. We also introduce an adaptive refinement process for the resultant smart cloud using an \textit{a posteriori} error indication. The scheme can be applied to any 2D or 3D geometry, to any PDE and can be applied to most point collocation approaches. We illustrate this with the meshfree Generalized Finite Difference (GFD) method applied to steady linear elasticity problems. We further show that each step of this process, from the initial discretization to the refinement strategy, is connected and is affected by the approach selected in the previous step, thus requiring an integrated scheme where the  whole solution process should be considered at once.
	\end{abstract}


\end{frontmatter}

\section{Introduction}

Computer Aided Design (CAD) software packages are used in many domains of engineering to design components or various nature. A new design is often proposed based on previous experience and knowledge. It can then be optimized using calculation and/or simulation  tools to increase its performance, lower the manufacturing costs or for many other reasons. Reaching a satisfactory design often requires iterations. To minimize the development costs, both industry and academia have been trying to speed-up this iteration process (design \verb|->| simulation \verb|->| design modification \verb|->| simulation \ldots) as much as possible.

Performing computer simulations directly from geometry can be a tedious task. We introduce in this article an integrated smart cloud collocation scheme. The scheme, based on point collocation, aims at simulating the behavior of components (the mechanical resistance of a solid part, for instance) directly from the designed CAD geometry. The proposed scheme minimizes the size of the solved problem by using \textit{a posteriori} error indication and adaptive refinement. The smart cloud is geometry-aware. Each point of the smart cloud has additional information related to the CAD geometry and to the boundary conditions applied to the domain.

Many collocation schemes use meshes to generate point clouds. The approach presented in this work does not require this step. Mesh generation poses several additional constraints on the aspect and shape of the elements which do not apply to point collocation methods. Therefore, skipping the mesh generation step could lead to a more robust domain discretization framework for point collocation methods.

Point collocation methods have been used for a long time and can be applied to smooth and non smooth solutions and domains \cite{Jacquemin2019}. The selection of collocation stencil is a key aspect of this family of methods. To ease their application to all types of problems and domains, a unified approach was introduced in reference \cite{Jacquemin2021}. The finite difference method was the first such collocation method. It was introduced by C. Runge in 1908 \cite{Runge1908}. A Cartesian grid was used to approximate the field derivatives, limiting the problems solved by this method to simple geometries. The method was then generalized in 1972 by Jensen \cite{Jensen1972}. He introduced the basis of the Generalized Finite Difference (GFD) method. The method was then successively developed by Liszka \cite{Liszka1980}, Orkisz \cite{Orkisz1998}, Benito \cite{Benito2001}, Milewski \cite{Milewski2012} and many other contributors.

Recently, numerical methods based on point collocation have regained interest with their use in meshfree frameworks. The GFD method shows good performance compared to other point collocation methods \cite{Jacquemin2019} and was used in this work. The scheme presented in this article can however be readily applied, without any modification, to other collocation schemes such as the Moving Least Squares (MLS) method \cite{Shepard1968,Lancaster1981,Onate1996} or the Radial Basis Finite Difference (RBF-FD) method \cite{Kansa1990a, Kansa1990b,Driscoll2002,Davydov2011}.

The idea of using the CAD geometry to solve problems defined by Partial Differential Equations (PDE) is at the heart of the Isogeometric Analysis (IGA) methods. These methods became rapidly popular after their first introduction in 2005 \cite{Hughes2005} and have been proven robust for problems of various nature. IGA methods use the functions of the CAD geometry as shape functions to approximate the field derivatives and solve PDEs over domains. The Isogeometric analysis boundary element method (IGABEM) is a popular IGA method introduced by Simpson et al. in 2012 \cite{Simpson2012}. It combines the benefits of the IGA method (direct use of the functions of the CAD geometry) and the benefits of the Boundary Element Method (use of a discretization of the surface boundaries only). The clearest advantage of IGABEM is the possibility to solve PDEs without any mesh generation, which is particularly convenient for shape optimization \cite{Lian2013,Simpson2014,Peng2017,Lian2017}. A collocation form of IGA has also been introduced \cite{Auricchio2010}.

Using the functions of the geometry to solve a given problem can be seen both as a strength and a weak point of IGA computational methods. This assumes that the CAD geometry is suitable for analysis, which is not always the case. Point collocation methods are deemed quite ``flexible" compared to element-based methods and IGA methods since there are no strong connections between the nodes. Furthermore, the discretization of the domain can be easily modified as needed during the simulation process to provide the best possible solution for a given computational cost.

Discretization adaptivity has always been of interest for point collocation methods. The performance of this approach in the framework of point collocation was shown in many articles for the GFD \cite{Milewski2012,Benito2008,Gavete2015,Suchde2019a} or the RBF-FD \cite{Davydov2011a,Oanh2017,slak2020adaptive} method.

In contrast to existing work, the smart cloud collocation scheme presented in this work uses the exact definition of a given geometry based on a CAD file. It can therefore be applied to most domains with the assurance that the exactness of the geometry is not lost as part of the refinement process. Such an approach was not present in previous publications, to the authors' knowledge. The proposed adaptive method uses new nodes, placed at key locations in the domain, to improve the solution. This approach is often referred to as $h$-adaptivity. The nodes of the initial point cloud are kept. This hierarchical approach implies that the point cloud does not need to be generated again before the next adaptive iteration step which saves computational effort.

The article is composed of two main sections. We present in the first one the method used to transform a CAD geometry into a smart cloud. We show in a second section how error indicators can be used to identify the zones of the domain where the error is the greatest. Once these zones are identified, we show how adaptive refinement can be used to converge efficiently to an accurate solution. Our work focuses on linear elasticity problems (2D and 3D) using the Generalized Finite Difference method. The smart cloud concept can also be applied to other types of elliptic problems and to most point collocation methods.

A key novelty in the present work is that the discretization of the domain and its adaptivity are based on the exact CAD geometry to minimize the user input in the simulation workflow and ensure that the domain is represented exactly throughout the refinement process.

\section{From CAD to smart cloud}\label{DiscretizationSec}
\subsection{General}

Most numerical methods used to solve partial differential equations use some sort of discretization of the otherwise infinite dimensional mathematical problem. Either the differential operator is discretized (e.g. GFD, RBF-FD) or the unknown field is discretized (e.g. MLS, Finite Element). Collocation methods use nodes spread in the domain, on the boundary of the domain and sometimes even outside of the domain to approximate the differential operator. Nodes placed outside of the domain, called ghost nodes, can be used as additional degrees of freedom to enforce boundary conditions or balance the stencils of the collocation nodes located on, or close to the boundary of the domain. The concept of ghost cells or ghost nodes was described by Fedkiw et al. in 1999 in reference \cite{Fedkiw1999}. The approach was applied to fluid flow solvers but can be applied to other types of problems such as linear elasticity as presented in reference \cite{Runnels2021}. In this work, we considered only discretizations featuring nodes placed in the domain and on the boundary of the domain.

The design of a new mechanical system very often requires a model of the geometry. CAD software packages are used for this purpose. Sophisticated user interfaces allow the design of complex geometries and the assembly of very large structures. Whilst many file formats are proprietary, the STEP file format (Standard for the Exchange of Product model) is a format defined by the international norm ISO 10303-21. Such a definition makes the STEP file format popular and most of the packages support it. It is mostly oriented toward 3D geometries but can also be used for 2D geometries. STEP files store the exact geometry of the domain using simple geometric features such as planes or conical surfaces but also B-spline surfaces and trimmed surfaces.

We present in this section a method to discretize a given geometry, provided in a STEP file format, with the aim of solving it using a point collocation method. We show how the smart cloud is generated from the discretization and from the CAD geometry. Finally, we analyze the impact of the selected parameters of the method.

\subsection{Domain discretization} \label{SecMethod}

We used the library Open CASCADE Technology \cite{opencascade} to communicate with the STEP files and to get information about the exact geometry. Our algorithm is composed of the following steps:

\begin{enumerate}
	\setlength{\itemsep}{0pt}
	      \setlength{\parskip}{0pt}
	\item Loading the information from the STEP file using Open CASCADE Technology;
	\item Discretization of the boundaries of the domain;
	\item Regular discretization of a rectangle or box comprising the geometry;
	\item Identification of the nodes of the rectangle or box included in the domain.
\end{enumerate}

\begin{figure}[h!]
	\centering
	\includegraphics[width=15cm]{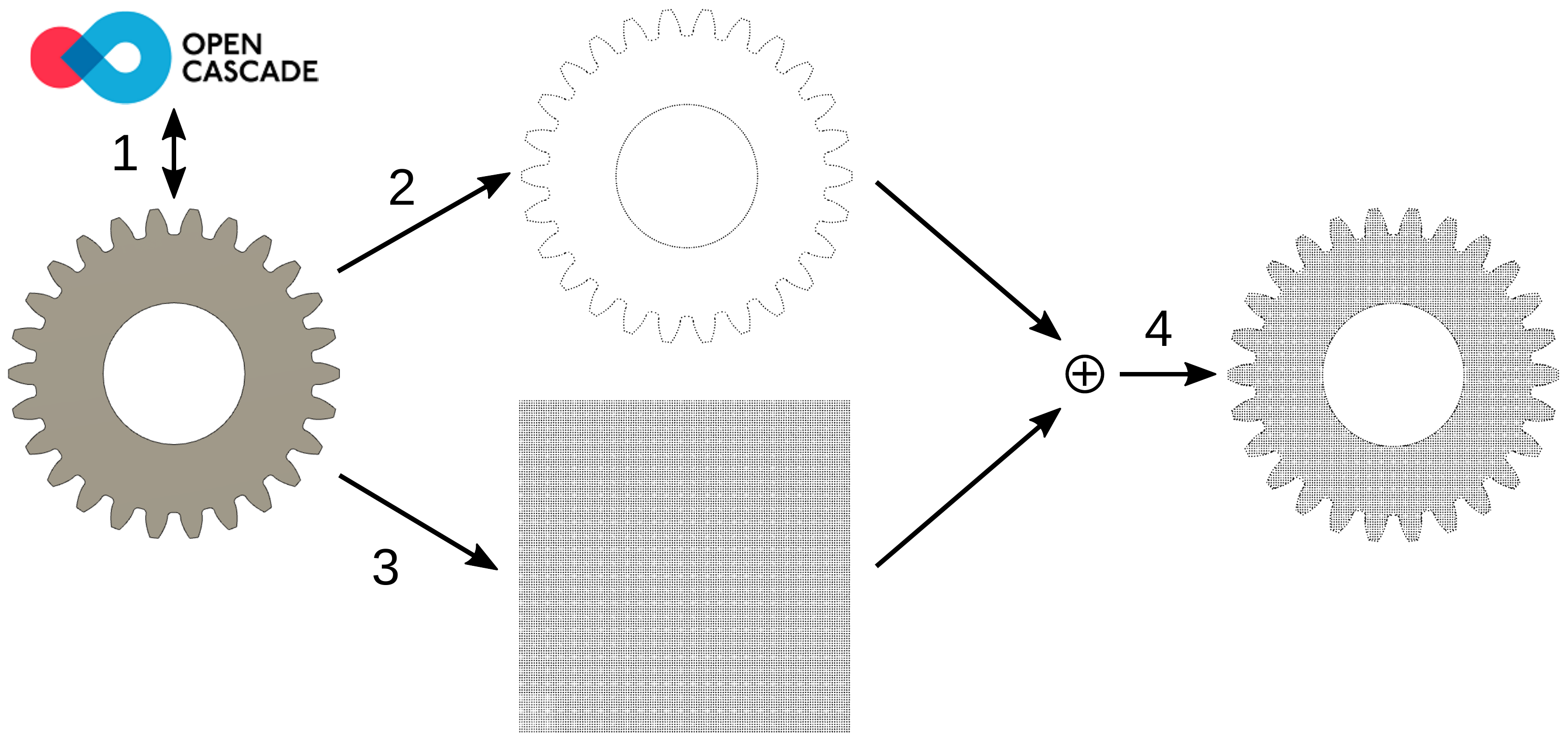}
	\caption{Steps of the discretization of a domain from a CAD file using the library Open CASCADE Technology \cite{opencascade}.}
	\label{DiscretizationSteps}
\end{figure}

We presented these steps in Figure \ref{DiscretizationSteps} for the case of a 2D gear. We give more details about each of these steps in the paragraphs below.

\paragraph{Step 1}

At the beginning of the discretization process key parameters, such as the bounding volume or the dimensions of the rectangle or box comprising the geometry, are computed from the CAD file.

A CAD geometry is composed of multiple topological entities. Those are:

\begin{itemize}
	\item solid;
	\item shell;
	\item face;
	\item edge loop;
	\item edge;
	\item vertex.
\end{itemize}

These entities are identified and used to discretize the boundaries of the domain and set the boundary conditions. Geometrical entities are associated with each topological entity. For instance, the geometrical entity associated with a face is a surface (e.g. plane, cylindrical surface, B-spline surface) and the geometrical entity associated with an edge is a curve (e.g. line, circle, B-spline curve).

The discretization process requires the selection of a characteristic length noted $h$. $h$ can be an input from the user or can be approximated based on a target number of nodes of the domain discretization. In this case, $h$ is computed in this first step.

\paragraph{Step 2}

The boundaries of the domain are discretized, based on the characteristic length $h$, ensuring that the distance between two adjacent boundary nodes is close to $h$. For 2D problems, all the edges of the domain are discretized using a fixed distance, close to $h$, between two consecutive nodes. The duplicated corner nodes are removed.

For 3D problems, we used a Delaunay triangulation of the boundary faces that compose the geometry. Generating such a mesh is robust since the boundary of the domain is composed of faces of simple geometry. The duplicated edge nodes are removed. We used the library Gmsh \cite{Geuzaine2009} to mesh the surfaces.

The exact normal vectors are computed, at each boundary collocation node, using the information about the exact geometry contained in the STEP file.

\paragraph{Step 3}

\begin{figure}[h!]
	\centering
	\begin{tabular}[b]{c c c}
		\subfloat[][]{\includegraphics[scale=0.5]{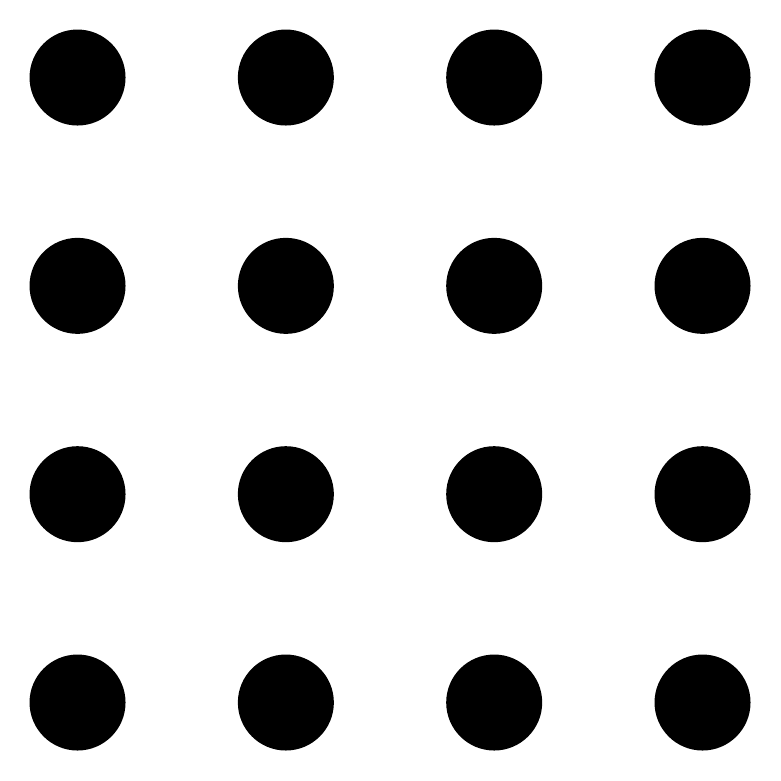}} & \quad &
		\subfloat[][]{\includegraphics[scale=0.5]{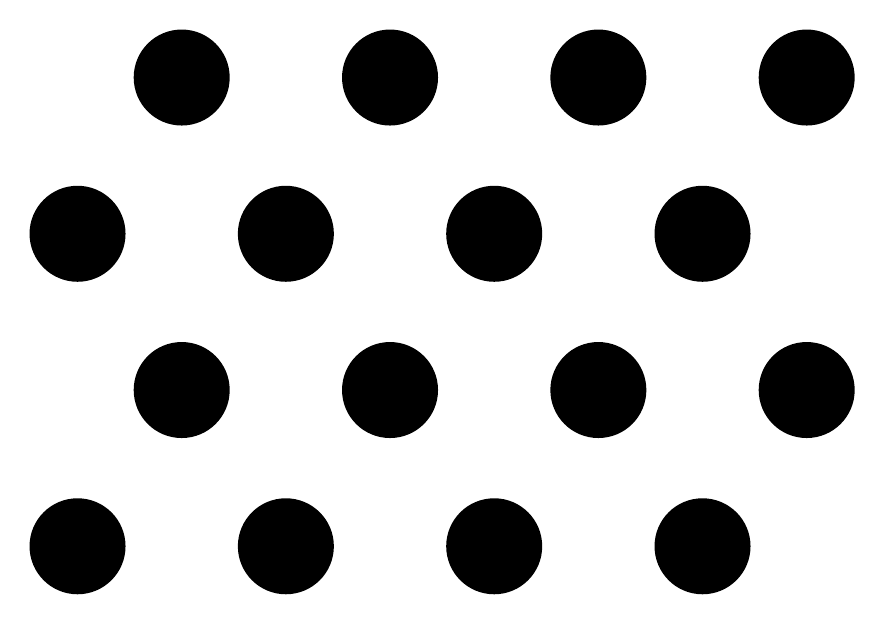}}               \\
		\subfloat[][]{\includegraphics[width=4cm]{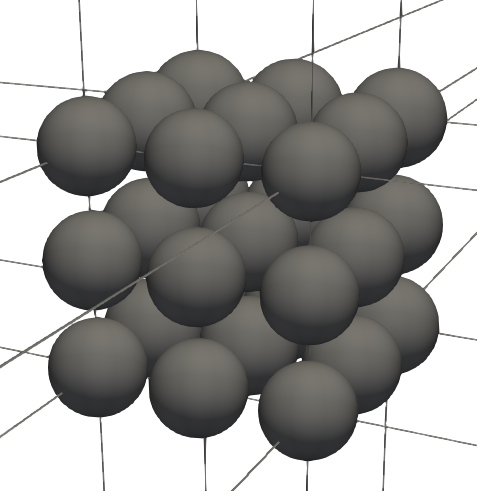}}   & \quad &
		\subfloat[][]{\begin{tabular}[b]{c}
				\includegraphics[width=3.5cm]{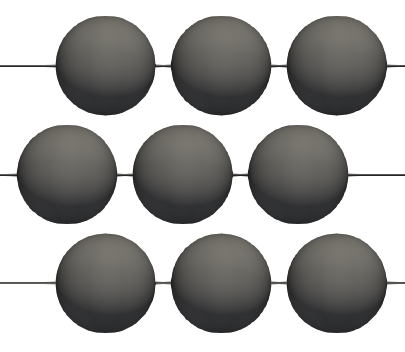}  \\
				Side view                                      \\
				\includegraphics[height=4cm]{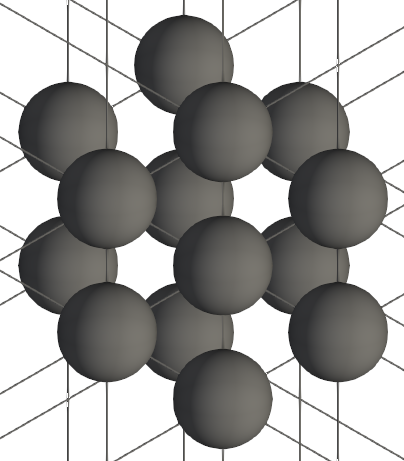} \\
				Top view                                       \\
			\end{tabular}}                                     \\
	\end{tabular}
	\caption{Node arrangement configurations in 2D (a), (b) and 3D (c), (d). The subfigures (a) and (b) show the node arrangements for square and triangular lattices, respectively. The subfigures (c) and (d) show the node arrangements for cubic and hexagonal close-packed lattices, respectively.}
	\label{2DLattices}
\end{figure}

The rectangle or box is discretized based on the characteristic length $h$. The nodes are placed regularly in the rectangle or box comprising the geometry. The nodes can be organized in different regular forms, also called lattices. In 2D and 3D, the nodes can be placed following the Cartesian grid. In 2D, the rectangle can also be discretized using equilateral triangles. This corresponds to a hexagonal close-packed lattice in 3D. We selected these two lattices for 2D and 3D problems as they lead to the most uniform discretizations. Figure \ref{2DLattices} shows a comparison of the mentioned 2D and 3D lattices. 
Many other lattices could be considered. The rectangle or box could also be filled using the advancing front method \cite{Lhner1998}.

\paragraph{Step 4}

The final step of the discretization process requires the identification of the nodes which are outside of the domain. Multiple algorithms can be used for this purpose. The CAD file can be used directly to assess the position of a node in a domain. We used this approach. The algorithms presented in reference \cite{Lhner1998} or in reference \cite{Jacquemin2021} are alternative algorithms that use boundary nodes and elements to decide upon the inclusion of nodes in the domain. Other algorithms such as the ``crossing number'' or the ``winding number'' methods, described in references \cite{Shimrat1962,Chinn1966,Orourke1998}, can be used in 2D. For 3D problems, the M\"oller-Trumbore algorithm \cite{Mller1997} or the AABB tree algorithm \cite{cgal_aabb_tree} can be used if the boundaries of the domain are triangulated surfaces. However, these alternative approaches are imprecise because they depend on an approximation of the boundary of the domain.

\begin{figure}[h!]
	\centering
	\includegraphics[width=12cm]{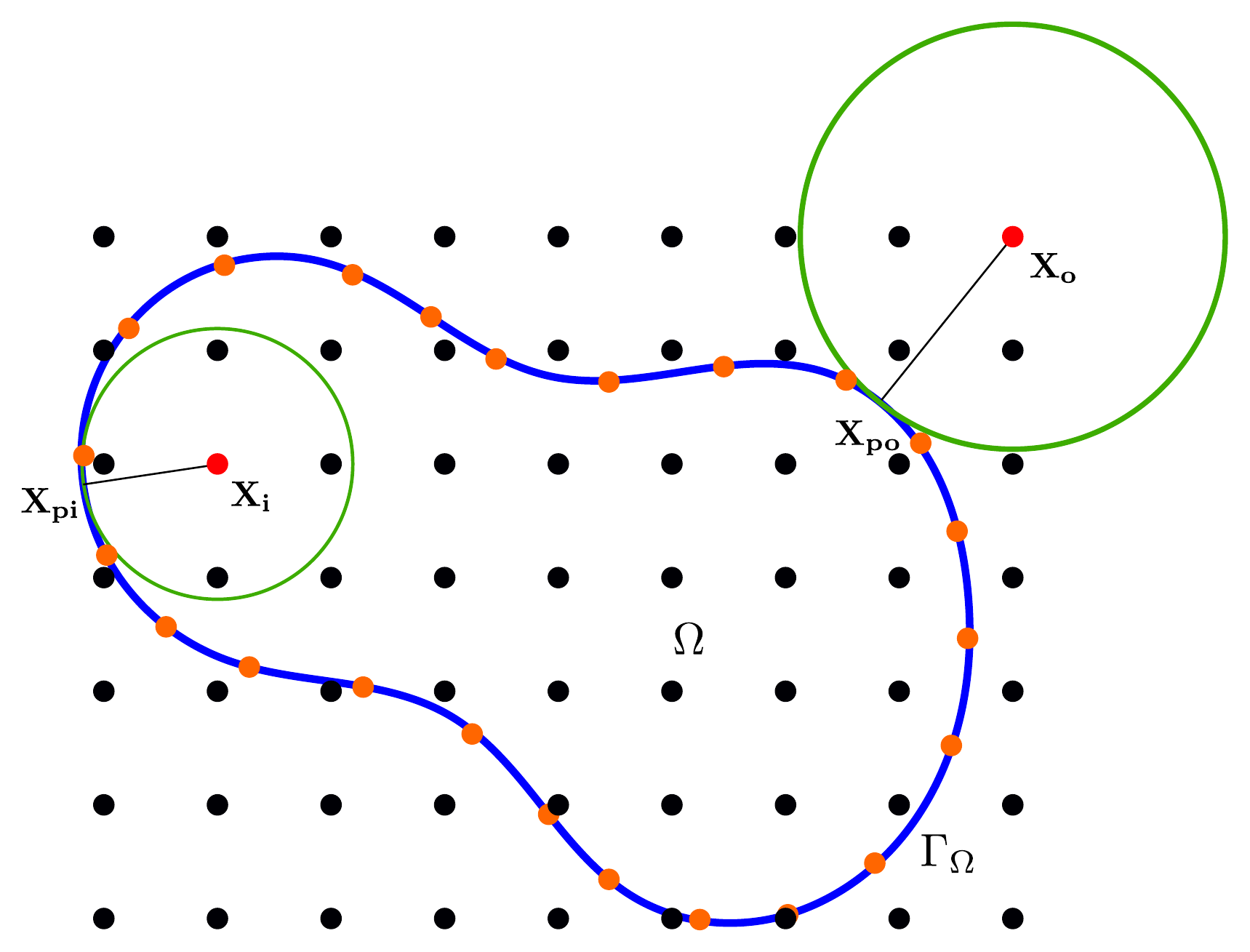};
	\caption{Identification of the position of interior nodes with respect to the domain $\Omega$. Considering a node $\mathbf{X_i}$ and its projection $\mathbf{X_{pi}}$ on the boundary of the domain $\Gamma_\Omega$, all the nodes located within a disc or sphere of radius $\left\lVert \mathbf{X_i}-\mathbf{X_{pi}} \right\rVert_2$ are located in the domain $\Omega$ if the node $\mathbf{X_i}$ is located in the domain $\Omega$. Similarly, all the nodes located within a disc or sphere of radius $\left\lVert \mathbf{X_o}-\mathbf{X_{po}} \right\rVert_2$ are located outside of the domain $\Omega$ if the node $\mathbf{X_o}$ is outside of the domain $\Omega$.}
	\label{2DNodePosition}
\end{figure}

The position of all the box nodes $\mathbf{X}$ with respect to a domain $\Omega$ do not need to be assessed. All the nodes located in a disc or sphere centered at a considered node $\mathbf{X_c}$ and of radius  $\left\lVert \mathbf{X_c}-\mathbf{X_{pc}} \right\rVert_2$, where $\mathbf{X_{pc}}$ is the the projection of $\mathbf{X_c}$ on the boundary of the domain $\Gamma_\Omega$, are located on the same side of the boundary as $\mathbf{X_c}$. This is illustrated by Figure \ref{2DNodePosition} for a node $\mathbf{X_i}$ located inside of the domain and a node $\mathbf{X_o}$ located outside of the domain.

The proximity of the interior nodes to nodes located on the boundary of the domain shall also be considered to avoid the ill conditioning of the system. This aspect is discussed in Subsection \ref{ThresholdSens}.

\subsection{From discretization to smart cloud}\label{CollocModelGen}

Transforming the discretization into a smart cloud, which contains all the required information for its solution using a collocation model and for model adaptivity, is the final step. It consists primarily in the enforcement of the boundary conditions and in the addition of additional information about the geometry useful to improve the solution and adaptive refinement. 

The topological entities are used to define the boundary conditions. The boundary conditions are most often defined on edges, for 2D problems, and on faces, for 3D problems. We defined in an input file the boundary condition associated to the topological entities of interest of the CAD geometry. The boundary conditions are transmitted from the topological entities to the collocation nodes during the discretization process. For nodes at the intersection between multiple CAD topological entities, we automatically select the most relevant boundary condition for each degree of freedom. We apply first non-zero Neumann boundary conditions, then Dirichlet boundary conditions and, finally, homogeneous Neumann boundary conditions.

We explained in Section \ref{SecMethod} that surface elements are used as part of the discretization of the process. In case of adaptive refinement, the smart nodes are used to carry the boundary condition information from the initial model to the refined models. This allows reducing the number of interactions between the collocation code and the geometry to the minimum, thus saving computational cost.

Each boundary node of the smart cloud has the following pieces of information:

\begin{itemize}
	\item the base CAD geometry;
	\item the exact normal vector;
	\item the boundary conditions;
	\item the parent CAD edge or surface(s);
	\item the boundary conditions applied to the parent CAD edge or surface(s);
	\item the connections to other boundary nodes if boundary elements are used to enforce the generalized visibility criterion \cite{Jacquemin2021} or to speed-up the refinement of the surface. 
\end{itemize}

\subsection{Threshold sensitivity analysis}\label{ThresholdSens}

Stencil nodes located too close to each other may lead to ill-conditioning of the linear system solved, at each collocation center, to obtain an approximation of the field derivatives as a function of the field itself. We used a threshold ratio, denoted by $t$, to determine if an interior node, obtained from the regular discretization of the rectangle or box, should be included in the point cloud. A node $\mathbf{X_i}$ located in the domain $\Omega$ is included in the point cloud only if the closest boundary node $\mathbf{X_b}$ is located at a distance larger than the product $t h$ (i.e. $\mathbf{X_i}$ $\in$ $\Omega\text{ if }\norm{\mathbf{X_i}-\mathbf{X_b}}_2>th$).

We used two benchmark problems from the field of linear elasticity to assess the impact of the threshold on the error. The problems considered are: a plate with an elliptical hole and an infinite body with a cylindrical hole. The exact solution is known for each of the considered benchmark problems. The governing equations of linear elasticity and the problems considered are presented in Subsections \ref{SubGovEqn} and \ref{SubPbConsidered}, respectively. 

To assess the impact of the threshold value $t$ on the solution of the benchmark problems, we selected threshold values ranging from 0.02 to 0.7. For the purpose of the sensitivity analysis, we generated coarse and fine discretizations of the considered benchmark problems (approximately 4,500 and 30,000 nodes, respectively) based on the discretization method presented in Section \ref{SecMethod}. We selected a square lattice to discretize the rectangles (bounding box) which contain the geometries. We compared the results in terms of the $L_2$ relative error norm of the von Mises stress. The calculation of this norm is described in Subsection \ref{SubPbConsidered}.

The results are presented in Figure \ref{TresholdSensitivity}. The results show that the threshold has little impact on the error for the problem of a plate with an elliptical hole. The error varies by less than 3\% for the coarse discretization. For the fine discretization, the error is lower by approximately 10\% for a threshold of 0.7. The variations are more important for the infinite body with a cylindrical hole. A threshold of 0.3 leads to the lowest error for the fine discretization. A threshold of 0.02 leads to the lowest error for the coarse discretization. For both node densities, the lowest error is approximately 75\% lower than the maximum error.

To better understand these results, we plot in Figure \ref{TresholdSensitivityCond} the condition number of the linear systems solved as part of the field derivatives approximation as a function of the threshold value. We see that the condition number increases as the threshold decreases. For the infinite body with a cylindrical hole, we observe a sharp increase in the condition number for the two lowest threshold values (i.e. 0.02 and 0.1). This increase is associated with an increase of the error. The condition number of the stencil is closely related to the node selection algorithm. In this work, we selected the stencil nodes based on the distance criterion \cite{Jensen1972}. Close to concave boundaries of the domain, we used the visibility criterion with a threshold angle of 5.0$^\circ$ as presented in reference \cite{Jacquemin2021}. Based on these results, we selected in this work a threshold value of 0.3 as it leads to the lowest error for the body with a cylindrical hole problem. We decided not to investigate other node selection algorithms which could be more suitable for the lowest threshold values.

\begin{figure}[!h]
	\centering
	\subfloat[][]{
		\begin{tikzpicture}[scale=1]
			\begin{axis}[height=7cm,width=8cm,title=Plate with an elliptical hole, ymin=0.03, ymax=1, ymode=log, xmin=0,xmax=0.8, legend entries={Coarse discretization,Fine discretization}, legend style={ at={(0.5,-0.2)}, anchor=south west,legend columns=1, cells={anchor=west},  font=\footnotesize, rounded corners=2pt,}, legend pos=north east,xlabel=Threshold, ylabel=$L_2$ Relative Error - $\sigma_{VM}$]
				\addplot+[Tblue,mark=triangle*,mark options={fill=Tblue}]   table [x=Threshold, y=L2R-VMS-Ell-Co
						, col sep=comma] {ThresholdSens.csv};
				\addplot+[Tred,mark=diamond*,mark options={fill=Tred}]   table [x=Threshold, y=L2R-VMS-Ell-Fi
						, col sep=comma] {ThresholdSens.csv};
			\end{axis}
		\end{tikzpicture}
	}
	\subfloat[][]{
		\begin{tikzpicture}[scale=1]
			\begin{axis}[height=7cm,width=8cm,title=Body with a cylindrical hole, ymin=0.0001, ymax=0.1, ymode=log, xmin=0,xmax=0.8, legend entries={Coarse discretization,Fine discretization}, legend style={ at={(0.5,-0.2)}, anchor=south west,legend columns=1, cells={anchor=west},  font=\footnotesize, rounded corners=2pt,}, legend pos=north east,xlabel=Threshold, ylabel=$L_2$ Relative Error - $\sigma_{VM}$]
				\addplot+[Tblue,mark=triangle*,mark options={fill=Tblue}]   table [x=Threshold, y=L2R-VMS-Cyl-Co
						, col sep=comma] {ThresholdSens.csv};
				\addplot+[Tred,mark=diamond*,mark options={fill=Tred}]   table [x=Threshold, y=L2R-VMS-Cyl-Fi
						, col sep=comma] {ThresholdSens.csv};
			\end{axis}
		\end{tikzpicture}
	}
	\caption{Error for threshold values ranging from $0.02$ to $0.7$ for coarse (approx. 4,500 nodes) and fine (approx. 30,000 nodes) discretizations for the plate with an elliptical hole problem (a) and for the body with a cylindrical hole problem (b). The error in terms of the $L_2$ relative error norm is presented for the von Mises stress noted $\sigma_{\text{VM}}$.}
	\label{TresholdSensitivity}

	\subfloat[][]{
		\begin{tikzpicture}[scale=1]
			\begin{axis}[height=7cm,width=8cm,title=Plate with an elliptical hole, ymin=1E3, ymax=5E7, ymode=log, xmin=0,xmax=0.8, legend entries={Coarse discretization,Fine discretization}, legend style={ at={(0.5,-0.2)}, anchor=south west,legend columns=1, cells={anchor=west},  font=\footnotesize, rounded corners=2pt,}, legend pos=north east,xlabel=Threshold, ylabel=Max. stencil condition number]
				\addplot+[Tblue,mark=triangle*,mark options={fill=Tblue}]   table [x=Threshold, y=Cond-Ell-Co
						, col sep=comma] {ThresholdSens.csv};
				\addplot+[Tred,mark=diamond*,mark options={fill=Tred}]   table [x=Threshold, y=Cond-Ell-Fi
						, col sep=comma] {ThresholdSens.csv};
			\end{axis}
		\end{tikzpicture}
	}
	\subfloat[][]{
		\begin{tikzpicture}[scale=1]
			\begin{axis}[height=7cm,width=8cm,title=Body with a cylindrical hole, ymin=1E3, ymax=5E7, ymode=log, xmin=0,xmax=0.8, legend entries={Coarse discretization,Fine discretization}, legend style={ at={(0.5,-0.2)}, anchor=south west,legend columns=1, cells={anchor=west},  font=\footnotesize, rounded corners=2pt,}, legend pos=north east,xlabel=Threshold, ylabel=Max. stencil condition number]
				\addplot+[Tblue,mark=triangle*,mark options={fill=Tblue}]   table [x=Threshold, y=Cond-Cyl-Co
						, col sep=comma] {ThresholdSens.csv};
				\addplot+[Tred,mark=diamond*,mark options={fill=Tred}]   table [x=Threshold, y=Cond-Cyl-Fi
						, col sep=comma] {ThresholdSens.csv};
			\end{axis}
		\end{tikzpicture}
	}
	\caption{Maximum stencil condition number for threshold values ranging from $0.02$ to $0.7$ for coarse (approx. 4,500 nodes) and fine (approx. 30,000 nodes) discretizations for the plate with an elliptical hole problem (a) and for the body with a cylindrical hole problem (b).}
	\label{TresholdSensitivityCond}
\end{figure}

\subsection{Discretization methods comparison}

\begin{figure}[!h]
	\centering
	\subfloat[][]{
		\begin{tikzpicture}[scale=1]
			\begin{axis}[height=7cm,width=8cm,title=Plate with an elliptical hole, ymin=1E-2, ymax=2, ymode=log, xmin=1000, xmax=1E6, xmode=log, legend entries={Square lattice discretization from CAD,Triangular lattice discretization from CAD,Discreatization from Delaunay trianglulation}, legend style={ at={(0.5,-0.2)}, anchor=south west,legend columns=1, cells={anchor=west},  font=\footnotesize, rounded corners=2pt,}, legend pos=north east,xlabel=Number of Nodes, ylabel=$L_2$ Relative Error - $\sigma_{VM}$]
				\addplot+[Tblue,mark=triangle*,mark options={fill=Tblue}]   table [x=NodeNum-Ell-CAD-Square, y=L2R-VMS-Ell-CAD-Square
						, col sep=comma] {DiscretizationComp.csv};
				\addplot+[Tred,mark=diamond*,mark options={fill=Tred}]   table [x=NodeNum-Ell-CAD-Triangle, y=L2R-VMS-Ell-CAD-Triangle
						, col sep=comma] {DiscretizationComp.csv};
				\addplot+[Tgreen,mark=*,mark options={fill=Tgreen}]   table [x=NodeNum-Ell-Gmsh, y=L2R-VMS-Ell-Gmsh
						, col sep=comma] {DiscretizationComp.csv};
				\logLogSlopeTriangle{0.85}{0.1}{0.15}{0.5}{black};
			\end{axis}
		\end{tikzpicture}
	}
	\subfloat[][]{
		\begin{tikzpicture}[scale=1]
			\begin{axis}[height=7cm,width=8cm,title=Body with a cylindrical hole, ymin=1E-5, ymax=1E-1, ymode=log, xmin=1000, xmax=1E6, xmode=log, legend entries={Square lattice discretization from CAD,Triangular lattice discretization from CAD,Discreatization from Delaunay trianglulation}, legend style={ at={(0.5,-0.2)}, anchor=south west,legend columns=1, cells={anchor=west},  font=\footnotesize, rounded corners=2pt,}, legend pos=north east,xlabel=Number of Nodes, ylabel=$L_2$ Relative Error - $\sigma_{VM}$]
				\addplot+[Tblue,mark=triangle*,mark options={fill=Tblue}]   table [x=NodeNum-Cyl-CAD-Square, y=L2R-VMS-Cyl-CAD-Square
						, col sep=comma] {DiscretizationComp.csv};
				\addplot+[Tred,mark=diamond*,mark options={fill=Tred}]   table [x=NodeNum-Cyl-CAD-Triangle, y=L2R-VMS-Cyl-CAD-Triangle
						, col sep=comma] {DiscretizationComp.csv};
				\addplot+[Tgreen,mark=*,mark options={fill=Tgreen}]   table [x=NodeNum-Cyl-Gmsh, y=L2R-VMS-Cyl-Gmsh
						, col sep=comma] {DiscretizationComp.csv};
				\logLogSlopeTriangle{0.85}{0.1}{0.35}{1}{black};
			\end{axis}
		\end{tikzpicture}
	}\\
	\caption{Comparison of the error in terms of the $L_2$ relative error norm obtained from different discretization techniques (i.e. square or triagular lattice discretization from CAD and Delaunay triangulation generated using Gmsh). The results are presented for the plate with an elliptical hole problem (a) and for the body with a cylindrical hole problem (b). The three discretization methods lead to similar errors for the plate with an elliptical hole problem. The square and triangular lattice discretizations from CAD lead to similar results for the body with a cylindrical hole problem. The error obtained with these discretizations is lower than the error obtained from the Delaunay triangulation of the domain.}
	\label{DiscretizationSensitivity}
\end{figure}
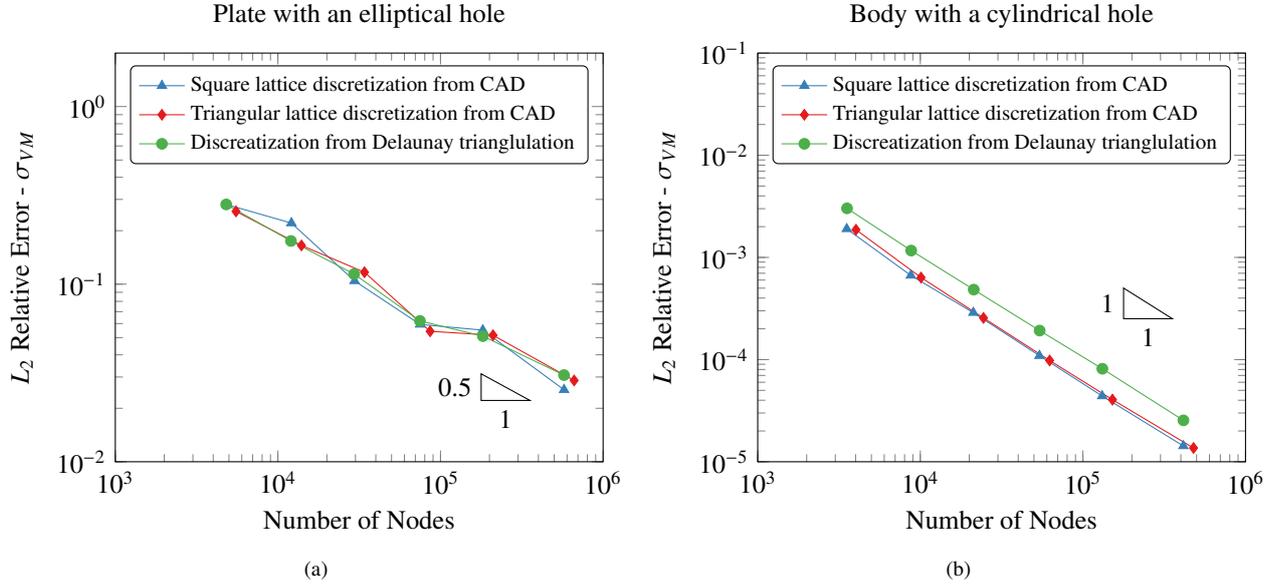

We compared the results obtained using the proposed discretization method to results obtained from discretizations generated using Gmsh \cite{Geuzaine2009}. Gmsh is a powerful finite element mesh generator which generates primarily meshes based on a Delaunay triangulation of the geometry. Gmsh is suitable for both 2D and 3D problems. We used this code to generate uniform discretizations of the geometry by using only the nodes of the generated meshes. We used the benchmark problems mentioned in Subsection \ref{ThresholdSens} for the purpose of this comparison. We compare results for a square lattice and a triangular lattice discretization of the rectangle comprising the geometry. We used a threshold value of 0.3. The results in terms of the $L_2$ relative error of the von Mises stress are shown in Figure \ref{DiscretizationSensitivity}.

The results show that the two discretization methods lead to similar errors for the plate with an elliptical hole. For the body with a cylindrical hole, the square and triangular lattice discretizations from CAD lead to errors approximately 40\% lower than those obtained from the discretization obtained from a Delaunay triangulation of the geometry.

These results give confidence in the proposed discretization methods as they lead to results not far from the ones obtained from a discretization method based on a triangulation of the domain.

\newpage

\section{Model adaptivity from CAD}\label{SecErrorIndicators}

We show in this section how the CAD geometry can be effectively used in an adaptivity scheme using smart cloud discretizations presented in Section \ref{DiscretizationSec}. We present results for problems from the field of linear elasticity. The governing equations for linear elasticity are introduced in Subsection \ref{SubGovEqn}. We used two 2D benchmark problems, for which analytical solutions are known, to assess the sensitivity of the parameters of the presented method on the quality of the error indicator and on the adaptive refinement scheme. The models considered are presented in Subsection \ref{SubPbConsidered}. Then, we present in Subsection \ref{SubErrorInd} two error indicators that we used to identify the zones where the error is the greatest. Finally, we show in Subsection \ref{AdaptSec} how error indicators are used to refine locally the domain and improve the convergence rate of the solution.

\subsection{Governing equations}\label{SubGovEqn}

We present in this section the governing equations for linear elasticity for the general case of a 3D problem. This section is based on Reference \cite{Jacquemin2021}.

The equilibrium of a domain $\Omega$ subject to body forces $\mathbf{b}$ is expressed as a function of the stress tensor $\boldsymbol{\sigma}$ by Newton's second law. For static problems, the equilibrium equation is:

\begin{equation} \label{EquilibriumEquation}
	\begin{aligned}
		\boldsymbol{\nabla} \cdot \boldsymbol{\sigma}+\mathbf{b}        & =0  \\
		\text{or} \quad \forall i \in \{1,2,3\} \quad \sigma_{ij,j}+b_i & =0.
	\end{aligned}
\end{equation}

The equilibrium equations are expressed as a function of the displacement field $\mathbf{u}$ at each node of the domain using:

\begin{itemize}
	\item the relationship between the displacement field and the strain field $\mathbf{\epsilon}$ (kinematics):
\end{itemize}
\begin{equation} \label{DispStrain}
	\begin{aligned}
		\mathbf{\epsilon}                                           & =\frac{1}{2}\left( \boldsymbol{\nabla}\mathbf{u} \otimes \boldsymbol{\nabla}\mathbf{u}^T \right) \\
		\text{or} \quad \forall i \in \{1,2,3\} \quad \epsilon_{ij} & =\frac{1}{2}\left(u_{i,j}+u_{j,i} \right),
	\end{aligned}
\end{equation}

\begin{itemize}
	\item Hooke's law which gives the relationship between the strain field and the stress field (presented here in Voigt form). This is the constitutive law:
\end{itemize}

\begin{equation} \label{HookLaw}
	\begin{bmatrix}
		\sigma_{11} \\
		\sigma_{22} \\
		\sigma_{33} \\
		\sigma_{23} \\
		\sigma_{13} \\
		\sigma_{12} \\
	\end{bmatrix} = \frac{E}{\left(1+\nu\right)\left(1-2\nu\right)}
	\begin{bmatrix}
		1-\nu & \nu   & \nu   & 0      & 0      & 0      \\
		\nu   & 1-\nu & \nu   & 0      & 0      & 0      \\
		\nu   & \nu   & 1-\nu & 0      & 0      & 0      \\
		0     & 0     & 0     & 1-2\nu & 0      & 0      \\
		0     & 0     & 0     & 0      & 1-2\nu & 0      \\
		0     & 0     & 0     & 0      & 0      & 1-2\nu \\
	\end{bmatrix}
	\begin{bmatrix}
		\epsilon_{11} \\
		\epsilon_{22} \\
		\epsilon_{33} \\
		\epsilon_{23} \\
		\epsilon_{13} \\
		\epsilon_{12} \\
	\end{bmatrix}.
\end{equation}

The above equations can be used for 2D problems using either the plane stress assumption (i.e. $\sigma_{33}=0$, $\sigma_{13}=0$ and $\sigma_{23}=0$) or the plane strain assumption (i.e. $\epsilon_{33}=0$, $\epsilon_{13}=0$ and $\epsilon_{33}=0$).

Dirichlet and Neumann boundary conditions are respectively applied to the degrees of freedom of the collocation nodes located on the boundaries $\Gamma_D$ and $\Gamma_N$. The known displacement field $\mathbf{u^e}$ is applied on $\Gamma_D$. An external pressure $\mathbf{f^e}$ is applied to the nodes located on $\Gamma_N$. The outer normal $\mathbf{n_N}$ allows the computation of the pressure at the nodes of $\Gamma_N$. Dirichlet and Neumann boundary conditions can be applied to different degrees of freedom of the same node.

\begin{equation} \label{BoundConditions}
	\begin{aligned}
		\mathbf{u}                      & =\mathbf{u^e}                                                                          & \quad & \text{on} \quad \Gamma_D  \\
		\boldsymbol{\sigma}\mathbf{n_N} & =\mathbf{f^e} \quad \text{or} \quad \forall i \in \{1,2,3\} \quad \sigma_{ij}n_j=f^e_i & \quad & \text{on} \quad \Gamma_N.
	\end{aligned}
\end{equation}

\subsection{Benchmark problems considered}\label{SubPbConsidered}

We present in this section the benchmark problems that we considered as part of our analysis. We selected 2D problems with known solutions from the field of linear elasticity. These are:

\begin{itemize}
	\item an infinite plate with an elliptical hole under biaxial loading;
	\item an infinite body with a cylindrical hole under remote stress loading.
\end{itemize}

\begin{figure}[!h]
	\centering
	\begin{tikzpicture}[x=0.7cm,y=0.7cm]
		\definecolor{LightBlue}{rgb}{0.733,0.878,0.890};

		\draw [-stealth,thick,red] (0.000,-3.000) -- (0.000,-5.045);
		\draw [-stealth,thick,red] (0.000,3.000) -- (0.000,5.045);
		\draw [-stealth,thick,red] (0.400,-3.000) -- (0.536,-5.032);
		\draw [-stealth,thick,red] (0.400,3.000) -- (0.536,5.032);
		\draw [-stealth,thick,red] (0.800,-3.000) -- (1.071,-4.993);
		\draw [-stealth,thick,red] (0.800,3.000) -- (1.071,4.993);
		\draw [-stealth,thick,red] (1.200,-3.000) -- (1.602,-4.929);
		\draw [-stealth,thick,red] (1.200,3.000) -- (1.602,4.929);
		\draw [-stealth,thick,red] (1.600,-3.000) -- (2.127,-4.843);
		\draw [-stealth,thick,red] (1.600,3.000) -- (2.127,4.843);
		\draw [-stealth,thick,red] (2.000,-3.000) -- (2.645,-4.739);
		\draw [-stealth,thick,red] (2.000,3.000) -- (2.645,4.739);
		\draw [-stealth,thick,red] (2.400,-3.000) -- (3.152,-4.623);
		\draw [-stealth,thick,red] (2.400,3.000) -- (3.152,4.623);
		\draw [-stealth,thick,red] (2.800,-3.000) -- (3.649,-4.503);
		\draw [-stealth,thick,red] (2.800,3.000) -- (3.649,4.503);
		\draw [-stealth,thick,red] (3.200,-3.000) -- (4.137,-4.386);
		\draw [-stealth,thick,red] (3.200,3.000) -- (4.137,4.386);
		\draw [-stealth,thick,red] (3.600,-3.000) -- (4.616,-4.278);
		\draw [-stealth,thick,red] (3.600,3.000) -- (4.616,4.278);
		\draw [-stealth,thick,red] (4.000,-3.000) -- (5.092,-4.185);
		\draw [-stealth,thick,red] (4.000,3.000) -- (5.092,4.185);
		\draw [-stealth,thick,red] (4.400,-3.000) -- (5.566,-4.106);
		\draw [-stealth,thick,red] (4.400,3.000) -- (5.566,4.106);
		\draw [-stealth,thick,red] (4.800,-3.000) -- (6.041,-4.041);
		\draw [-stealth,thick,red] (4.800,3.000) -- (6.041,4.041);
		\draw [-stealth,thick,red] (5.200,-3.000) -- (6.519,-3.990);
		\draw [-stealth,thick,red] (5.200,3.000) -- (6.519,3.990);
		\draw [-stealth,thick,red] (5.600,-3.000) -- (7.000,-3.949);
		\draw [-stealth,thick,red] (5.600,3.000) -- (7.000,3.949);
		\draw [-stealth,thick,red] (6.000,-3.000) -- (7.484,-3.916);
		\draw [-stealth,thick,red] (6.000,-2.625) -- (7.459,-3.424);
		\draw [-stealth,thick,red] (6.000,-2.250) -- (7.431,-2.931);
		\draw [-stealth,thick,red] (6.000,-1.875) -- (7.402,-2.438);
		\draw [-stealth,thick,red] (6.000,-1.500) -- (7.373,-1.946);
		\draw [-stealth,thick,red] (6.000,-1.125) -- (7.346,-1.456);
		\draw [-stealth,thick,red] (6.000,-0.750) -- (7.324,-0.969);
		\draw [-stealth,thick,red] (6.000,-0.375) -- (7.309,-0.484);
		\draw [-stealth,thick,red] (6.000,0.000) -- (7.304,0.000);
		\draw [-stealth,thick,red] (6.000,0.375) -- (7.309,0.484);
		\draw [-stealth,thick,red] (6.000,0.750) -- (7.324,0.969);
		\draw [-stealth,thick,red] (6.000,1.125) -- (7.346,1.456);
		\draw [-stealth,thick,red] (6.000,1.500) -- (7.373,1.946);
		\draw [-stealth,thick,red] (6.000,1.875) -- (7.402,2.438);
		\draw [-stealth,thick,red] (6.000,2.250) -- (7.431,2.931);
		\draw [-stealth,thick,red] (6.000,2.625) -- (7.459,3.424);
		\draw [-stealth,thick,red] (6.000,3.000) -- (7.484,3.916);

		\draw [draw=none,fill=LightBlue] (0,-3) rectangle (6,3);
		\draw [draw=none,fill=white] (0,-0.2) arc(-91:91:2.95 and 0.2);
		\draw [thick] (0,-0.2) arc(-90:90:3 and 0.2);
		\draw [thick] (0,0.2) -- (0,3) -- (6,3) -- (6,-3) -- (0,-3) -- (0,-0.2);
		\draw plot [only marks, mark=square*, mark size=0.2] coordinates {(0,-0.2) (0,0.2)};

		\draw (-0.75,0) edge[out=0,in=-90,-stealth] (0.378,0.198);
		\node at (-0.75,0) [left] {Stress-free surface};
		\draw (-0.75,1.7) edge[-stealth] (0,1.7);
		\node at (-2.7,2) [right] {$u_1$=0};
		\node at (-2.7,1.4) [right] {$u_2$=free};

		\node [red,align=right] at (-0.6,3.9) [left] {Applied displacement};
		\node [red,align=right] at (-0.6,3.5) [left] {field};

		\draw (-4,-3) edge[-stealth] (-3,-3);
		\node at (-3,-3) [right] {$X_1$};
		\draw (-4,-3) edge[-stealth] (-4,-2);
		\node at (-4,-2) [above] {$X_2$};

	\end{tikzpicture}
	\caption{2D model of an infinite plate with an elliptical hole under biaxial loading. Symmetry boundary conditions are applied to the vertical edge on the left. Considering a displacement field denoted $u$, this boundary condition corresponds to $u_1$=0, $u_2$=free. Stress-free surface boundary conditions are applied to the boundary of the elliptical hole. The displacement field of the exact solution is applied to the other boundaries of the domain. To improve the quality of the solution, three boundary nodes at the tip are considered as interior nodes.}
	\label{EllipticalHolePlate}
\end{figure}
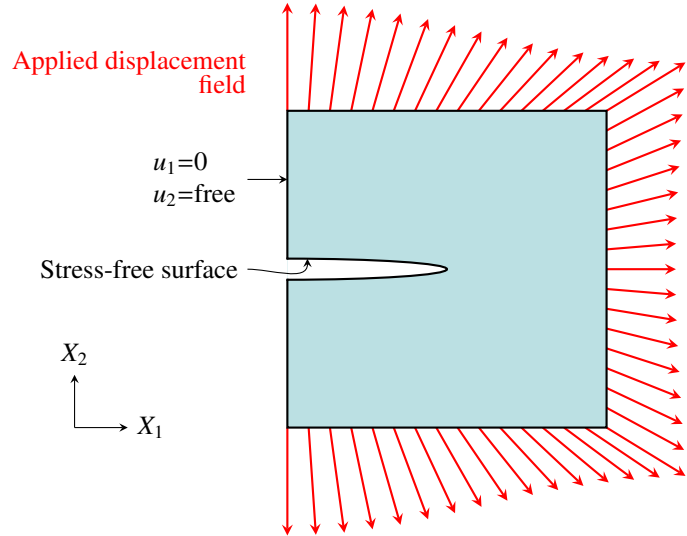

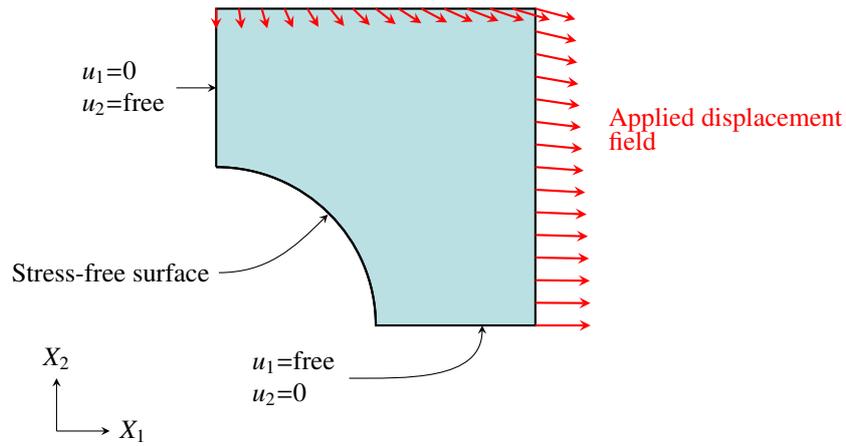
\begin{figure}[!h]
	\centering
	\begin{tikzpicture}[x=0.7cm,y=0.7cm]
		\definecolor{LightBlue}{rgb}{0.733,0.878,0.890};

		\draw [draw=none,fill=LightBlue] (0,0) rectangle (6,6);
		\draw [draw=none,fill=white] (3.0,0) arc(0:181:3);
		\draw [thick] (3,0) arc(0:90:3);
		\draw [thick] (0,3) -- (0,6) -- (6,6) -- (6,0) -- (3,0) arc(0:90:3);

		\draw (-0.75,4.5) edge[-stealth] (0,4.5);
		\node at (-2.7,4.8) [right] {$u_1$=0};
		\node at (-2.7,4.2) [right] {$u_2$=free};

		\draw (0.05,1) edge[out=0,in=-135,-stealth] (2.1,2.1);
		\node at (0.05,1) [left] {Stress-free surface};

		\draw (2.5,-1) edge[out=0,in=-90,-stealth] (5,0);
		\node at (0.5,-0.7) [right] {$u_1$=free};
		\node at (0.5,-1.3) [right] {$u_2$=0};

		\node [red,align=right] at (7.2,3.9) [right] {Applied displacement};
		\node [red,align=right] at (7.2,3.5) [right] {field};

		\draw [-stealth,thick,red] (6.000,0.000) -- (7.038,0.000);
		\draw [-stealth,thick,red] (6.000,6.000) -- (6.754,5.798);
		\draw [-stealth,thick,red] (0.000,6.000) -- (0.000,5.618);
		\draw [-stealth,thick,red] (6.000,0.429) -- (7.035,0.423);
		\draw [-stealth,thick,red] (6.000,0.857) -- (7.024,0.846);
		\draw [-stealth,thick,red] (6.000,1.286) -- (7.008,1.267);
		\draw [-stealth,thick,red] (6.000,1.714) -- (6.987,1.687);
		\draw [-stealth,thick,red] (6.000,2.143) -- (6.963,2.104);
		\draw [-stealth,thick,red] (6.000,2.571) -- (6.937,2.519);
		\draw [-stealth,thick,red] (6.000,3.000) -- (6.909,2.932);
		\draw [-stealth,thick,red] (6.000,3.429) -- (6.882,3.344);
		\draw [-stealth,thick,red] (6.000,3.857) -- (6.856,3.754);
		\draw [-stealth,thick,red] (6.000,4.286) -- (6.832,4.163);
		\draw [-stealth,thick,red] (6.000,4.714) -- (6.809,4.572);
		\draw [-stealth,thick,red] (6.000,5.143) -- (6.789,4.980);
		\draw [-stealth,thick,red] (6.000,5.571) -- (6.771,5.389);

		\draw [-stealth,thick,red] (5.571,6.000) -- (6.274,5.789);
		\draw [-stealth,thick,red] (5.143,6.000) -- (5.792,5.779);
		\draw [-stealth,thick,red] (4.714,6.000) -- (5.310,5.767);
		\draw [-stealth,thick,red] (4.286,6.000) -- (4.827,5.754);
		\draw [-stealth,thick,red] (3.857,6.000) -- (4.343,5.739);
		\draw [-stealth,thick,red] (3.429,6.000) -- (3.859,5.723);
		\draw [-stealth,thick,red] (3.000,6.000) -- (3.375,5.705);
		\draw [-stealth,thick,red] (2.571,6.000) -- (2.891,5.687);
		\draw [-stealth,thick,red] (2.143,6.000) -- (2.407,5.669);
		\draw [-stealth,thick,red] (1.714,6.000) -- (1.925,5.653);
		\draw [-stealth,thick,red] (1.286,6.000) -- (1.443,5.639);
		\draw [-stealth,thick,red] (0.857,6.000) -- (0.961,5.627);
		\draw [-stealth,thick,red] (0.429,6.000) -- (0.480,5.620);

		\draw (-3,-2) edge[-stealth] (-2,-2);
		\node at (-2,-2) [right] {$X_1$};
		\draw (-3,-2) edge[-stealth] (-3,-1);
		\node at (-3,-1) [above] {$X_2$};

	\end{tikzpicture}
	\caption{2D model of an infinite body with a cylindrical hole under remote stress loading. Symmetry boundary conditions are applied to the vertical edge on the left and to the horizontal edge at the bottom. Considering a displacement field denoted $u$, these boundary conditions correspond respectively to $u_1$=0, $u_2$=free and $u_1$=free, $u_2$=0. Stress-free surface boundary conditions are applied to the boundary of the hole. The displacement field of the exact solution is applied to the other boundaries of the domain.}
	\label{CylindricalHolePlate}
\end{figure}

\begin{figure}[h!]
	\centering
	\subfloat[][]
	{
		\begin{tabular}{@{}r@{}l@{}}
			\begin{tabular}{@{}l@{}} \includegraphics[width=7cm]{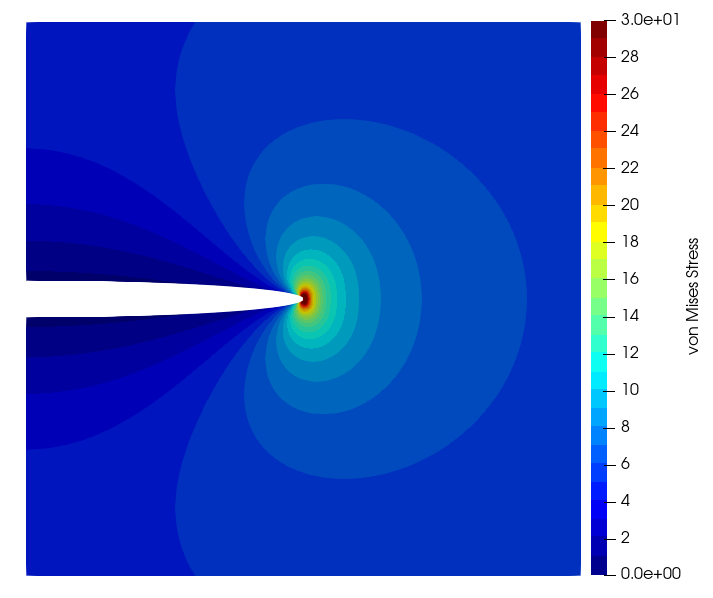} \\
			\end{tabular}
		\end{tabular}
	}
	\subfloat[][]
	{
		\begin{tabular}{@{}r@{}l@{}}
			\begin{tabular}{@{}l@{}} \quad \includegraphics[width=7cm]{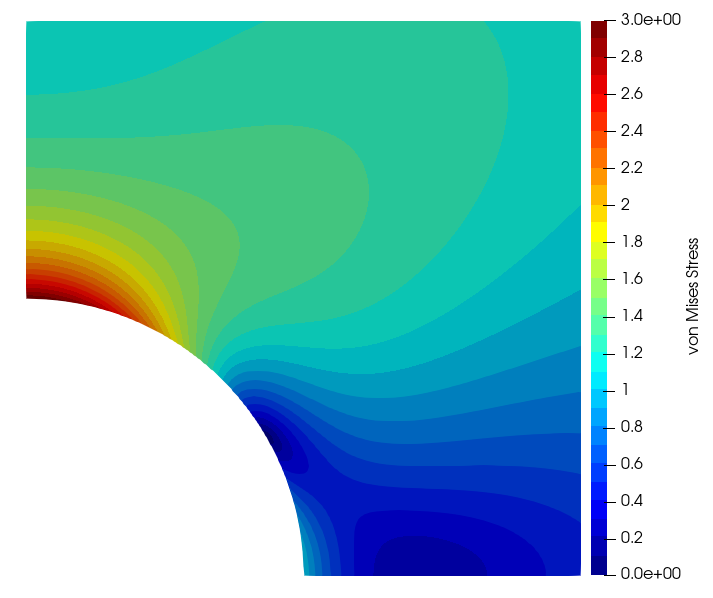} \\
			\end{tabular}
		\end{tabular}
	}
	\caption{Solutions in terms of von Mises stress for the problem of an infinite plate with an elliptical hole under biaxial loading (a) and for the problem of an infinite body with a cylindrical hole under remote stress loading (b). The maximum von Mises stress is 150 at the point of highest curvature of the elliptical hole for the problem of an infinite plate with an elliptical hole. The von Mises stress decreases rapidly as the distance from the tip increases. We truncated the colorbar of the von Mises stress field to 30 to show the variations of the stress over the domain.}
	\label{2DProbemsSolution}
\end{figure}

The domain considered and the boundary conditions applied are presented in Figure \ref{EllipticalHolePlate} and \ref{CylindricalHolePlate}, respectively, for the first and and second problems. We show the exact solution of these problems in terms of von Mises stress in Figure \ref{2DProbemsSolution}. The von Mises stress solution of the first problem varies rapidly at the point of highest curvature of the ellipse. The solution in terms of von Mises stress is smoother for the second problem.

We used the $L_2$ relative error norm (denoted by $L_2R$) and the $L_2$ weighted error norm (denoted by $L_2W$) in this work to compare the results obtained to the reference solutions. 

At a collocation node $\mathbf{X_k}$ the exact stress and approximated stress solutions are denoted $\sigma_{ij}^{e}(\mathbf{X_k})$ and $\sigma_{ij}^{h}(\mathbf{X_k})$, respectively. Considering a domain $\Omega$ discretized by $n$ collocation nodes, the $L_2$ relative error norm is calculated as per Equation \ref{L2RelativeNorm}. The $L_2$ weighted error norm is calculated as per Equation \ref{L2WeightedNorm}.

\begin{equation} \label{L2RelativeNorm}
	L_2R(\sigma_{ij})= \frac{\sqrt{\sum_{k=1}^{n}{\left(\sigma_{ij}^{e}(\mathbf{X_k}) - \sigma_{ij}^{h}(\mathbf{X_k})\right)^2}}}{\sqrt{\sum_{k=1}^{n}{{\sigma_{ij}^{e}(\mathbf{X_k})}^2}}}.
\end{equation}

\begin{equation} \label{L2WeightedNorm}
	L_2W(\sigma_{ij})= \frac{\sqrt{\sum_{k=1}^{n}{\left(\sigma_{ij}^{e}(\mathbf{X_k}) - \sigma_{ij}^{h}(\mathbf{X_k})\right)^2}}}{n}.
\end{equation}

\begin{samepage}
	\subsection{Error indicators}\label{SubErrorInd}

	We describe in this section two types of error indicators that we used to assess the need for local refinement of the discretization:

	\begin{itemize}
		\item a ZZ-type error indicator;
		\item a residual-type error indicator.
	\end{itemize}
\end{samepage}

We use the term ``error indicator" in this article rather than the term ``error estimator". Our methods only give an indication of the zones of the domain where the solution is expected to be the most imprecise rather than an estimation of the exact error. Therefore, the computed error should be considered relatively to the error computed at other locations of the domain rather than as an estimation of the true error. The considered indicators are described in the subsections below.

\subsubsection{ZZ-type error indicator for the GFD method}\label{ZZIndicatorSection}

The ZZ-type indicator refers to the class of error estimators introduced by Zienkiewicz and Zhu in 1987 \cite{Zienkiewicz1987} and extended by Bordas and Duflot \cite{Bordas2007b,Duflot2008} and Rodenas et al. \cite{Rodenas2008} to enriched approximations. Zienkiewicz and Zhu used a moving least square approximation of the stress field (for linear elastic problems) computed at superconvergent points to estimate the error. The moving least square approximation is used to extrapolate the stress computed at the superconvergent points at nodes of a selected patch. ZZ-type error estimators can be understood as an indication of the smoothness of the computed stress field over the selected patch. If the stress field is smooth, the difference in terms of stress at the recovery nodes will be small. A sharp variation of the stress field leads to a large error. We build on this idea to define an error indicator in the framework of the GFD method.

In many engineering problems, the stress field and in particular the von Mises stress, is the parameter of primary interest as it is a criterion associated with the onset of material yielding (von Mises plasticity). We considered this stress criterion to compute a ZZ-type error indicator. The von Mises criterion (noted $\sigma_{\text{vM}}$) is calculated using the equation:

\begin{equation} \label{vonMisesCalc}
	\sigma_{\text{vM}}=\sqrt{\frac{\left( \sigma_{11}-\sigma_{22} \right)^2+ \left( \sigma_{22}-\sigma_{33} \right)^2 + \left( \sigma_{33}-\sigma_{11} \right)^2 + 6\left( \sigma_{12}^2+\sigma_{23}^2+\sigma_{31}^2 \right)}{2}}.
\end{equation}

The solution of a linear elastic problem using the GFD methods requires the estimation of the first and second derivatives of the displacement field $\mathbf{u}$ to solve Equation (\ref{EquilibriumEquation}) at the collocation nodes located in the domain and Equation (\ref{BoundConditions}) at the collocation nodes located on the boundary of the domain. A detailed description of all the steps involved in the solution of a problem using the GFD method is presented in reference \cite{Jacquemin2019}. Once the system of equations is solved, the solution in terms of displacement field (usually at the collocation nodes) is used to compute the stress field using Equation (\ref{DispStrain}) and Equation (\ref{HookLaw}).

In this work, collocation is done at the nodes. This means that the partial differential equations and the boundary conditions are solved at the nodes. Therefore, the von Mises stress can be computed at each node of the domain.

We computed an indication of the error at each node of the domain. For this, we used the von Mises stress calculated at each node based on the classical GFD method and based on a smoothed (recovered) von Mises field computed using a moving least square approximation of the von Mises stress values obtained at each collocation node.

\begin{figure}[h!]
	\centering
	\includegraphics[width=8cm]{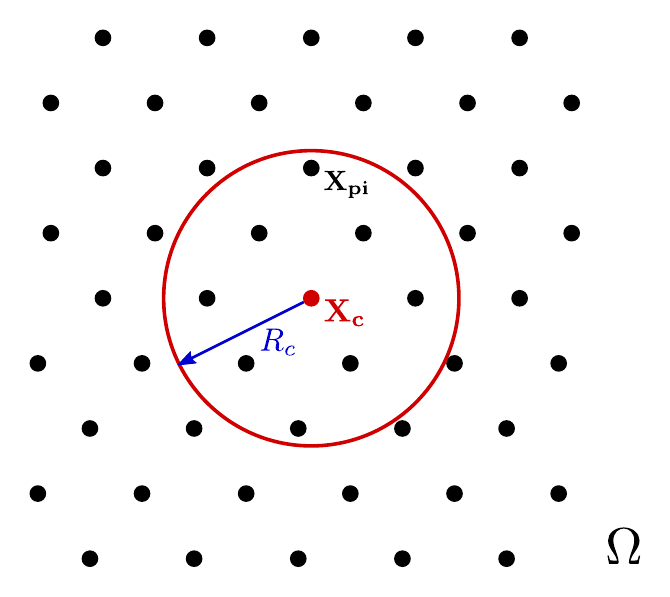};
	\caption{Discretization of a portion of a domain $\Omega$ and identification of the nodes $\mathbf{X_{pi}}$, located within a support of radius $R_c$ of a collocation node $\mathbf{X_c}$, involved in the computation of the smooth von Mises field $\sigma_{\text{vM}}^s\left(\mathbf{X_c}\right)$.}
	\label{CollocDomain}
\end{figure}

More specifically, considering a domain $\Omega$, the von Mises stress obtained at a collocation node $\mathbf{X_c}$ (noted $\sigma_{\text{vM}}^c\left(\mathbf{X_c}\right)$) is compared to a moving least square approximation of the von Mises stress field at the collocation node $\mathbf{X_c}$ (noted $\sigma_{\text{vM}}^s\left(\mathbf{X_c}\right)$). The moving least square approximation is calculated based on the von Mises stress $\sigma_{\text{vM}}^c\left(\mathbf{X_{pi}}\right)$ computed at the support nodes $\mathbf{X_{pi}}$ of the collocation node $\mathbf{X_c}$ (see Figure \ref{CollocDomain}).

The moving least square approximation is computed using a second order polynomial basis $\mathbf{p}$ (or a polynomial of the same order as the GFD approximation) and a vector of coefficients $\mathbf{a}$ determined for each collocation node. For the case of a 2D problem, the polynomial basis $\mathbf{p}$ at a point $\mathbf{X}=\left[ x,y \right]^T$ in the vicinity of $\mathbf{X_c}=\left[ x_c,y_c \right]^T$ is:

\begin{equation}\label{ZZ_MLS_PBasis}
	\mathbf{p\left(\mathbf{X},\mathbf{X_c}\right)}=
	\begin{bmatrix}
		1                  \\
		(x - x_c)          \\	(y - y_c) \\	(x - x_c)^2 \\
		(x - x_c)(y - y_c) \\	(y - y_c)^2
	\end{bmatrix}.
\end{equation}

For a collocation node $\mathbf{X_c}$, the smooth von Mises stress field $\sigma_{\text{vM}}^s$ is written:

\begin{equation} \label{ZZ_MLS}
	\sigma_{\text{vM}}^s\left(\mathbf{X},\mathbf{X_c}\right) = \mathbf{p\left(\mathbf{X},\mathbf{X_c}\right)}^T \mathbf{a}\left(\mathbf{X_c}\right).
\end{equation}

The coefficients $\mathbf{a}\left(\mathbf{X_c}\right)$ are computed to minimize the error between $\sigma_{\text{vM}}^c$ and $\sigma_{\text{vM}}^s$. For a collocation node $\mathbf{X_c}$ which has $m$ support nodes $\mathbf{X_{pi}}$ ($\mathbf{X_c}$ is not considered as a support node), we write the functional $B\left(\mathbf{X_c}\right)$ presented in Equation (\ref{ZZ_FunctionalB}). The error is weighted by a function $w\mathbf{\left(\mathbf{X},\mathbf{X_c}\right)}$ which depends on the support radius of the collocation node $\mathbf{X_c}$ and on a selected radial basis function.

\begin{equation} \label{ZZ_FunctionalB}
	B\left(\mathbf{X_c}\right)= \sum_{i=1}^m {w\mathbf{\left(\mathbf{X_{pi}},\mathbf{X_c}\right)} \left(\mathbf{p\left(\mathbf{X_{pi}},\mathbf{X_c}\right)}^T \mathbf{a\left(X_c\right)} - \sigma^c\left(\mathbf{X_{pi}}\right) \right)^2}.
\end{equation}

The error is minimized at each collocation node $\mathbf{X_c}$ when:

\begin{equation} \label{ZZ_FunctionalB_Deriv}
	\frac{\partial {B\left(\mathbf{X_c}\right)}}{\partial \mathbf{a\left(X_c\right)}}=0.
\end{equation}

This minimization problem can be transformed into a linear problem of the form:

\begin{equation} \label{ZZ_FunctionalB_LinearProb}
	\mathbf{A\left(X_c\right)} \mathbf{a\left(X_c\right)} = \mathbf{E\left(X_c\right)} \mathbf{f\left(X_c\right)}.
\end{equation}

For a polynomial basis of size $q$ ($q=6$ for a 2D second order basis), the matrices $\mathbf{A\left(X_c\right)}$, $\mathbf{E\left(X_c\right)}$ and $\mathbf{f\left(X_c\right)}$ are:

\begin{align}
	\mathbf{A\left(X_c\right)} & =\begin{bmatrix}
		m_{11} & m_{12} & \dots & m_{1q} \\
		m_{21} & m_{22} & \dots & m_{2q} \\
		\vdots &        &       & \vdots \\
		m_{q1} & m_{q2} & \dots & m_{qq} \\
	\end{bmatrix} \in \R^{q \times q}, \label{ZZ_MLS_MatA} \\[2ex]
	\mathbf{E\left(X_c\right)} & =\begin{bmatrix}
		w\mathbf{\left(\mathbf{X_{p1}},\mathbf{X_c}\right)} \mathbf{p\left(\mathbf{X_{p1}},\mathbf{X_c}\right)}_1 & \dots & w\mathbf{\left(\mathbf{X_{pm}},\mathbf{X_c}\right)} \mathbf{p\left(\mathbf{X_{pm}},\mathbf{X_c}\right)}_1 \\
		w\mathbf{\left(\mathbf{X_{p1}},\mathbf{X_c}\right)} \mathbf{p\left(\mathbf{X_{p1}},\mathbf{X_c}\right)}_2 & \dots & w\mathbf{\left(\mathbf{X_{pm}},\mathbf{X_c}\right)} \mathbf{p\left(\mathbf{X_{pm}},\mathbf{X_c}\right)}_2 \\
		\vdots                                                                                                    &       & \vdots                                                                                                    \\
		w\mathbf{\left(\mathbf{X_{p1}},\mathbf{X_c}\right)} \mathbf{p\left(\mathbf{X_{p1}},\mathbf{X_c}\right)}_q & \dots & w\mathbf{\left(\mathbf{X_{pm}},\mathbf{X_c}\right)} \mathbf{p\left(\mathbf{X_{pm}},\mathbf{X_c}\right)}_q \\
	\end{bmatrix} \in \R^{q \times m}, \label{ZZ_MLS_MatC} \\[2ex]
	\mathbf{f\left(X_c\right)} & =\begin{bmatrix}
		\sigma^c\left(\mathbf{X_{p1}}\right) & \sigma^c\left(\mathbf{X_{p2}}\right) & \dots & \sigma^c\left(\mathbf{X_{pm}}\right)
	\end{bmatrix}^T, \label{ZZ_MLS_MatF}
\end{align}
where
\begin{equation}\label{ZZ_MLS_Moments}
	m_{ij}= \sum_{k=1}^m {w\mathbf{\left(\mathbf{X_{pk}},\mathbf{X_c}\right)} \mathbf{p\left(\mathbf{X_{pk}},\mathbf{X_c}\right)}_i \mathbf{p\left(\mathbf{X_{pk}},\mathbf{X_c}\right)}_j}.
\end{equation}

$\mathbf{p\left(\mathbf{X_{pk}},\mathbf{X_c}\right)}_i$ refers to the $i^\text{th}$ component of the vector $\mathbf{p\left(\mathbf{X_{pk}},\mathbf{X_c}\right)}$.

The solution of Equation (\ref{ZZ_FunctionalB_LinearProb}) at $\mathbf{X_c}$ allows the computation of $\sigma^s\left(\mathbf{X_c}\right)$ and of an error indicator $e\left(\mathbf{X_c}\right)$ as follows:

\begin{equation}\label{ZZ_ErrorIndicator}
	e\left(\mathbf{X_c}\right)=\left|\sigma^c\left(\mathbf{X_c}\right)-\sigma^s\left(\mathbf{X_c}\right)\right|.
\end{equation}

We presented in this section a ZZ-type error indicator based on an assessment of the smoothness of the solution of the PDE. We present in the next section of the article a residual-type error indicator.

\subsubsection{Residual-type error indicator}\label{ResIndicatorSection}

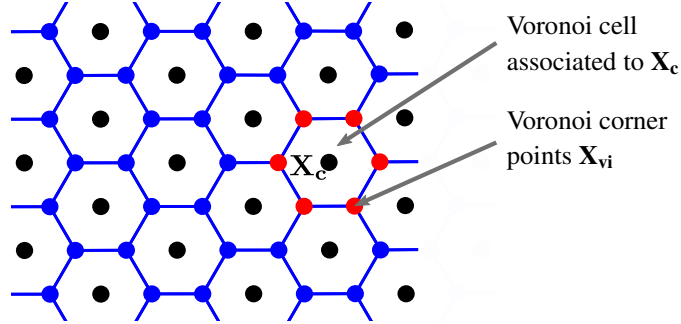
\begin{figure}[h!]
	\centering
	\begin{tikzpicture}
		\def\svgwidth{7cm}
		\node at (0,0) {\includegraphics[width=7cm]{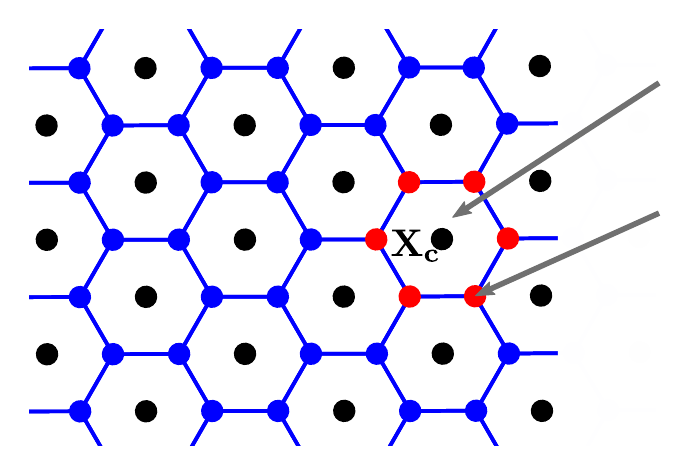}};
		\node[color=black] at (3.2,1.8) [right] {Voronoi cell};
		\node[color=black] at (3.2,1.3) [right] {associated to $\mathbf{X_c}$};
		\node[color=black] at (3.2,0.55) [right] {Voronoi corner};
		\node[color=black] at (3.2,0.05) [right] {points $\mathbf{X_{vi}}$};
	\end{tikzpicture}
	\caption{Voronoi diagram for a set of collocation nodes. The residual-type error indicator at a collocation node $\mathbf{X_c}$ is computed based on the residual of the PDE at the corner points $\mathbf{X_{vi}}$ of the Voronoi cell associated to $\mathbf{X_c}$.}
	\label{VoroCell}
\end{figure}

Residual-type error indicators are based on an estimation of the residual of the PDE at locations where it is not enforced as part of the solution process. This type of estimator has been widely used in literature: for example in reference \cite{Driscoll2007,Oh2018}. In this work, we selected the corners of the Voronoi cell surrounding a considered collocation node $\mathbf{X_c}$ to compute the residual of the PDE. The residual error at a collocation node $\mathbf{X_c}$ is calculated as the average of the residual of the PDE at each Voronoi center $\mathbf{X_{vi}}$ of the Voronoi cell associated to $\mathbf{X_c}$ (see Figure \ref{VoroCell}).

The residual at each Voronoi center $\mathbf{X_{vi}}$ is approximated using the solution at the neighboring collocation nodes using the GFD method. For a collocation node $\mathbf{X_c}$ for which the associated Voronoi cell has $q$ corner points $\mathbf{X_{vi}}$, we calculate the error indicator $e\left(\mathbf{X_c}\right)$ as follows:

\begin{equation}\label{Res_ErrorIndicator}
	e\left(\mathbf{X_c}\right)=\frac{1}{q}\sum_{i=1}^{q}{\left| \boldsymbol{\nabla} \cdot \boldsymbol{\sigma}\left( \mathbf{X_{vi}} \right)+\mathbf{b}\left( \mathbf{X_{vi}} \right) \right|}.
\end{equation}

\subsubsection{Parameter Variation and Indicators Comparison}

We compare in this section the error obtained from the indicators presented in Subsections \ref{ZZIndicatorSection} and \ref{ResIndicatorSection}. We show the impact of some parameters of the methods on the calculated error, compare the error indicators to the true error, and also compare the convergence rates.

We start by presenting a comparison of the spatial pattern of the error for the two benchmark problems considered. For this, we selected discretizations composed of approximately 200,000 nodes for the plate with an elliptical hole and 140,000 nodes for the body with a cylindrical hole.

\begin{figure}[h!]
	\centering
	\begin{tabular}{|M{4cm}|c|c|}
		\hline
		                                           & \textbf{Square lattice}                                      & \textbf{Triangular lattice}                                            \\
		\hline
		\multirow{2}{*}{Exact Error $\sigma_{vM}$} &                                                              &                                                                        \\[-1ex]
		                                           & \includegraphics[width=4cm]{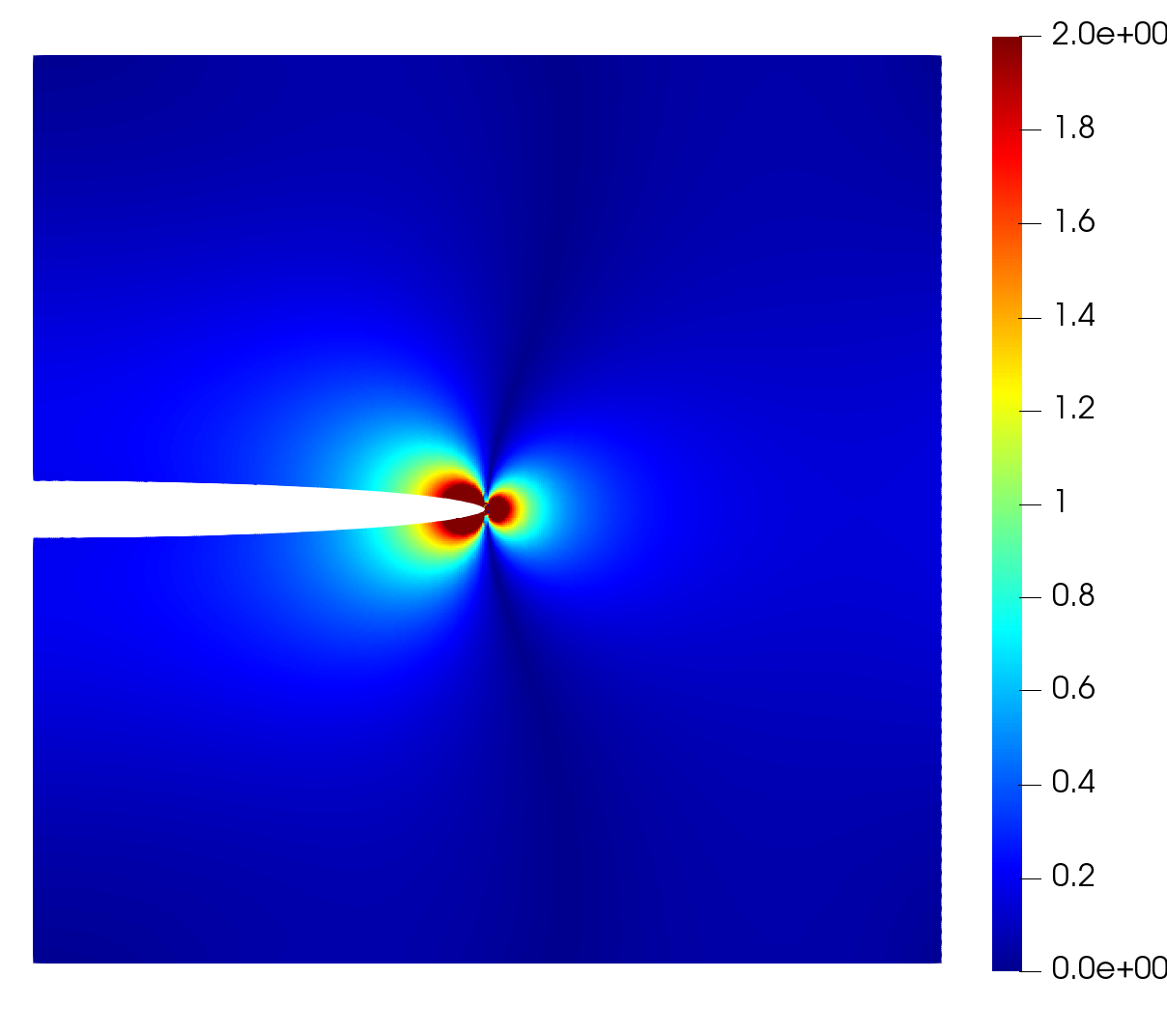} & \includegraphics[width=4cm]{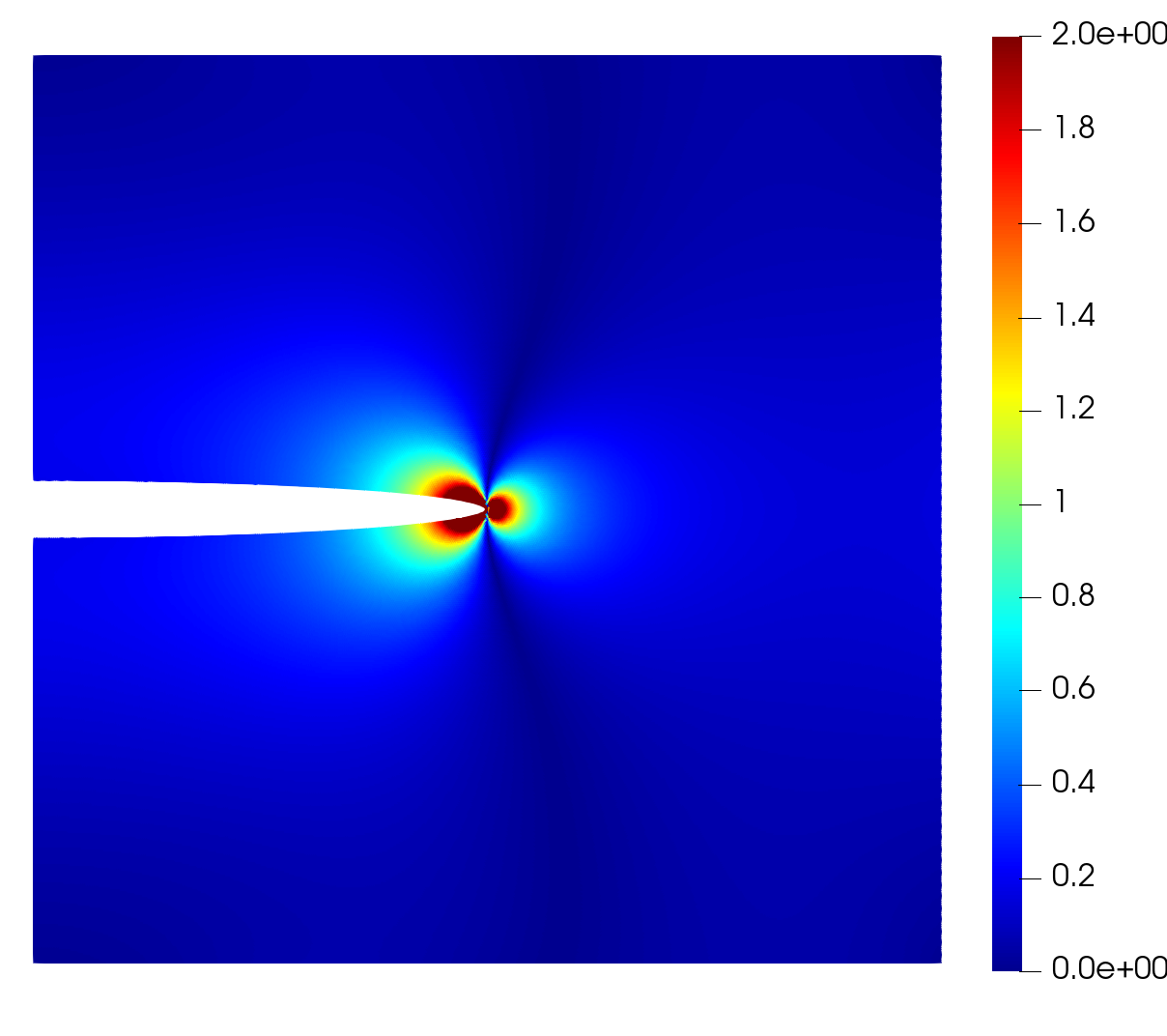} \\
		\hline
	\end{tabular}

	\caption{Exact von Mises stress error for the problems of an infinite plate with an elliptical hole under biaxial loading. The exact error is shown for a square lattice discretization of the domain (left) and for a triangular lattice discretization of the domain (right). We see that both node arrangements lead to the same pattern of the exact error.}
	\label{ExactErrorEllipticalHole}

	\centering
	\begin{tabular}{|M{4cm}|c|c|}
		\hline
		                                           & \textbf{Square lattice}                                        & \textbf{Triangular lattice}                                       \\
		\hline
		\multirow{2}{*}{Exact Error $\sigma_{vM}$} &                                                                &                                                                   \\[-1ex]
		                                           & \includegraphics[width=4cm]{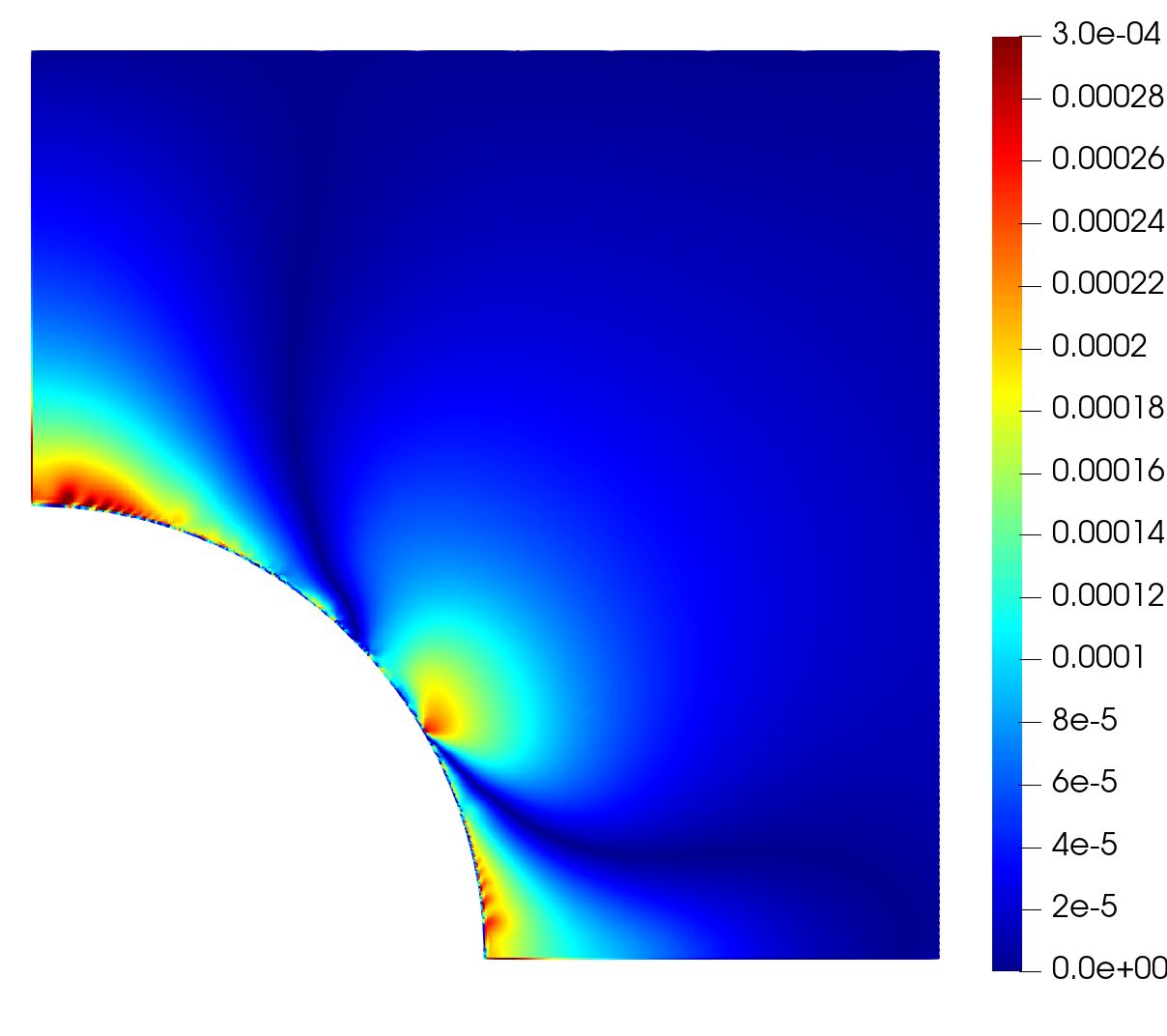} & \includegraphics[width=4cm]{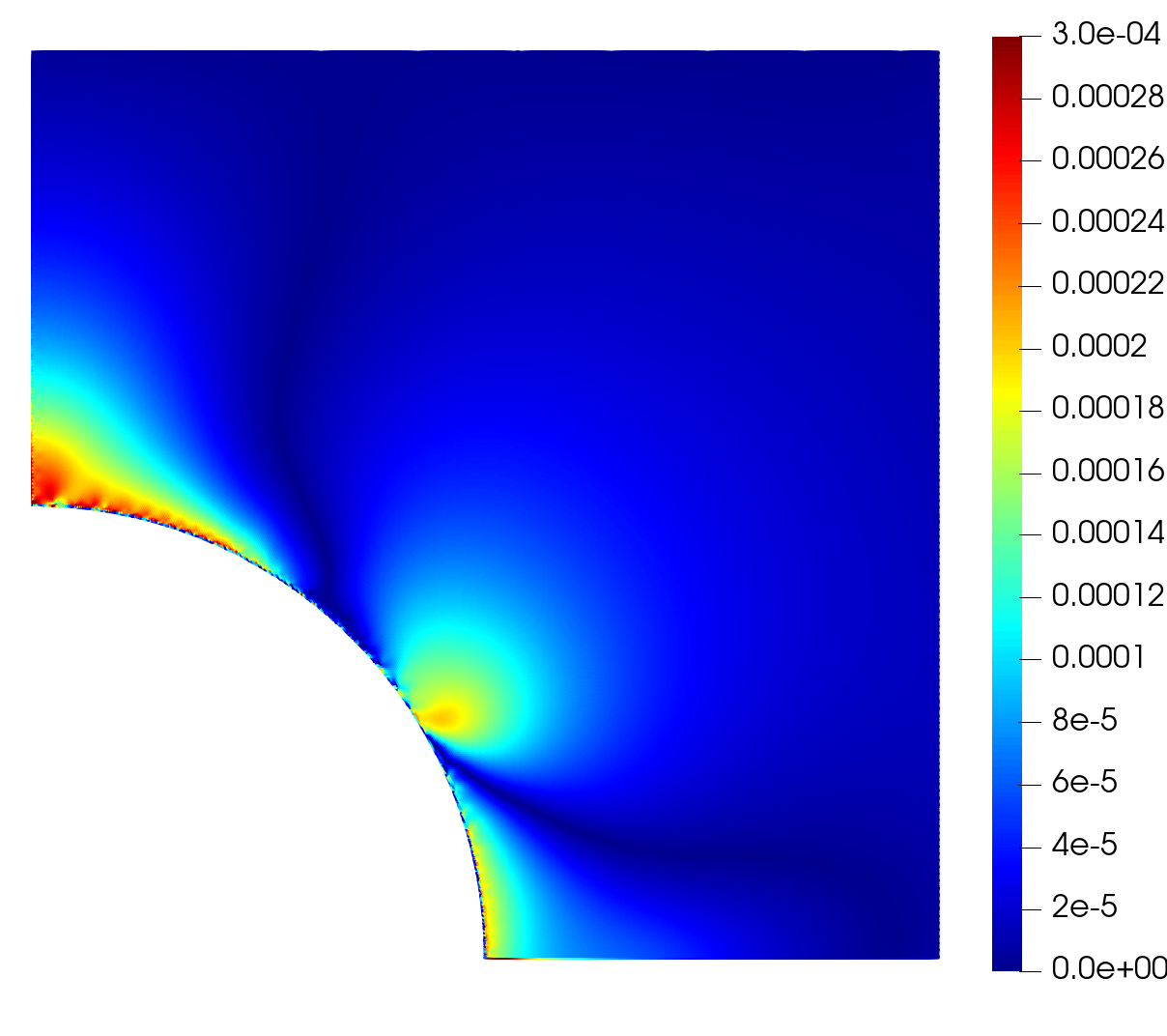} \\
		\hline
	\end{tabular}

	\caption{Exact von Mises stress error for the problems of an infinite body with a cylindrical hole under remote stress loading. The exact error is shown for a square lattice discretization of the domain (left) and for a triangular lattice discretization of the domain (right). We see that both node arrangements lead to the same pattern of the exact error.}
	\label{ExactErrorCylindricalHole}
\end{figure}

The exact error for the regular nodes arrangements presented in Section \ref{DiscretizationSec} (i.e. square lattice and triangular lattice for 2D problems) are presented in Figure \ref{ExactErrorEllipticalHole} for the problem of a plate with an elliptical hole and in Figure \ref{ExactErrorCylindricalHole} for the problem of an infinite body with a cylindrical hole.

We observe that, for both problems, the pattern of the exact error is the same for both discretization techniques. For the first problem, we see that the error is the highest close to the elliptical hole and in the region where the stress is the largest (see Figure \ref{2DProbemsSolution} (a))). For the second problem, two regions can be identified, both close to the hole. A region at the top of the hole which corresponds to the region where the stress is the largest and a region in the middle of the considered portion of the hole. This second region corresponds to a region where the stress solution in terms of von Mises stress is the lowest and changes rapidly along the hole (see Figure \ref{2DProbemsSolution} (b)).

\begin{figure}[h!]
	\centering
	\begin{tabular}{|M{4cm}|c|c|}
		\hline
		                                            & \textbf{Square lattice}                                              & \textbf{Triangular lattice}                                                    \\
		\hline
		\multirow{2}{*}{\begin{minipage}[t]{0.8\columnwidth} ZZ-type indicator \\ Unweighted - Indirect $\sigma_{vM}^s$ \end{minipage}} &                                                                      &                                                                                \\[-3ex]
		                                            & \includegraphics[width=4cm]{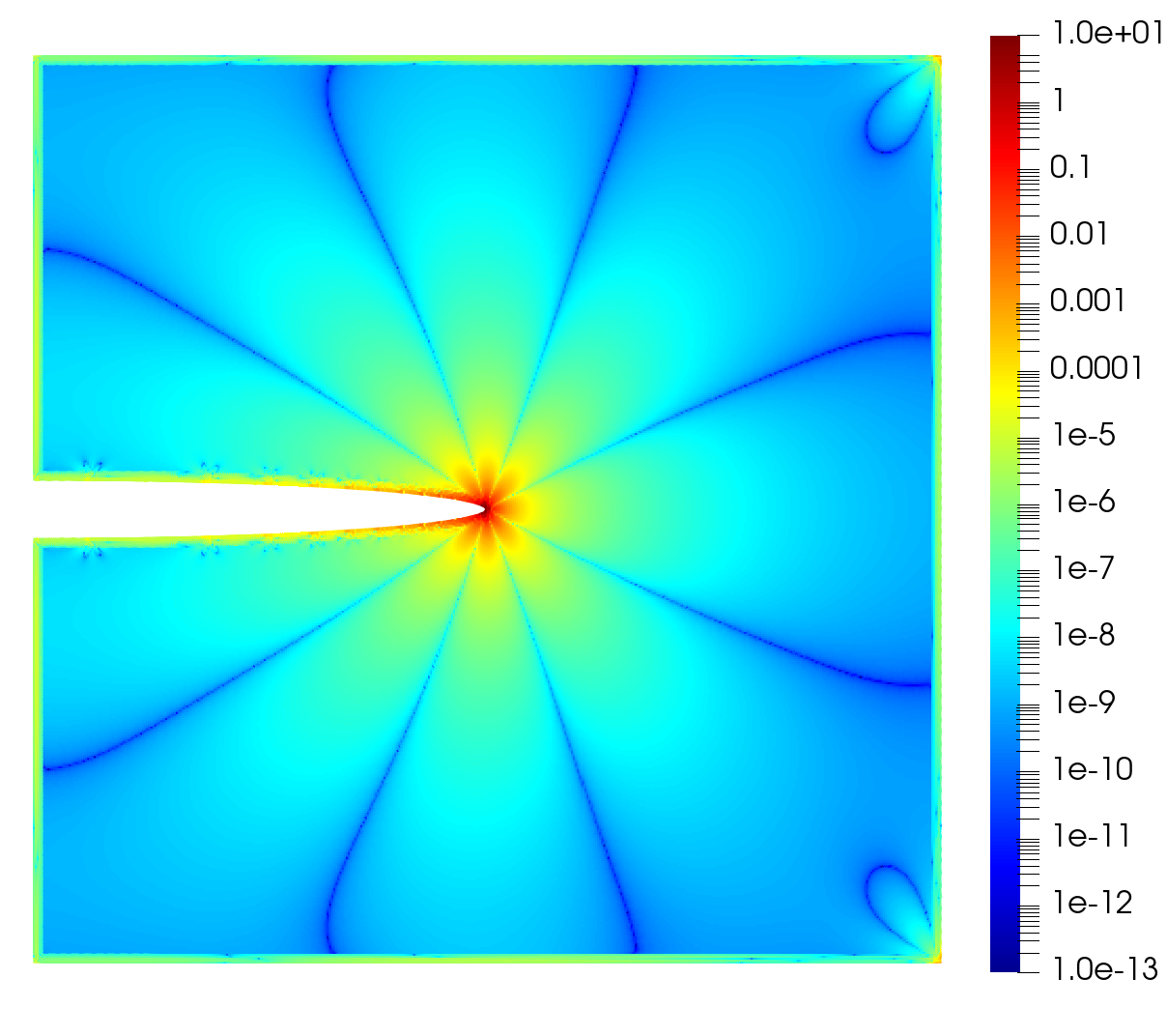} & \includegraphics[width=4cm]{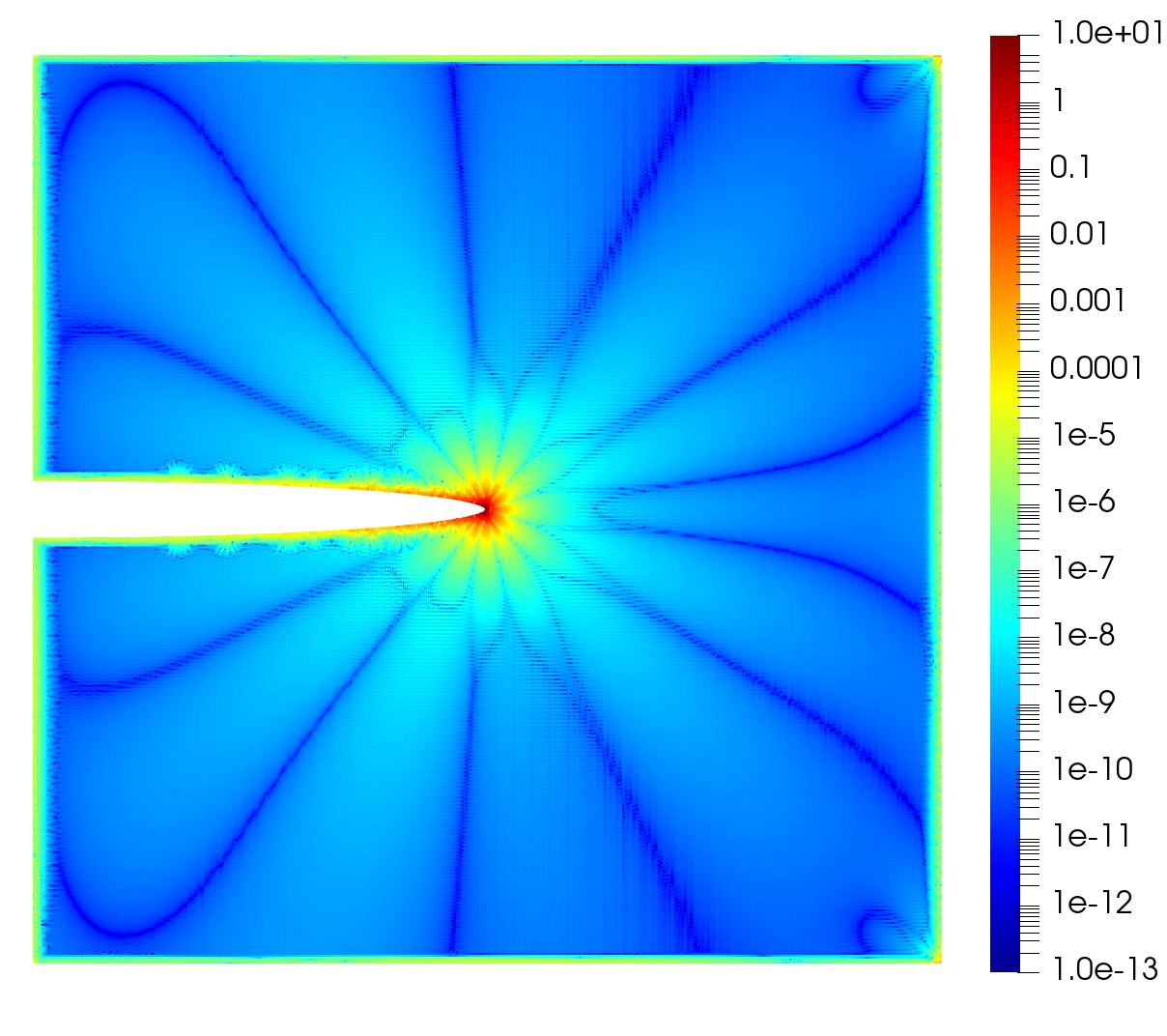} \\
		\hline
		\multirow{2}{*}{\begin{minipage}[t]{0.8\columnwidth} ZZ-type indicator \\ Weighted - Indirect $\sigma_{vM}^s$ \end{minipage}} &                                                                      &                                                                                \\[-3ex]
		                                            & \includegraphics[width=4cm]{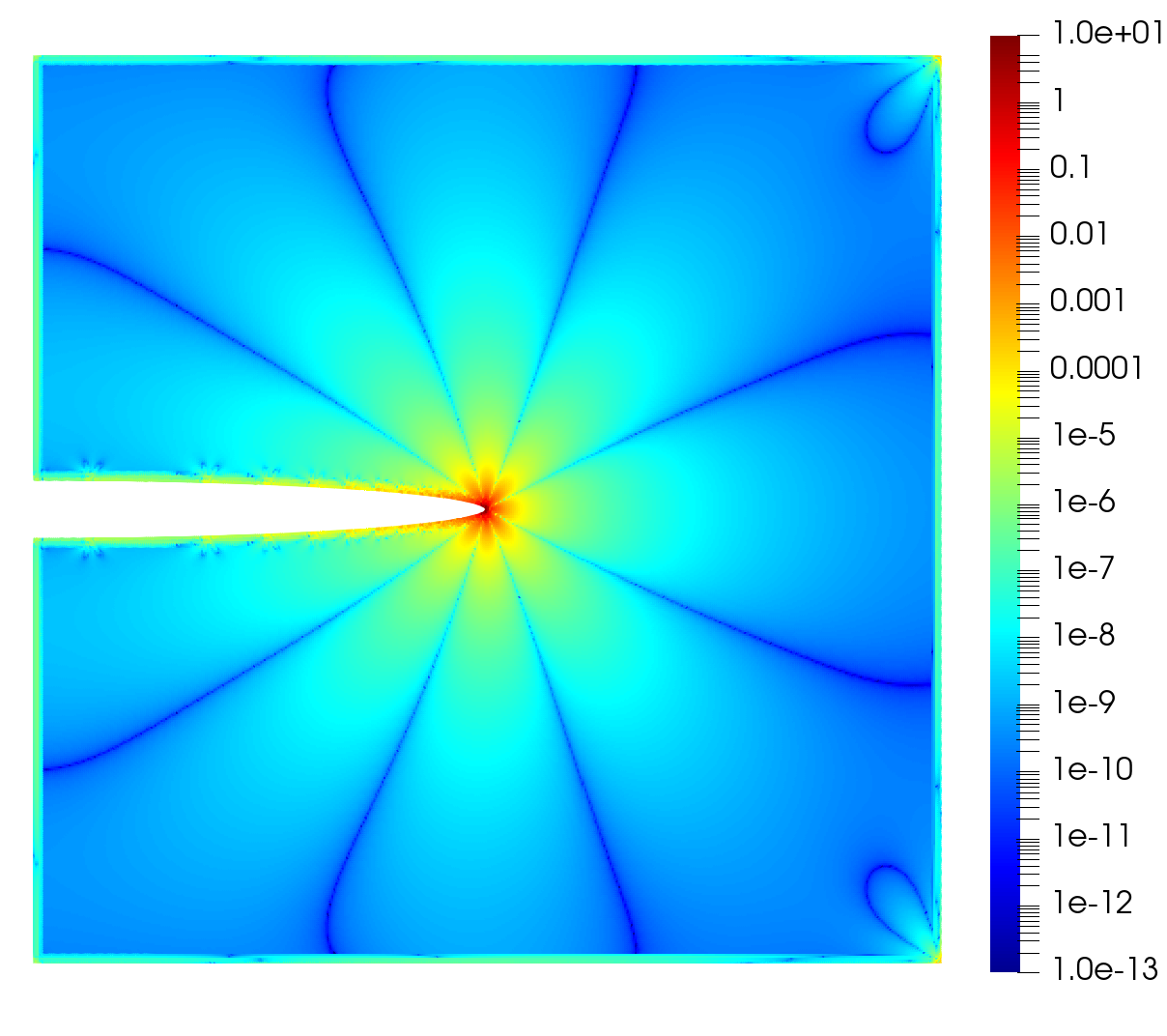}       & \includegraphics[width=4cm]{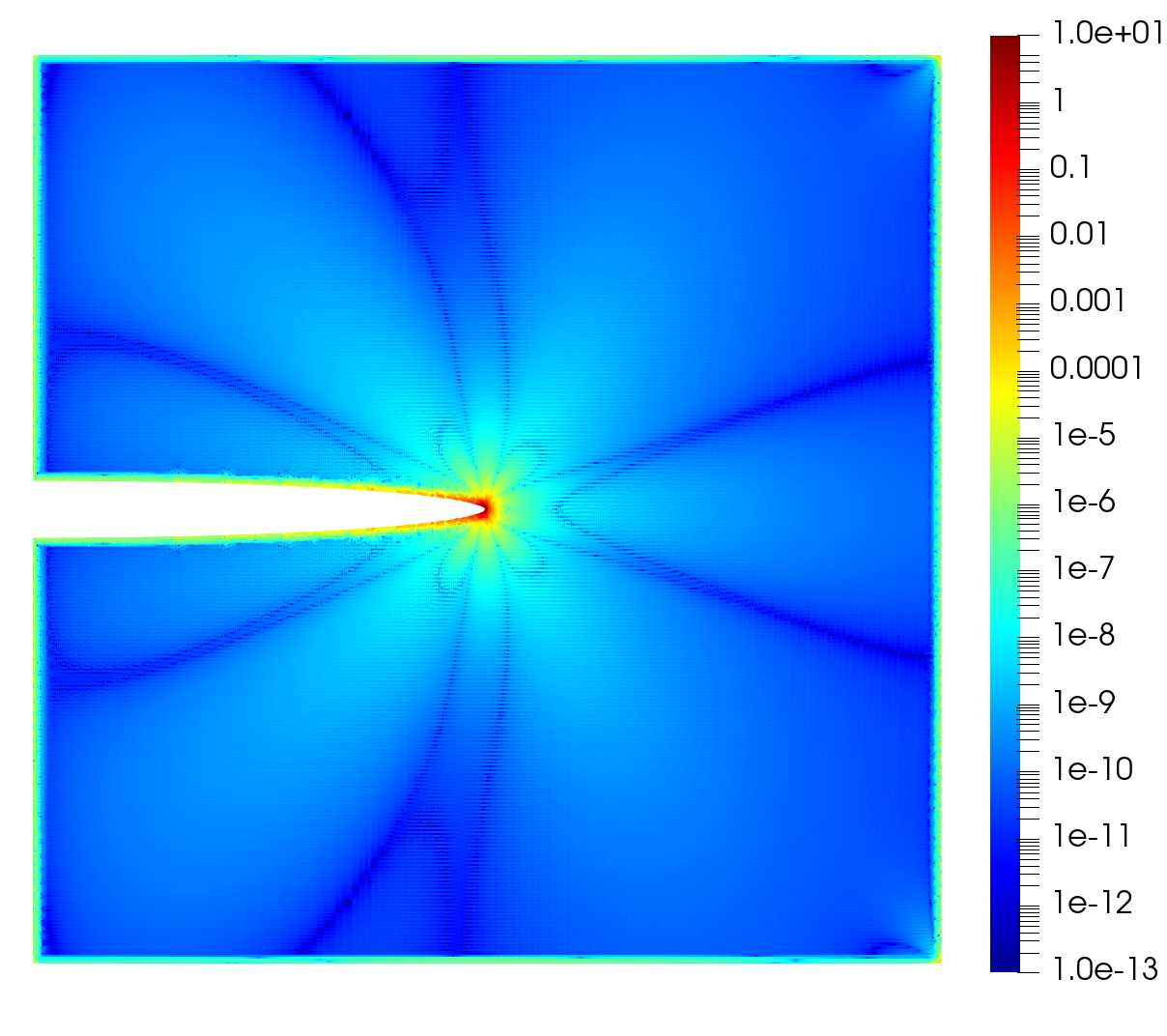}       \\
		\hline
		\multirow{2}{*}{\begin{minipage}[t]{0.8\columnwidth} ZZ-type indicator \\ Unweighted - Direct $\sigma_{vM}^s$ \end{minipage}} &                                                                      &                                                                                \\[-3ex]
		                                            & \includegraphics[width=4cm]{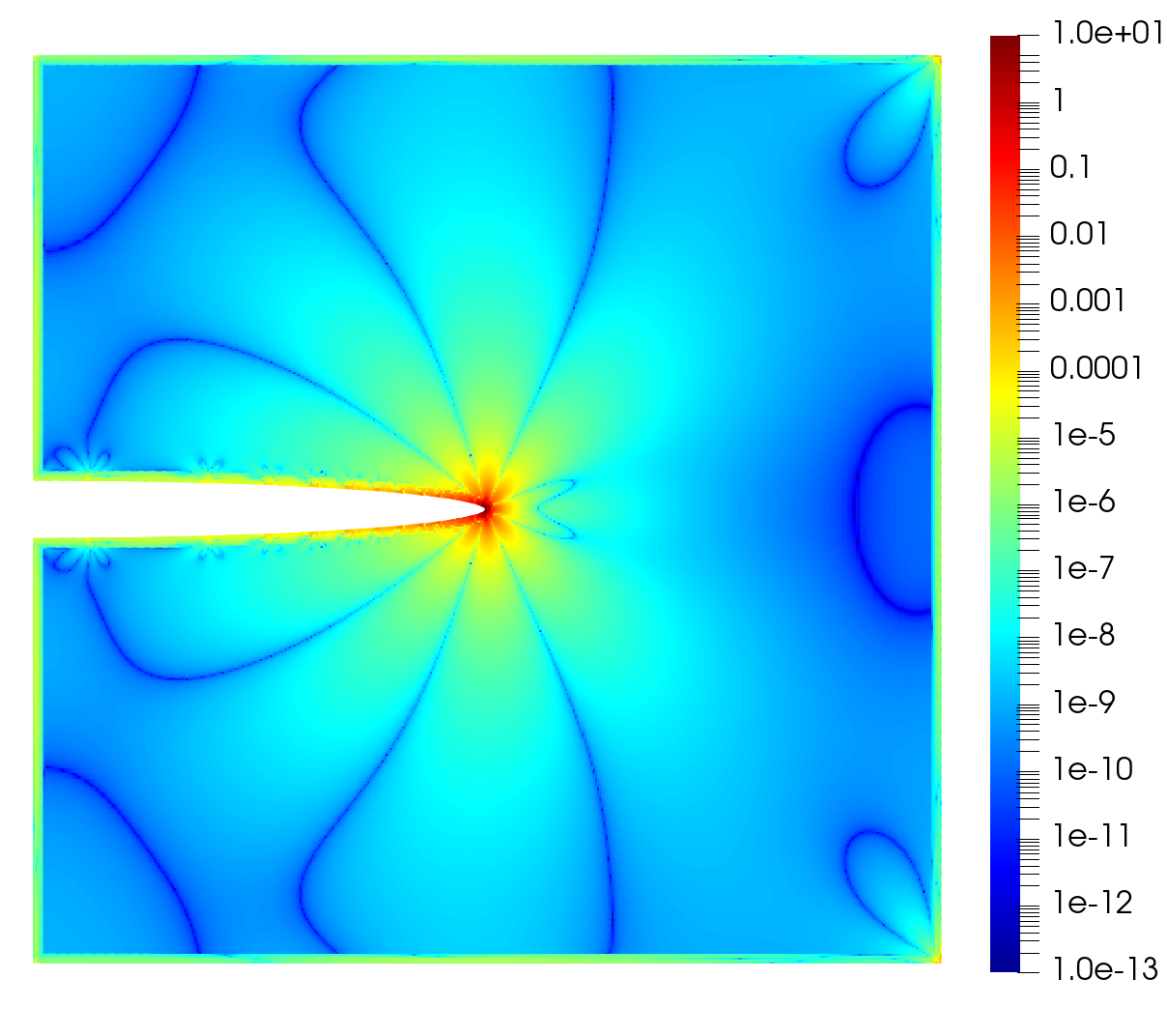}       & \includegraphics[width=4cm]{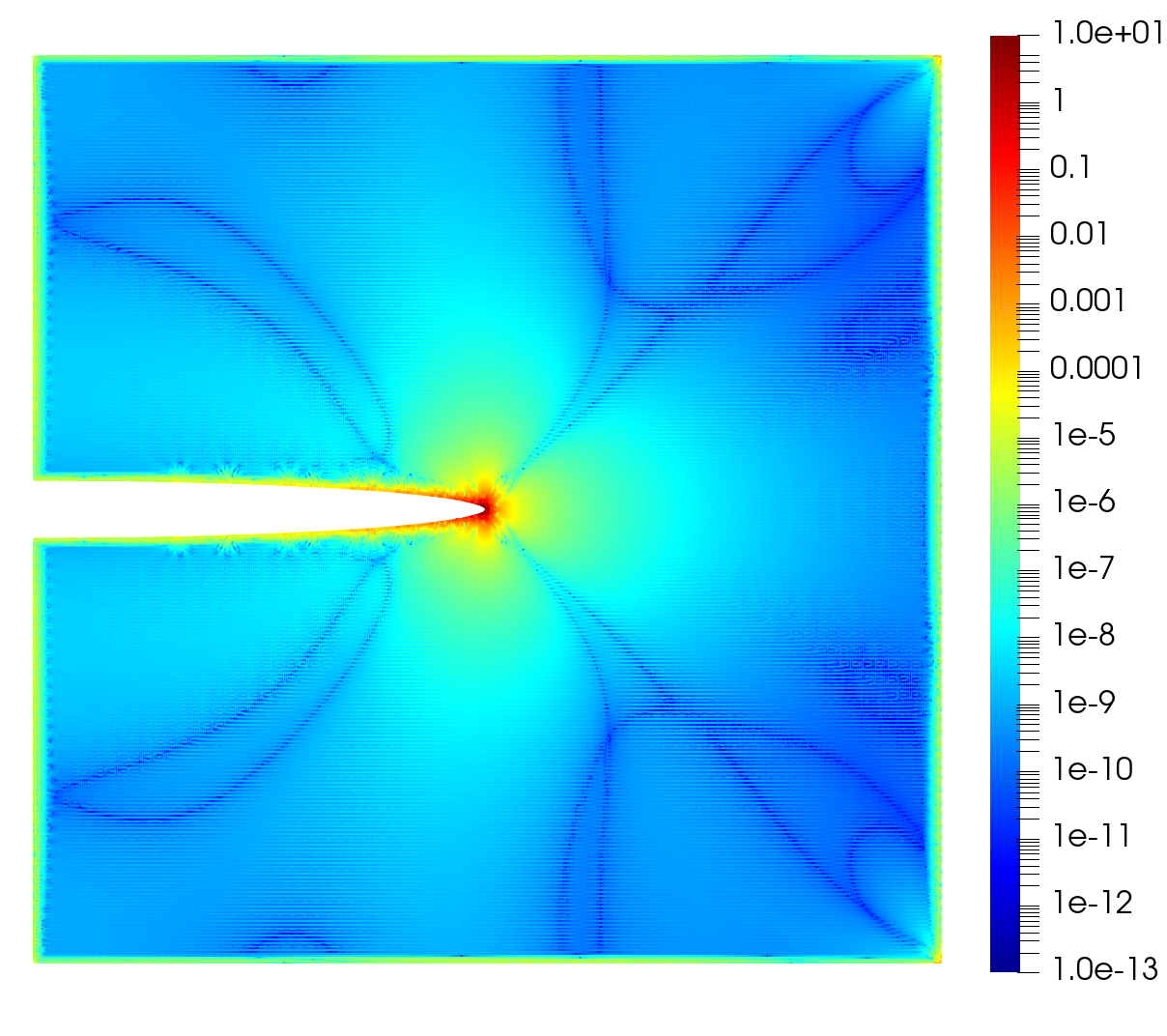}       \\
		\hline
		\multirow{2}{*}{\begin{minipage}[t]{0.8\columnwidth} ZZ-type indicator \\ Weighted - Direct $\sigma_{vM}^s$ \end{minipage}} &                                                                      &                                                                                \\[-3ex]
		                                            & \includegraphics[width=4cm]{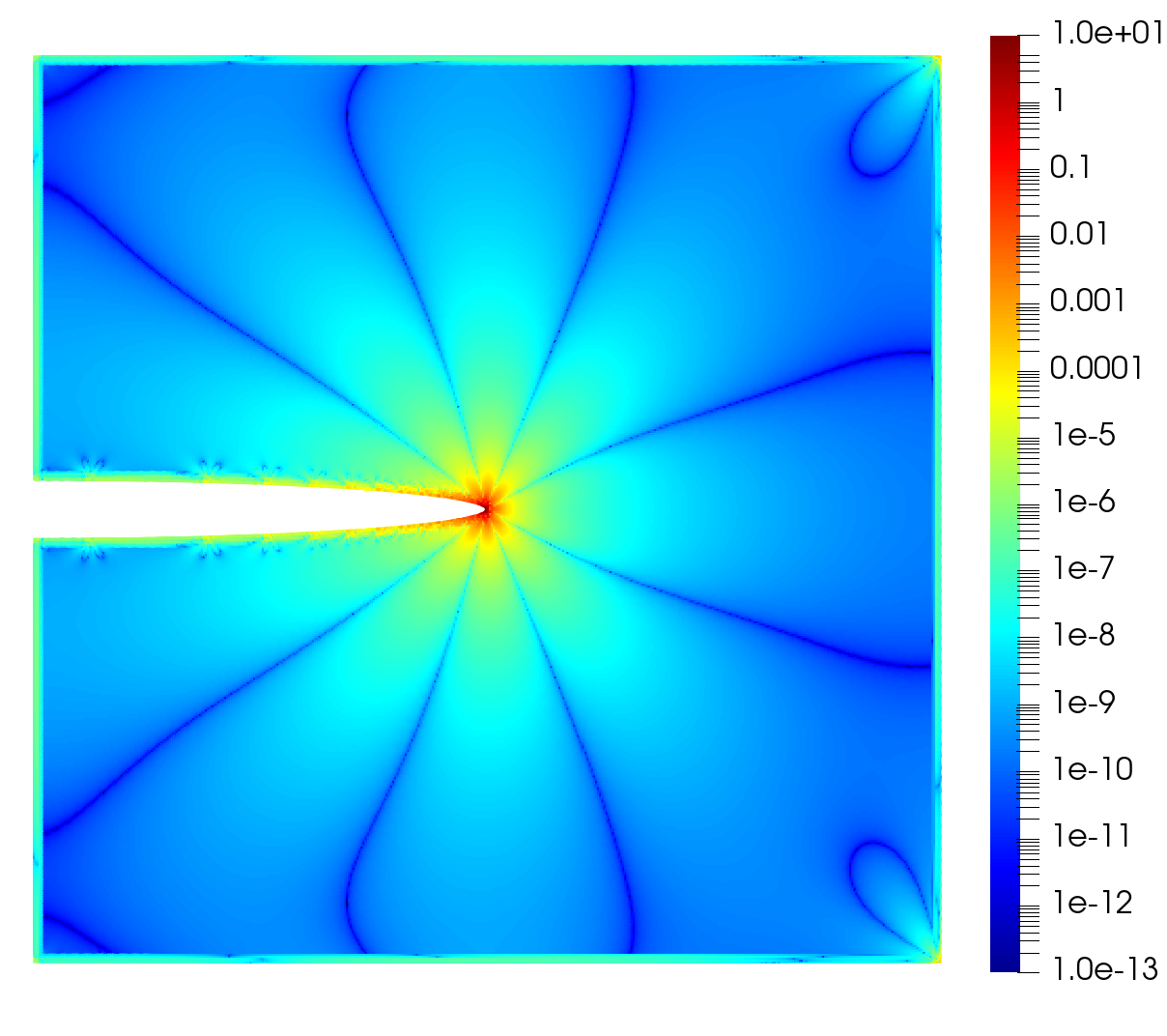}     & \includegraphics[width=4cm]{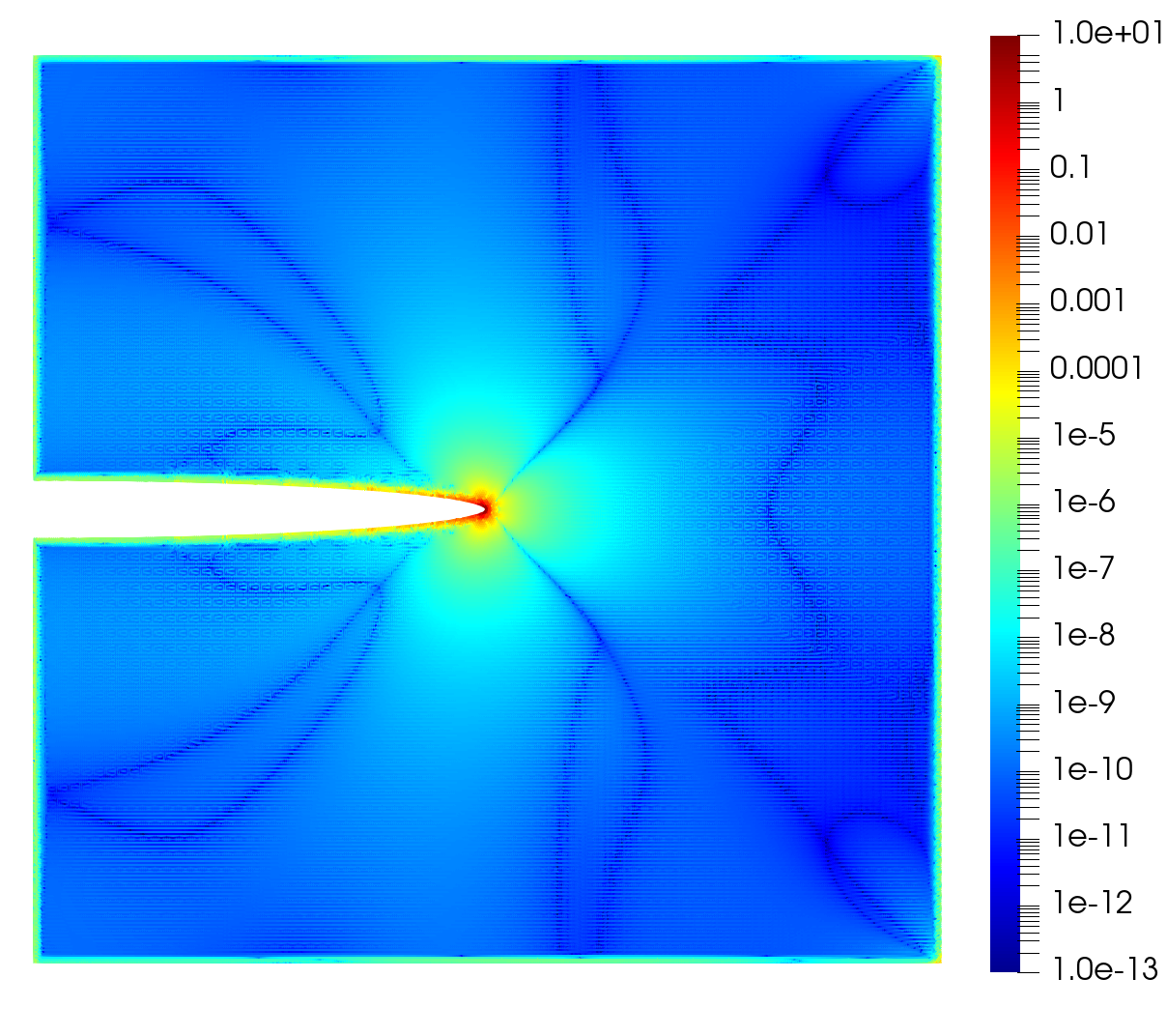}     \\
		\hline
		\multirow{2}{*}{\begin{minipage}[t]{0.8\columnwidth} Residual-type indicator \end{minipage}} &                                                                      &                                                                                \\[-3ex]
		                                            & \includegraphics[width=4cm]{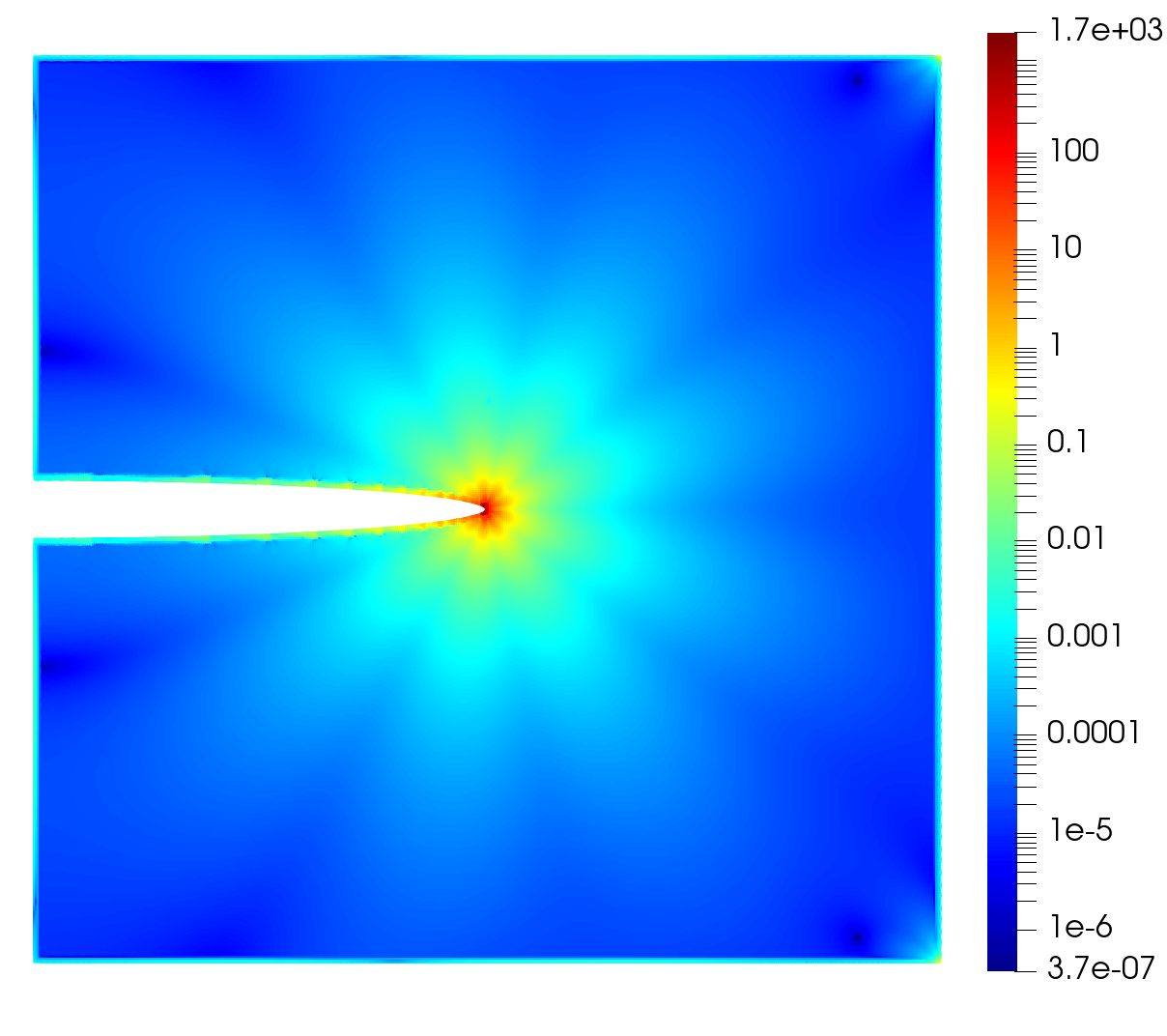}       & \includegraphics[width=4cm]{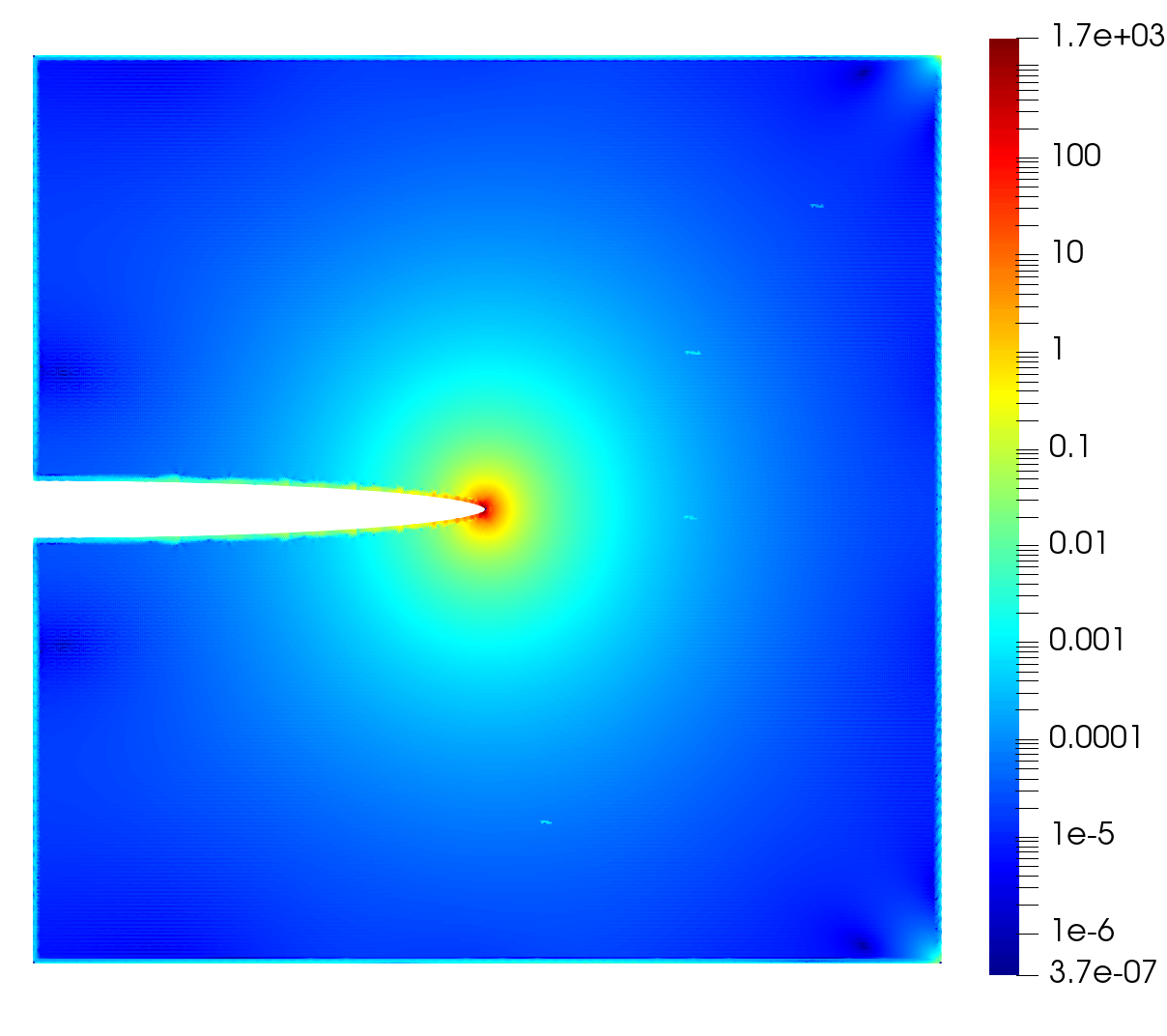}                 \\
		\hline
	\end{tabular}

	\caption{Comparison of the error pattern for ZZ-type error indicators computed with various parameters and for a residual-type error indicator for the problem of a plate with an elliptical hole. The results are shown for square and triangular lattice discretizations of the interior of the domain.}
	\label{EstimErrorPatternEllipticalHole}
\end{figure}

\begin{figure}[h!]
	\centering
	\begin{tabular}{|M{4cm}|c|c|}
		\hline
		                                            & \textbf{Square lattice}                                                & \textbf{Triangular lattice}                                               \\
		\hline
		\multirow{2}{*}{\begin{minipage}[t]{0.8\columnwidth} ZZ-type indicator \\ Unweighted - Indirect $\sigma_{vM}^s$ \end{minipage}} &                                                                        &                                                                           \\[-3ex]
		                                            & \includegraphics[width=4cm]{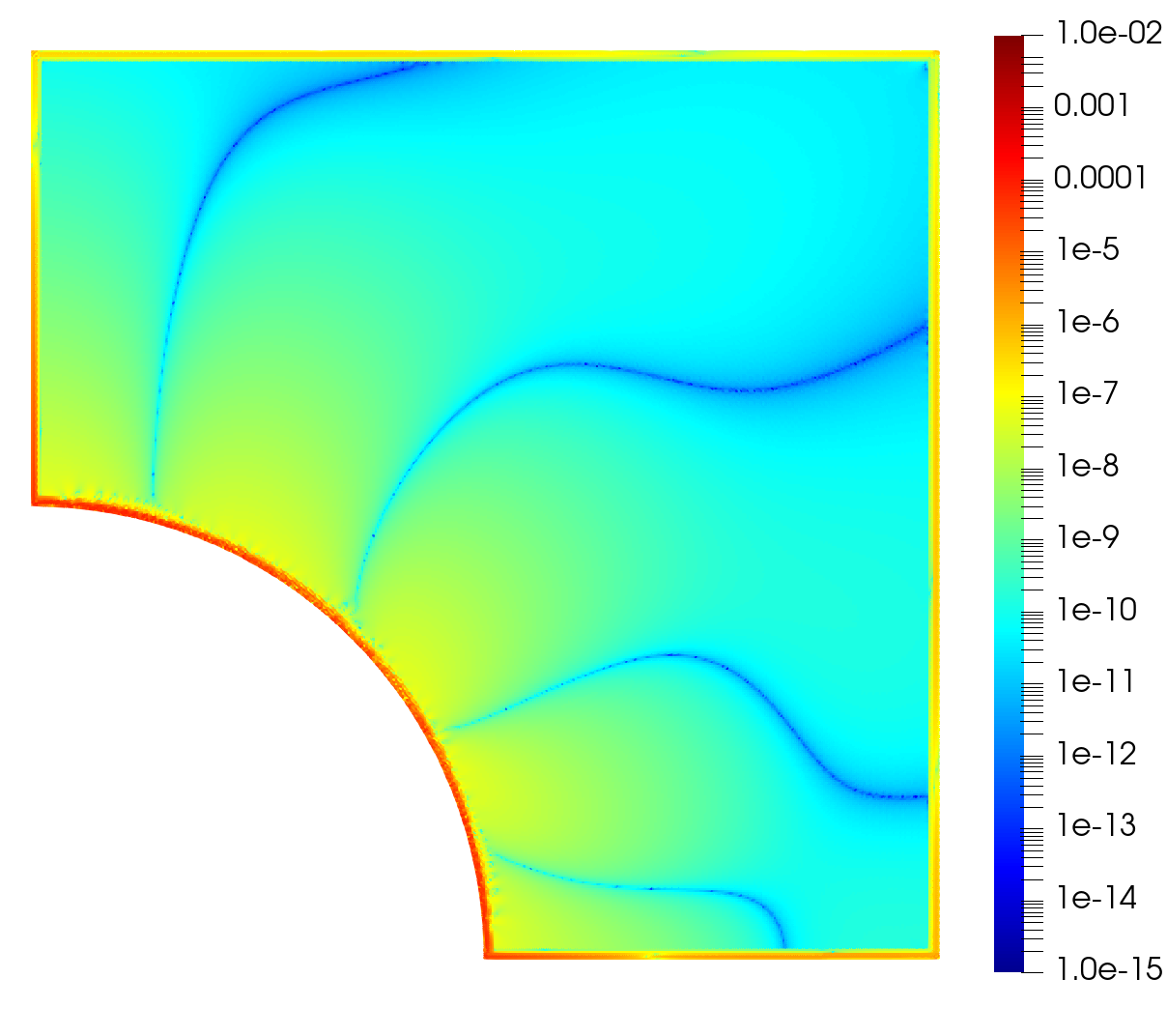} & \includegraphics[width=4cm]{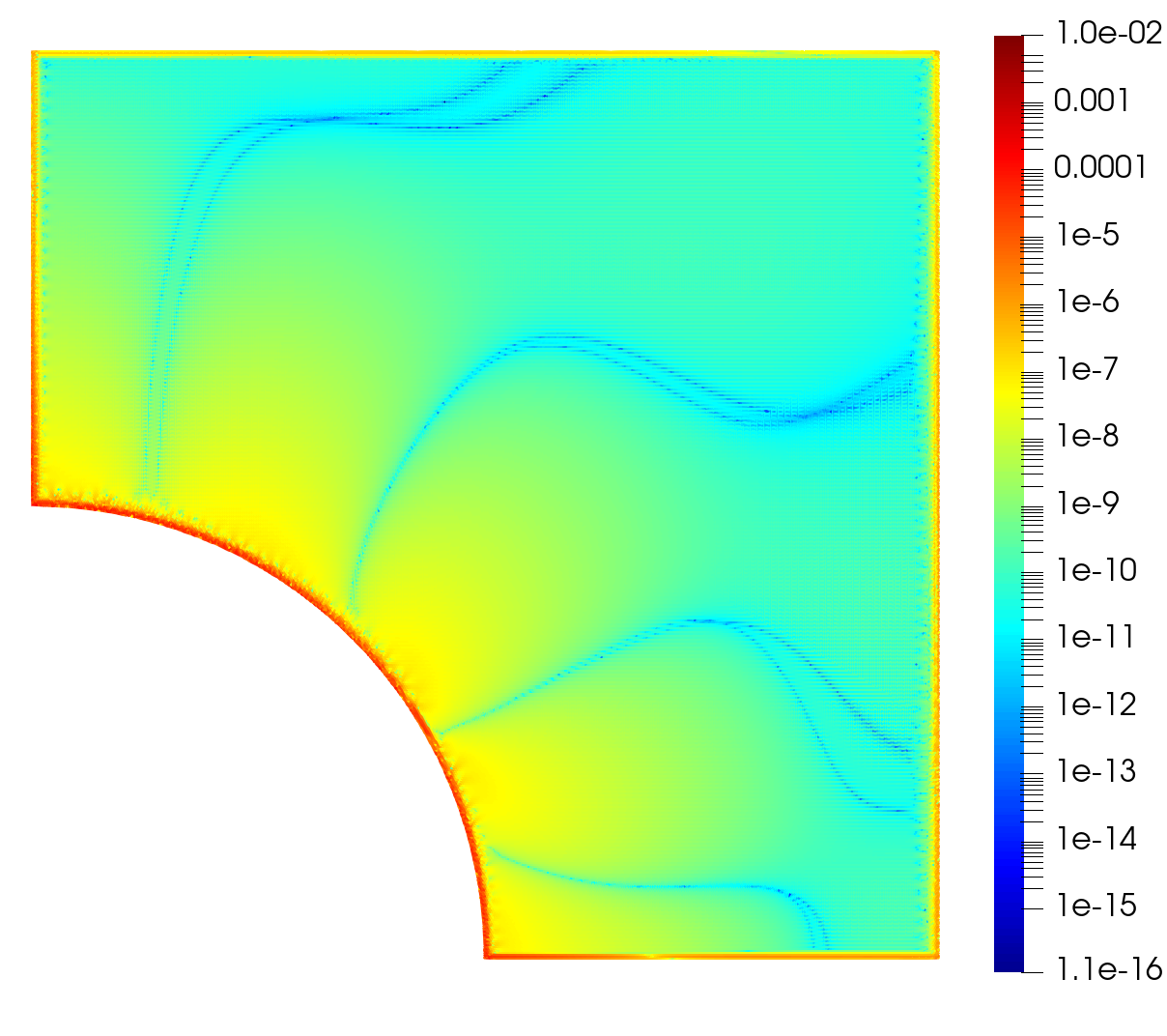} \\
		\hline
		\multirow{2}{*}{\begin{minipage}[t]{0.8\columnwidth} ZZ-type indicator \\ Weighted - Indirect $\sigma_{vM}^s$ \end{minipage}} &                                                                        &                                                                           \\[-3ex]
		                                            & \includegraphics[width=4cm]{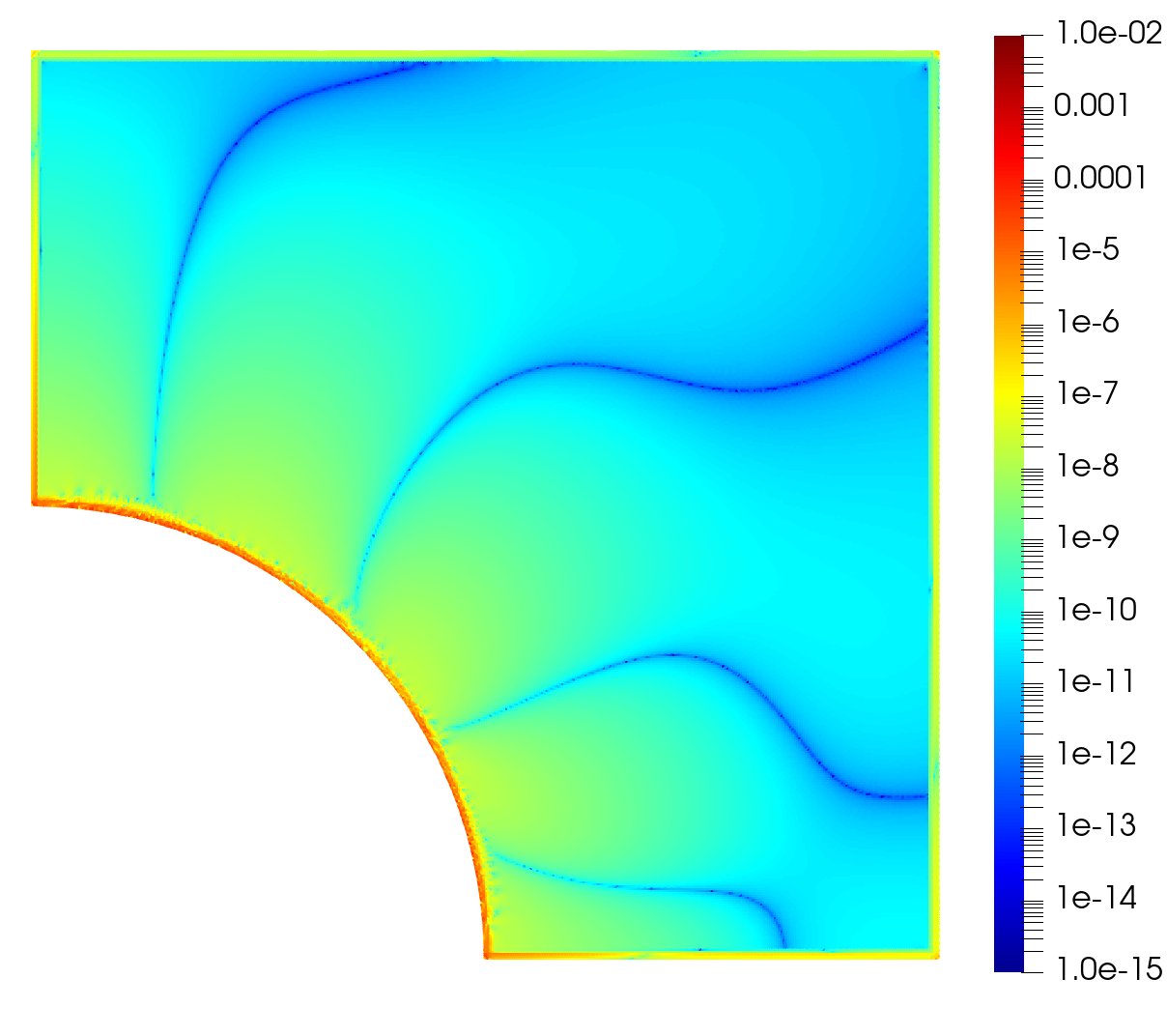}   & \includegraphics[width=4cm]{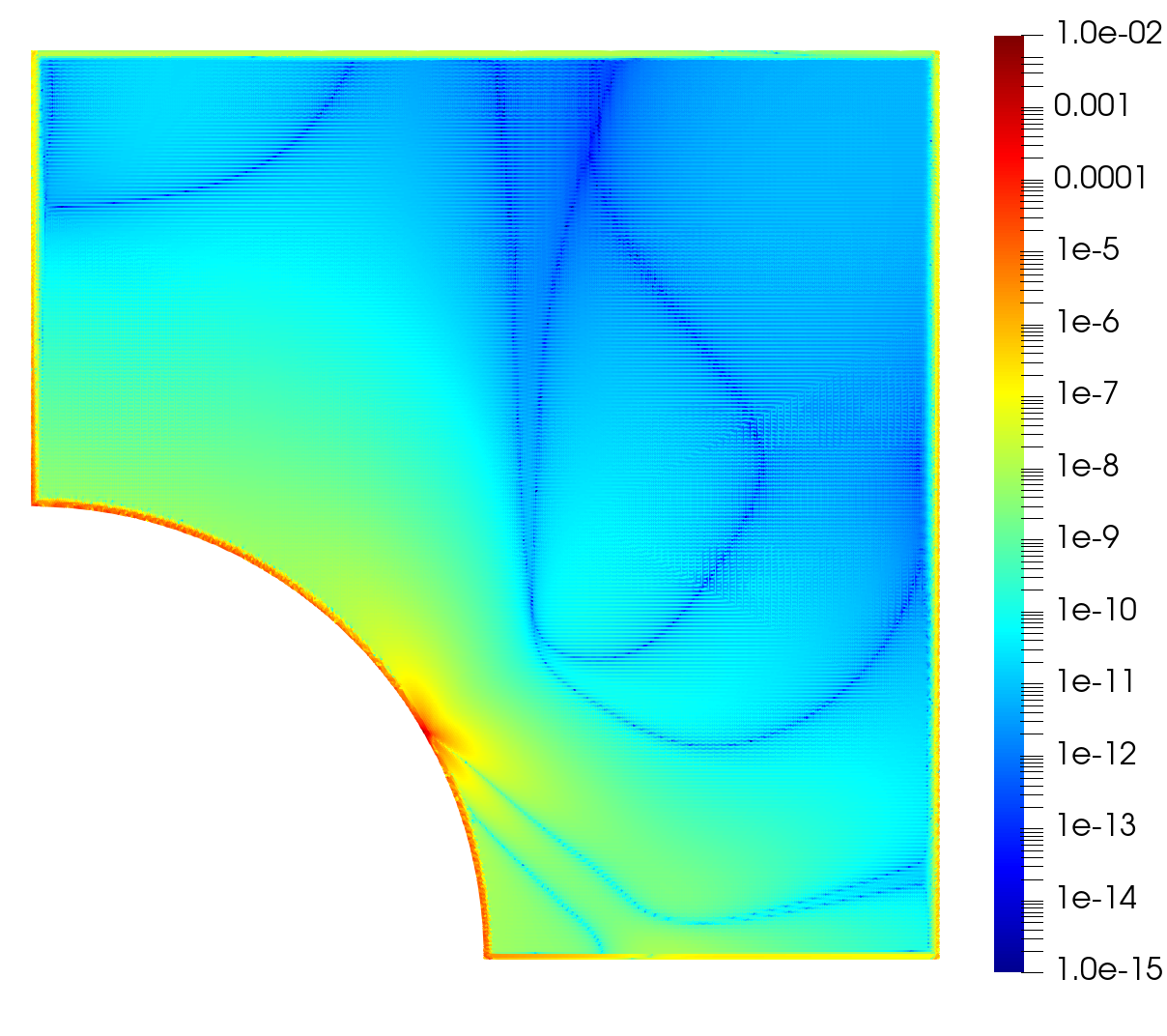}   \\
		\hline
		\multirow{2}{*}{\begin{minipage}[t]{0.8\columnwidth} ZZ-type indicator \\ Unweighted - Direct $\sigma_{vM}^s$ \end{minipage}} &                                                                        &                                                                           \\[-3ex]
		                                            & \includegraphics[width=4cm]{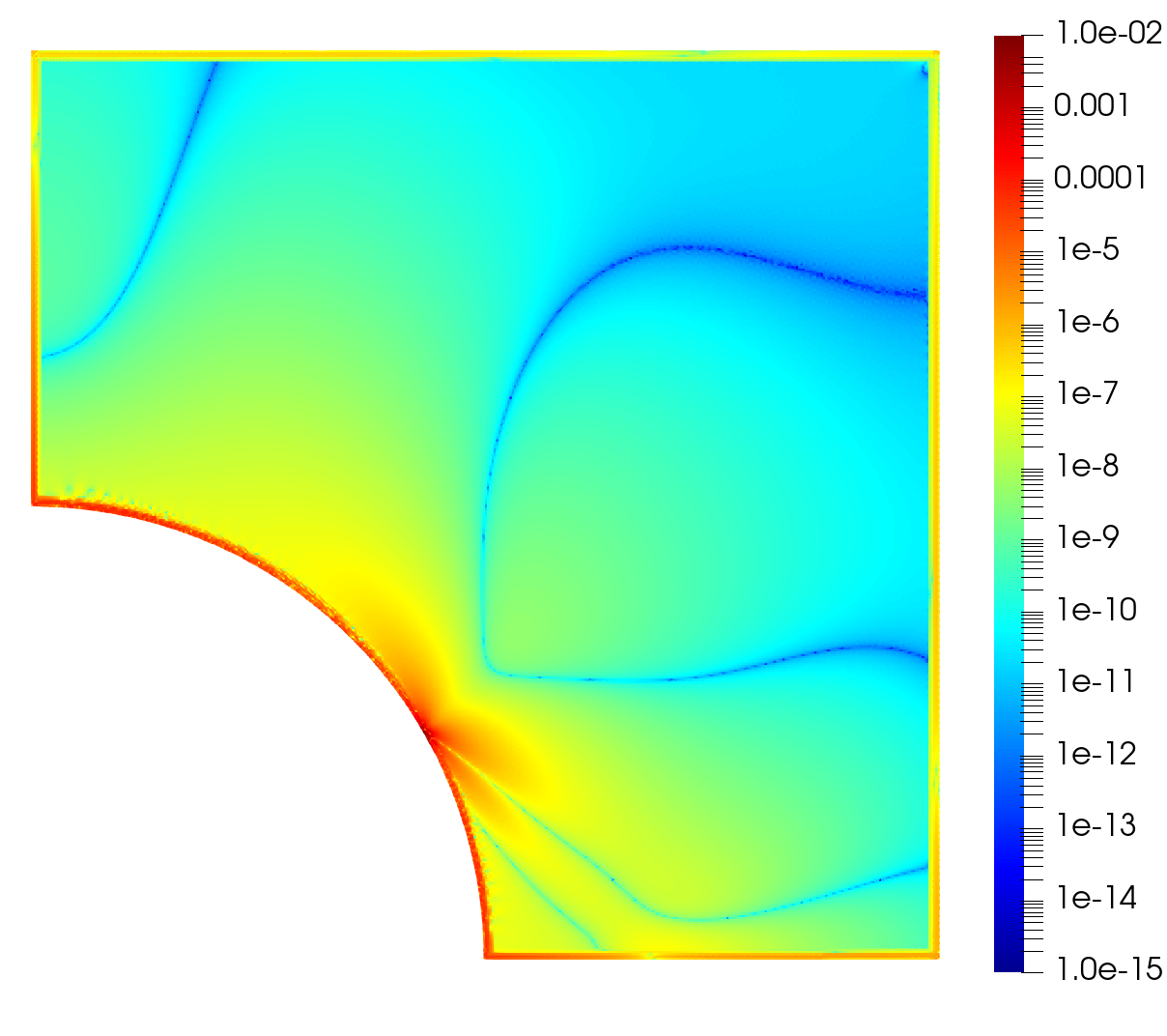}   & \includegraphics[width=4cm]{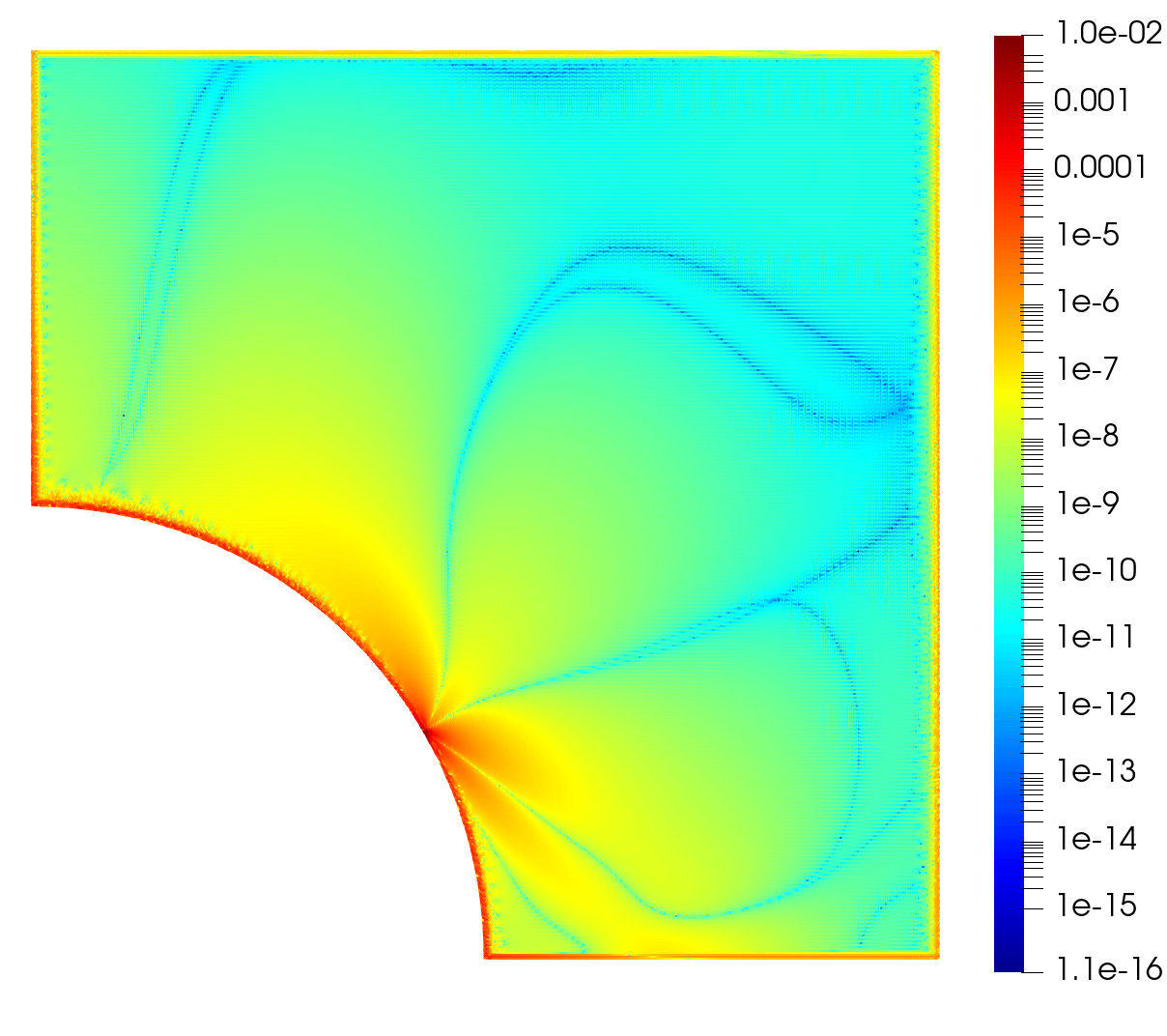}   \\
		\hline
		\multirow{2}{*}{\begin{minipage}[t]{0.8\columnwidth} ZZ-type indicator \\ Weighted - Direct $\sigma_{vM}^s$ \end{minipage}} &                                                                        &                                                                           \\[-3ex]
		                                            & \includegraphics[width=4cm]{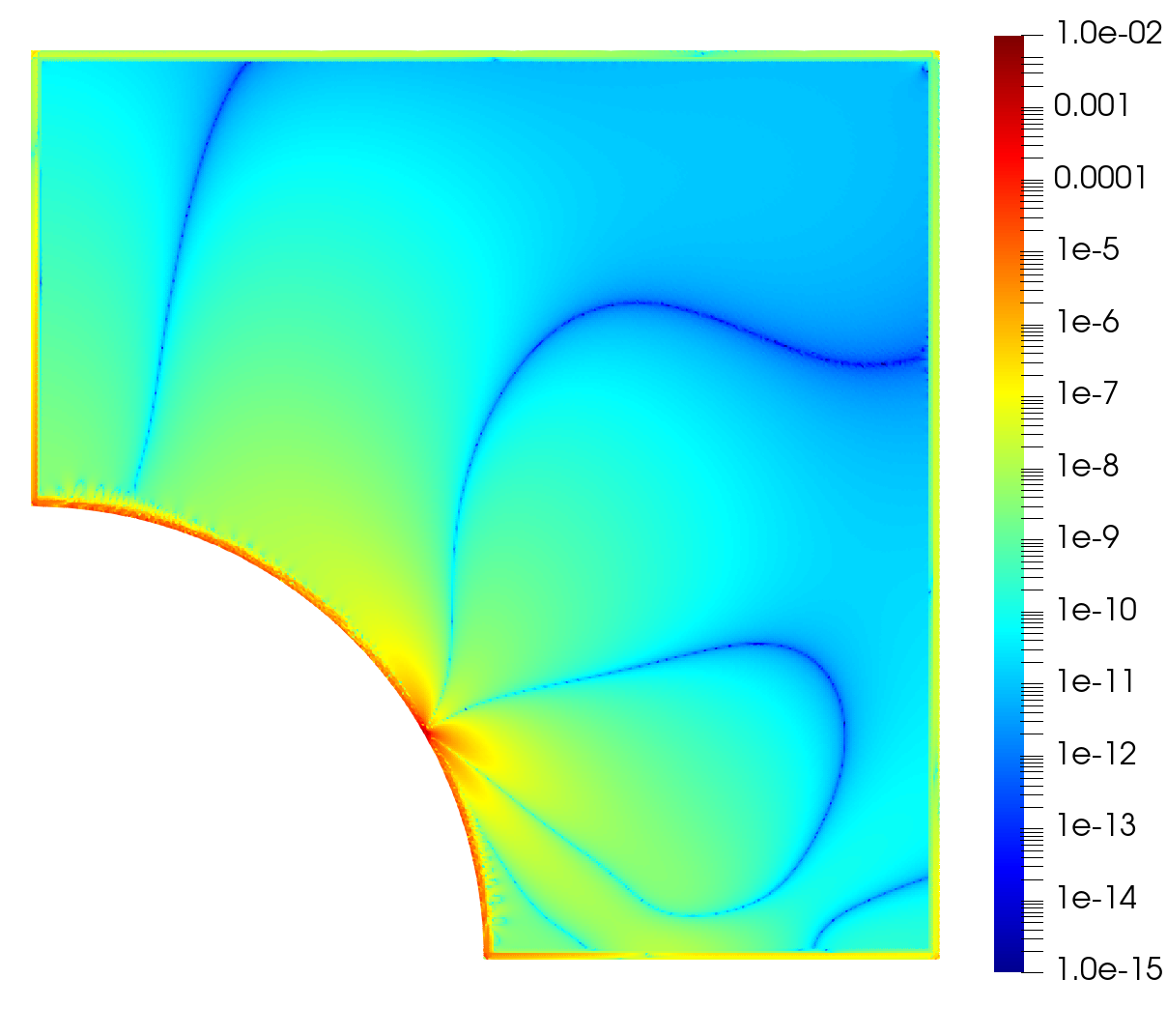}     & \includegraphics[width=4cm]{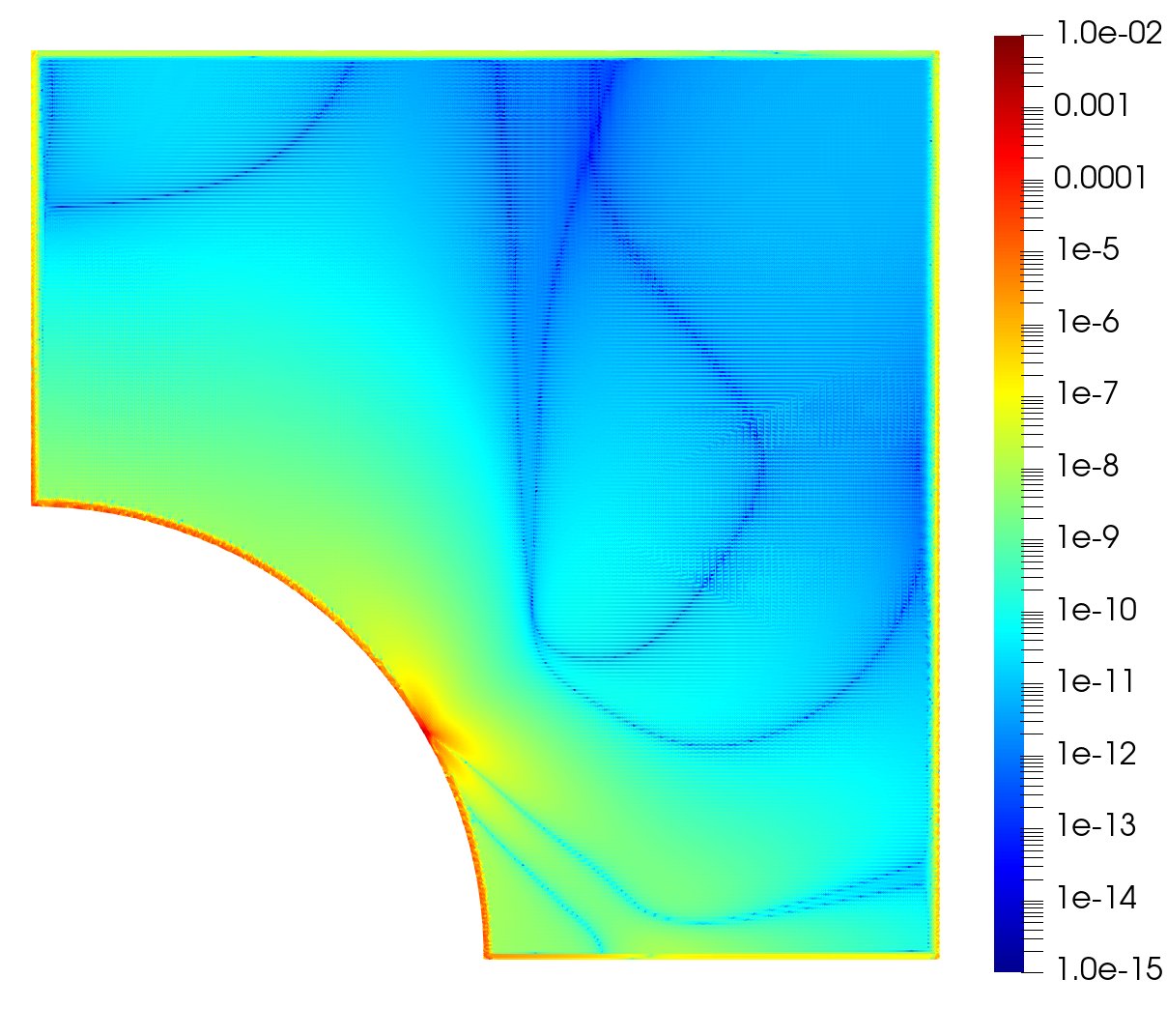}     \\
		\hline
		\multirow{2}{*}{\begin{minipage}[t]{0.8\columnwidth} Residual-type indicator \end{minipage}}    &                                                                        &                                                                           \\[-3ex]
		                                            & \includegraphics[width=4cm]{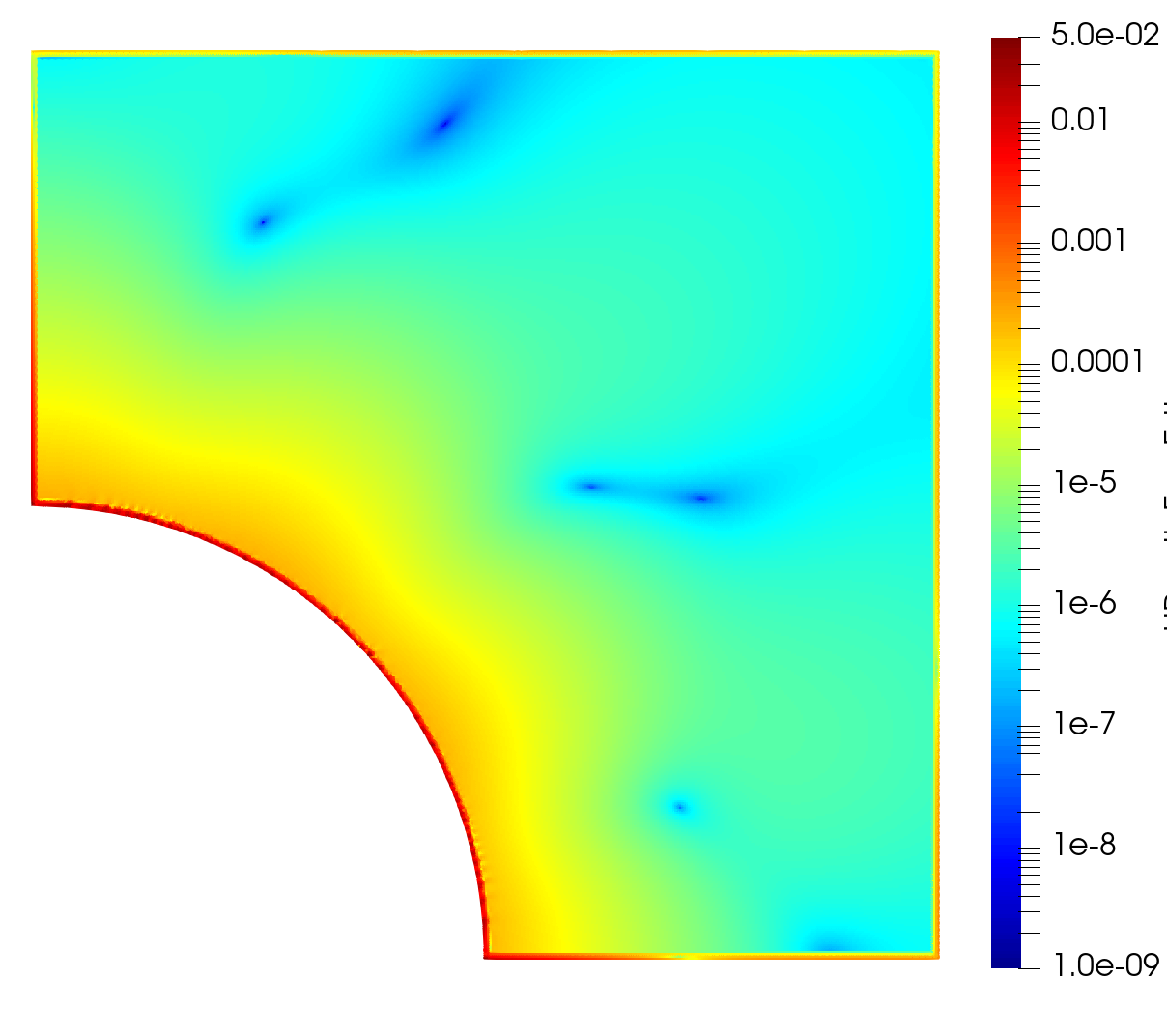}           & \includegraphics[width=4cm]{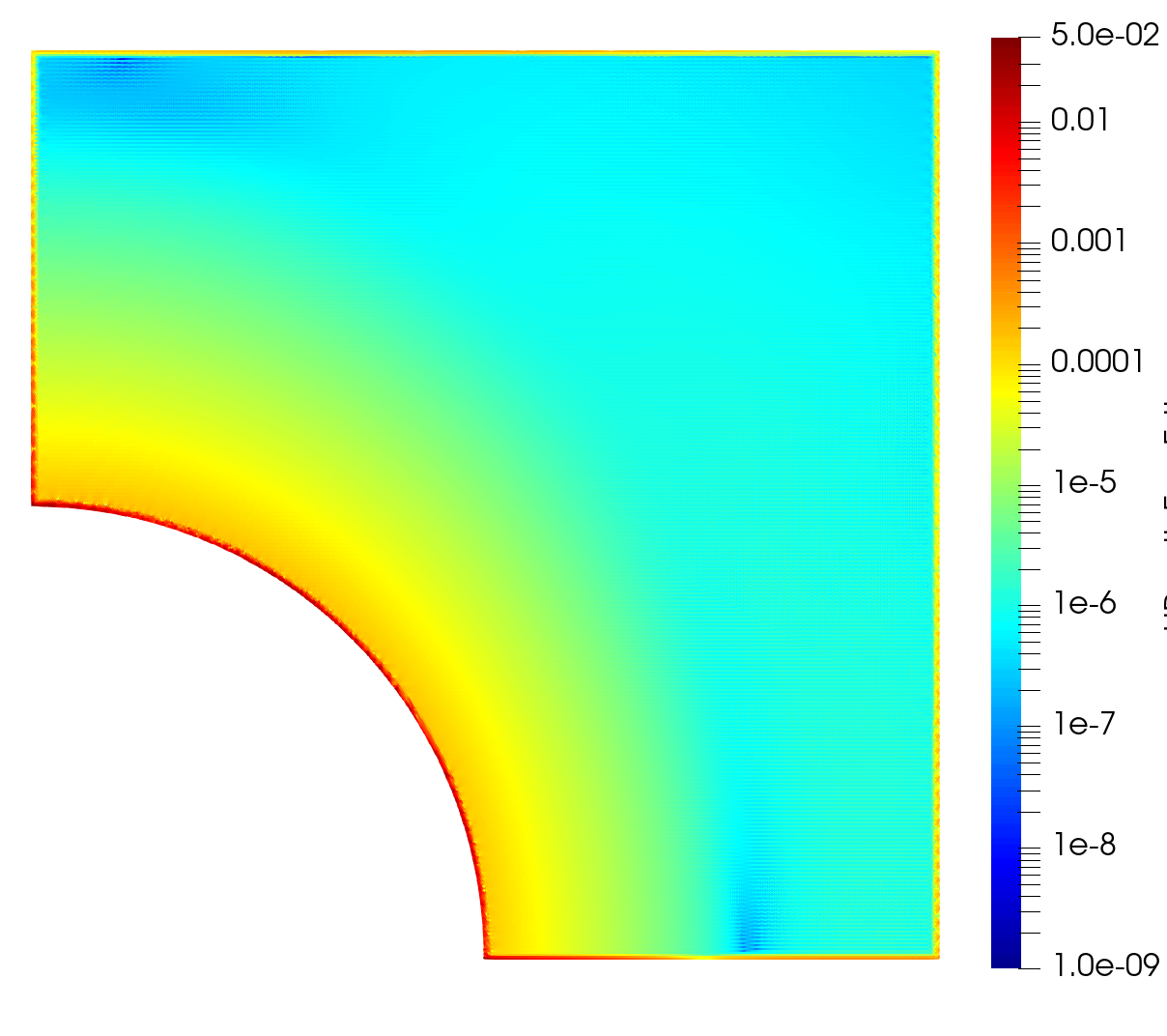}           \\
		\hline
	\end{tabular}

	\caption{Comparison of the error pattern for ZZ-type error indicators computed with various parameters and for a residual-type error indicator for the problem of a body with a cylindrical hole. The results are shown for square and triangular lattice discretizations of the interior of the domain.}
	\label{EstimErrorPatternCylindricalHole}
\end{figure}

We compare in Figure \ref{EstimErrorPatternEllipticalHole} and Figure \ref{EstimErrorPatternCylindricalHole} the error patterns for various parameters of the ZZ-type indicator and for the residual error indicator, respectively, for the two benchmark problems. We compare the results obtained for an ``unweighted'' indicator (i.e. $w\mathbf{\left(\mathbf{X},\mathbf{X_c}\right)}=1$ in Equation (\ref{ZZ_FunctionalB})) and for a ``weighted'' indicator. We selected a 4$^{\text{th}}$ order spline for the radial basis function $w$ in Equation (\ref{ZZ_FunctionalB}). The equation of the spline is presented in Equation \ref{Spline4Weight}.

\begin{equation}\label{Spline4Weight}
	w(s)=
	\begin{cases}
		1 - 6 s^2 +	8 s^3 - 3 s^4 & \text{ \quad if } s \leq 1 \\
		0                        & \text{ \quad if } s > 1.   \\
	\end{cases}\\
\end{equation}

\begin{samepage}

	\noindent
	For a node $\mathbf{X}$ within a support of radius $R_c$ of a collocation node $\mathbf{X_c}$, the weight function $w$ based on the 4$^{\text{th}}$ order spline is:

	\begin{equation}\label{Spline4WeightFunction}
		w\mathbf{\left(\mathbf{X},\mathbf{X_c}\right)}=w(s) \quad \text{with} \quad s=\frac{\lVert \mathbf{X}-\mathbf{X_c} \rVert_2}{R_c}.
	\end{equation}
\end{samepage}

In Figures \ref{EstimErrorPatternEllipticalHole} and \ref{EstimErrorPatternCylindricalHole} we also compare the results obtained from a ``indirect'' computation of the indicator from the von Mises stress field to results obtained from an ``direct'' computation of the indicator. The ``indirect'' computation of the indicator uses the computation of the individual components of the stress tensor $\sigma^s_{11}\left(\mathbf{X_c}\right)$, $\sigma^s_{12}\left(\mathbf{X_c}\right)$ and $\sigma^s_{22}\left(\mathbf{X_c}\right)$ to compute $\sigma^s_{vM}\left(\mathbf{X_c}\right)$. The ``direct'' computation of the indicator uses the computation of the von Mises stress components $\sigma^s_{vM}\left(\mathbf{X_{vi}}\right)$ at each support nodes $\mathbf{X_{vi}}$ of $\mathbf{X_c}$ to compute $\sigma^s_{vM}\left(\mathbf{X_c}\right)$. The figures show the result for both the unweighted and weighted configurations, and both the direct and indirect computation methods. The results for each of these configurations are shown for both the square and triangular lattice discretizations. We presented the error indicator using a logarithmic color scale as it allows a better identification of the different error zones of the solution. The amplitude of the scale has been set constant for the ZZ-type error indicator and for the residual-type error indicator to facilitate the analysis of the results.

The results obtained for the plate with an elliptical hole show that the discretization method selected for the interior of the domain impacts significantly the pattern of the error indicator. For both discretization techniques, we can observe lines where the computed error indicator is lower. This phenomenon is the most significant for the square lattice discretization techniques, for the case of a ZZ-type error indicator computed using an ``indirect'' computation of $\sigma_{vM}^s$. Such lines are not observed for the residual-type error indicator computed for the domain discretized using a triangular lattice. We also see that the presence of weights in the computation of the ZZ-type error indicator has little impact on the results. The ``direct'' computation of $\sigma_{vM}^s$ leads to better results than an ``indirect'' computation of $\sigma_{vM}^s$ in the sense that the error indicator appears less affected by the discretization strategy.

The trend of the results obtained for the problem of a body with a cylindrical hole is the same as the trend of the results obtained for the problem of a plate with an elliptical hole. The discretization technique impacts the pattern of the error. The error pattern is the most uniform for the residual-type error indicator. It can be noticed however that the ZZ-type error indicator allows for the identification of a zone, close to the cavity, where the error is greater. This zone is not identified with the residual-type error indicator. The exact error presented in Figure \ref{ExactErrorCylindricalHole} shows that this zone corresponds to a zone where the error is significant.

The trend of the results presented in figures \ref{EstimErrorPatternEllipticalHole} and  \ref{EstimErrorPatternCylindricalHole} is in favor of the residual-type indicator because this indicator appears to be the least affected by the discretization of the geometry. Among the different settings of the proposed ZZ-type indicator, the weighted - direct computation of $\sigma_{vM}^s$, combined with a discretization of the domain based on a triangular lattice, is the configuration that appears to be the least affected by the discretization of the geometry. We further investigated the impact of other parameters of the method for this configuration of the indicator.

\begin{figure}[h!]
	\centering
	\begin{tabular}{|M{4cm}|c|c|}
		\hline
		                                              & \textbf{Stencil Scaling Factor = 0.8}                                           & \textbf{Stencil Scaling Factor = 1.5}                                           \\
		\hline
		\multirow{2}{4cm}{\begin{minipage}[t]{0.8\columnwidth} Plate with an elliptical hole \\ ZZ-type indicator \\ Weighted - Direct $\sigma_{vM}^s$ \end{minipage}} &                                                                                 &                                                                                 \\[-1ex]
		                                              & \includegraphics[width=3.8cm]{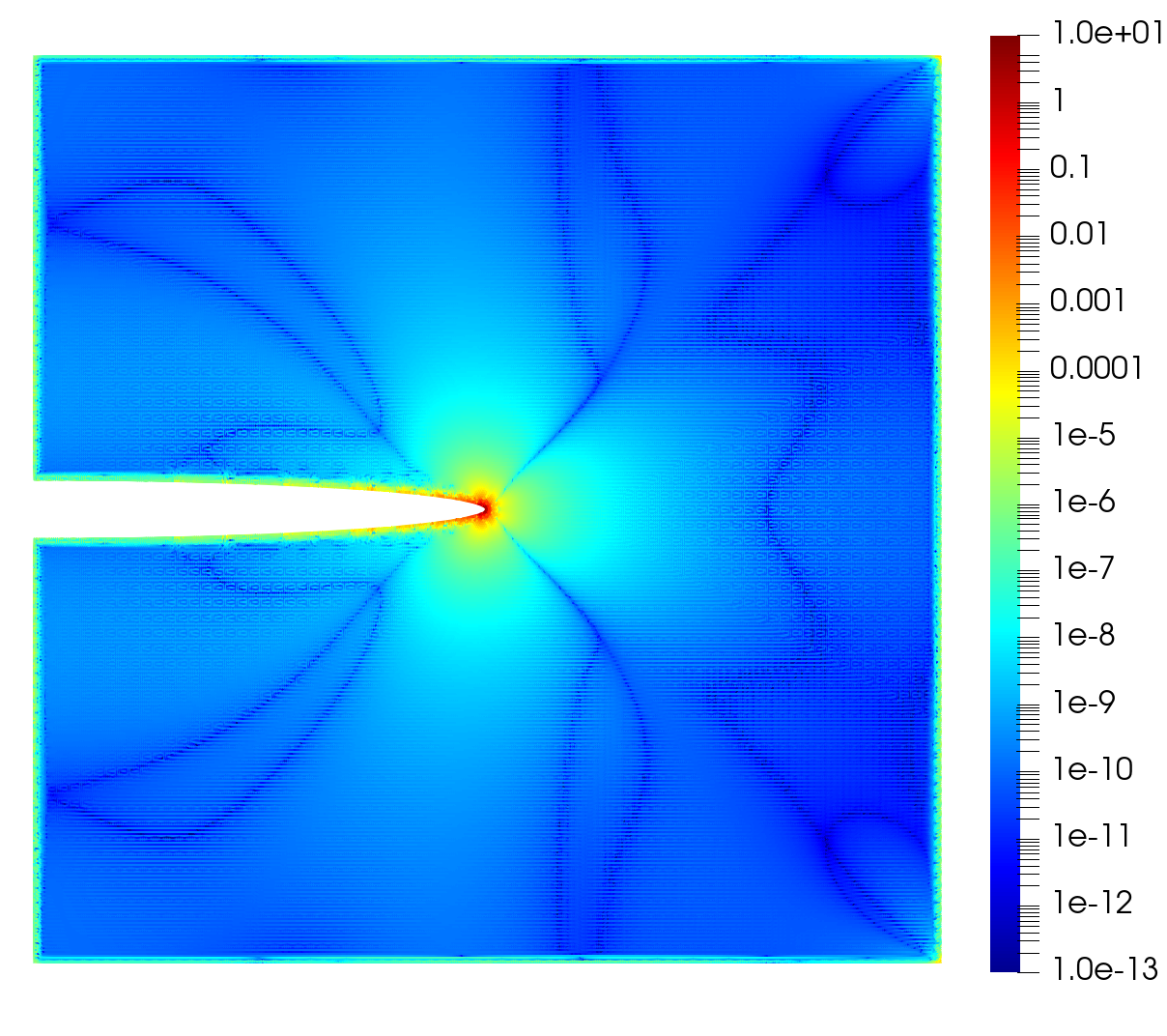} & \includegraphics[width=3.8cm]{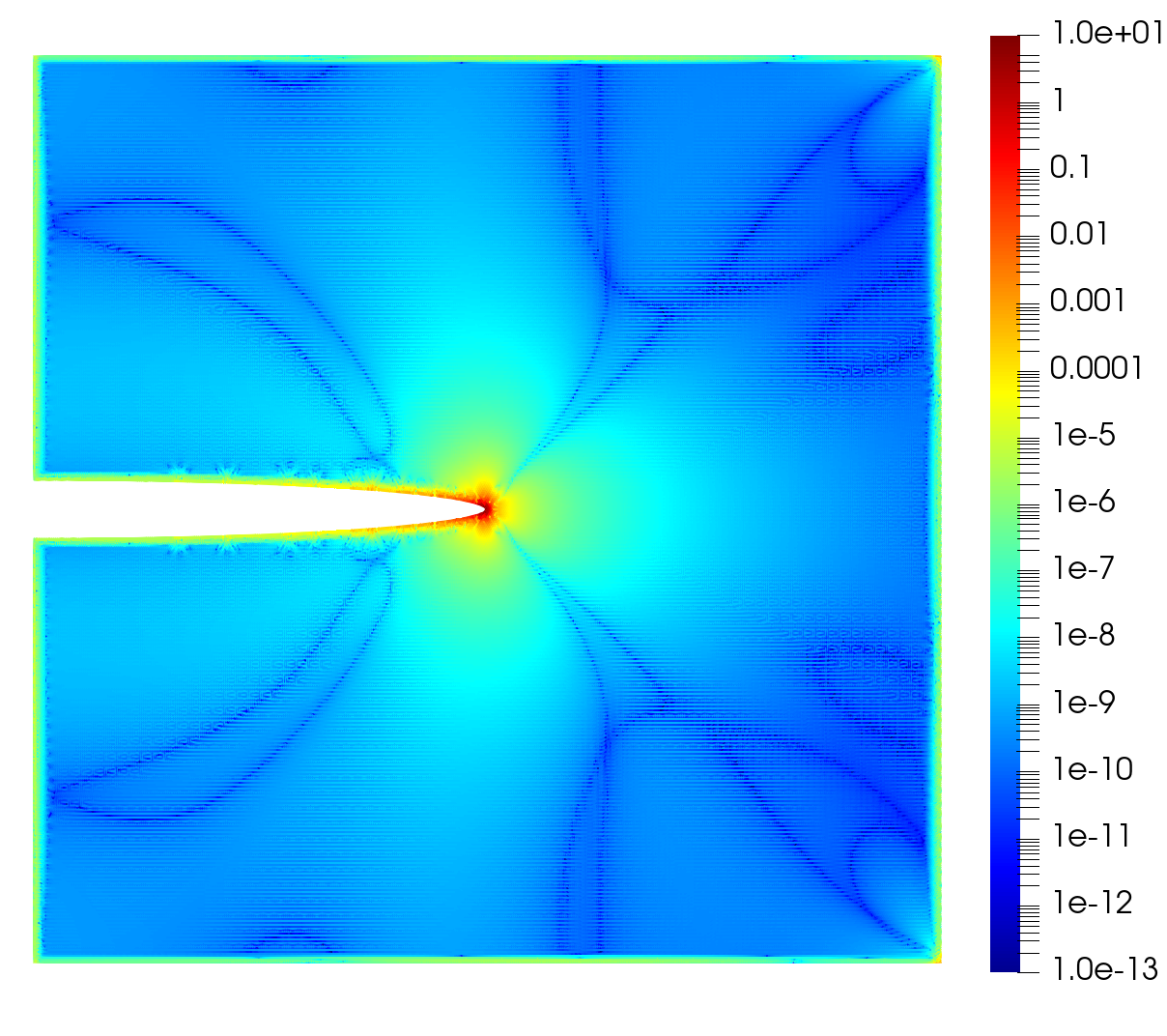} \\
		\hline
		\multirow{2}{4cm}{\begin{minipage}[t]{0.8\columnwidth} Body with a cylindrical hole \\ ZZ-type indicator \\ Weighted - Direct $\sigma_{vM}^s$ \end{minipage}} &                                                                                 &                                                                                 \\[-1ex]
		                                              & \includegraphics[width=3.8cm]{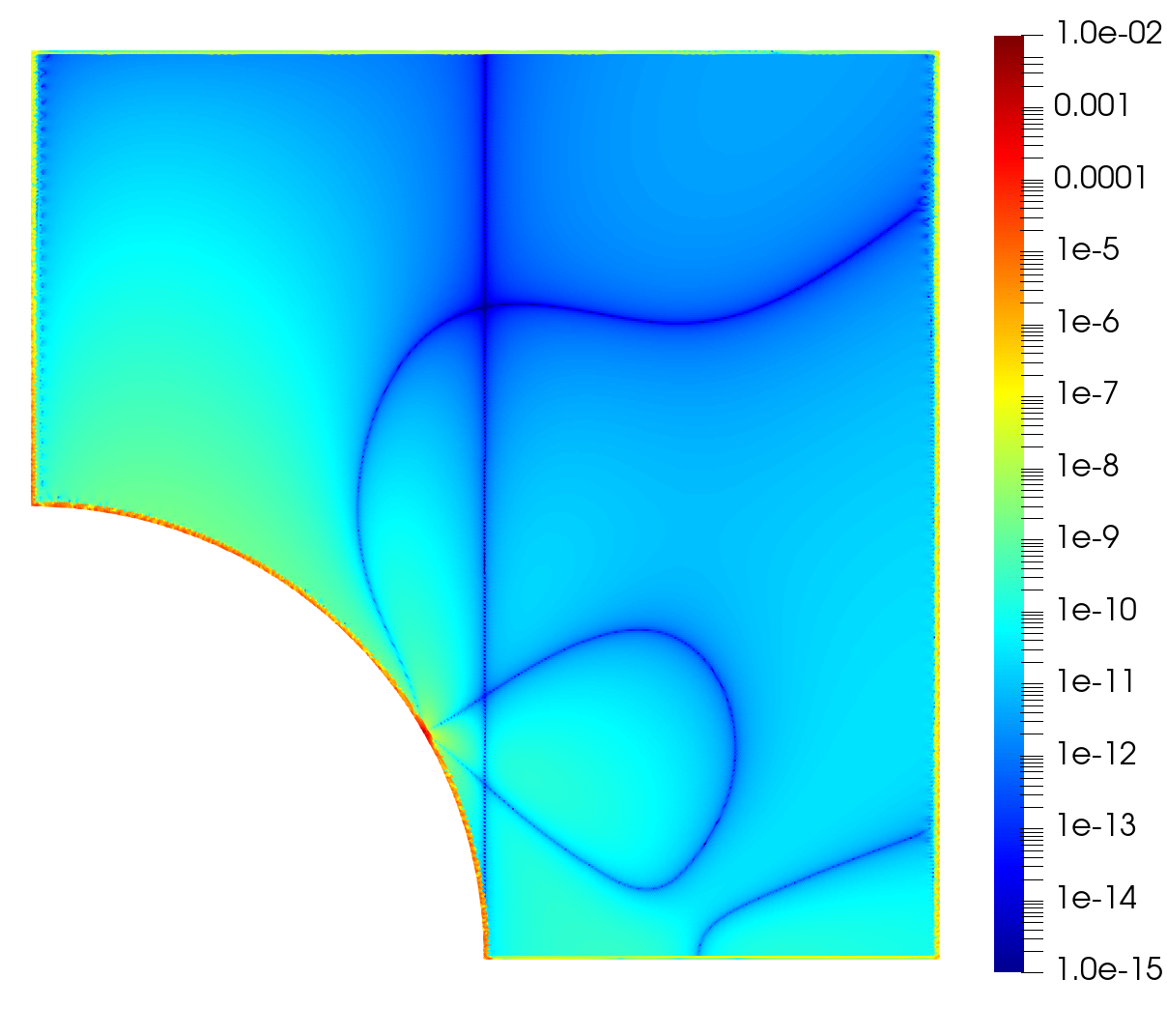}      & \includegraphics[width=3.8cm]{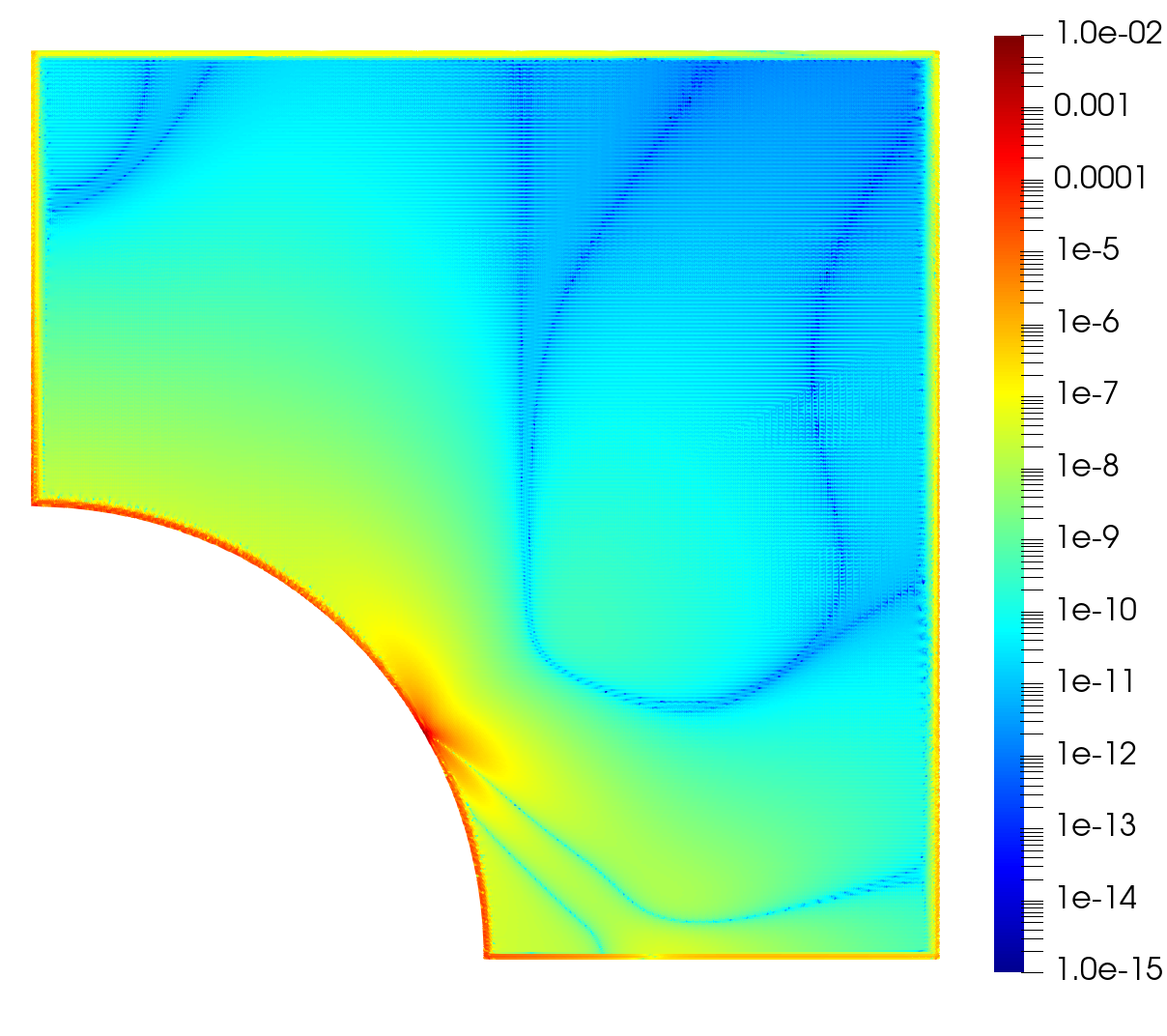}      \\
		\hline
	\end{tabular}
	\caption{Impact of the size of the stencil considered in the computation of the ZZ-type error indicator. Results are shown for both benchmark problems for two scaling factors applied to the selected stencil size considered in the solution of the collocation problem.}
	\label{SensitivityZZSupSize}
\end{figure}

The computation of the ZZ-type error indicator necessitates the selection of a stencil size or radius since we considered the distance criterion for the selection of the stencil nodes. The results presented in Figure \ref{EstimErrorPatternEllipticalHole} and Figure \ref{EstimErrorPatternCylindricalHole} were computed based on the same stencil as the stencil considered for the solution of the collocation problem. We present in Figure \ref{SensitivityZZSupSize} results showing the impact of the size of the stencil on the pattern of the error indicator. We considered two scaling factors applied to the size of the stencil selected for the solution of the collocation problem. We selected domains discretized based on a triangular lattice and a direct computation of $\sigma_{vM}^s$ considering a weighted moving least square approximation.

\clearpage

\clearpage

The results presented in Figure \ref{SensitivityZZSupSize} show, for the problem of a plate with an elliptical hole, that the scaling factor has little impact on the pattern of the error indicator. The impact of the scaling factor is more significant for the problem of a body with a cylindrical hole. We observe that the intensity of some zones where the error is the low is more significant for a scaling factor of 0.8. In these zones, the error indicator does not represent well the pattern of the exact error and is expected to lead to an incorrect refinement of the domain. A scaling factor of 1.5 leads to similar results than the base case (i.e. scaling factor of 1.0 presented in Figure \ref{EstimErrorPatternCylindricalHole}).

\subsubsection{Discussion}

Both the ZZ-type error indicator and the residual-type error indicator enable the identification of zones of the point cloud where the error is the greatest. We observe that a discretization of the interior of the domain based on a triangular lattice leads to an error indicator which is less dependent on the geometry than one based on a square lattice. This result is observed even though both discretization methods lead to similar exact error patterns as per the results presented in Figures \ref{ExactErrorEllipticalHole} and \ref{ExactErrorCylindricalHole}.

The weighted - direct computation of $\sigma_{vM}^s$ is the configuration of the ZZ-type error indicator that leads to the best results although little difference is observed between the weighted and the unweighted indicators. The selection of a stencil larger or smaller than the one used as part of the solution of the global collocation model does not improve the indicator. The residual-type error indicator leads to smoother results which are not affected by the discretization of the geometry. However, the computational cost of this indicator is much greater than for the ZZ-type error indicator. Voronoi corner nodes are computed for all the nodes of the domain. The stencil of each Voronoi corner node needs to be determined and the derivatives approximated.

\begin{figure}[!h]
	\centering
	\begin{tikzpicture}[scale=1]
		\begin{axis}[height=8cm,width=10.5cm, ymin=0.01, ymax=3, ymode=log, xmin=2000, xmax=1E6, xmode=log, legend style={text width=3.5cm, at={(0.5,-0.2)}, anchor=south west,legend columns=1, cells={anchor=west},  font=\footnotesize, rounded corners=2pt,}, legend pos=outer north east,xlabel=Number of Nodes, ylabel=\scalebox{1.35}{$\frac{\text{Indicator computation time}}{\text{Assembly + solution time of collocation problem}}$}]
			\addplot+[Tblue,mark=triangle*,mark options={fill=Tblue}]   table [x=NodeNum-Ell-CAD-Tri, y=L2R-VMS-Ell-Est-PP-Frac
					, col sep=comma] {ZZEstimatorResults.csv};
			\addlegendentry[minimum height=0.5in]{Plate w/ an elliptical hole - ZZ-type indicator};
			\addplot+[Tred,mark=diamond*,mark options={fill=Tred}]   table [x=NodeNum-Ell-CAD-Tri, y=L2R-U-Ell-ResEst-PP-Frac
					, col sep=comma] {ZZEstimatorResults.csv};
			\addlegendentry[minimum height=0.5in]{Plate w/ an elliptical hole - residual-type indicator};
			\addplot+[Tgreen,mark=*,mark options={fill=Tgreen}]   table [x=NodeNum-Hole-CAD-Tri, y=L2R-VMS-Hole-Est-PP-Frac
					, col sep=comma] {ZZEstimatorResults.csv};
			\addlegendentry[minimum height=0.5in]{Body w/ a cylindrical hole - ZZ-type indicator};
			\addplot+[Torange,mark=square*,mark options={fill=Torange}]   table [x=NodeNum-Hole-CAD-Tri, y=L2R-U-Hole-ResEst-PP-Frac
					, col sep=comma] {ZZEstimatorResults.csv};
			\addlegendentry[minimum height=0.5in]{Body w/ a cylindrical hole - residual-type indicator};
		\end{axis}
	\end{tikzpicture}
	\caption{Comparison of the ratio of the indicator computation time to the time needed to assemble and solve the collocation problem. Both problems lead to similar results. The computation time of the residual-type error indicator is approximately 10 times the computation time of the ZZ-type error indicator. The computation of the residual-type error indicator corresponds to between 165\% and 23\% of the assembly and solution time of the collocation problem. The computation of the ZZ-type error indicator corresponds to between 11\% and 3\% of the assembly and solution time of the collocation problem.}
	\label{EstimatorComputationTime}
\end{figure}
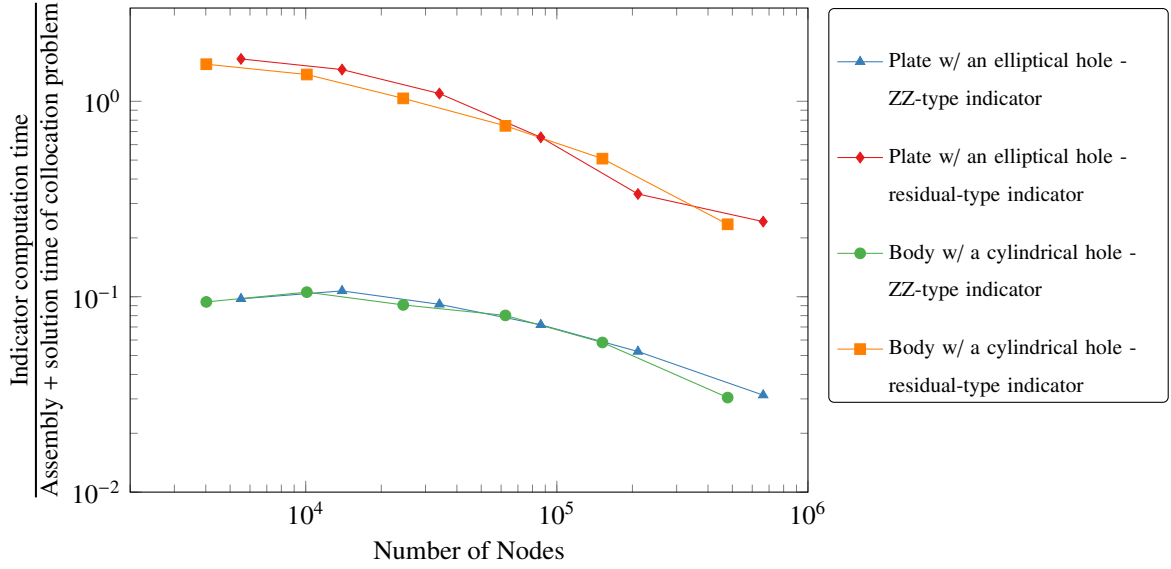

We show a comparison of the time required to compute the error indicators for the two considered problems in Figure \ref{EstimatorComputationTime}. We present these results in the form of a ratio of the indicator computation time to the time needed to assemble and solve the collocation problem. The results are indicative as they depend heavily on the method selected to solve the linear system and on the number of threads/processes involved in each step of the solution process. The computation of the indicators at each node of the domain is independent from the computation of the indicator at other nodes of the domain. Therefore, both error indicators can be parallelized with no extra effort. We used two threads to assemble the linear system and compute the error indicator and solved the linear system using a LU factorization (one thread, one process). We see from Figure \ref{EstimatorComputationTime} that the computation time of the residual-type error indicator is approximately 10 times the computation time of the ZZ-type error indicator. For both indicators, the ratio mostly decreases as the number of nodes increases. Such a result is expected since the computation of the error indicators scales linearly with the number nodes while the solution time of the collocation problem increases at an increasing rate when solved with LU factorization. The problem of the plate with an elliptical hole and the problem of a body with a cylindrical hole lead to similar results. The computation of the residual-type error indicator corresponds to between 165\% and 23\% of the assembly and solution time of the collocation problem. The computation of the ZZ-type error indicator corresponds to between 11\% and 3\% of the assembly and solution time of the collocation problem for the considered discretizations.

Based on these results, the ZZ-type indicator appears to be a better choice in terms of computational cost.

To complete the comparison of both error indicators, we computed ZZ-type and residual-type error indicators for the problem of a plate with an elliptical hole and the problem of a body with a cylindrical hole for various node densities obtained with a global refinement of the domains. We selected for both problems a discretization of the domain based on a triangular lattice. The results are presented in Figure \ref{ConvResults_Ell} and Figure \ref{ConvResults_Hole}. The ZZ-type error indicator is based on a weighted - direct computation of $\sigma_{vM}^s$. We computed the exact error in terms of the $L_2$ relative error norm for the von Mises stress component and compared it to a $L_2$ relative error norm where the ``smooth'' von Mises stress field is considered as the reference solution. We also computed the $L_2$ weighted error norm of the residual-type indicator computed using Equation \ref{L2WeightedNorm_Estim}. We selected this norm since a reference solution is not computed as part of the indicator.

\begin{equation} \label{L2WeightedNorm_Estim}
	L_2W(e)= \frac{\sqrt{\sum_{k=1}^{n}{\left( e\left(\mathbf{X_c}\right) \right)^2}}}{n}.
\end{equation}

We observe from the results presented in Figure \ref{ConvResults_Ell} and in Figure \ref{ConvResults_Hole} a convergence of both error indicators for both problems. The convergence rate is higher for the problem of a body with a cylindrical hole. The convergence rate of the $L_2$ weighted error norm computed for the residual-type indicator is 0.1 for the problem of a plate with an elliptical hole and 1.2 for the problem of a body with a cylindrical hole. This result is expected since the solution is smooter for this problem than for the problem of a plate with an elliptical hole. For both problems, the convergence rate of the $L_2$ relative error norm is similar when the exact or smooth von Mises stress components are considered.

\begin{figure}[!h]
	\centering
	{
		\begin{tikzpicture}[scale=1]
			\begin{axis}[height=5cm,width=7cm, ymin=0.01, ymax=1, ymode=log, xmin=2000, xmax=1E6, xmode=log, legend entries={Exact Error - $\sigma_{vM}$}, legend style={ at={(0.5,-0.2)}, anchor=south west,legend columns=1, cells={anchor=west},  font=\footnotesize, rounded corners=2pt,}, legend pos=north east,xlabel=Number of Nodes, ylabel=$L_2$ Relative Error]
				\addplot+[Tblue,mark=triangle*,mark options={fill=Tblue}]   table [x=NodeNum-Ell-CAD-Tri, y=L2R-VMS-Ell-Ex
						, col sep=comma] {ZZEstimatorResults.csv};
				\logLogSlopeTriangle{0.85}{0.1}{0.15}{0.51}{black};
			\end{axis}
		\end{tikzpicture}
	}
	{
		\begin{tikzpicture}[scale=1]
			\begin{axis}[height=5cm,width=7cm, ymin=0.001, ymax=0.1, ymode=log, xmin=2000, xmax=1E6, xmode=log, legend entries={ZZ-type error indicator - $\sigma_{vM}$}, legend style={ at={(0.5,-0.2)}, anchor=south west,legend columns=1, cells={anchor=west},  font=\footnotesize, rounded corners=2pt,}, legend pos=north east,xlabel=Number of Nodes, ylabel=$L_2$ Relative Error]
				\addplot+[Tblue,mark=triangle*,mark options={fill=Tblue}]   table [x=NodeNum-Ell-CAD-Tri, y=L2R-VMS-Ell-Est-S4-DirectVMS
						, col sep=comma] {ZZEstimatorResults.csv};
				\logLogSlopeTriangle{0.85}{0.1}{0.15}{0.45}{black};
			\end{axis}
		\end{tikzpicture}
	}\\
	{
	\begin{tikzpicture}[scale=1]
		\begin{axis}[height=5cm,width=7cm, ymin=0.005, ymax=0.1, ymode=log, xmin=2000, xmax=1E6, xmode=log, legend entries={Residual-type indicator}, legend style={ at={(0.5,-0.2)}, anchor=south west,legend columns=1, cells={anchor=west},  font=\footnotesize, rounded corners=2pt,}, legend pos=north east,xlabel=Number of Nodes, ylabel=$L_2$ Weighted Error]
			\addplot+[Tblue,mark=triangle*,mark options={fill=Tblue}]   table [x=NodeNum-Ell-CAD-Tri, y=L2R-U-Ell-Est-ResU-L2W
					, col sep=comma] {ZZEstimatorResults.csv};
			\logLogSlopeTriangle{0.85}{0.1}{0.15}{0.10}{black};
		\end{axis}
	\end{tikzpicture}
	}\\
	\caption{Plate with an elliptical hole - Comparison of the exact and estimated error. The results are presented in terms of the $L_2$ relative error norm for the exact error on the von Mises stress and for the ZZ-type indicator of the error on the von Mises stress. The results are presented in terms of the $L_2$ weighted error norm for the residual-type indicator.}
	\label{ConvResults_Ell}
\end{figure}
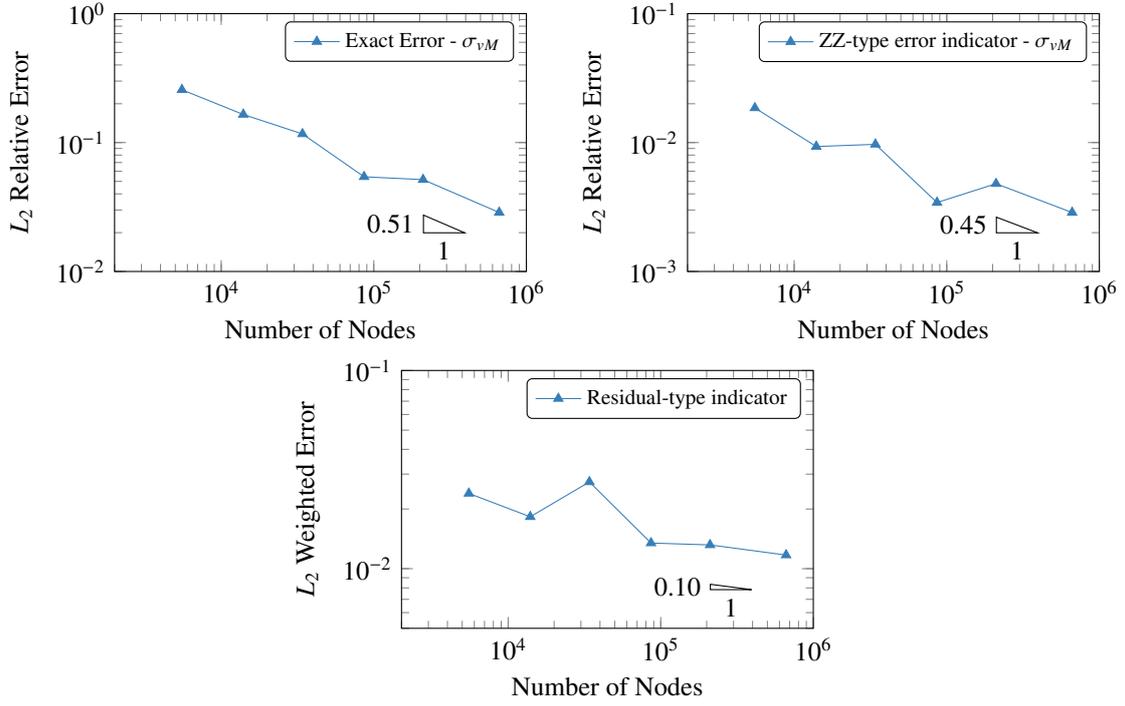

\begin{figure}[!h]
	\centering
	{
		\begin{tikzpicture}[scale=1]
			\begin{axis}[height=5cm,width=7cm, ymin=1e-6, ymax=1e-2, ymode=log, xmin=2000, xmax=1E6, xmode=log, legend entries={Exact Error - $\sigma_{vM}$}, legend style={ at={(0.5,-0.2)}, anchor=south west,legend columns=1, cells={anchor=west},  font=\footnotesize, rounded corners=2pt,}, legend pos=north east,xlabel=Number of Nodes, ylabel=$L_2$ Relative Error]
				\addplot+[Tblue,mark=triangle*,mark options={fill=Tblue}]   table [x=NodeNum-Hole-CAD-Tri, y=L2R-VMS-Hole-Ex
						, col sep=comma] {ZZEstimatorResults.csv};
				\logLogSlopeTriangle{0.85}{0.1}{0.15}{1.02}{black};
			\end{axis}
		\end{tikzpicture}
	}
	{
		\begin{tikzpicture}[scale=1]
			\begin{axis}[height=5cm,width=7cm, ymin=1e-7, ymax=5e-3, ymode=log, xmin=2000, xmax=1E6, xmode=log, legend entries={ZZ-type error indicator - $\sigma_{vM}$}, legend style={ at={(0.5,-0.2)}, anchor=south west,legend columns=1, cells={anchor=west},  font=\footnotesize, rounded corners=2pt,}, legend pos=north east,xlabel=Number of Nodes, ylabel=$L_2$ Relative Error]
				\addplot+[Tblue,mark=triangle*,mark options={fill=Tblue}]   table [x=NodeNum-Hole-CAD-Tri, y=L2R-VMS-Hole-Est-S4-DirectVMS
						, col sep=comma] {ZZEstimatorResults.csv};
				\logLogSlopeTriangle{0.85}{0.1}{0.15}{0.95}{black};
			\end{axis}
		\end{tikzpicture}
	}\\
	{
	\begin{tikzpicture}[scale=1]
		\begin{axis}[height=5cm,width=7cm, ymin=1e-7, ymax=5e-4, ymode=log, xmin=2000, xmax=1E6, xmode=log, legend entries={Residual Residual}, legend style={ at={(0.5,-0.2)}, anchor=south west,legend columns=1, cells={anchor=west},  font=\footnotesize, rounded corners=2pt,}, legend pos=north east,xlabel=Number of Nodes, ylabel=$L_2$ Weighted Error]
			\addplot+[Tblue,mark=triangle*,mark options={fill=Tblue}]   table [x=NodeNum-Hole-CAD-Tri, y=L2R-U-Hole-Est-ResU-L2W
					, col sep=comma] {ZZEstimatorResults.csv};
			\logLogSlopeTriangle{0.85}{0.1}{0.15}{1.21}{black};
		\end{axis}
	\end{tikzpicture}
	}\\
	\caption{Body with a cylindrical hole - Comparison of the exact and indicative error. The results are presented in terms of the $L_2$ relative error norm for the exact error on the von Mises stress and for the ZZ-type indicator of the error on the von Mises stress. The results are presented in terms of the $L_2$ weighted error norm for the residual-type indicator.}
	\label{ConvResults_Hole}
\end{figure}

Based on the results presented in this subsection, we selected the ZZ-type error indicator, using a weighted - direct computation of $\sigma_{vM}^s$, as the main input to the discretization refinement scheme presented in the next section. We preferred this indicator rather than the residual-type error indicator because of its reduced computational expense and because of its simplicity.

\subsection{Discretization refinement}\label{AdaptSec}

We describe in this subsection a scheme to refine a point cloud based on the results obtained from \textit{a posteriori} error indicators. We use the information of the smart cloud to place new nodes on the exact geometry of the domain and to set the boundary conditions of the updated collocation model.

\textit{A posteriori} error indicators allow the identification of the areas of the domain where the error is expected to be the greatest. Several techniques can then be used to improve the solution in these zones. The most common ones are $h$- and $p$- adaptivity. $h$-adaptivity consists of an increase of the node density in the zone where the error is the greatest. Such type of adaptivity is the most commonly used in the literature. For instance, it was used by Benito \cite{Benito2003}, Davydov \cite{Davydov2011a}, Gavete \cite{Gavete2015} or Slak and Kosec \cite{Slak2019} within the framework of point collocation. $h$-adaptivity has also be considered in literature based on geometric indicators \cite{Suchde2019b, Suchde2017}. $p$-adaptivity is another technique which consists in an increase of the order of the approximation. Liszka et al. and Duarte et al. made use of $p$-adaptivity, in the framework of meshless methods, as part of the hp-meshless cloud method \cite{Liszka1996,Duarte1996}. Jancic et al. showed the benefits of $p$-refinement for a Poisson problem with a strong source within the domain \cite{Jancic2021} in the framework of the radial basis function-generated finite difference method (RBF-FD).

Our scheme is based on $h$-adaptivity. We selected this technique to be able to perform successive refinement iterations and reduce the observed error as much as possible. An $h$-adaptive scheme is based on the repetition of a succession of steps. The main steps that we followed are presented in the form of a pseudo code in Figure \ref{AlgorithmAdaptivity}.

\begin{figure}[!h]
	\begin{tabular}{ | p{15.5cm} |}
		\hline
		\textbf{$h$-adaptivity algorithm}                                                                    \\
		\hline
		Initial discretization of the domain;                                                                \\
		Application of boundary conditions and generation of the initial collocation input file;             \\
		\textbf{While} (Continue adaptive refinement = true)                                                 \\
		\quad Solution of the collocation problem;                                                           \\
		\quad Computation of the \textit{a posteriori} error indicator;                                        \\
		\quad Assessment of the need for a new refinement iteration;                                         \\
		\quad \textbf{if} (Refine the discretization = true)                                                 \\
		\quad \quad Identification of the nodes of highest error for discretization refinement;              \\
		\quad \quad Placement of new nodes in the areas of high error;                                       \\
		\quad \quad Application of boundary conditions and generation of the updated collocation input file; \\
		\quad \textbf{else}                                                                                  \\
		\quad \quad Continue adaptive refinement = false.                                                    \\
		\textbf{loop}                                                                                        \\
		\hline
	\end{tabular}
	\caption{$h$-adaptivity algorithm considered for the presented adaptive refinement method.}
	\label{AlgorithmAdaptivity}
\end{figure}

\begin{samepage}
	The algorithm is composed of three main steps:
	\begin{enumerate}
		\setlength{\itemsep}{0pt}
		      \setlength{\parskip}{0pt}
		\item Identification of the refinement areas;
		\item Placement of new nodes;
		\item Generation of the updated collocation model
	\end{enumerate}
	We describe these steps in the subsections below.
\end{samepage}

\subsubsection{Identification of the refinement areas} \label{RefinementIdSec}

The error indicators presented in Section \ref{SecErrorIndicators} are computed at the collocation nodes. Therefore, we decided to identify the areas of the domain to be refined based on a selection of marked collocation nodes. The selection of these nodes is based on the computed error indicator. Different criteria can be used to determine the collocation nodes to be marked for local refinement of the domain.

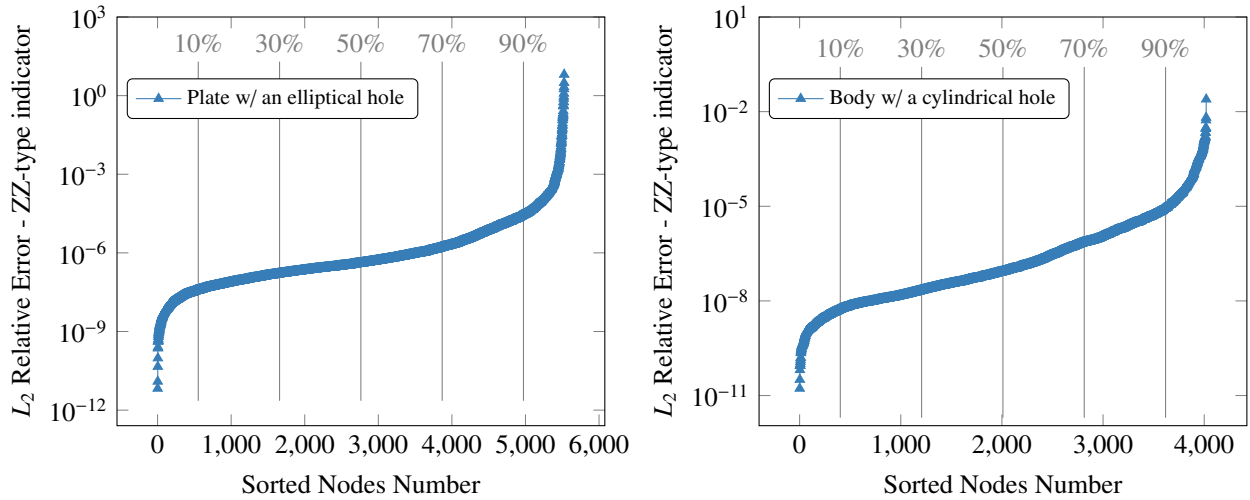
\begin{figure}[!h]
	\centering
	\begin{tikzpicture}[scale=1]
		\begin{axis}[height=7cm,width=8cm, ymax=1000, ymode=log, legend entries={Plate w/ an elliptical hole}, legend style={ at={(0.02,0.85)}, anchor=north west,legend columns=1, cells={anchor=west},  font=\footnotesize, rounded corners=2pt,},xlabel=Sorted Nodes Number, ylabel=$L_2$ Relative Error - ZZ-type indicator]
			\addplot+[Tblue,mark=triangle*,mark options={fill=Tblue}]   table [x=N, y=ZZ-Error
					, col sep=comma] {EllipticalHole-SortedZZError.csv};
			\node(s1) at (axis cs: 553,1E-12){};\node[text=gray](d1) at (axis cs: 553,100) {10\%};
			\draw[gray](s1)--(d1);
			\node(s2) at (axis cs: 1658,1E-12){};\node[text=gray](d2) at (axis cs: 1658,100) {30\%};
			\draw[gray](s2)--(d2);
			\node(s3) at (axis cs: 2763,1E-12){};\node[text=gray](d3) at (axis cs: 2763,100) {50\%};
			\draw[gray](s3)--(d3);
			\node(s4) at (axis cs: 3868,1E-12){};\node[text=gray](d4) at (axis cs: 3868,100) {70\%};
			\draw[gray] (s4)--(d4);
			\node(s5) at (axis cs: 4973,1E-12){};\node[text=gray](d5) at (axis cs: 4973,100) {90\%};
			\draw[gray](s5)--(d5);
		\end{axis}
	\end{tikzpicture}
	\begin{tikzpicture}[scale=1]
		\begin{axis}[height=7cm,width=8cm, ymax=10, ymode=log, legend entries={Body w/ a cylindrical hole}, legend style={ at={(0.02,0.85)}, anchor=north west,legend columns=1, cells={anchor=west},  font=\footnotesize, rounded corners=2pt,},xlabel=Sorted Nodes Number, ylabel=$L_2$ Relative Error - ZZ-type indicator]
			\addplot+[Tblue,mark=triangle*,mark options={fill=Tblue}]   table [x=N, y=ZZ-Error
					, col sep=comma] {RoundHole-SortedZZError.csv};
			\node(s1) at (axis cs: 402,1E-12){};\node[text=gray](d1) at (axis cs: 402,1) {10\%};
			\draw[gray](s1)--(d1);
			\node(s2) at (axis cs: 1206,1E-12){};\node[text=gray](d2) at (axis cs: 1206,1) {30\%};
			\draw[gray](s2)--(d2);
			\node(s3) at (axis cs: 2010,1E-12){};\node[text=gray](d3) at (axis cs: 2010,1) {50\%};
			\draw[gray](s3)--(d3);
			\node(s4) at (axis cs: 2814,1E-12){};\node[text=gray](d4) at (axis cs: 2814,1) {70\%};
			\draw[gray] (s4)--(d4);
			\node(s5) at (axis cs: 3618,1E-12){};\node[text=gray](d5) at (axis cs: 3618,1) {90\%};
			\draw[gray](s5)--(d5);
		\end{axis}
	\end{tikzpicture}
	\caption{Distribution of the error in terms of ZZ-type error indicator for the problem of a plate with an elliptical hole (left) and for the problem of a body with a cylindrical hole (right). Three distinct zones are observed on these graphs. Less than 10\% of the nodes have a very low error (nearly zero). Approximately 80-90\% of the nodes have an error in the similar range (10$^{\text{-8}}$-10$^{\text{-4}}$
		for the plate with an elliptical hole, 10$^{\text{-8}}$-10$^{\text{-5}}$ for the body with a cylindrical hole). Less than 10\% of the nodes have an error much larger than the other nodes.}
	\label{ZZSortedErrorRound}
\end{figure}

An error indicator threshold could be selected by the user. However such a threshold is problem specific and cannot be easily generalized to all problems. To help the definition of a node selection criterion, we presented the computed error indicator, sorted from the lowest error to the highest error in Figure \ref{ZZSortedErrorRound} for the considered benchmark problems. The presented error indicators are ZZ-type error indicators, computed based on the parameters presented in Section \ref{SubErrorInd}. Vertical bars are used to visualize the percentage of nodes in the different error zones.

We observe three distinct zones from the results presented in Figure \ref{ZZSortedErrorRound}. A limited number of nodes have a very low error (approximately below 10$^{\text{-8}}$ for both problems). A majority of the nodes have an error between 10$^{\text{-8}}$ and 10$^{\text{-4}}$ for the plate with an elliptical hole an between 10$^{\text{-8}}$ and 10$^{\text{-5}}$ for the body with a cylindrical hole. Finally, a limited number of nodes have an error larger than 10$^{\text{-4}}$ and 10$^{\text{-5}}$ for the plate with an elliptical hole and the body with a cylindrical hole, respectively. The zones of low and high error represent each approximately 5\% of the total number of nodes. Based on the these observations, we decided to select the nodes of highest error based on a defined fraction $f$ of the total number of nodes $n$. Therefore, the $fn$ nodes having the highest error are selected for local refinement. The impact of the selected threshold on the convergence rate is analyzed in a later section (see Subsection \ref{ThresholdSensSec}). This approach has similarities with the D\"{o}rfler marking strategy \cite{Dorfler1996} used in the framework of adaptive finite element method \cite{Bulle2021}.

\begin{figure}[h!]
	\centering
	\begin{tabular}{|M{4cm}|c|c|}
		\hline
		                                              & \textbf{Plate with an elliptical hole}                 & \textbf{Body with a cylindrical hole}             \\
		\hline
		\multirow{2}{4cm}{\begin{minipage}[t]{0.8\columnwidth} ZZ-type indicator \\ Weighted - Direct $\sigma_{vM}^s$ \end{minipage}} &                                                        &                                                   \\[-1ex]
		                                              & \includegraphics[width=5cm]{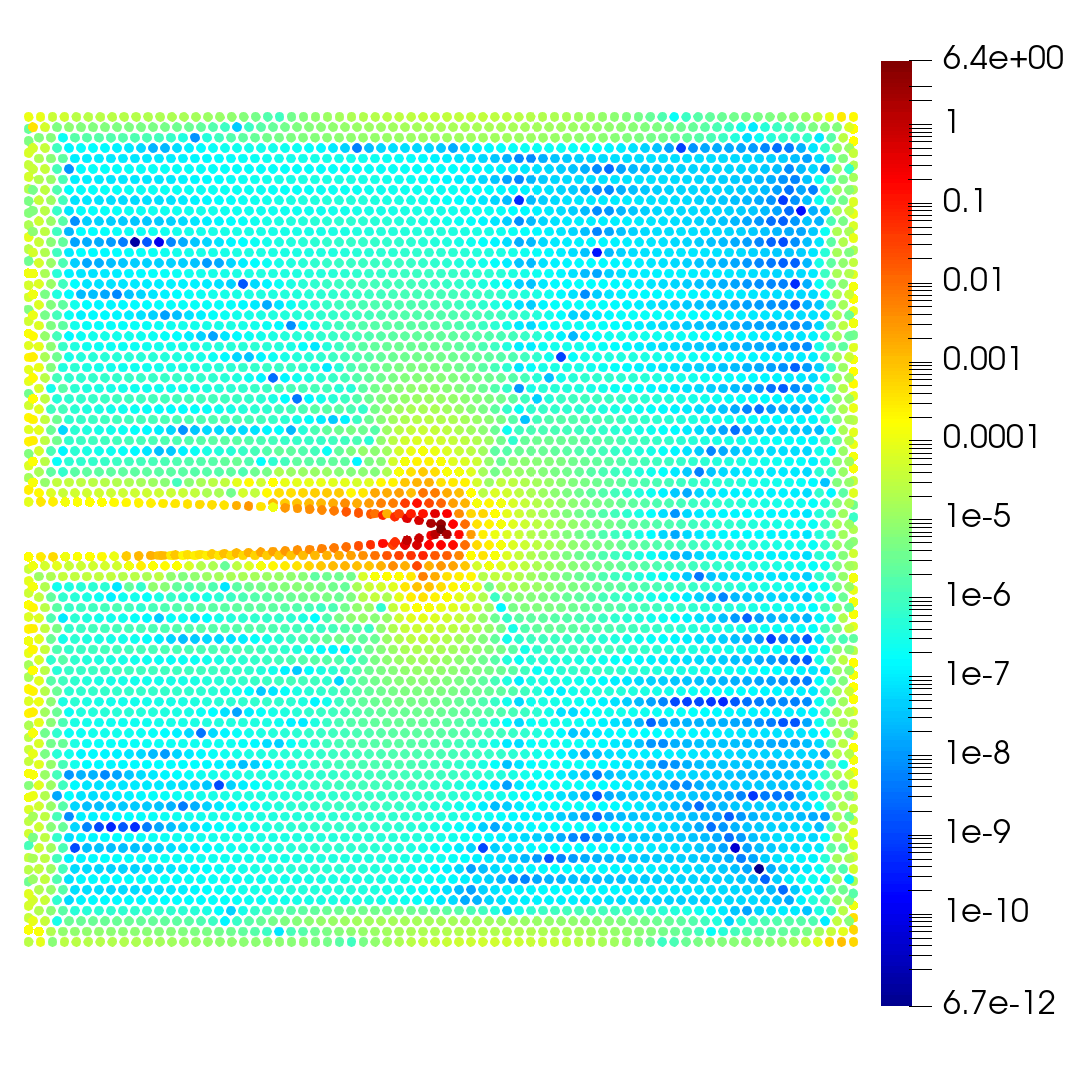}  & \includegraphics[width=5cm]{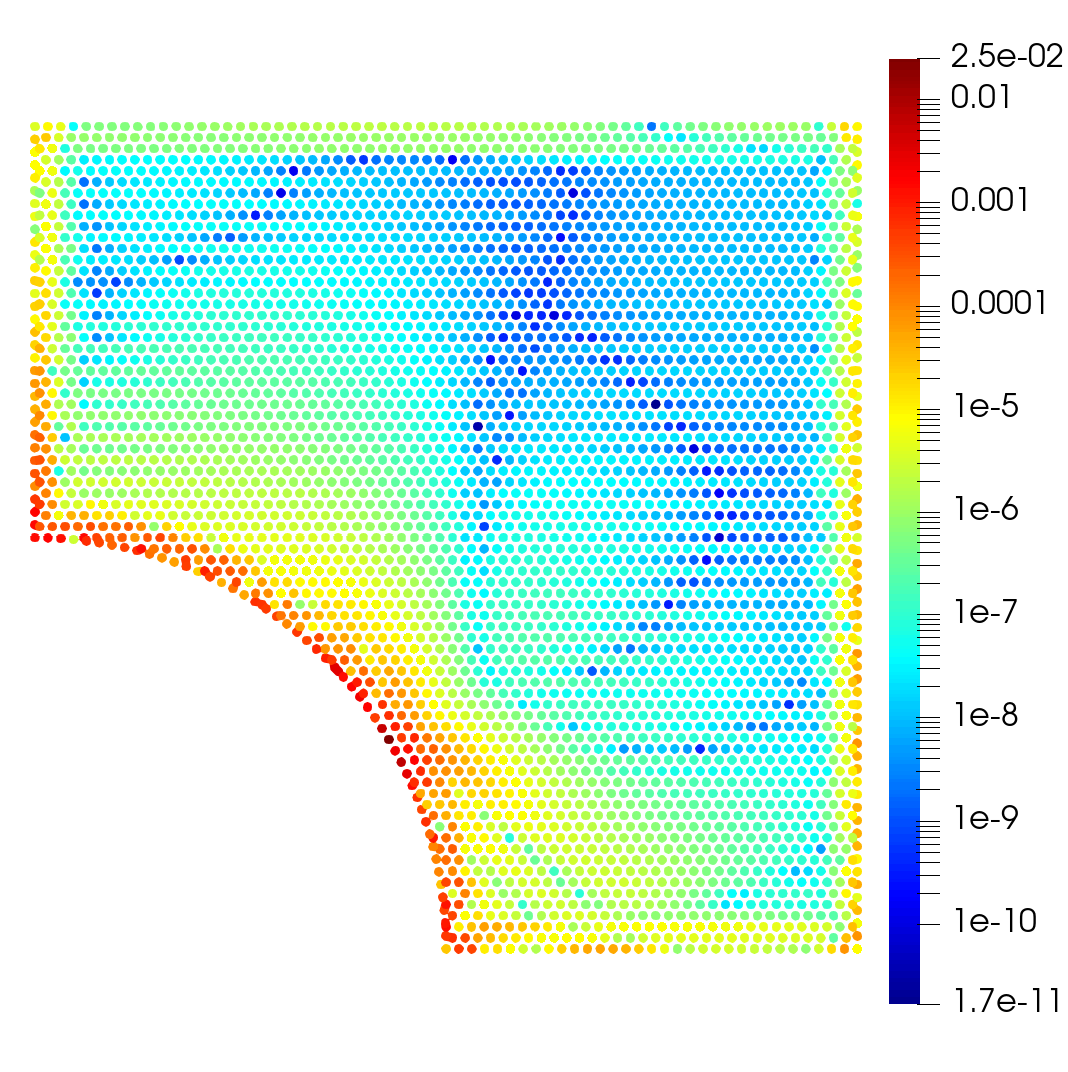}  \\
		\hline
		\multirow{2}{4cm}{\begin{minipage}[t]{0.8\columnwidth} Marked collocation nodes \\ (in red) based on a \\ threshold criteria \\ corresponding to 10\% of the \\ nodes of highest error \end{minipage}} &                                                        &                                                   \\[-1ex]
		                                              & \includegraphics[width=4cm]{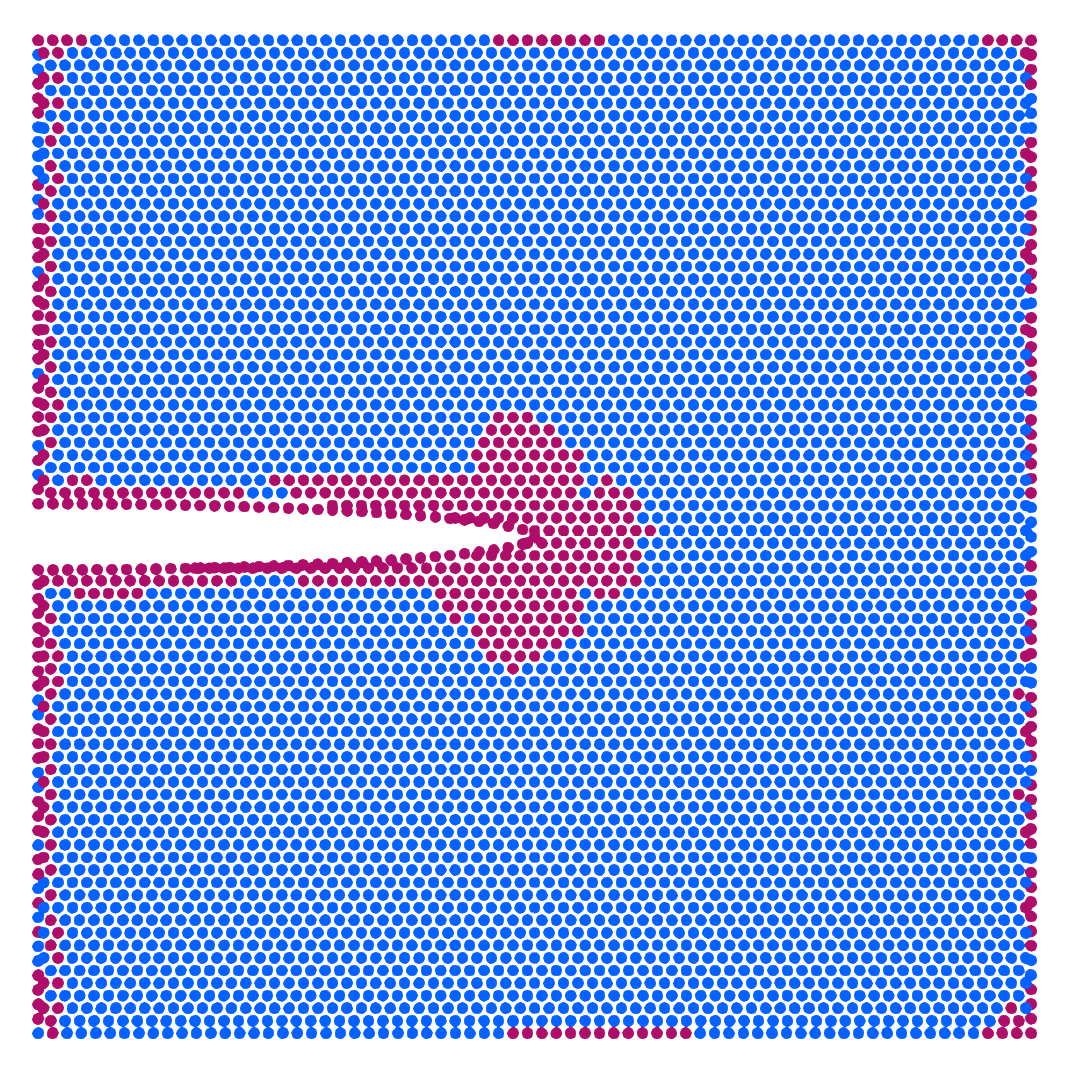}    & \includegraphics[width=4cm]{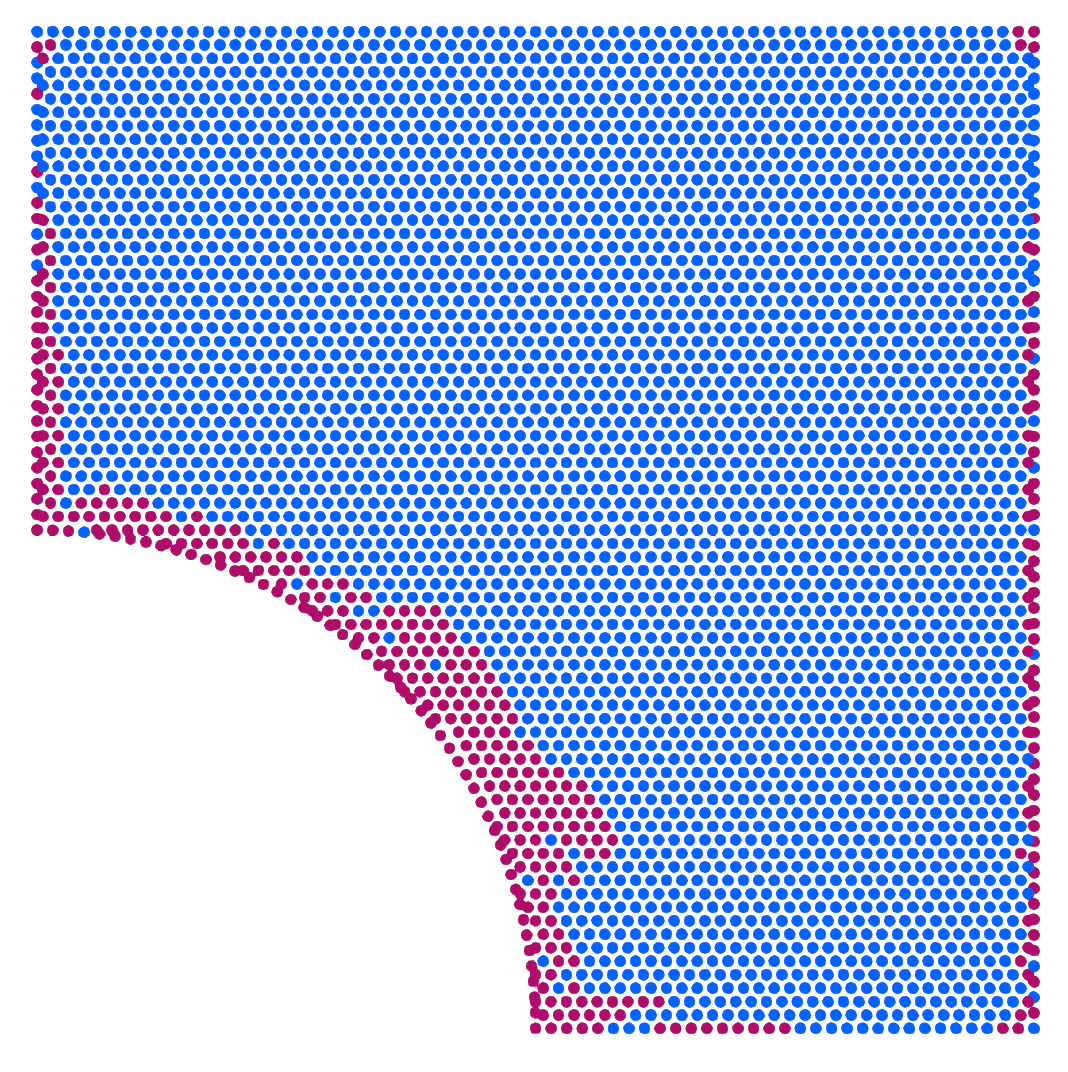}    \\
		\hline
		\multirow{2}{4cm}{\begin{minipage}[t]{0.8\columnwidth} Marked collocation nodes \\ (in red) based on a \\ threshold criteria \\ corresponding to 10\% of the \\ nodes of highest error and \\ the associated stencil nodes \end{minipage}} &                                                        &                                                   \\[-1ex]
		                                              & \includegraphics[width=4cm]{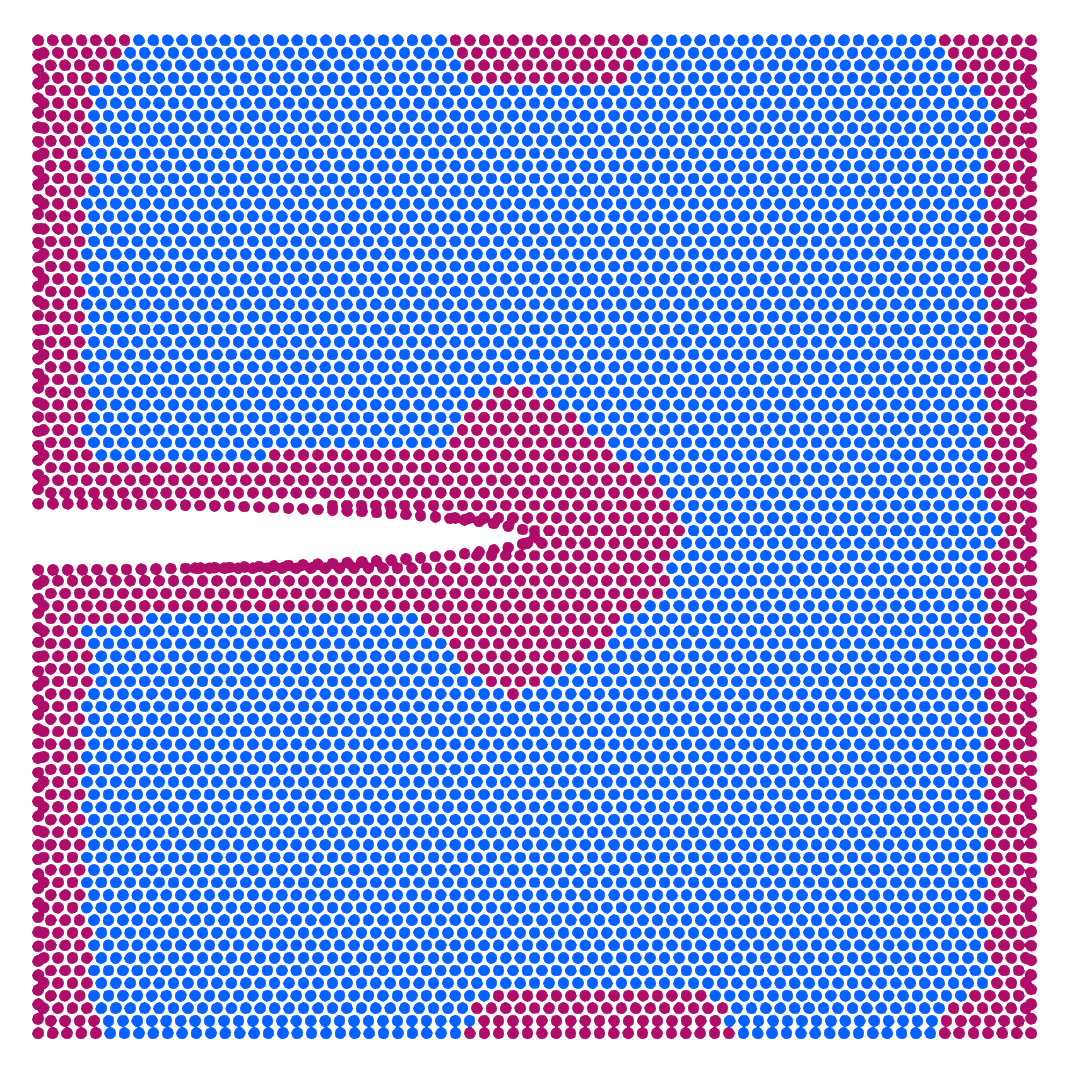} & \includegraphics[width=4cm]{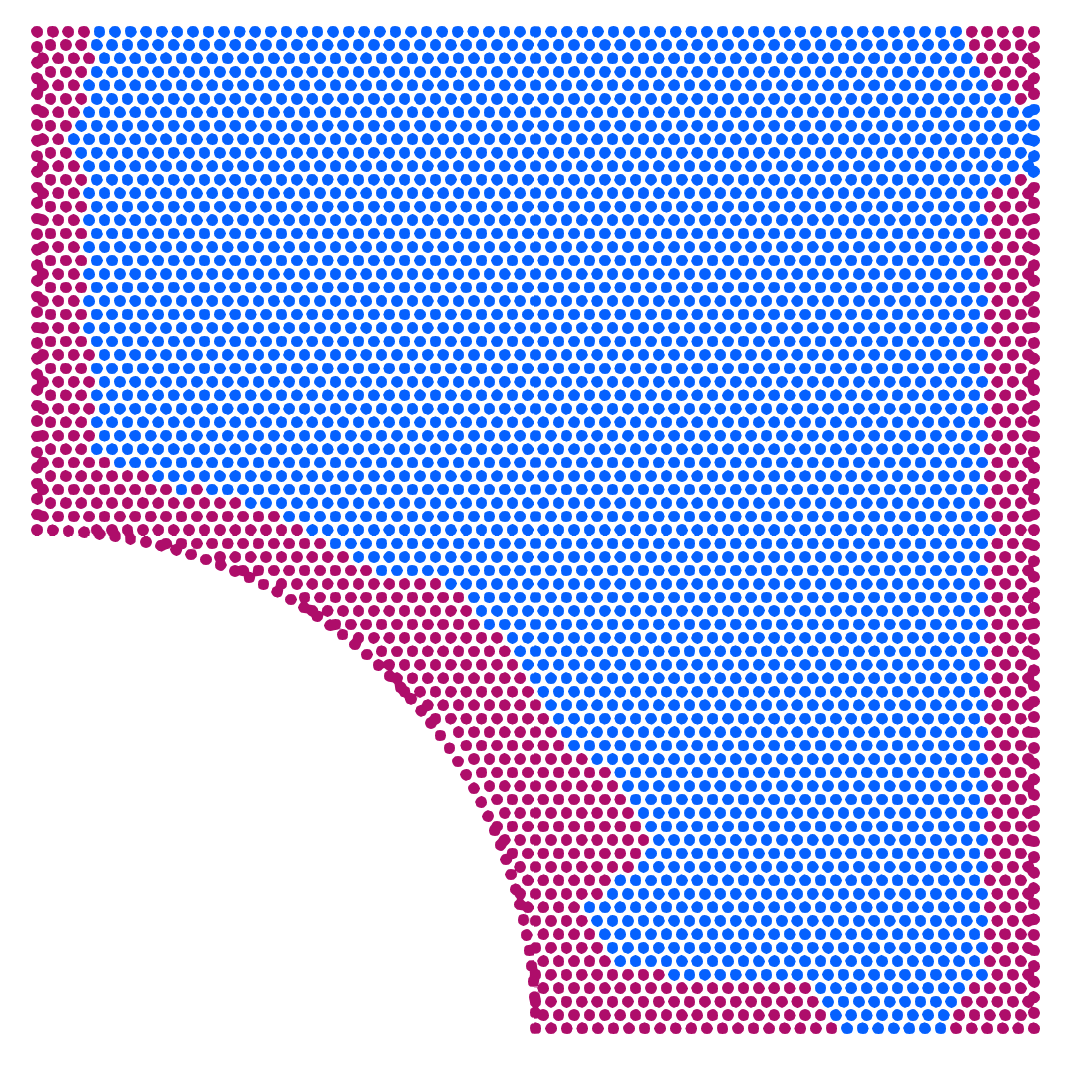} \\
		\hline
	\end{tabular}
	\caption{Pattern of the ZZ-type error indicator and selection of the nodes for local refinement based on a fraction of the nodes of highest error and based on a fraction of the nodes of highest error and all the stencil nodes of these collocation nodes. The results are presented for the problems of a plate with an elliptical hole and for the body with a cylindrical hole.}
	\label{NodeSelectionForRefinement}
\end{figure}

In order to obtain as smooth a refinement pattern as possible, we also selected, for local refinement, all the stencil nodes associated with the selected nodes of highest error. The results presented in Figure \ref{NodeSelectionForRefinement} show the benefits of this approach for the problems of a plate with an elliptical hole and for the body with a cylindrical hole. We presented first the pattern of the error indicator, based on a ZZ-type error indicator. Then, we present the nodes selected based on a fraction of the nodes of highest error. 10\%  of the nodes showing the largest computed error indicator are marked in red. Finally, we show all the nodes marked for local refinement based on the method described above (i.e. the nodes of highest error and their corresponding stencil nodes). We see that the boundaries of the zones marked for local refinement are smooth and correspond to the zones of highest computed error indicator.

 \clearpage

\subsubsection{Placement of new nodes}\label{SecNodePacement}

We presented in Section \ref{RefinementIdSec} how we select the zones for local refinement. Once these zones are identified, we add new nodes to the domain with the aim of placing them as far as possible from existing nodes to avoid ill-conditioning of the refined discretization. The refinement process consists in three steps:
\begin{enumerate}
	\setlength{\itemsep}{0pt}
	      \setlength{\parskip}{0pt}
	\item Refinement of the boundaries of the domain;
	\item Refinement of the interior of the domain;
	\item Deletion of the nodes which are too close to other nodes.
\end{enumerate}
\paragraph{Step 1: Refinement of the boundaries of the domain} \

The boundaries of the domain are refined first. We use boundary elements, to discretize surfaces for 3D problems and facilitate the implementation of the visibility criterion as described in reference \cite{Jacquemin2021}. Those elements are used as part of the refinement process to facilitate the placement of new nodes far from existing nodes. The approach for 2D and 3D problems is different. We present both approaches below.

In 2D, the new nodes are added in the center of all the elements, which have not already been refined, connected to the boundary nodes marked for local refinement. The new nodes are then projected onto the boundary of the domain using the exact geometry from the CAD file. The refined surface elements are split into two new surface elements. The element normal vectors are then computed.

In 3D, the new nodes are added in the center of the edges of all the elements, which have not already been refined, connected to the boundary nodes marked for local refinement. The new nodes are then projected onto the boundary faces or edges of the domain. The refined surface elements are split into four new surface elements. The element normal vectors are computed.

\begin{figure}[h!]
	\centering
	\begin{tikzpicture}
		\node at (0,0) {\includegraphics[width=16cm]{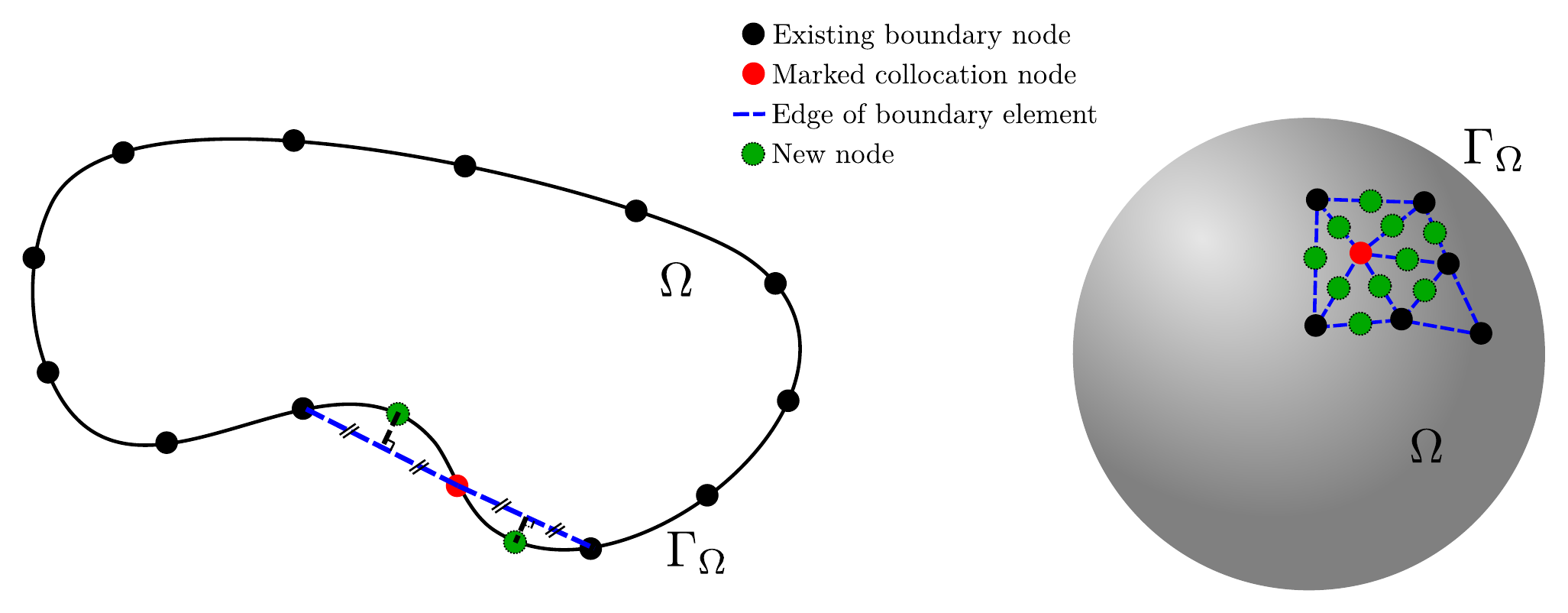}};
	\end{tikzpicture}
	\caption{Computation of the position of new boundary nodes. Collocation nodes marked for local refinement are shown in red color. New boundary nodes, shown in green color bounded by dotted lines, are added in the middle of the edges of the elements connected to the collocation node and projected onto the surface or edge of the domain. The obtained nodes are part of the refined discretization.}
	\label{NewBCNode}
\end{figure}

Figure \ref{NewBCNode} shows a domain $\Omega$ with boundary $\Gamma_\Omega$ for the case of a 2D problem (left) and a 3D problem (right). A collocation node marked for local refinement is shown in red color. The edges of surface elements connected to the marked collocation node are shown by dashed lines in blue color. New boundary nodes are added in the vicinity of the marked collocation node and are shown in green color bounded by dotted lines. The nodes obtained are part of the refined discretization.

\begin{figure}[h!]
	\centering
	\begin{tikzpicture}
		\node at (0,0) {\includegraphics[width=16cm]{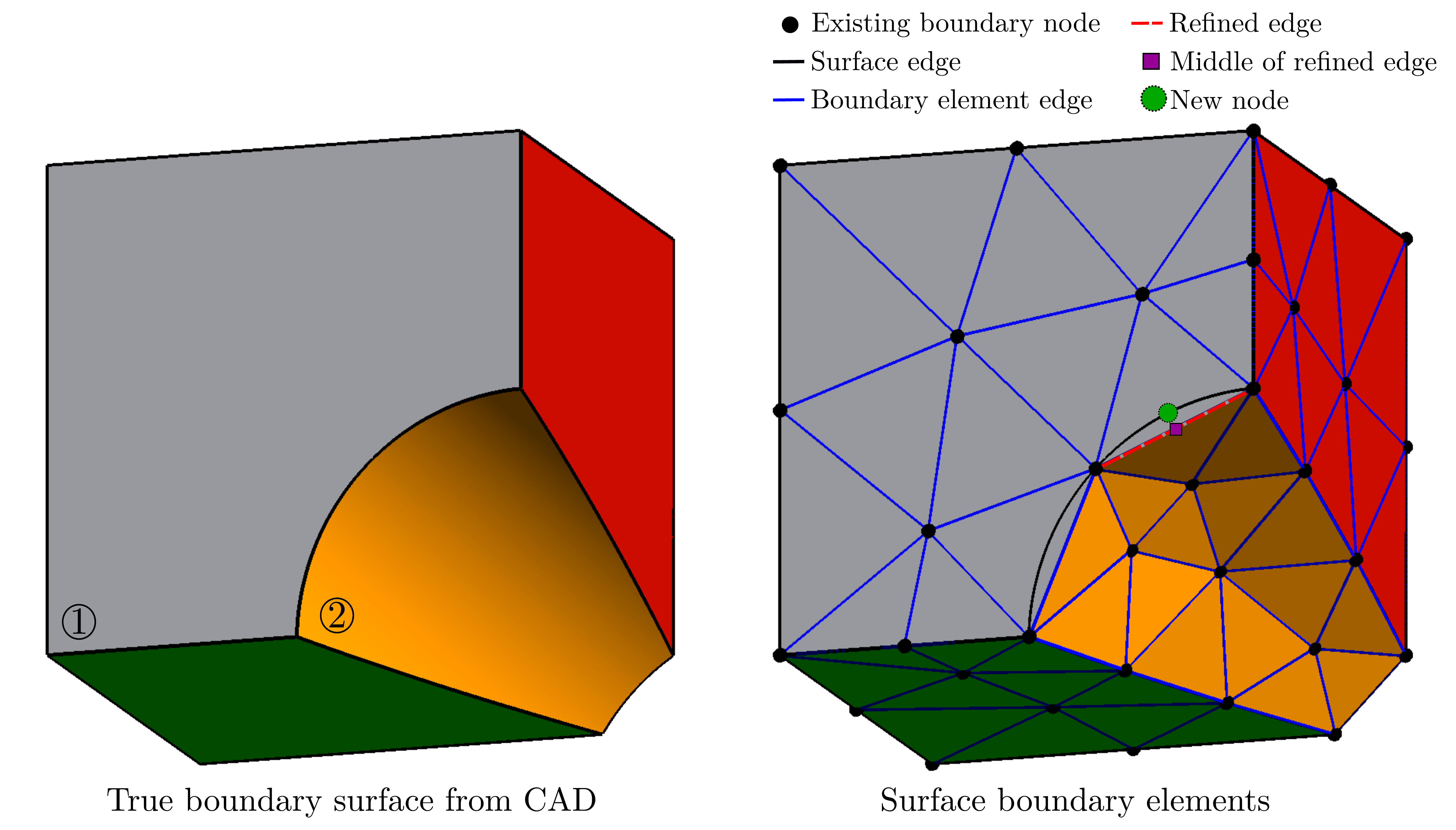}};
	\end{tikzpicture}
	\caption{3D edge refinement. The true faces of the CAD geometry are shown on the left. The faces are discretized using boundary elements. The edge  shown in the red dotted line is located at the intersection between two discretized surfaces. This edge is refined. The middle point of the edge, shown as a purple square, should be projected onto the edge at the intersection between the gray and orange surfaces, marked \textcircled{\raisebox{-.9pt} {1}} and \textcircled{\raisebox{-.9pt} {2}} respectively since a projection onto one of these surfaces would lead to a node outside of the domain.}
	\label{EdgeProjection}
\end{figure}

In 2D, the projection of new nodes onto the boundary of the domain consists only in a projection of the new node onto the edge of the domain. In 3D, the new boundary nodes, not located at the intersection between two faces, are projected onto the surface, parent to the boundary element. If some edges of the boundary elements are located at the intersection of faces of the domain, the new nodes placed in the middle of those edges are projected onto the intersection of faces. Those intersections often correspond to edges of the domain. Therefore, such a projection allows for an accurate refinement of the edges of the domain. This projection is illustrated by Figure \ref{EdgeProjection}. The new node, in the middle of the refined edge, shown in red dotted line should be projected onto the edge at the intersection between the gray and orange surfaces, marked \textcircled{\raisebox{-.9pt} {1}} and \textcircled{\raisebox{-.9pt} {2}} respectively. A projection on either surface would lead to a node outside of the domain.

To facilitate the projection process, the reference of the CAD topological entities (edge or face) should be associated with each node of the smart cloud. For 3D problems, the reference of the parent edge should also be associated with the node of the smart cloud node located at the intersection between two surfaces.

\paragraph{Step 2: Refinement of the interior of the domain} \

\begin{figure}[h!]
	\centering
	\begin{tikzpicture}
		\def\svgwidth{7cm}
		\node at (0,0) {\includegraphics{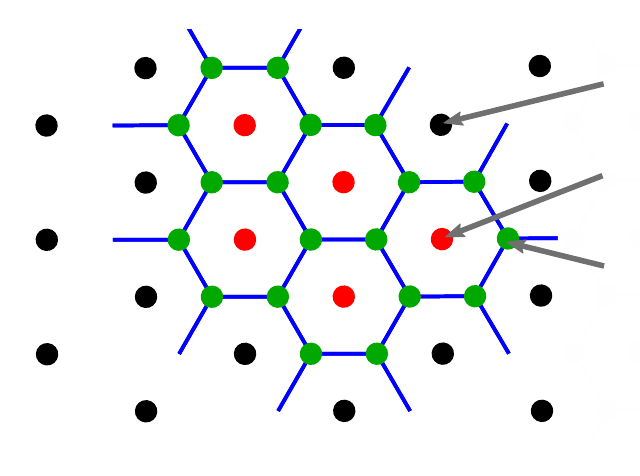}};
		\node[color=black] at (3,1.6) [right] {Collocation node};
		\node[color=black] at (3,0.7) [right] {Collocation node marked for local refinement};
		\node[color=black] at (3,-0.2) [right] {Voronoi corner selected for the refined discretization};
	\end{tikzpicture}
	\caption{Computed Voronoi diagram for a set of collocation nodes - Use for new node selection.}
	\label{NewNodesVoro}
\end{figure}
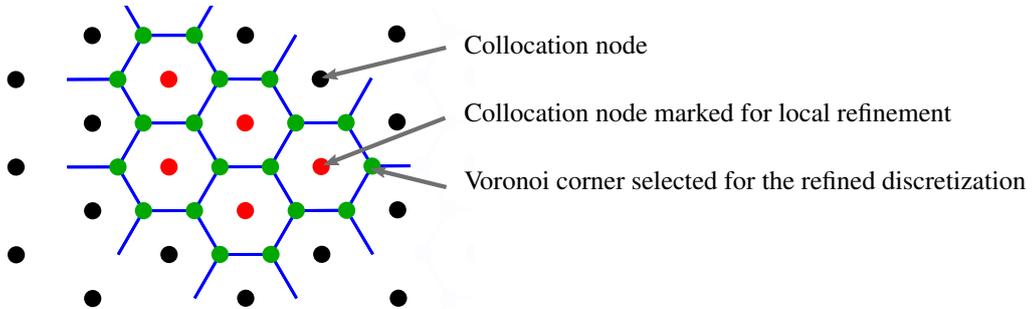

The zones around the interior nodes marked for local refinement are considered once the new boundary nodes have been added to the refined discretization. Voronoi diagrams are used to place the new nodes as far as possible from existing nodes and new boundary nodes. In two dimensions, all the nodes located on the edges of a Vornoi cell are located at equal distance from two adjacent nodes. The corners of a Voronoi cell are located at equal distance from three adjacent nodes as shown in Figure \ref{NewNodesVoro}.

Voronoi diagrams are computed at all the collocation nodes marked for local refinement. The existing nodes of the discretization and the new boundary nodes are considered in the computation of the Voronoi diagrams. Only the Voronoi cells around the marked collocation nodes are of interest. All the corners of the selected cells are added to a list of new candidate nodes. The duplicated nodes are deleted from this list.

\paragraph{Step 3: Deletion of the nodes too close to another node} \

The final step of the node placement process consists in the deletion of the nodes too close from other nodes. For this, a characteristic length is computed at each node of the initial discretization. It corresponds to the minimum distance between the considered node and its adjacent closest node (also based on the initial discretization). The characteristic length associated with a new node of the domain is the characteristic length of its closest node from the initial discretization. This minimum distance is noted $h_{\text{Loc}}$. The minimum allowed inter-node distance, noted $h_{\text{Min}}$, is the ratio of $h_{\text{Loc}}$ to a factor $\alpha$:

\begin{equation} \label{hminEqn}
	h_{\text{Min}}=\frac{h_{\text{Loc}}}{\alpha}.
\end{equation}

We selected $\alpha=3$ in our work.

\subsubsection{Generation of the updated collocation model}

We presented in Subsection \ref{SecNodePacement} a method to place new nodes on the boundary of the domain and in the domain based on a selection of marked collocation nodes. The generation of the updated collocation model is the next step.

The process is the same as the one described in Subsection \ref{CollocModelGen}. The boundary nodes of the smart cloud have references to their parent edge or surface and also references to the boundary conditions associated to these entities. New boundary nodes can be located at the intersection of multiple surfaces. In this case, we apply first the boundary conditions of the surface associated to stress loading boundary conditions, if any. Then, we set the Dirichlet boundary conditions and finally the stress-free boundary conditions.

\subsubsection{Node Selection Threshold Sensitivity}\label{ThresholdSensSec}

We assessed in this section the impact of the node selection threshold on the convergence of the error using model adaptivity. We considered node selection threshold values between 0.02 and 0.2 based on the trend of the results presented in Figure \ref{ZZSortedErrorRound}. The results are presented in Figure \ref{ThresholdSensitivityEllipticalHole} and Figure \ref{ThresholdSensitivityRoundHole} for the benchmark problems introduced in Section \ref{SubPbConsidered}. We compared the results in terms of ``exact'' $L_2$ relative error norm and in terms of ``indicative'' $L_2$ relative error norm. We also compared the results obtained from the adaptive refinement technique presented in this article to results obtained using global refinement (uniform refinement of the domain).

\begin{figure}[!h]
	\centering
	\subfloat[][]{
		\begin{tikzpicture}[scale=1]
			\begin{axis}[height=8cm,width=8cm, ymode=log, ymin=0.005, ymax=2, xmin=3000, xmax=2E6, xmode=log, legend entries={Global refinement,Adapt. refinement. Threshold=0.02,Adapt. refinement. Threshold=0.05,Adapt. refinement. Threshold=0.10,Adapt. refinement. Threshold=0.20}, legend style={ at={(0.5,-0.2)}, anchor=south west,legend columns=1, cells={anchor=west},  font=\footnotesize, rounded corners=2pt,}, legend pos=north east,xlabel=Number of Nodes, ylabel=Exact $L_2$ Relative Error - $\sigma_{VM}$]
				\addplot+[Tblue,mark=triangle*,mark options={fill=Tblue}]   table [x=X-Glob, y=L2R-VMS-Ex-Glob
						, col sep=comma] {ZZAdaptThresholdSensitivity-EllipticalHole.csv};
				\addplot+[Tred,mark=diamond*,mark options={fill=Tred}]   table [x=X-0.02, y=L2R-VMS-Ex-0.02
						, col sep=comma] {ZZAdaptThresholdSensitivity-EllipticalHole.csv};
				\addplot+[Tgreen,mark=*,mark options={fill=Tgreen}]   table [x=X-0.05, y=L2R-VMS-Ex-0.05
						, col sep=comma] {ZZAdaptThresholdSensitivity-EllipticalHole.csv};
				\addplot+[Torange,mark=square*,mark options={fill=Torange}]   table [x=X-0.10, y=L2R-VMS-Ex-0.10
						, col sep=comma] {ZZAdaptThresholdSensitivity-EllipticalHole.csv};
				\addplot+[Tpurple,mark=+,mark options={fill=Tpurple}]   table [x=X-0.20, y=L2R-VMS-Ex-0.20
						, col sep=comma] {ZZAdaptThresholdSensitivity-EllipticalHole.csv};
				\logLogSlopeTriangle{0.85}{0.1}{0.45}{0.51}{Tblue};
				\logLogSlopeTriangle{0.5}{0.05}{0.15}{2.24}{Tgreen};
			\end{axis}
		\end{tikzpicture}
	}
	\subfloat[][]{
		\begin{tikzpicture}[scale=1]
			\begin{axis}[height=8cm,width=8cm, ymode=log, ymin=0.0005, ymax=0.5, xmin=3000, xmax=2E6, xmode=log, legend entries={Global refinement,Adapt. refinement. Threshold=0.02,Adapt. refinement. Threshold=0.05,Adapt. refinement. Threshold=0.10,Adapt. refinement. Threshold=0.20}, legend style={ at={(0.5,-0.2)}, anchor=south west,legend columns=1, cells={anchor=west},  font=\footnotesize, rounded corners=2pt,}, legend pos=north east,xlabel=Number of Nodes, ylabel=Indicative $L_2$ Relative Error - $\sigma_{VM}$]
				\addplot+[Tblue,mark=triangle*,mark options={fill=Tblue}]   table [x=X-Glob, y=L2R-VMS-Est-Glob
						, col sep=comma] {ZZAdaptThresholdSensitivity-EllipticalHole.csv};
				\addplot+[Tred,mark=diamond*,mark options={fill=Tred}]   table [x=X-0.02, y=L2R-VMS-Est-0.02
						, col sep=comma] {ZZAdaptThresholdSensitivity-EllipticalHole.csv};
				\addplot+[Tgreen,mark=*,mark options={fill=Tgreen}]   table [x=X-0.05, y=L2R-VMS-Est-0.05
						, col sep=comma] {ZZAdaptThresholdSensitivity-EllipticalHole.csv};
				\addplot+[Torange,mark=square*,mark options={fill=Torange}]   table [x=X-0.10, y=L2R-VMS-Est-0.10
						, col sep=comma] {ZZAdaptThresholdSensitivity-EllipticalHole.csv};
				\addplot+[Tpurple,mark=+,mark options={fill=Tpurple}]   table [x=X-0.20, y=L2R-VMS-Est-0.20
						, col sep=comma] {ZZAdaptThresholdSensitivity-EllipticalHole.csv};
				\logLogSlopeTriangle{0.95}{0.1}{0.55}{0.45}{Tblue};
				\logLogSlopeTriangle{0.5}{0.05}{0.10}{3.88}{Tgreen};
			\end{axis}
		\end{tikzpicture}
	}\\
	\caption{Exact and indicative $L_2$ relative error norm for the problem of a plate with an elliptical hole. Results from adaptive refinement based on a ZZ-type error indicator for adaptive threshold ratio between 0.02 and 0.20 are compared results obtained from a global refinement of the domain discretization.}
	\label{ThresholdSensitivityEllipticalHole}

	\centering
	\subfloat[][]{
		\begin{tikzpicture}[scale=1]
			\begin{axis}[height=8cm,width=8cm, ymode=log, ymax=0.04, xmin=3000, xmax=1E6, xmode=log, legend entries={Global refinement,Adapt. refinement. Threshold=0.02,Adapt. refinement. Threshold=0.05,Adapt. refinement. Threshold=0.10,Adapt. refinement. Threshold=0.20}, legend style={ at={(0.5,-0.2)}, anchor=south west,legend columns=1, cells={anchor=west},  font=\footnotesize, rounded corners=2pt,}, legend pos=north east,xlabel=Number of Nodes, ylabel=$L_2$ Relative Error - $\sigma_{VM}$ - Exact]
				\addplot+[Tblue,mark=triangle*,mark options={fill=Tblue}]   table [x=X-Glob, y=L2R-VMS-Ex-Glob
						, col sep=comma] {ZZAdaptThresholdSensitivity-RoundHole.csv};
				\addplot+[Tred,mark=diamond*,mark options={fill=Tred}]   table [x=X-0.02, y=L2R-VMS-Ex-0.02
						, col sep=comma] {ZZAdaptThresholdSensitivity-RoundHole.csv};
				\addplot+[Tgreen,mark=*,mark options={fill=Tgreen}]   table [x=X-0.05, y=L2R-VMS-Ex-0.05
						, col sep=comma] {ZZAdaptThresholdSensitivity-RoundHole.csv};
				\addplot+[Torange,mark=square*,mark options={fill=Torange}]   table [x=X-0.10, y=L2R-VMS-Ex-0.10
						, col sep=comma] {ZZAdaptThresholdSensitivity-RoundHole.csv};
				\addplot+[Tpurple,mark=+,mark options={fill=Tpurple}]   table [x=X-0.20, y=L2R-VMS-Ex-0.20
						, col sep=comma] {ZZAdaptThresholdSensitivity-RoundHole.csv};
				\logLogSlopeTriangle{0.9}{0.1}{0.35}{0.95}{Tblue};
				\logLogSlopeTriangle{0.25}{0.1}{0.3}{1.61}{Tgreen};
			\end{axis}
		\end{tikzpicture}
	}
	\subfloat[][]{
		\begin{tikzpicture}[scale=1]
			\begin{axis}[height=8cm,width=8cm, ymode=log, ymax=0.02, xmin=3000, xmax=1E6, xmode=log, legend entries={Global refinement,Adapt. refinement. Threshold=0.02,Adapt. refinement. Threshold=0.05,Adapt. refinement. Threshold=0.10,Adapt. refinement. Threshold=0.20}, legend style={ at={(0.5,-0.2)}, anchor=south west,legend columns=1, cells={anchor=west},  font=\footnotesize, rounded corners=2pt,}, legend pos=north east,xlabel=Number of Nodes, ylabel=$L_2$ Relative Error - $\sigma_{VM}$ - Indicator]
				\addplot+[Tblue,mark=triangle*,mark options={fill=Tblue}]   table [x=X-Glob, y=L2R-VMS-Est-Glob
						, col sep=comma] {ZZAdaptThresholdSensitivity-RoundHole.csv};
				\addplot+[Tred,mark=diamond*,mark options={fill=Tred}]   table [x=X-0.02, y=L2R-VMS-Est-0.02
						, col sep=comma] {ZZAdaptThresholdSensitivity-RoundHole.csv};
				\addplot+[Tgreen,mark=*,mark options={fill=Tgreen}]   table [x=X-0.05, y=L2R-VMS-Est-0.05
						, col sep=comma] {ZZAdaptThresholdSensitivity-RoundHole.csv};
				\addplot+[Torange,mark=square*,mark options={fill=Torange}]   table [x=X-0.10, y=L2R-VMS-Est-0.10
						, col sep=comma] {ZZAdaptThresholdSensitivity-RoundHole.csv};
				\addplot+[Tpurple,mark=+,mark options={fill=Tpurple}]   table [x=X-0.20, y=L2R-VMS-Est-0.20
						, col sep=comma] {ZZAdaptThresholdSensitivity-RoundHole.csv};
				\logLogSlopeTriangle{0.9}{0.1}{0.35}{1.02}{Tblue};
				\logLogSlopeTriangle{0.25}{0.1}{0.3}{1.72}{Tgreen};
			\end{axis}
		\end{tikzpicture}
	}\\
	\caption{Exact and indicative $L_2$ relative error norm for the problem of a body with a cylindrical hole. Results from adaptive refinement based on a ZZ-type error indicator for adaptive threshold ratio between 0.02 and 0.20 are compared results obtained from a global refinement of the domain discretization.}
	\label{ThresholdSensitivityRoundHole}
\end{figure}

The results show that the convergence rate is the largest for a threshold of 0.02 for the problem of a plate with an elliptical hole. After six iterations of iterative refinement, the error obtained is more than three times as low as those obtained with the largest node density considered for the case of global refinement. Such a solution is obtained with approximately 40 times fewer nodes. For the problem of a body with a cylindrical hole, a threshold of 0.02 leads to a non monotonic error reduction. A node selection threshold of 0.2 leads to results similar to the results obtained using a global refinement strategy. We selected a threshold of 0.05 for the problems solved in the next sections of the article. This value is a good compromise between rapid convergence of the solution and non monotonic error reduction.


\subsubsection{Node relaxation}

The addition of new nodes in the domain leads to regions of relatively high node density and regions of relatively low node density. Node relaxation can be used to smooth-out transitions between these regions. It consists in the application of a Laplace smoothing operator to nodes of the point cloud to obtain a uniform discretization. We used the relaxation method of the library Medusa for this purpose \cite{Medusa}.

Rather than applying node relaxation to the entire domain, we preferred a local approach to smooth-out only the transitions between fine and coarse regions. At each node of the refined domain we select the 30 nearest neighbor nodes for 2D problems. Then, we count the number of nodes in the ``fine'' region of the domain. The nodes marked for adaptive refinement and the new nodes are considered located in the ``fine'' region of the domain. Node relaxation is performed if between 20\% and 80\% of the nearest neighbor nodes belong to the ``fine'' region of the domain.

\begin{figure}[h!]
	\centering
	\begin{tikzpicture}
		\node at (0,0) {\includegraphics[height=5.5cm]{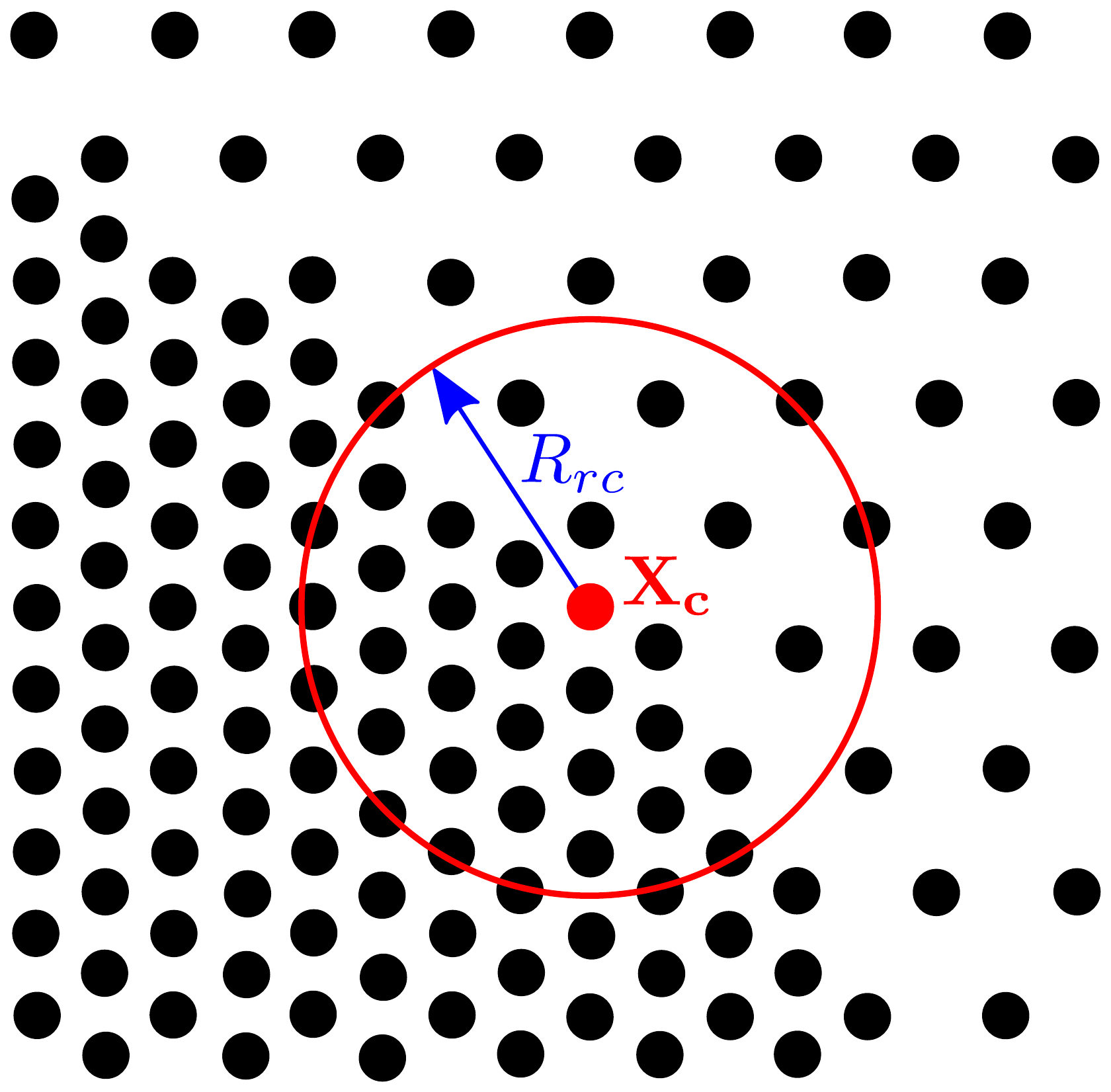}};
	\end{tikzpicture}
	\begin{tikzpicture}
		\node at (0,0) {\includegraphics[height=5.5cm]{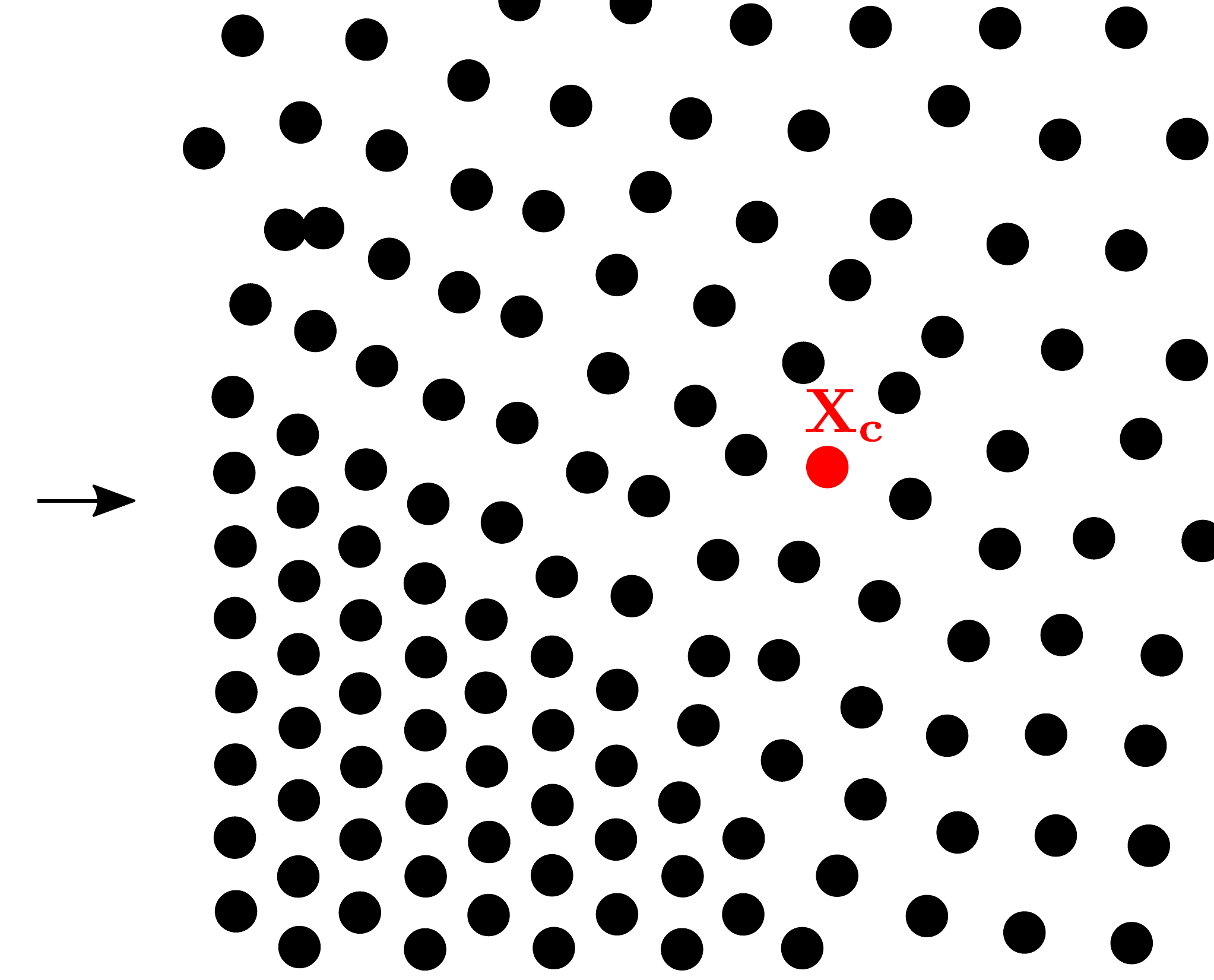}};
	\end{tikzpicture}
	\caption{Point cloud before (left) and after (right) node relaxation. Node relaxation is performed locally around collocation nodes such as $X_c$ over the domain of radius $R_{rc}$. Relaxation is performed when between 20\% and 80\% of the nodes within $R_{rc}$ are located in ``fine'' discretization regions.}
	\label{NodesRelax}
\end{figure}

The radius $R_{rc}$ is the distance between the farthest selected neighbor node and the considered collocation node $\mathbf{X_c}$. All the neighbor nodes located at a distance between $R_{rc}$ and $\frac{2}{3}R_{rc}$ are considered as boundary nodes, and fixed in position, during the relaxation process. The position of the nodes is updated once relaxation is performed. The nodes moved are no longer considered part of the ``fine'' region of the domain. This process is repeated for all the nodes of the discretization. The relaxation process is illustrated by Figure \ref{NodesRelax}.

We show in Figure \ref{RelaxEllipticalHole} and Figure \ref{RelaxRoundHole} the impact of node relaxation on the convergence of the adaptivity scheme and on the error indicator.

We observe from these figures that node relaxation improves the convergence rate of the error for both benchmark problems considered. This result is expected since the node discretization is more uniform. The improvement is however not significant and its usefulness may be considered debatable given the computational effort it necessitates. We did not investigate this approach further in the context of the problems solved in the next section of the article.

\begin{figure}[!h]
	\centering
	\subfloat[][]{
		\begin{tikzpicture}[scale=1]
			\begin{axis}[height=8cm,width=8cm, ymode=log, ymin=0.007, ymax=2, xmin=3000, xmax=2E5, xmode=log, legend entries={Adapt. refinement w/o node relaxation,Adapt. refinement w/ node relaxation}, legend style={ at={(0.5,-0.2)}, anchor=south west,legend columns=1, cells={anchor=west},  font=\footnotesize, rounded corners=2pt,}, legend pos=north east,xlabel=Number of Nodes, ylabel=Exact $L_2$ Relative Error - $\sigma_{VM}$]
				\addplot+[Tblue,mark=triangle*,mark options={fill=Tblue}]   table [x=X-0.05, y=L2R-VMS-Ex-0.05
						, col sep=comma] {ZZAdaptThresholdSensitivity-EllipticalHole.csv};
				\addplot+[Tred,mark=diamond*,mark options={fill=Tred}]   table [x=X-0.05-Relax, y=L2R-VMS-Ex-0.05-Relax
						, col sep=comma] {ZZAdaptThresholdSensitivity-EllipticalHole.csv};
			\end{axis}
		\end{tikzpicture}
	}
	\subfloat[][]{
		\begin{tikzpicture}[scale=1]
			\begin{axis}[height=8cm,width=8cm, ymode=log, ymin=0.0007, ymax=0.2, xmin=3000, xmax=2E5, xmode=log, legend entries={Adapt. refinement w/o node relaxation,Adapt. refinement w/ node relaxation}, legend style={ at={(0.5,-0.2)}, anchor=south west,legend columns=1, cells={anchor=west},  font=\footnotesize, rounded corners=2pt,}, legend pos=north east,xlabel=Number of Nodes, ylabel=Indicative $L_2$ Relative Error - $\sigma_{VM}$]
				\addplot+[Tblue,mark=triangle*,mark options={fill=Tblue}]   table [x=X-0.05, y=L2R-VMS-Est-0.05, col sep=comma] {ZZAdaptThresholdSensitivity-EllipticalHole.csv};
				\addplot+[Tred,mark=diamond*,mark options={fill=Tred}]   table [x=X-0.05-Relax, y=L2R-VMS-Est-0.05-Relax
						, col sep=comma] {ZZAdaptThresholdSensitivity-EllipticalHole.csv};
			\end{axis}
		\end{tikzpicture}
	}\\
	\caption{Comparison of the error with and without node relaxation for the problem of a plate with an elliptical hole. The exact and indicative $L_2$ relative error norms are shown in subfigure (a) and (b), respectively. A node selection threshold ratio of 0.05 is selected. The error obtained with node relaxation is lower than the one obtained without node relaxation. The trend of the results is similar for both cases.}
	\label{RelaxEllipticalHole}
\end{figure}
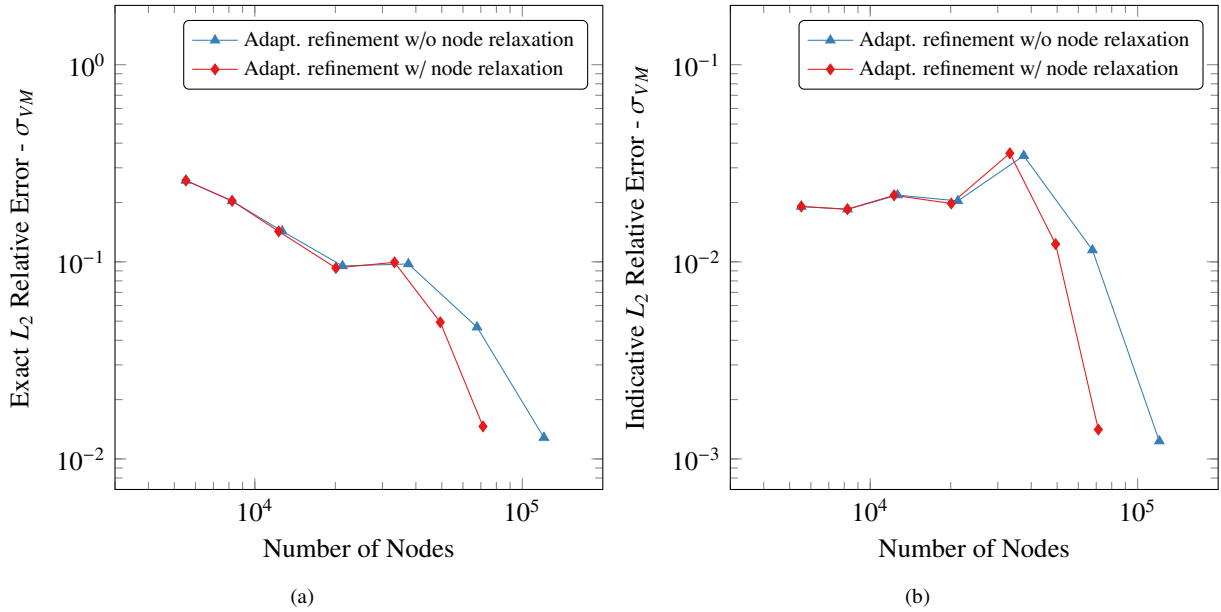
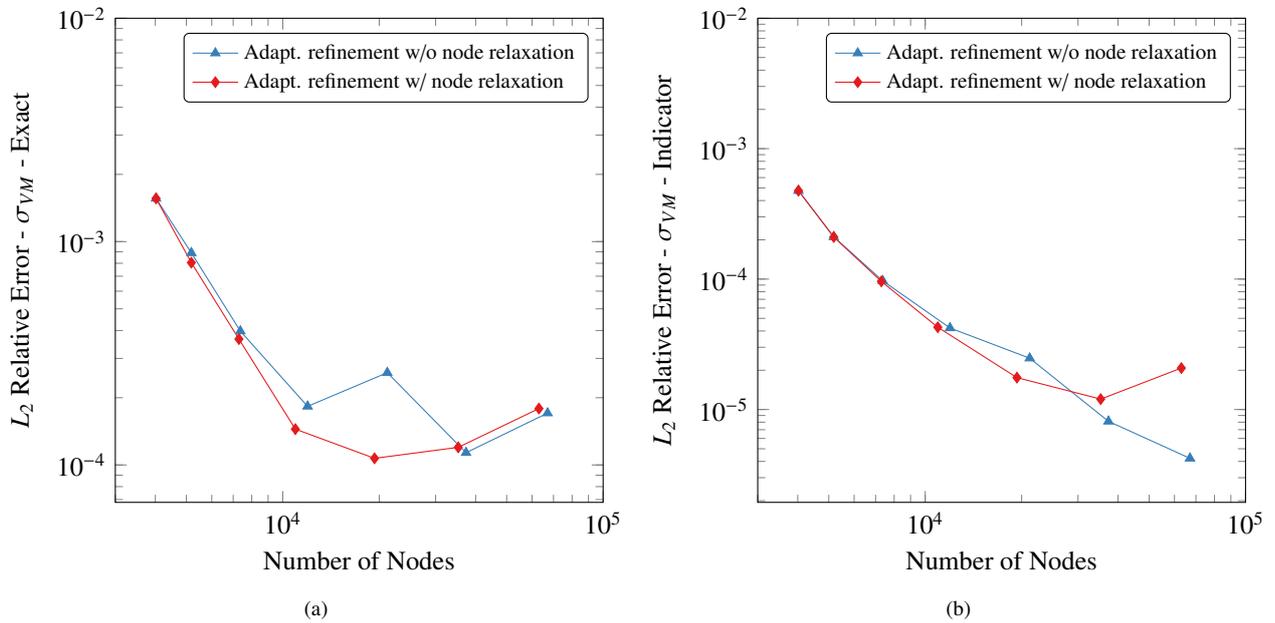
\begin{figure}[!h]
	\centering
	\subfloat[][]{
		\begin{tikzpicture}[scale=1]
			\begin{axis}[height=8cm,width=8cm, ymode=log, ymax=0.01, xmin=3000, xmax=1E5, xmode=log, legend entries={Adapt. refinement w/o node relaxation,Adapt. refinement w/ node relaxation}, legend style={ at={(0.5,-0.2)}, anchor=south west,legend columns=1, cells={anchor=west},  font=\footnotesize, rounded corners=2pt,}, legend pos=north east,xlabel=Number of Nodes, ylabel=$L_2$ Relative Error - $\sigma_{VM}$ - Exact]
				\addplot+[Tblue,mark=triangle*,mark options={fill=Tblue}]   table [x=X-0.05, y=L2R-VMS-Ex-0.05
						, col sep=comma] {ZZAdaptThresholdSensitivity-RoundHole.csv};
				\addplot+[Tred,mark=diamond*,mark options={fill=Tred}]   table [x=X-0.05-Relax, y=L2R-VMS-Ex-0.05-Relax
						, col sep=comma] {ZZAdaptThresholdSensitivity-RoundHole.csv};
			\end{axis}
		\end{tikzpicture}
	}
	\subfloat[][]{
		\begin{tikzpicture}[scale=1]
			\begin{axis}[height=8cm,width=8cm, ymode=log, ymax=0.01, xmin=3000, xmax=1E5, xmode=log, legend entries={Adapt. refinement w/o node relaxation,Adapt. refinement w/ node relaxation}, legend style={ at={(0.5,-0.2)}, anchor=south west,legend columns=1, cells={anchor=west},  font=\footnotesize, rounded corners=2pt,}, legend pos=north east,xlabel=Number of Nodes, ylabel=$L_2$ Relative Error - $\sigma_{VM}$ - Indicator]
				\addplot+[Tblue,mark=triangle*,mark options={fill=Tblue}]   table [x=X-0.05, y=L2R-VMS-Est-0.05
						, col sep=comma] {ZZAdaptThresholdSensitivity-RoundHole.csv};
				\addplot+[Tred,mark=diamond*,mark options={fill=Tred}]   table [x=X-0.05-Relax, y=L2R-VMS-Est-0.05-Relax
						, col sep=comma] {ZZAdaptThresholdSensitivity-RoundHole.csv};
			\end{axis}
		\end{tikzpicture}
	}\\
	\caption{Comparison of the error with and without node relaxation for the problem of a body with a cylindrical hole. The exact and indicative $L_2$ relative error norms are shown in subfigure (a) and (b), respectively. A node selection threshold ratio of 0.05 is selected. The error obtained with node relaxation is lower than the one obtained without node relaxation for the four first iteration steps. The trend of the error indicator is closer to the trend of the exact error for the case of adaptive refinement with node relaxation.}
	\label{RelaxRoundHole}
\end{figure}

\clearpage

\subsection{Practical applications}

To show the applicability of our method to more complicated test cases, we present in this section additional results for practical 2D and 3D problems. The selected problems are:

\begin{itemize}
	\item a gear coupled to a shaft (2D);
	\item a closed cylinder subject to pressure (3D).
\end{itemize}

We generated initial discretizations of the domains, directly from the CAD geometry, using a triangular lattice and a hexagonal close-packed lattice for the 2D and the 3D problems, respectively. We used a threshold of 0.3 to select the nodes close to the boundary of the domain that shall be deleted.

We performed several model adaptivity iterations. We used the ZZ-type error indicator, based on the parameters presented in Section \ref{SubErrorInd}, to determine the areas of the domain where the error is likely to be high. We selected the areas to be refined based on a threshold ratio of 0.05. Finally, we placed new nodes on the boundaries of the domain and in the domain based on the method presented in Section \ref{AdaptSec}.

\subsubsection{Gear coupled to a shaft}

\begin{figure}[!h]
	\centering
	\begin{tikzpicture}
		\def\svgwidth{7cm};
		\node at (0,0) {\includegraphics[width=7cm]{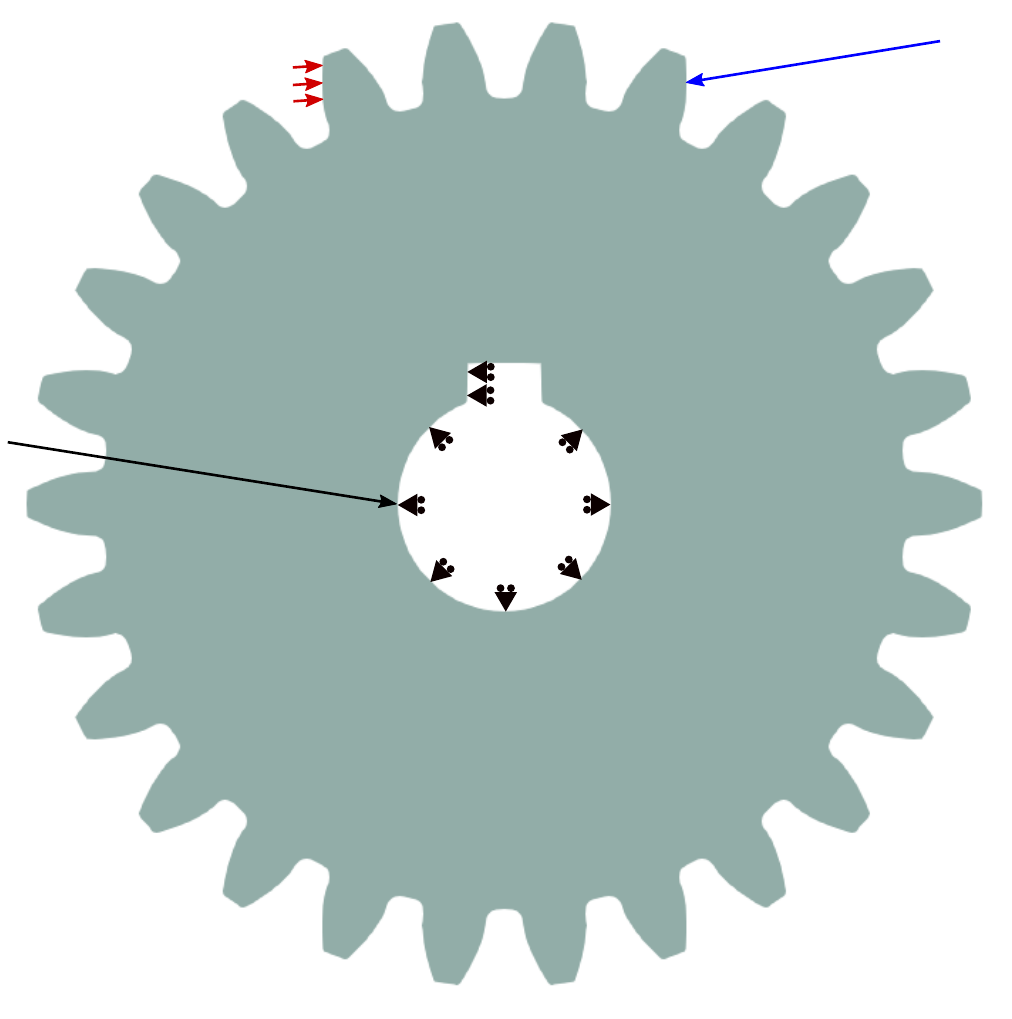}};
		\node[color=red] at (-5,3.8) [right] {Uniform pressure loading};
		\node[color=red] at (-5,3.4) [right] {on one tooth};
		\node[color=blue] at (3.0,3.4) [right] {Homogeneous Neumann};
		\node[color=blue] at (3.0,3.0) [right] {boundary conditions};
		\node[color=black] at (-8.3,0.6) [right] {Sliding boundary conditions};
	\end{tikzpicture}
	\caption{Boundary conditions applied to model a gear coupled to a shaft by a key.}
	\label{2DGearLoading}
\end{figure}
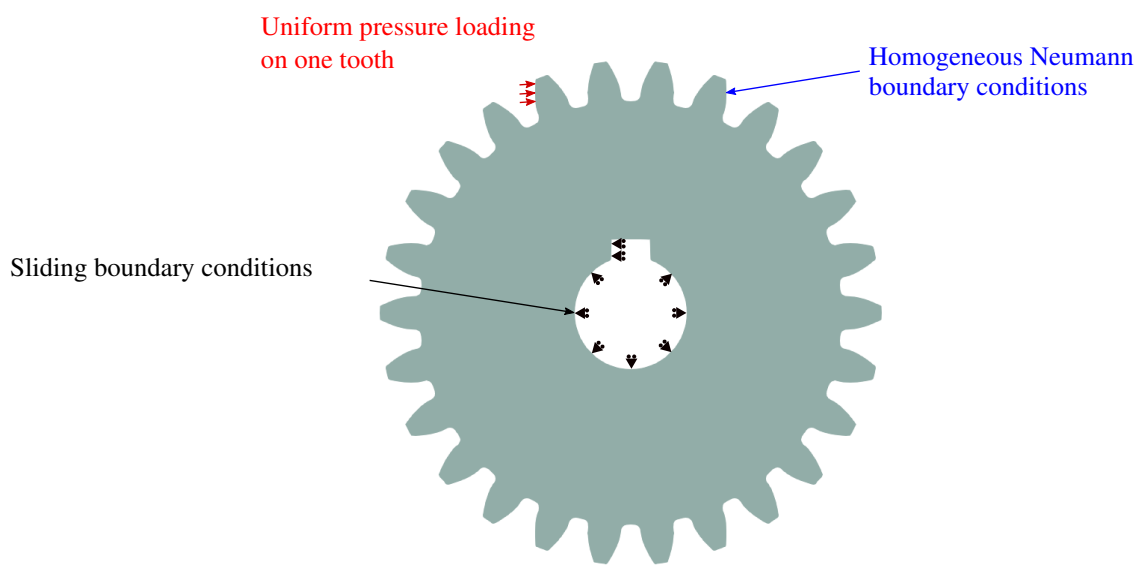

\begin{figure}[!h]
	\centering
	\subfloat[][]{%
		\begin{tabular}{@{}c@{}}
			\includegraphics[width=6.5cm]{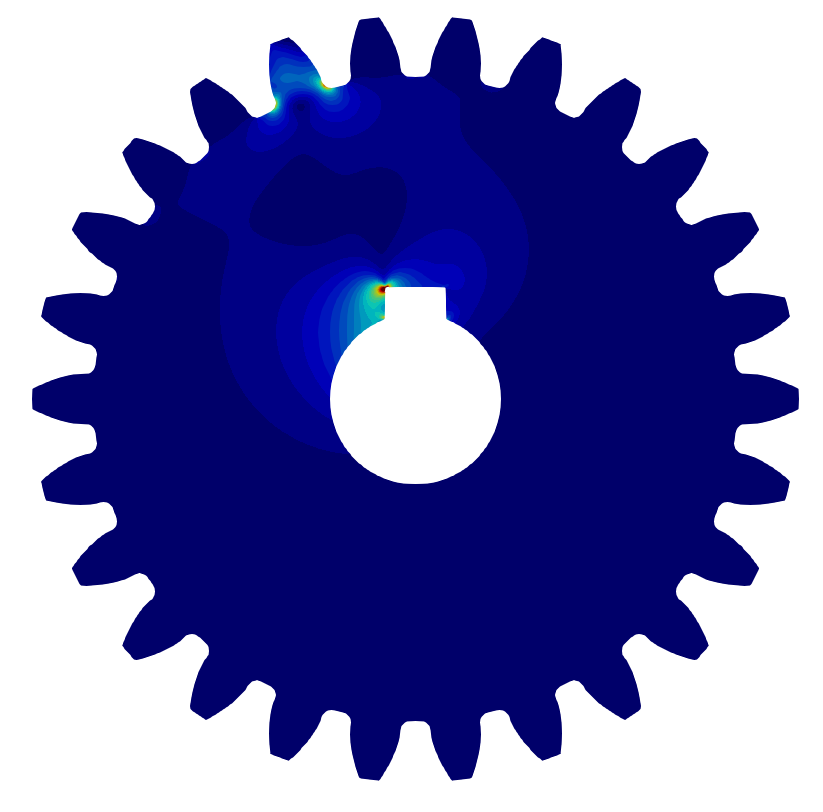}
			\includegraphics[height=6.5cm]{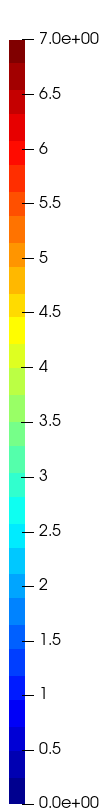} \\
		\end{tabular}
	}
	\subfloat[][]{%
		\begin{tabular}{@{}c@{}}
			\includegraphics[width=6.5cm]{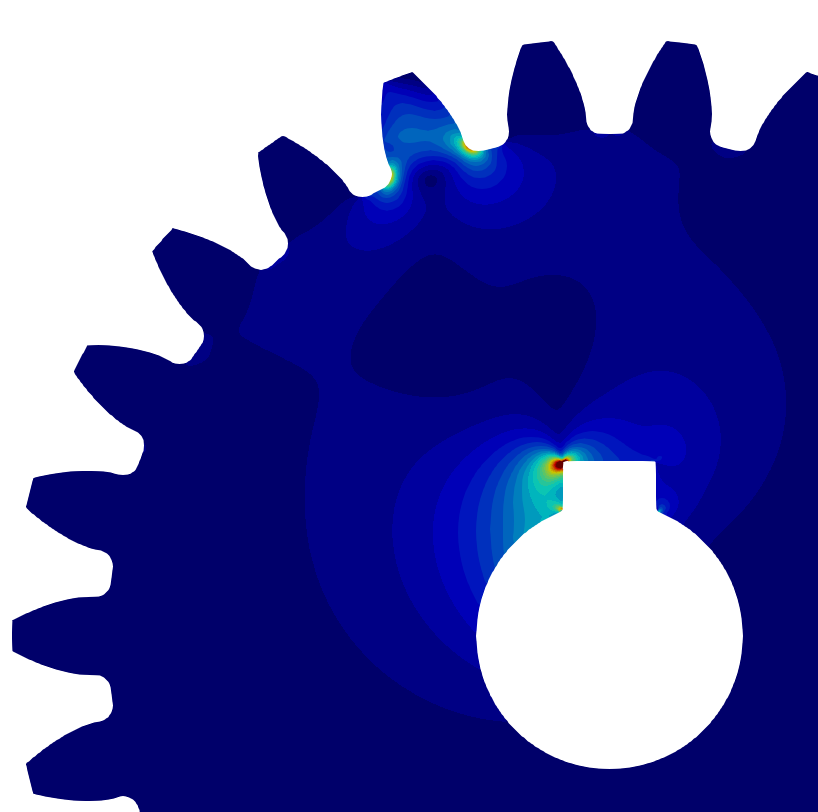}
			\includegraphics[height=6.5cm]{Gear2D_Scale.png} \\
		\end{tabular}
	}\\
	\caption{Gear coupled to a shaft - Solution in terms of von Mises stress obtained from a finite element model composed of 132,665 nodes and 262,193 linear triangular elements. The results are shown for the stress range 0-7 for comparison purposes.}
	\label{2DGearFEASolution}
\end{figure}

We present in Figure \ref{2DGearLoading} the geometry of a gear, coupled to a shaft by a key that we considered. The gear is loaded by uniform pressure on a tooth. We used a finite element solution of the test case as a reference solution which we computed using code\_aster \cite{CodeASTER} for this purpose. The finite element reference solution to this problem is shown in Figure \ref{2DGearFEASolution}. We observed three areas of stress concentration: the top left corner of the groove and both roots of the loaded tooth.

\begin{figure}[h!]
	\centering
	\begin{tabular}{|M{4cm}|c|c|}
		\hline
		                                            & \textbf{Overall view}                                                                           & \textbf{Groove detail view}                                                                     \\
		\hline
		\multirow{2}{*}{\begin{minipage}[t]{0.8\columnwidth} Initial discretization \\ 46,868 nodes \end{minipage}} &                                                                                                 &                                                                                                 \\[-3ex]
		                                            & \includegraphics[width=3.5cm]{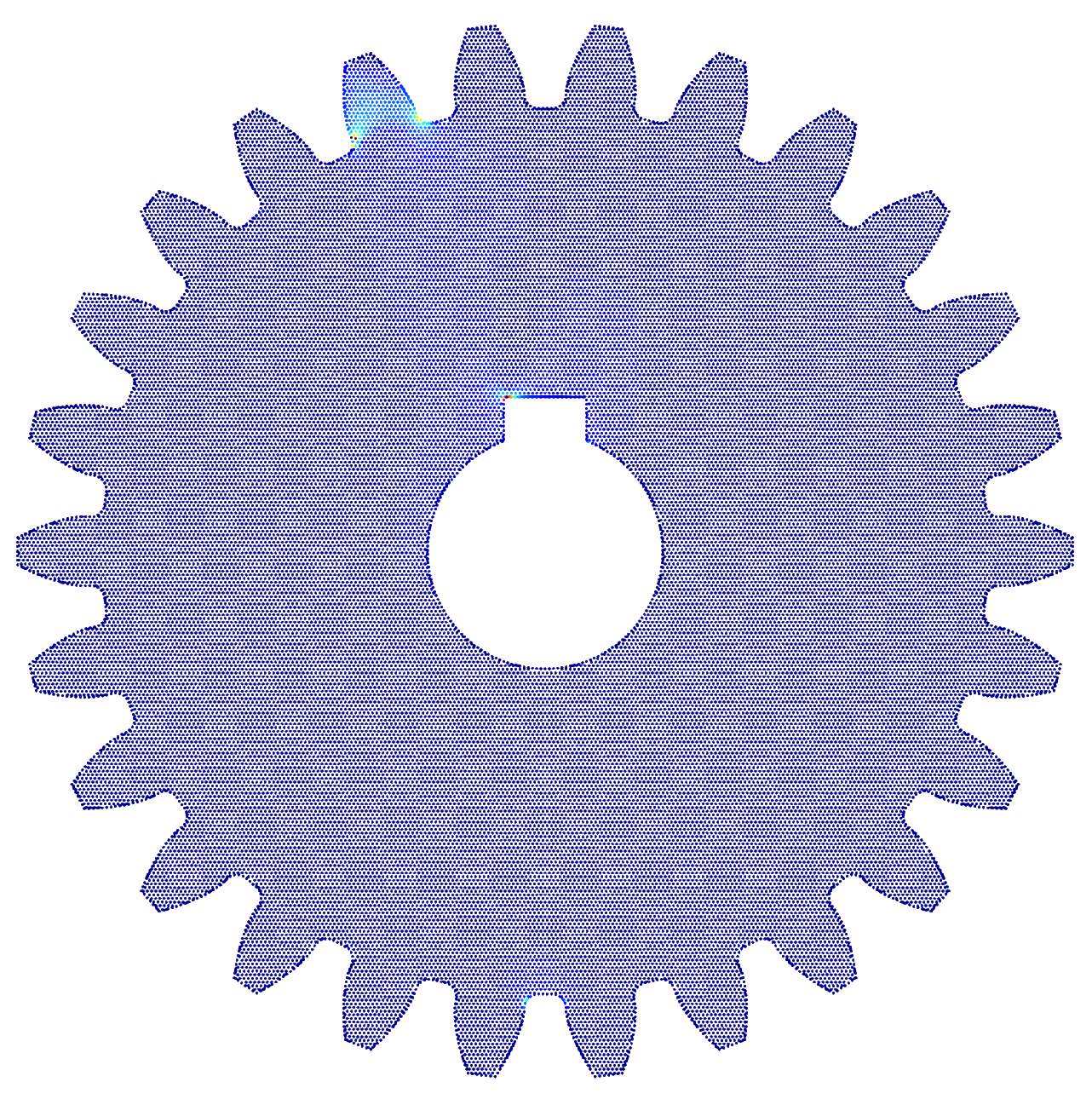} \includegraphics[height=3.5cm]{Gear2D_Scale.png} & \includegraphics[width=3.5cm]{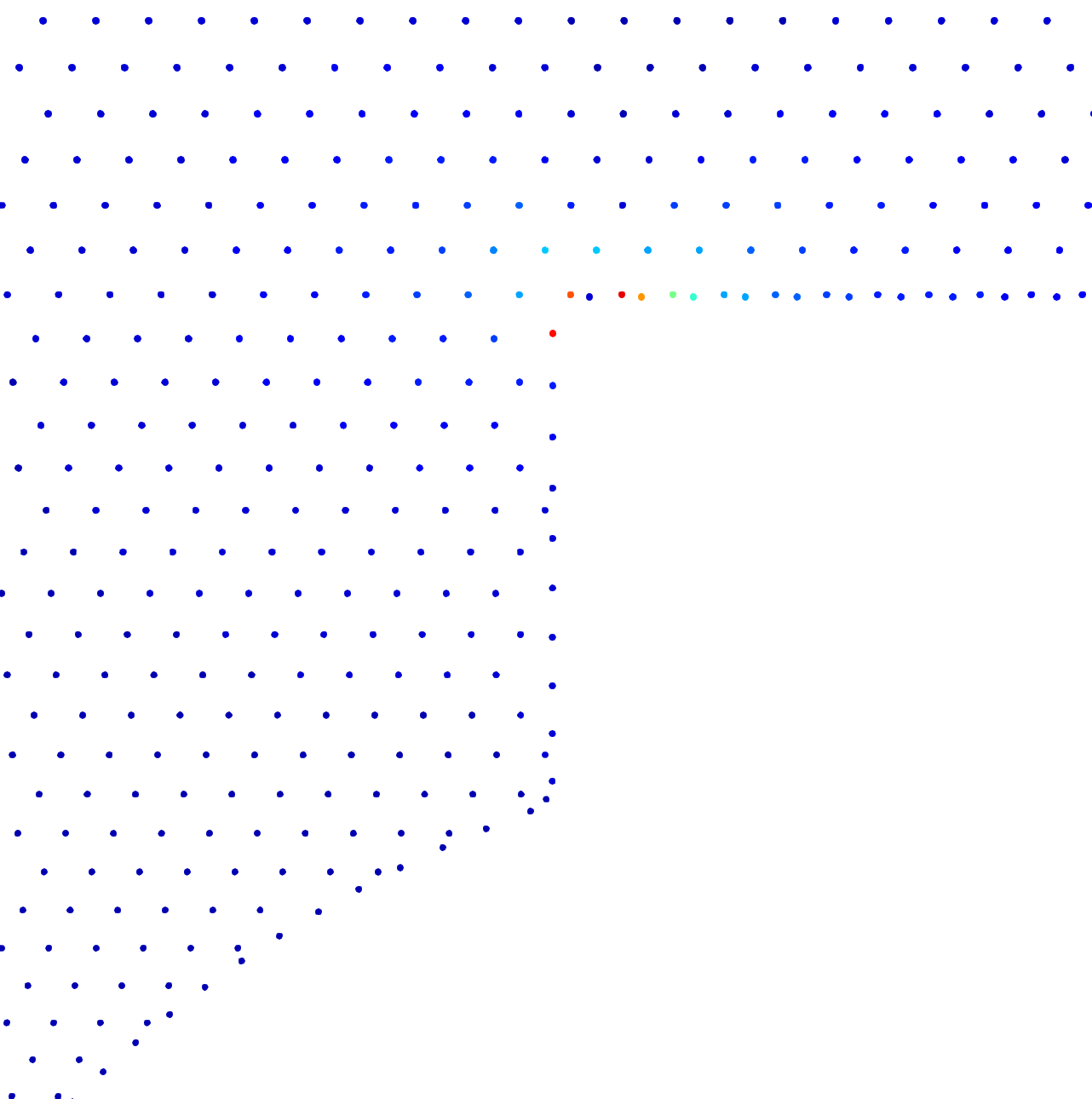} \includegraphics[height=3.5cm]{Gear2D_Scale.png} \\
		\hline
		\multirow{2}{*}{\begin{minipage}[t]{0.8\columnwidth} Refinement iteration \#1 \\ 59,544 nodes \end{minipage}} &                                                                                                 &                                                                                                 \\[-3ex]
		                                            & \includegraphics[width=3.5cm]{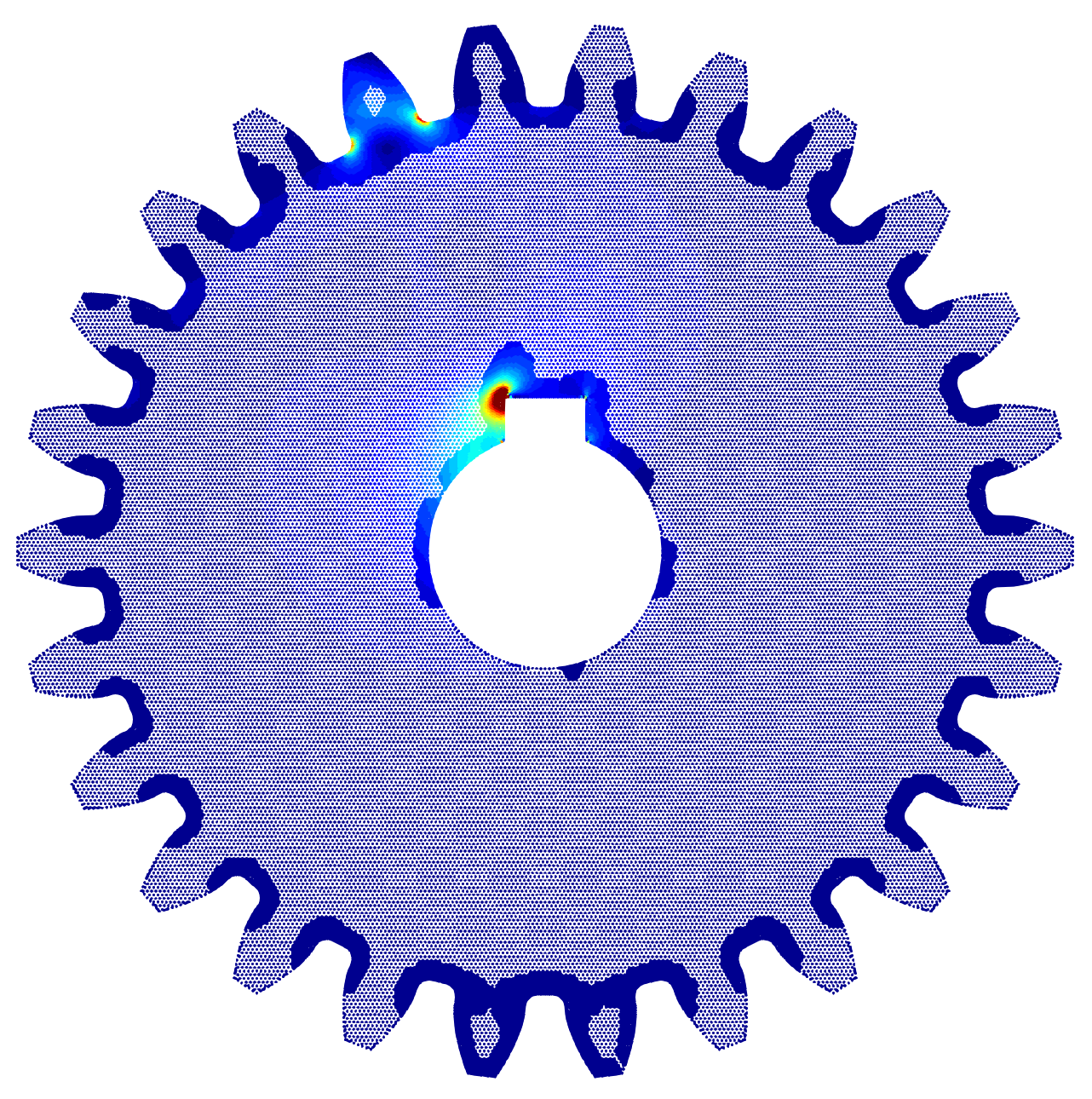} \includegraphics[height=3.5cm]{Gear2D_Scale.png} & \includegraphics[width=3.5cm]{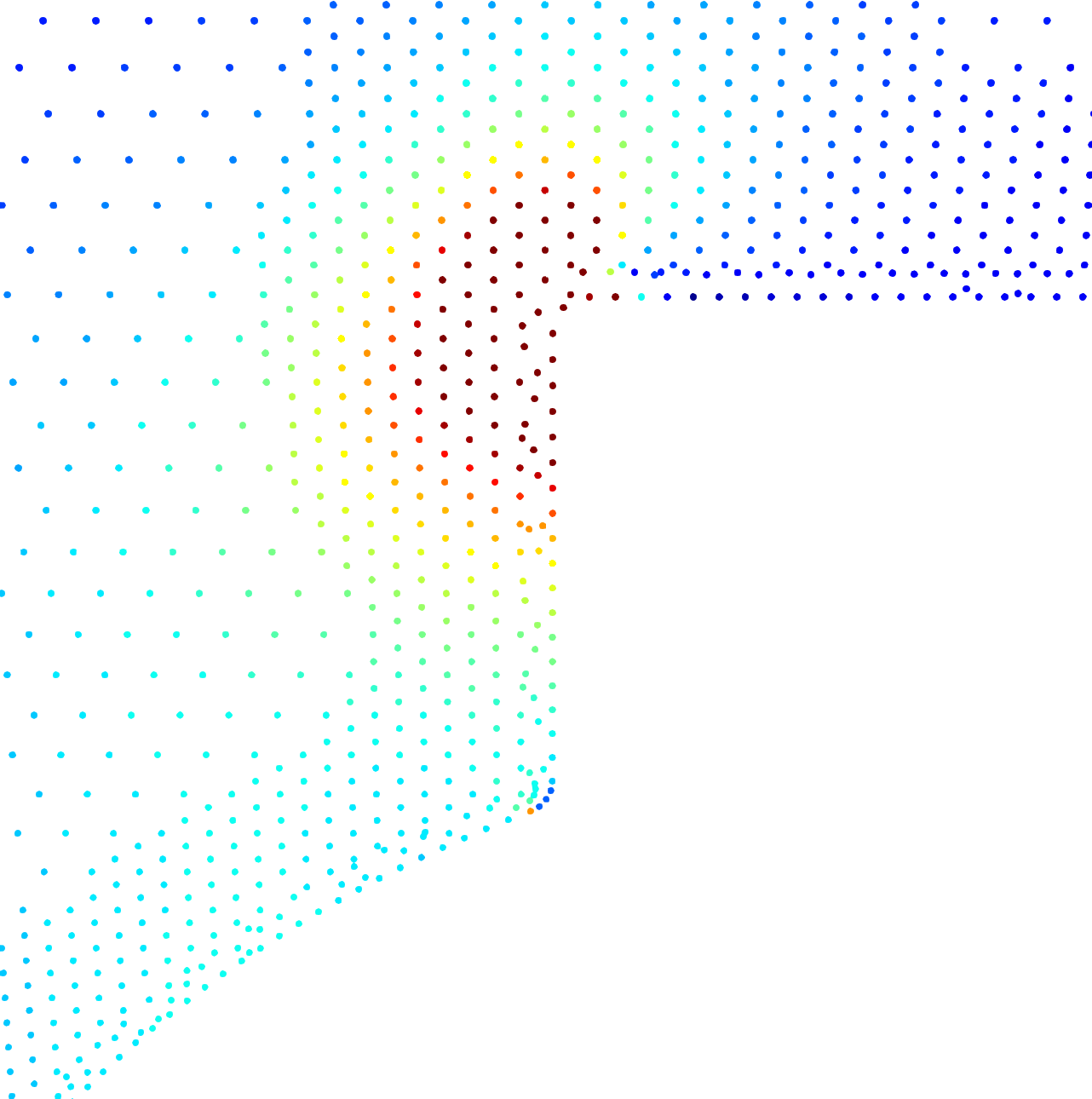} \includegraphics[height=3.5cm]{Gear2D_Scale.png} \\
		\hline
		\multirow{2}{*}{\begin{minipage}[t]{0.8\columnwidth} Refinement iteration \#2 \\ 78,004 nodes \end{minipage}} &                                                                                                 &                                                                                                 \\[-3ex]
		                                            & \includegraphics[width=3.5cm]{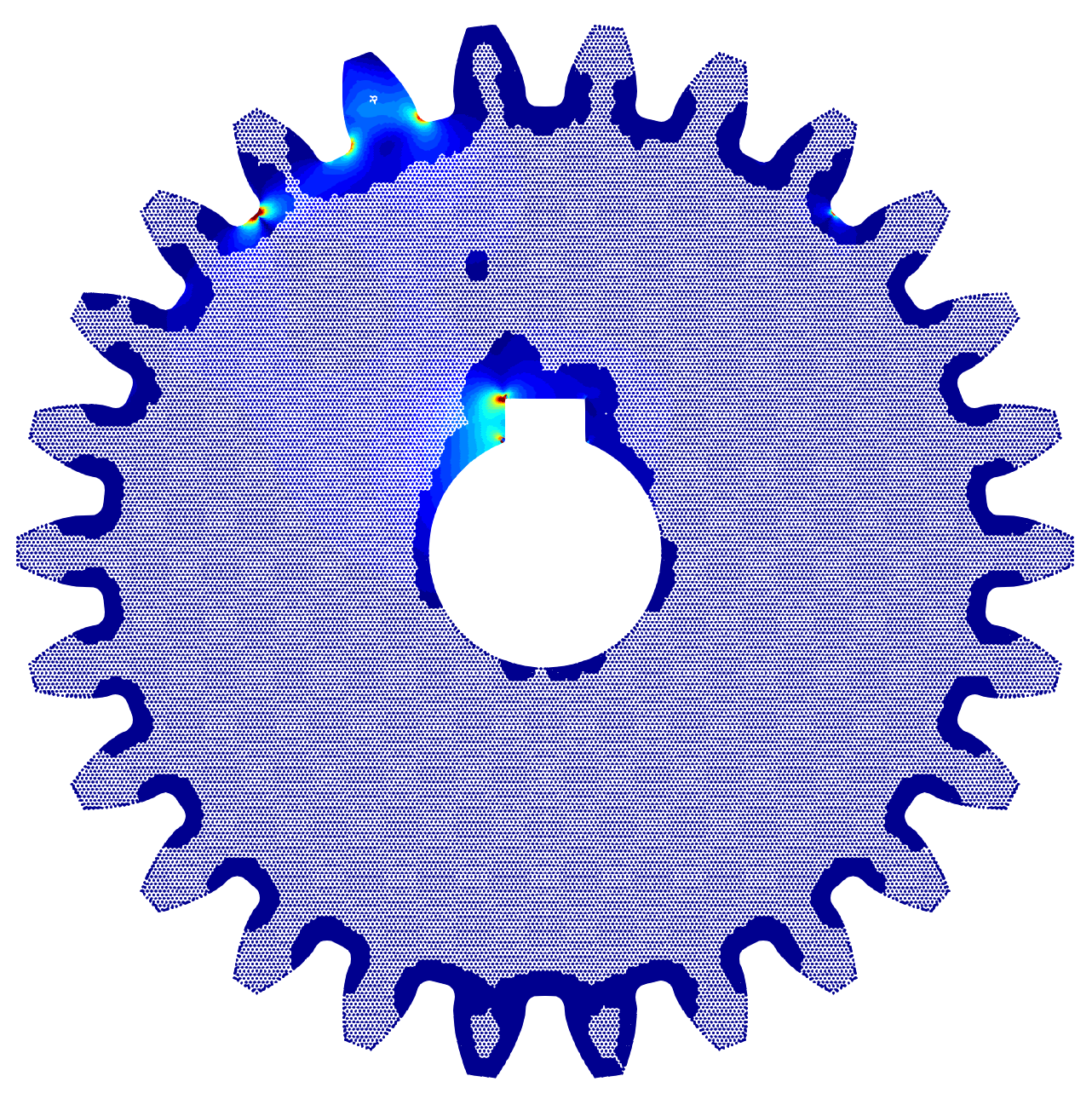} \includegraphics[height=3.5cm]{Gear2D_Scale.png} & \includegraphics[width=3.5cm]{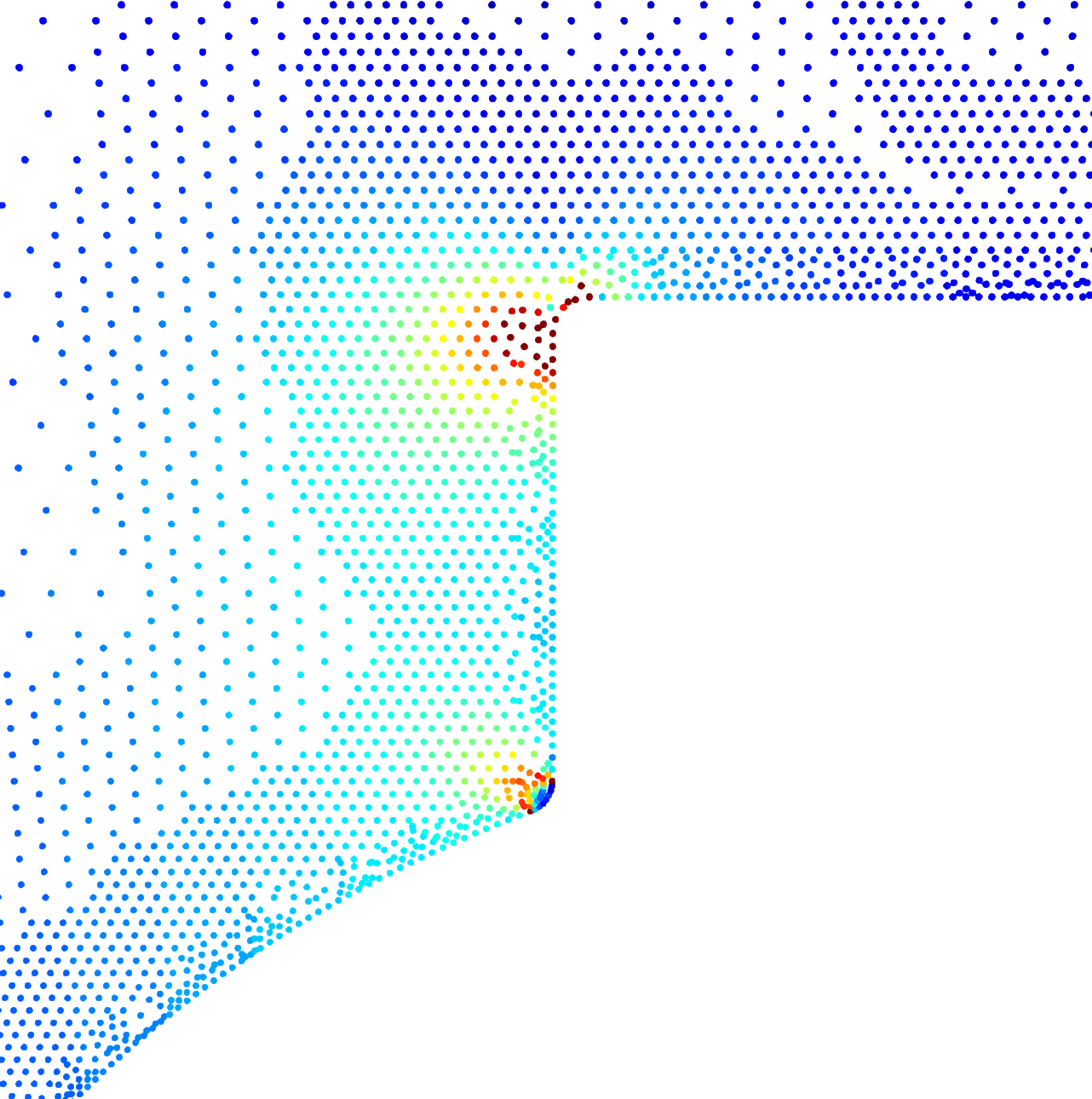} \includegraphics[height=3.5cm]{Gear2D_Scale.png} \\
		\hline
		\multirow{2}{*}{\begin{minipage}[t]{0.8\columnwidth} Refinement iteration \#3 \\ 104,527 nodes \end{minipage}} &                                                                                                 &                                                                                                 \\[-3ex]
		                                            & \includegraphics[width=3.5cm]{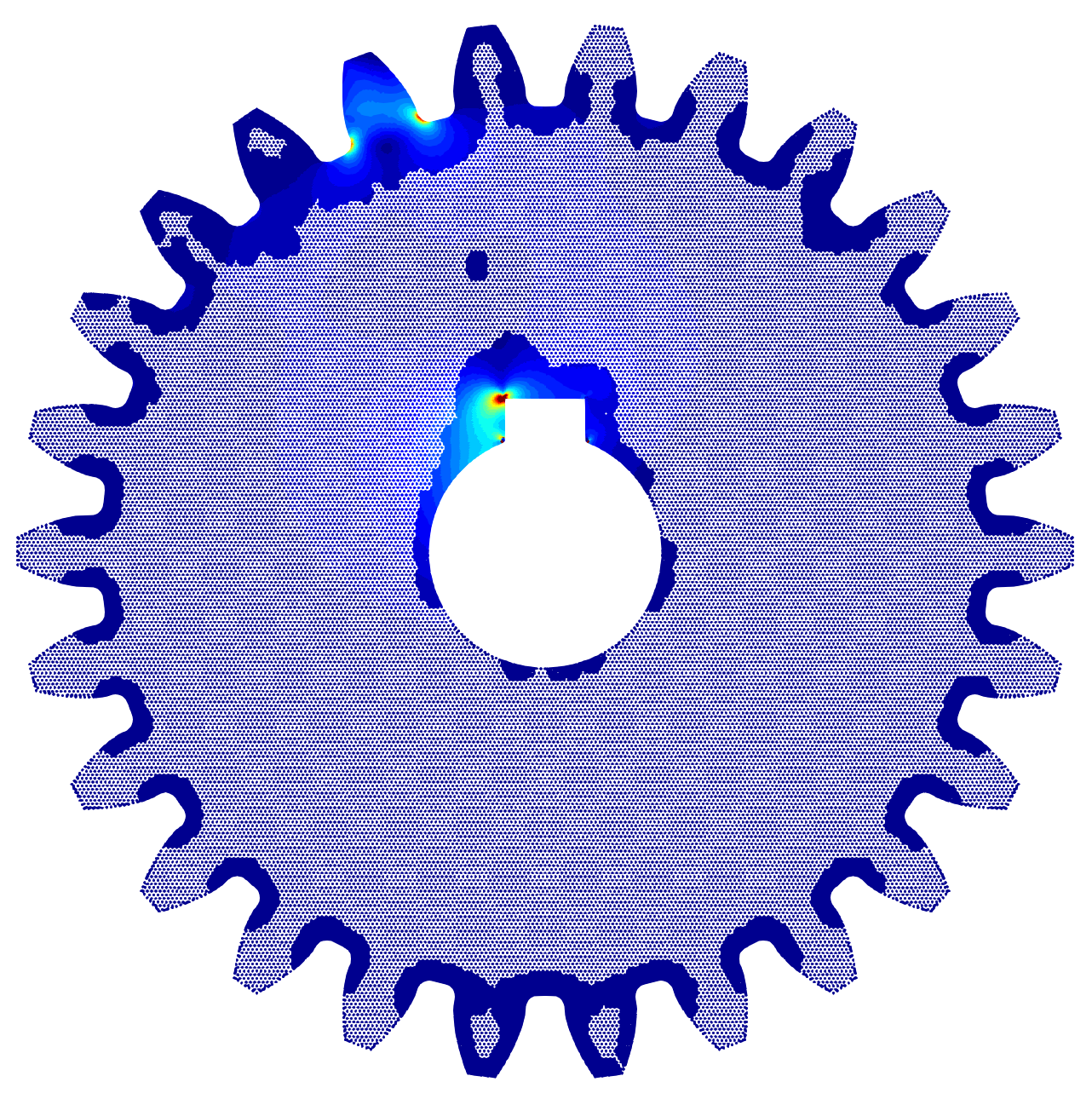} \includegraphics[height=3.5cm]{Gear2D_Scale.png} & \includegraphics[width=3.5cm]{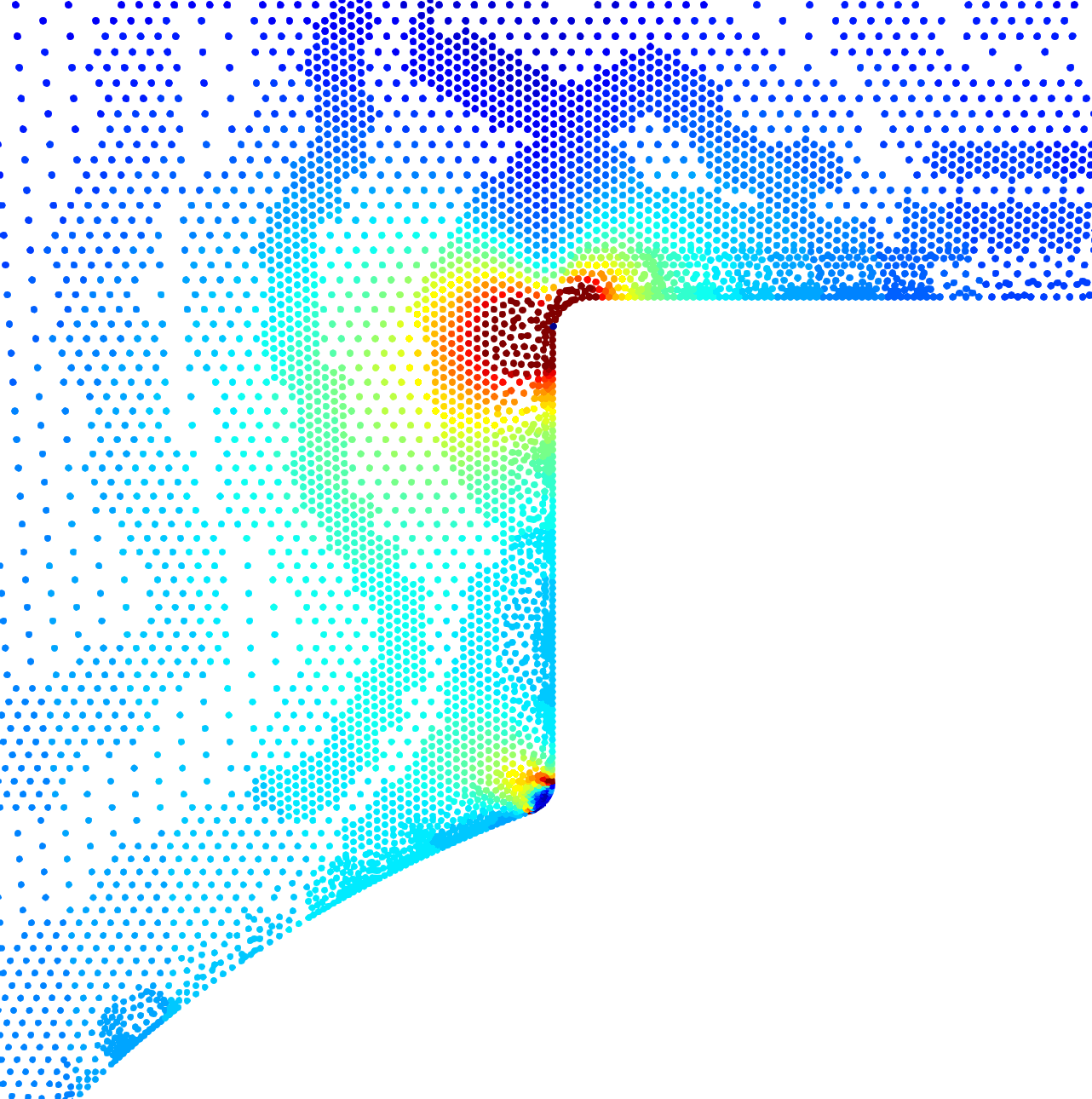} \includegraphics[height=3.5cm]{Gear2D_Scale.png} \\
		\hline
		\multirow{2}{*}{\begin{minipage}[t]{0.8\columnwidth} Refinement iteration \#4 \\ 143,763 nodes \end{minipage}} &                                                                                                 &                                                                                                 \\[-3ex]
		                                            & \includegraphics[width=3.5cm]{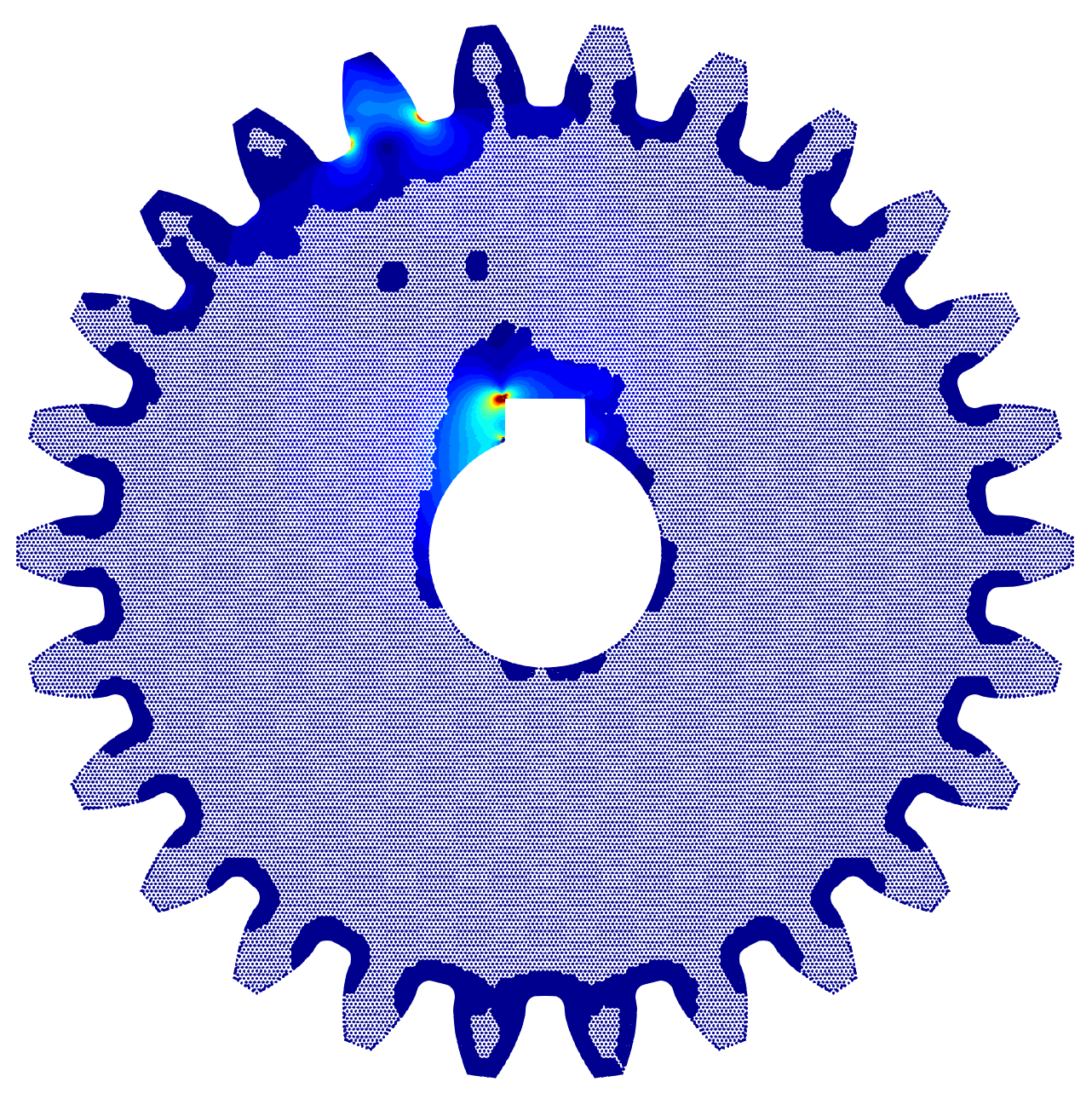} \includegraphics[height=3.5cm]{Gear2D_Scale.png} & \includegraphics[width=3.5cm]{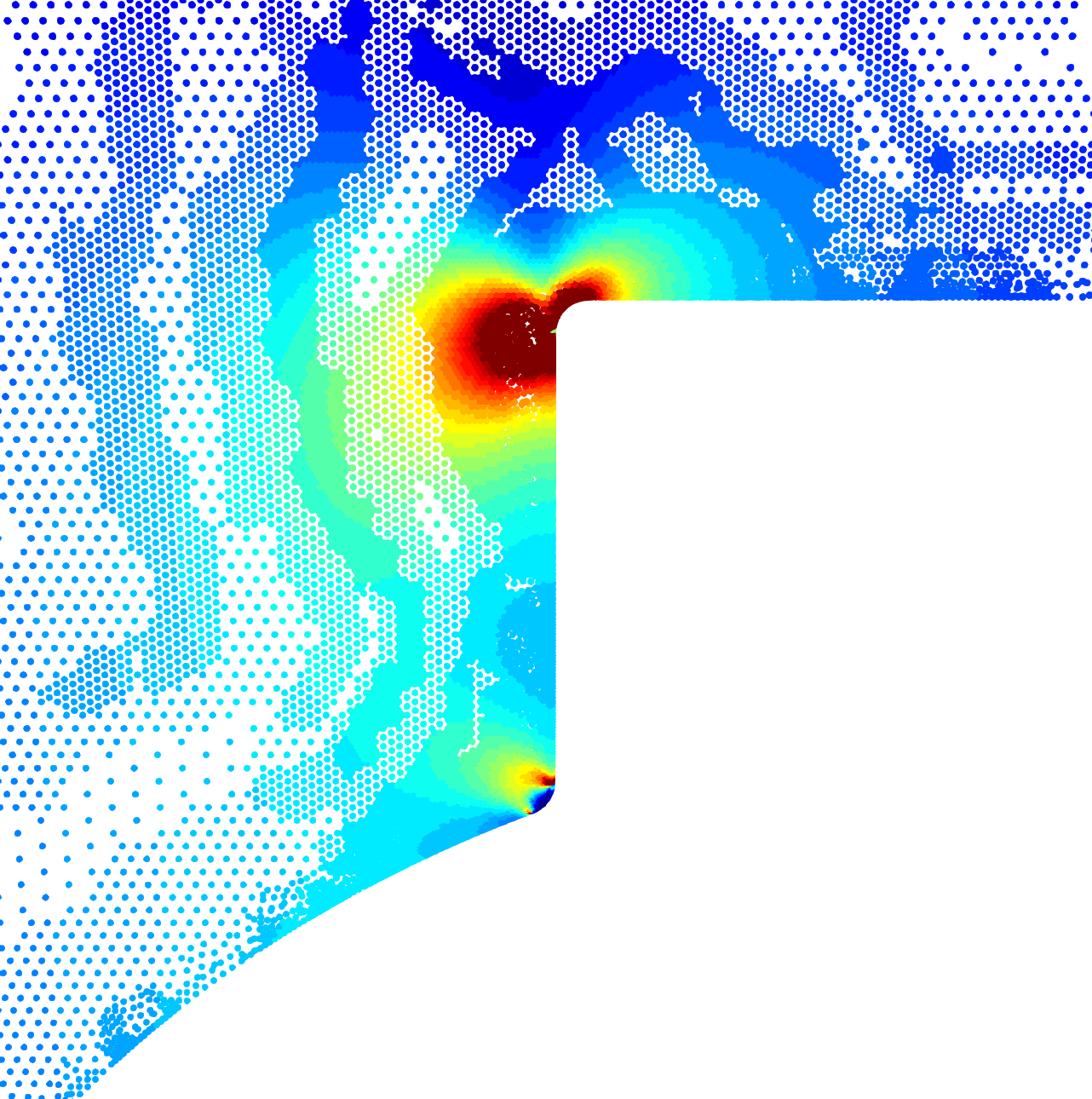} \includegraphics[height=3.5cm]{Gear2D_Scale.png} \\
		\hline
	\end{tabular}

	\caption{Gear coupled to a shaft - Evolution of the discretization of the domain and of the solution in terms of von Mises stress through four iterations of adaptive refinement. The results are shown for the stress range 0-7 for comparison purposes.}
	\label{RefinementStepsGearKey}
\end{figure}

We show in Figure \ref{RefinementStepsGearKey} the evolution of the point cloud and of the results in terms of von Mises stress through four adaptive refinement iterations. We show a general view of the domain and a detailed view of the top left corner of the groove where the stress concentration is the greatest. The results are shown for the stress range 0-7 for comparison purposes.

We observe from this figure that new nodes are placed at the roots of all the teeth at the first refinement iteration. The tooth subject to pressure loading is the tooth being the most refined. The area around the groove is also refined successively at each refinement iteration. We see from the third and fourth iterations that areas at the interface between coarse and fine discretization zones are also being refined. The zones correspond to areas of the domain where the discretization is not uniform leading to larger values of the computed error indicator.

\begin{figure}[h!]
	\centering
	\begin{tabular}{|M{5cm}|c|c|}
		\hline
		                                            & \textbf{Groove}                                                                                 & \textbf{Loaded tooth}                                                                             \\
		\hline
		\multirow{2}{*}{\begin{minipage}[t]{0.8\columnwidth} \underline{Smart cloud collocation solution} \\ $\blacktriangleright$ Refinement iteration \#4 \\ $\blacktriangleright$ 143,763 nodes  \\ $\blacktriangleright$ Approx. inter-node distance: \\ \hspace*{2ex} 0.018 in dense regions \\ \hspace*{2ex} 0.29 in coarse regions \end{minipage}} &                                                                                                 &                                                                                                   \\[-3ex]
		                                            & \includegraphics[width=3.5cm]{KeyView_It4.png} \includegraphics[height=3.5cm]{Gear2D_Scale.png} & \includegraphics[width=3.5cm]{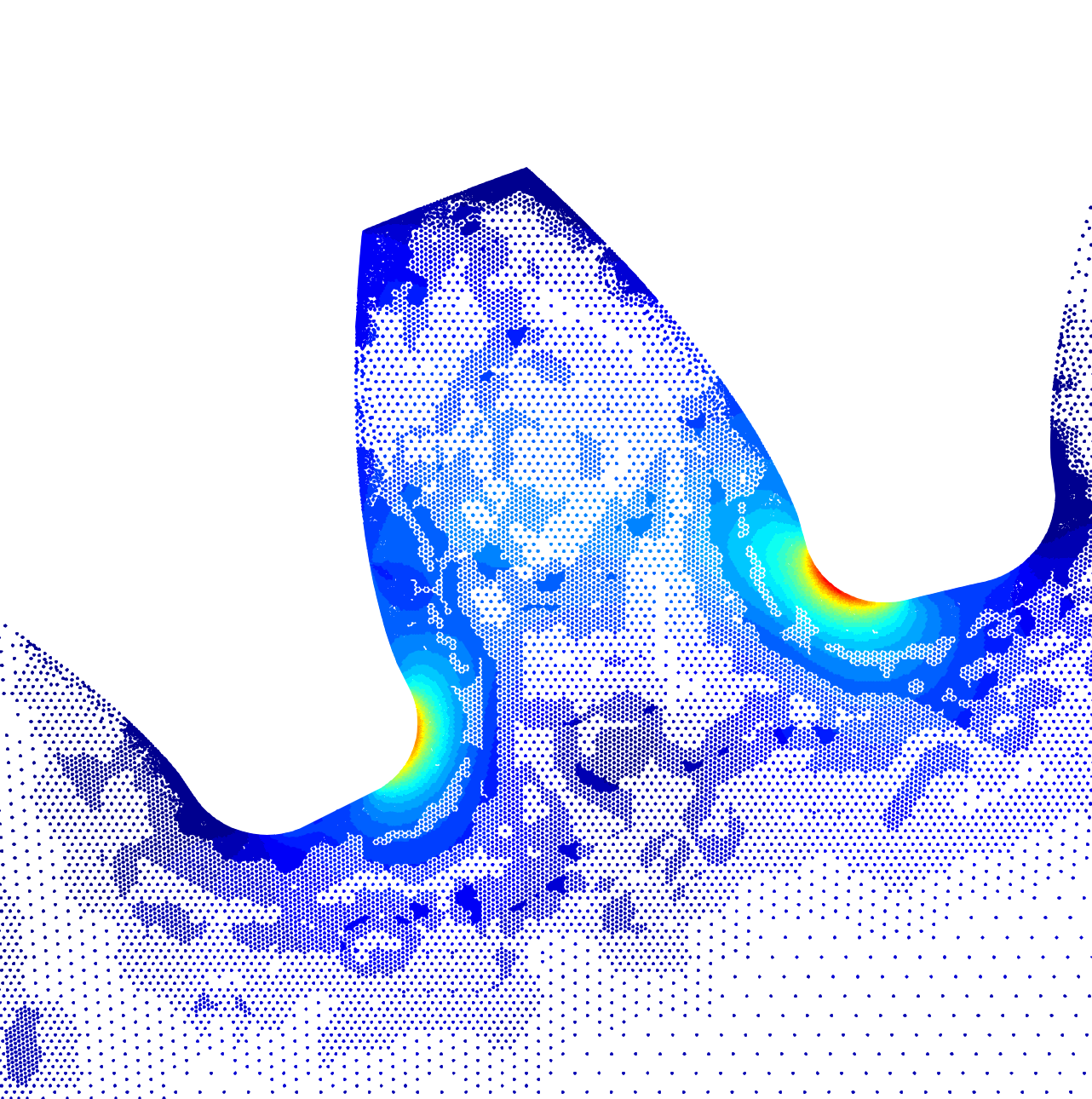} \includegraphics[height=3.5cm]{Gear2D_Scale.png} \\
		\hline
		\multirow{2}{*}{\begin{minipage}[t]{0.8\columnwidth}  \underline{FEA solution} \\ $\blacktriangleright$ 132,665 nodes  \\ $\blacktriangleright$ Approx. inter-node distance:\\ \hspace*{2ex} 0.15 \end{minipage}} &                                                                                                 &                                                                                                   \\[-3ex]
		                                            & \includegraphics[width=3.5cm]{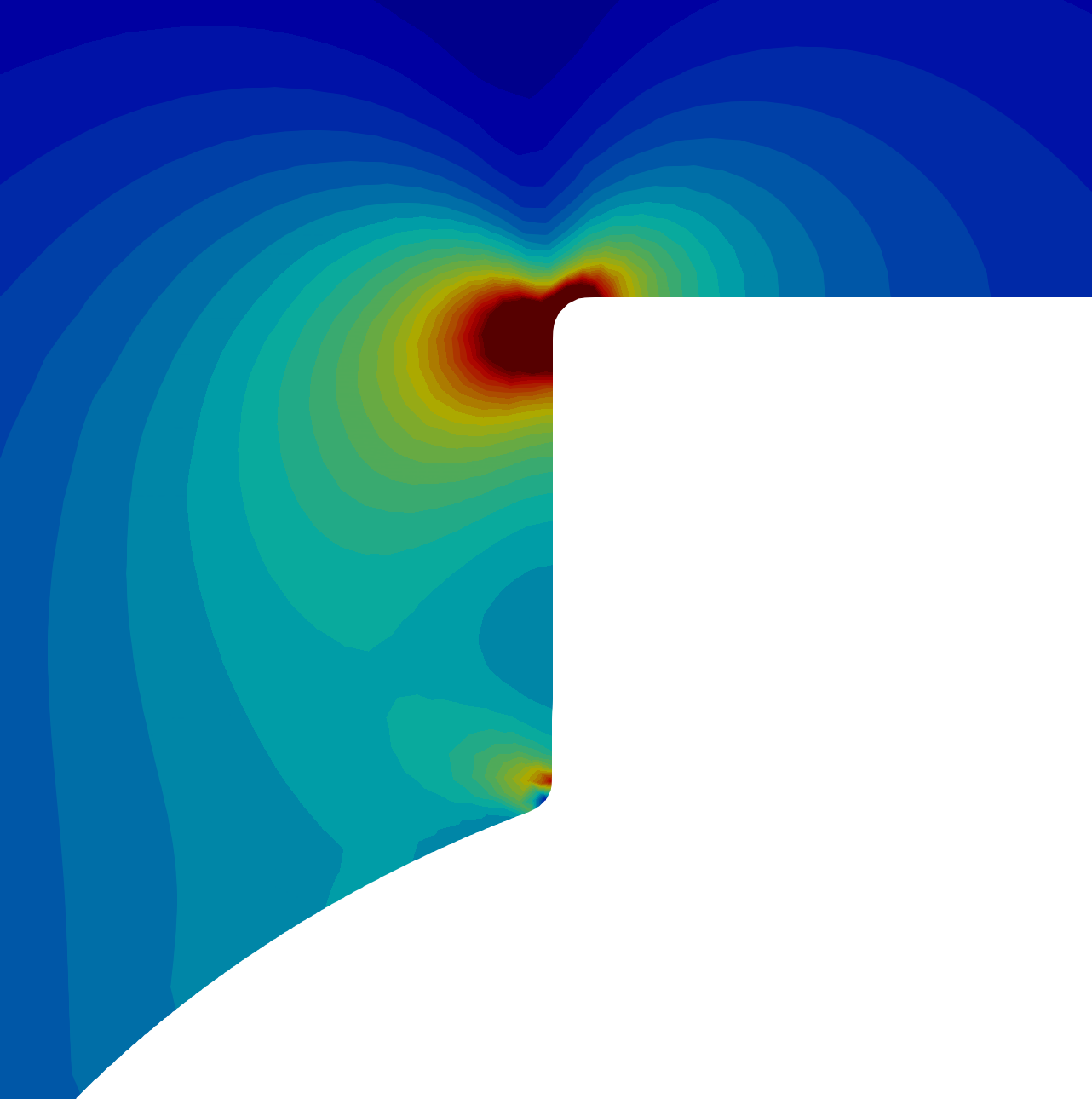} \includegraphics[height=3.5cm]{Gear2D_Scale.png} & \includegraphics[width=3.5cm]{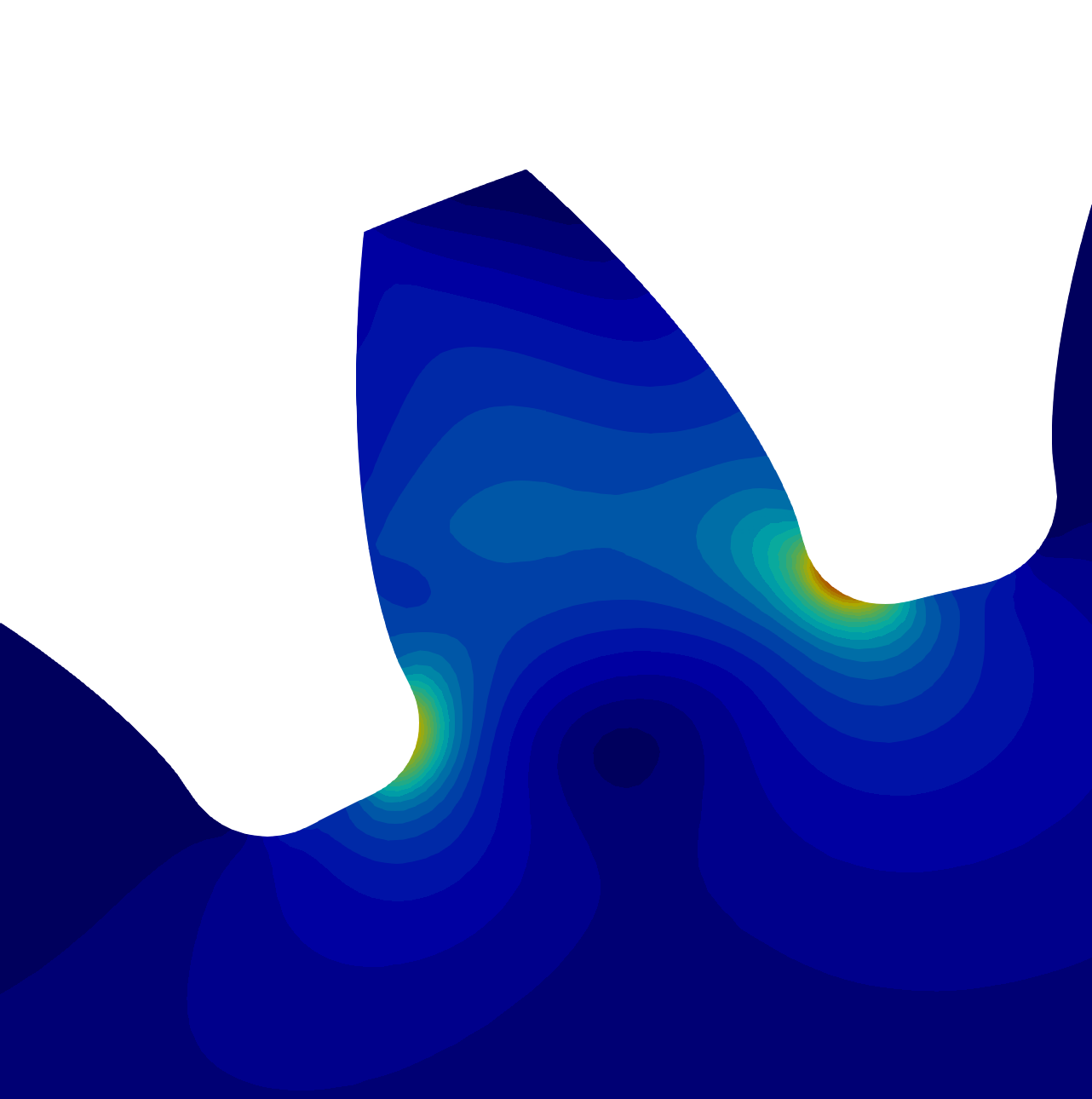} \includegraphics[height=3.5cm]{Gear2D_Scale.png} \\
		\hline
	\end{tabular}

	\caption{Comparison of the results, in terms of von Mises stress, based on the smart cloud adaptive collocation scheme after four adaptive refinement iterations to results from the reference finite element solution. The comparison focuses on the groove and on the loaded tooth which are the areas where the stress concentration is the most important.}
	\label{GearKeyVSFEA}
\end{figure}

The von Mises stress field obtained after four iterations of adaptive refinement is smooth and close to the von Mises stress field of the reference finite element solution presented in Figure \ref{2DGearFEASolution}. We compare both solutions in greater detail in Figure \ref{GearKeyVSFEA}. We focus on the stress field on the left side of the groove and in the loaded tooth. In the area around the groove, we see no difference in the stress field obtained from both methods. At the root of the loaded tooth, we see a slightly larger stress concentration for the  smart cloud collocation solution. The  smart cloud has a much higher node density at the root of the loaded tooth than the finite element discretization. The inter-node distance is approximately 0.018 in the dense region of the  smart cloud. For the finite element discretization, the inter-node distance is approximately 0.15. This tends to explain the larger computed von Mises stress values.

\clearpage

\begin{figure}[!h]
	\centering
	\begin{tikzpicture}[scale=1]
		\begin{axis}[height=7cm,width=11cm, ymode=log, ymin=0.01, ymax=1, xmin=39000, xmode=log, legend entries={ZZ-type indicator}, legend style={ at={(0.5,-0.2)}, anchor=south west,legend columns=1, cells={anchor=west},  font=\footnotesize, rounded corners=2pt,}, legend pos=north west,xlabel=Number of Nodes, ylabel=$L_2$ Relative Error - $\sigma_{VM}$ - Indicator]
			\addplot+[Tblue,mark=triangle*,mark options={fill=Tblue}]   table [x=X, y=L2R-VMS-Est
					, col sep=comma] {ZZIndicator-GearKey.csv};
			\draw [gray,<-] (axis cs:46868,0.1617)+(0,-1mm) -- (axis cs:46868,0.03);
			\node[text=gray](d1) at (axis cs:46868,0.025) {Initial};
			\node[text=gray](d2) [below=-1mm of d1] {discretization};
			\draw [gray,<-] (axis cs:59544,0.1995)+(0,-1mm) -- (axis cs:59544,0.03);
			\node[text=gray](d1) at (axis cs:59544,0.025) {It. \#1};
			\draw [gray,<-] (axis cs:78004,0.0970)+(0,-1mm) -- (axis cs:78004,0.03);
			\node[text=gray](d1) at (axis cs:78004,0.025) {It. \#2};
			\draw [gray,<-] (axis cs:104527,0.0746)+(0,-1mm) -- (axis cs:104527,0.03);
			\node[text=gray](d1) at (axis cs:104527,0.025) {It. \#3};
			\draw [gray,<-] (axis cs:143763,0.0598)+(0,-1mm) -- (axis cs:143763,0.03);
			\node[text=gray](d1) at (axis cs:143763,0.025) {It. \#4};
		\end{axis}
	\end{tikzpicture}
	\caption{Evolution of the error indicator for the four iterations of adaptive refinement for the problem of the gear coupled to a shaft.}
	\label{GearKeyIndicator}
\end{figure}
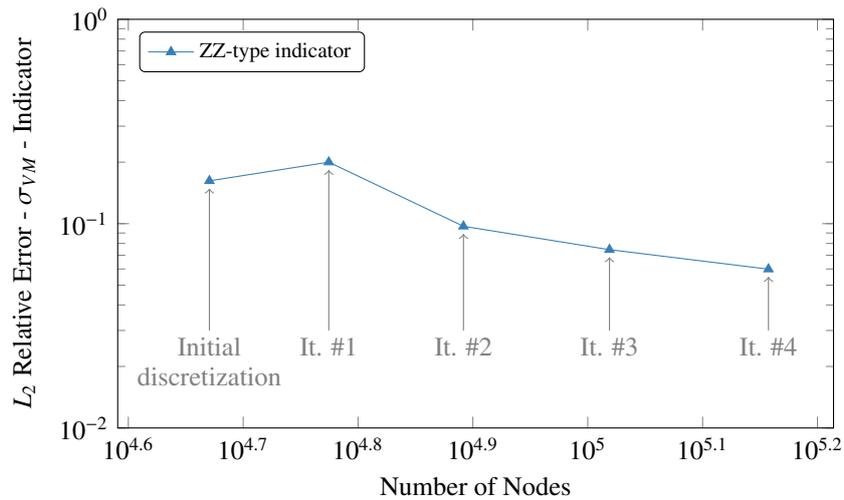

We present in Figure \ref{GearKeyIndicator} the computed error indicator for the results obtained from the initial discretization and the four refined discretizations. The error indicator is shown in terms of the $L_2$ relative error. We see that the value of the error increases at the first iteration and then reduces relatively steadily. The error reduction between the initial discretization and the final refinement iteration is approximately 63\%.

\subsubsection{Closed cylinder subject to pressure}

\begin{figure}[!h]
	\centering
	\begin{tikzpicture}
		\node at (0,0) {\includegraphics[width=12cm]{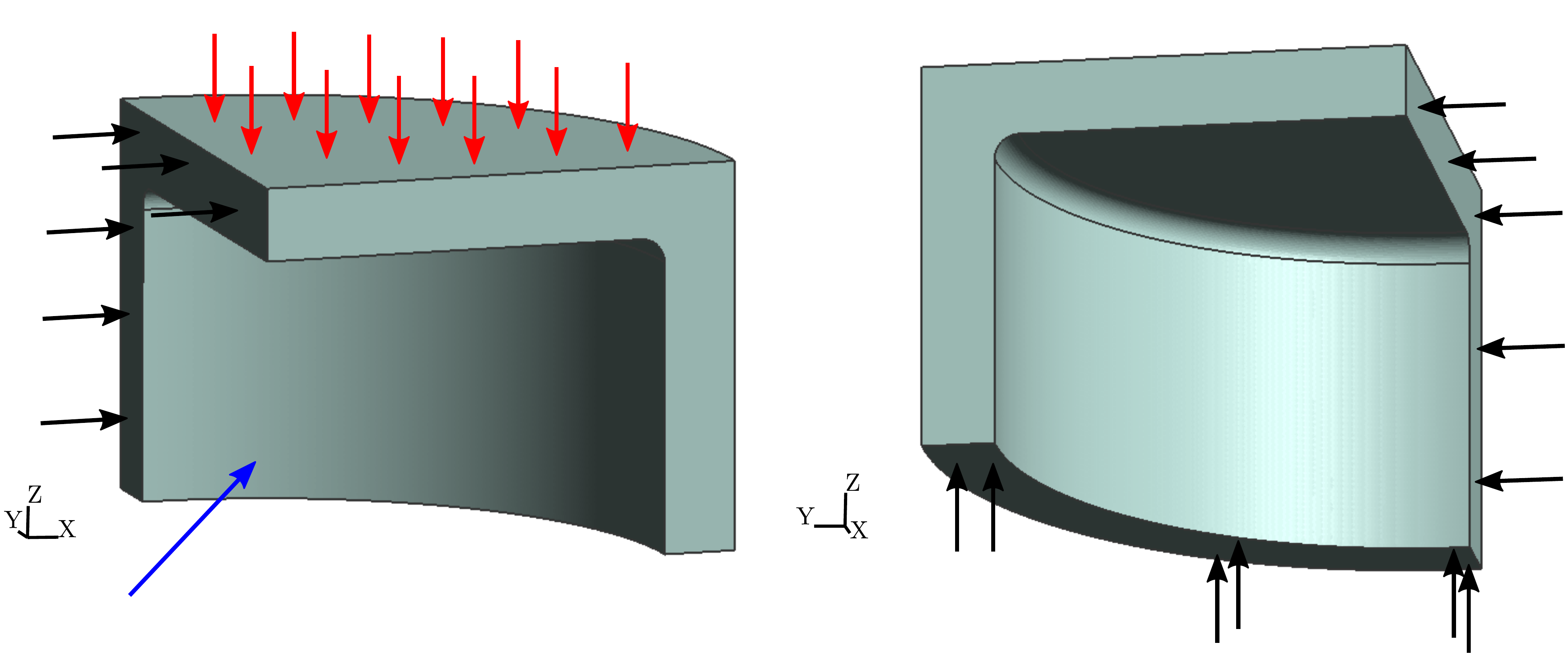}};
		\node[color=red] at (-4.7,2.5) [right] {Uniform pressure loading};
		\node[color=black] at (5.5,-2.0) [right] {Sliding boundary};
		\node[color=black] at (5.5,-2.4) [right] {conditions};
		\node[color=black] at (-8.3,0.6) [right] {Sliding boundary};
		\node[color=black] at (-8.3,0.2) [right] {conditions};
		\node[color=blue] at (-8.3,-2.4) [right] {Homogeneous Neumann};
		\node[color=blue] at (-8.3,-2.8) [right] {boundary conditions};
	\end{tikzpicture}
	\caption{Closed cylinder subject to pressure loading on the top. Uniform pressure loading is applied on the top surface. The displacement is limited in the directions normal to the surface on the symmetry planes XZ and YZ on the bottom surface. Stress-free boundary conditions are applied to the internal and external surface of the cylinder.}
	\label{ModelClosedCylinder}
\end{figure}
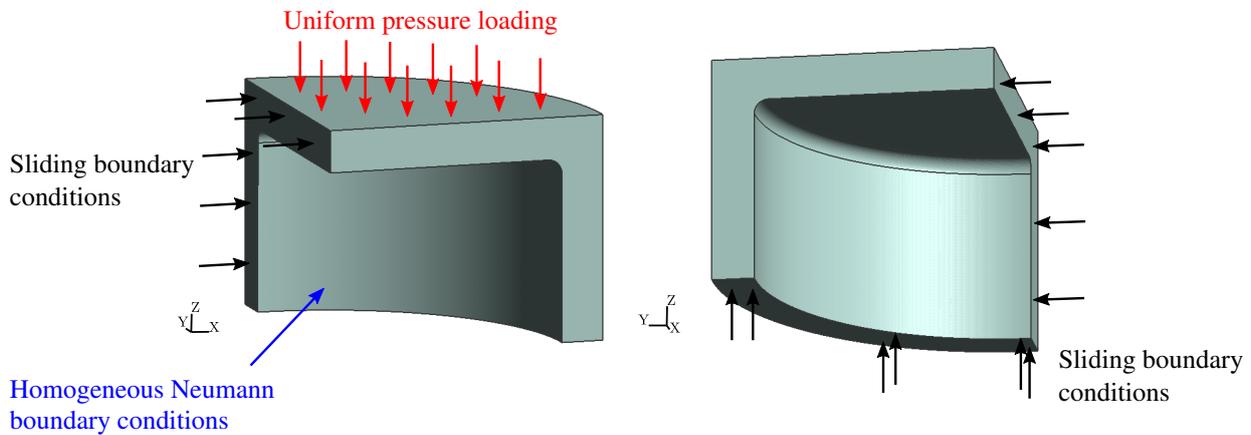

We present in Figure \ref{ModelClosedCylinder} the geometry of a closed cylinder subject to uniform pressure loading on the top of its closed end and the associated boundary condition. The bottom of the closed cylinder is fixed in the direction normal to the closed end of the cylinder. We considered only a quarter of the cylinder to reduce the computational cost. We applied Dirichlet boundary conditions in the directions normal to the surface on the symmetry planes XZ and YZ.

\begin{figure}[h!]
	\centering
	\begin{tabular}{|M{4cm}|c|c|}
		\hline
		                                             & \textbf{Bottom view}                                                                                                  & \textbf{Slice}                                                                                                       \\
		\hline
		\multirow{2}{*}{\begin{minipage}[t]{0.8\columnwidth} Initial discretization \\ 14,827 nodes \end{minipage}} &                                                                                                                       &                                                                                                                      \\[-3ex]
		                                             & \includegraphics[width=3.5cm]{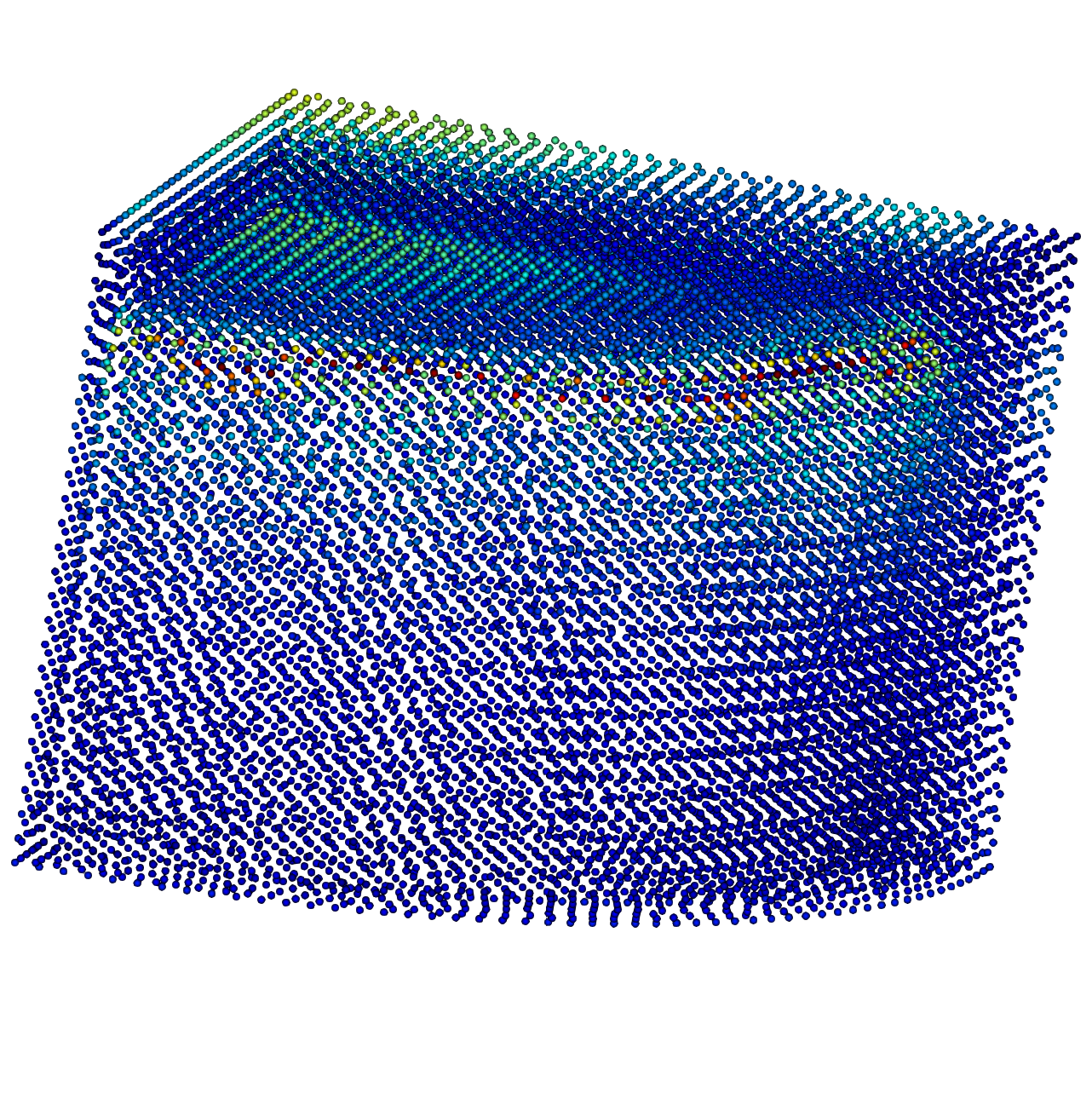} \includegraphics[height=3.5cm]{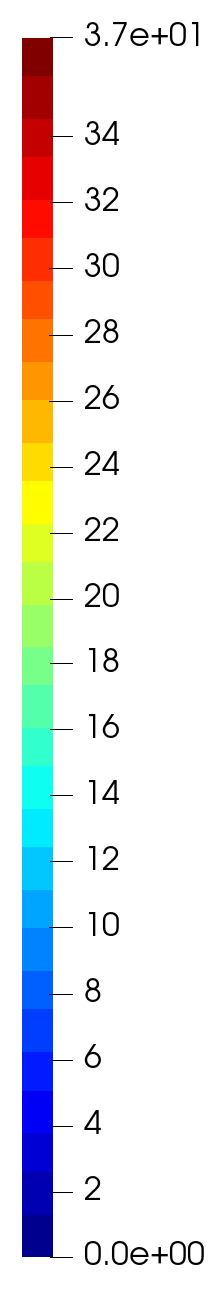} & \includegraphics[width=3.5cm]{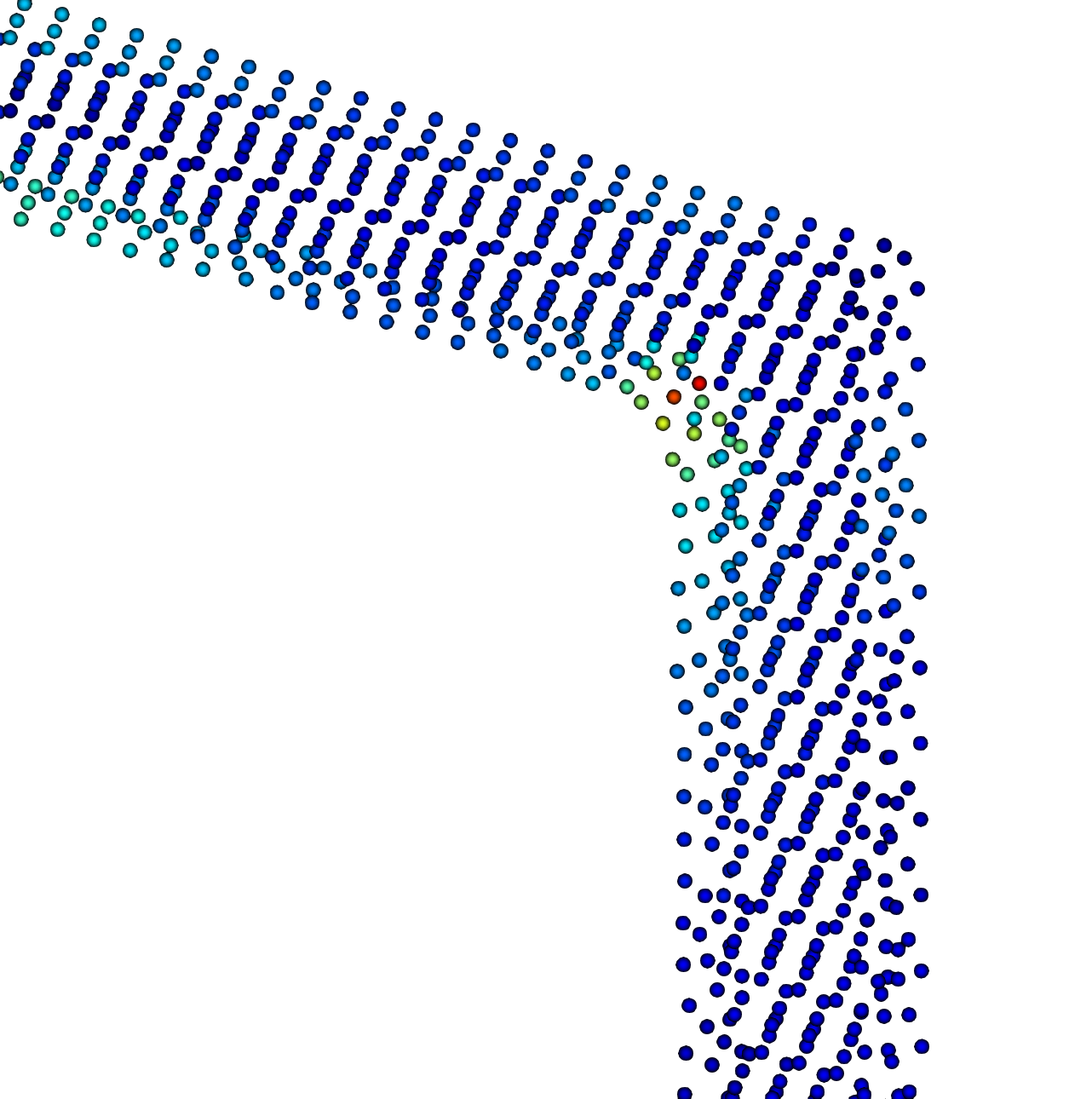} \includegraphics[height=3.5cm]{ClosedCylinder_Scale.png} \\
		\hline
		\multirow{2}{*}{\begin{minipage}[t]{0.8\columnwidth} Refinement iteration \#1 \\ 39,517 nodes \end{minipage}} &                                                                                                                       &                                                                                                                      \\[-3ex]
		                                             & \includegraphics[width=3.5cm]{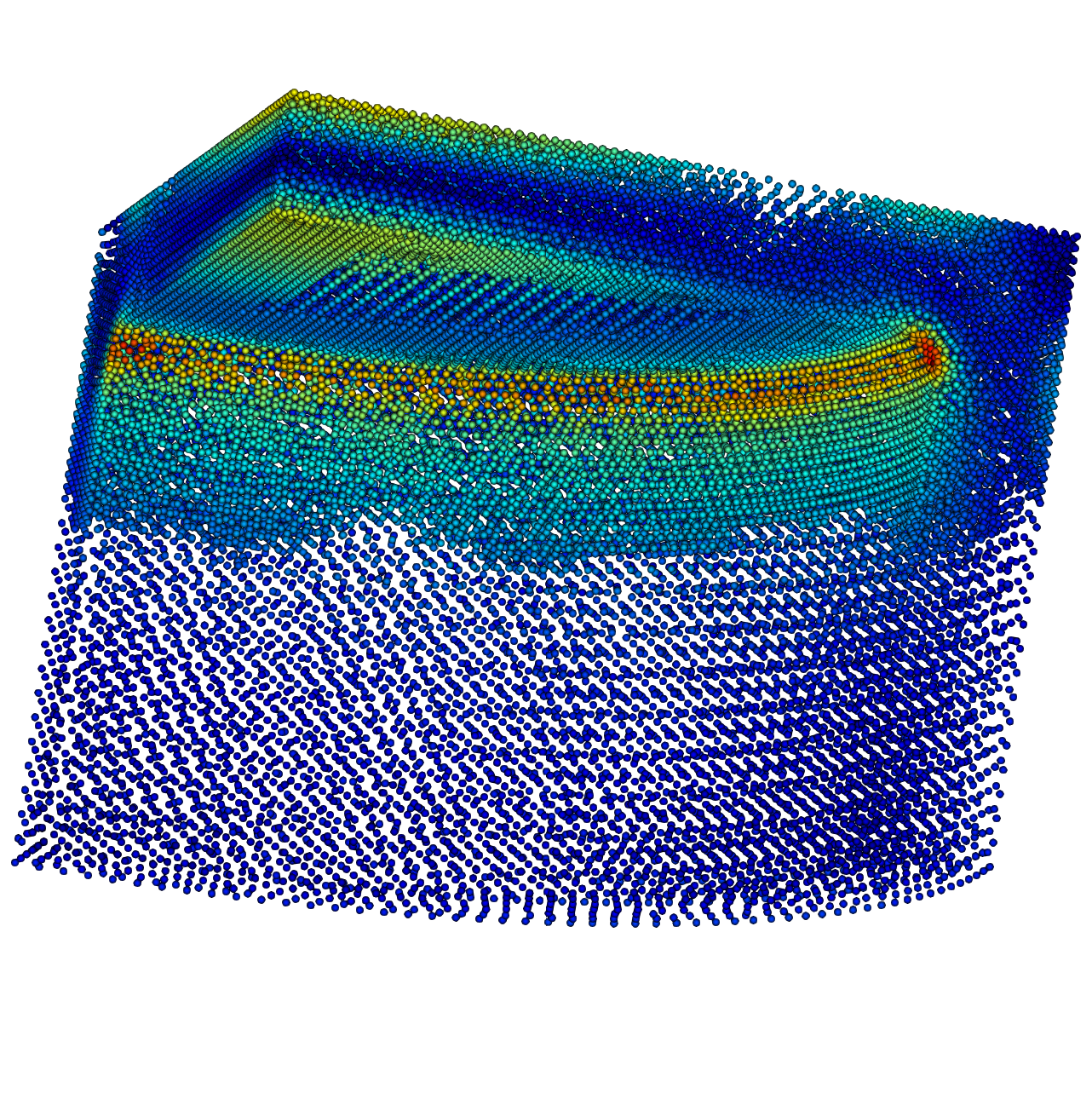} \includegraphics[height=3.5cm]{ClosedCylinder_Scale.png} & \includegraphics[width=3.5cm]{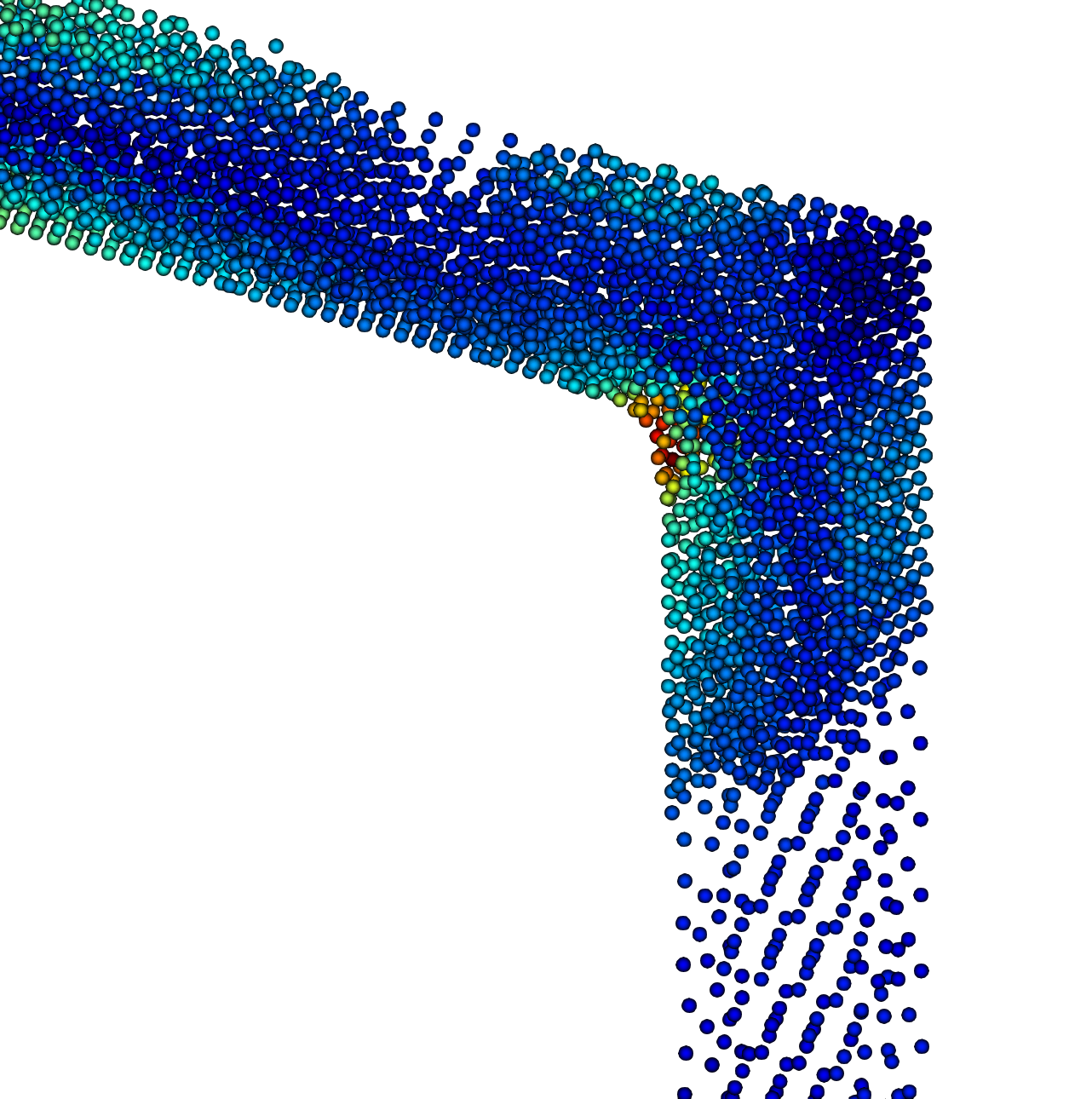} \includegraphics[height=3.5cm]{ClosedCylinder_Scale.png} \\
		\hline
		\multirow{2}{*}{\begin{minipage}[t]{0.8\columnwidth} Refinement iteration \#2 \\ 104,314 nodes \end{minipage}} &                                                                                                                       &                                                                                                                      \\[-3ex]
		                                             & \includegraphics[width=3.5cm]{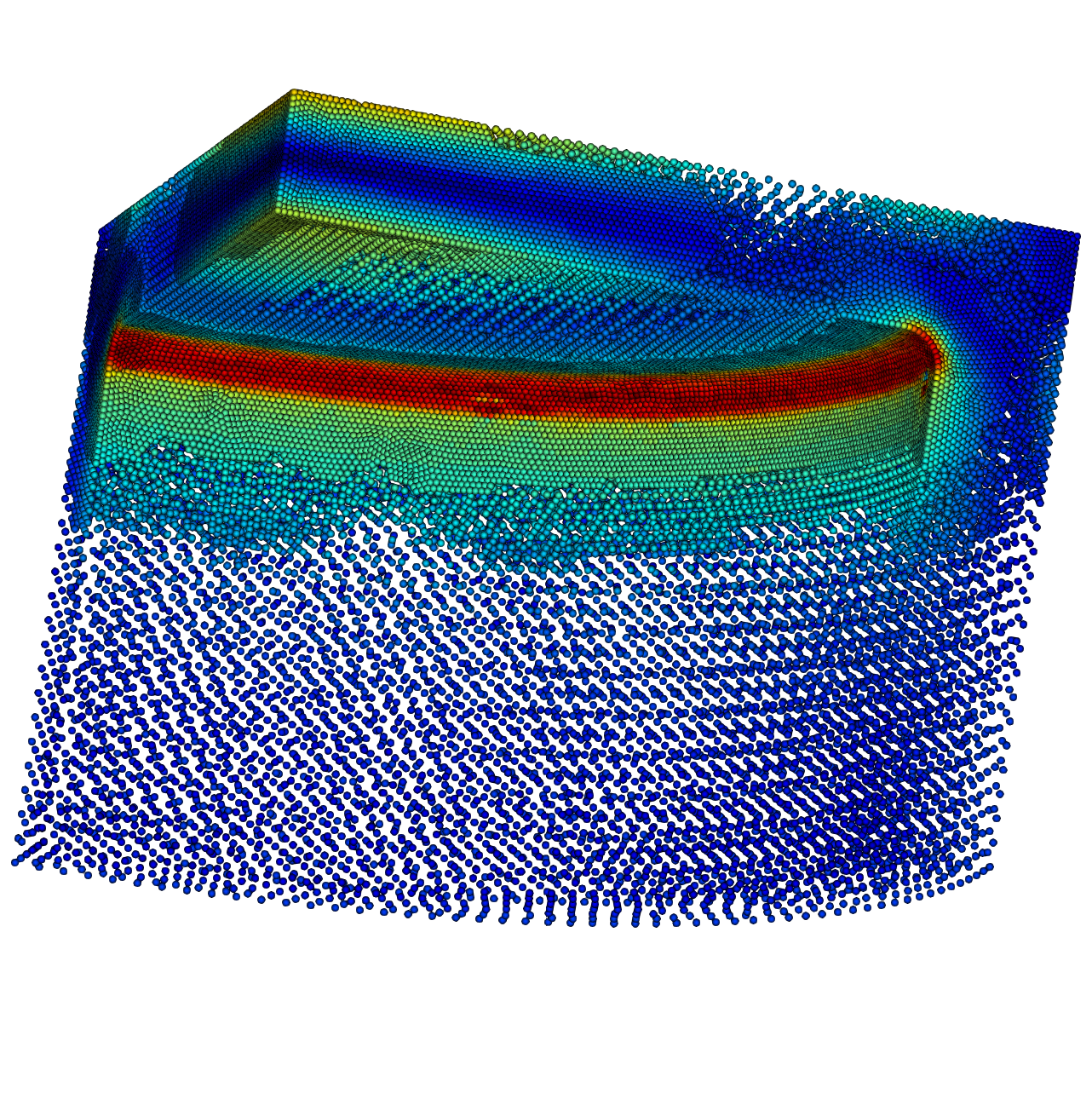} \includegraphics[height=3.5cm]{ClosedCylinder_Scale.png} & \includegraphics[width=3.5cm]{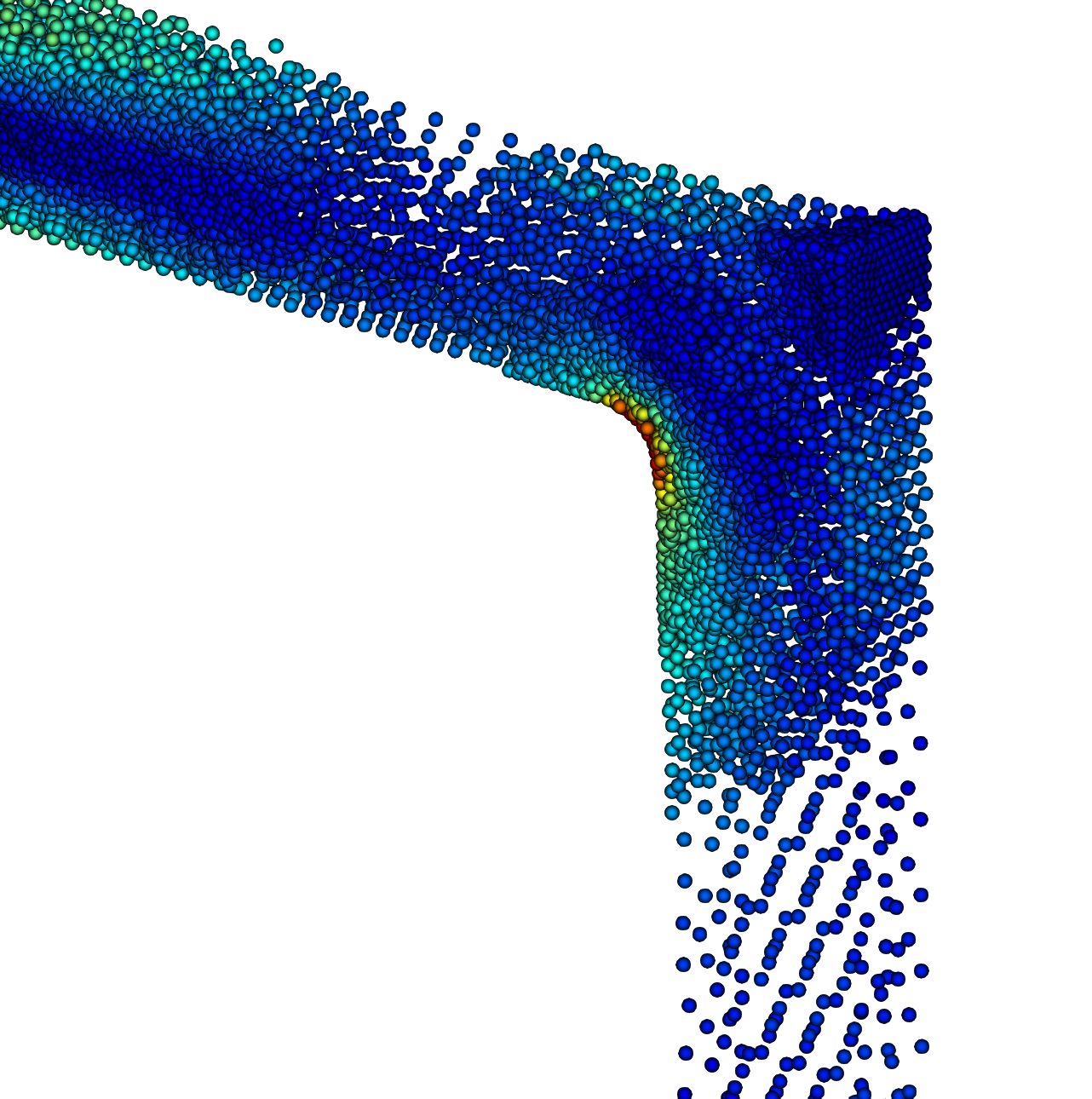} \includegraphics[height=3.5cm]{ClosedCylinder_Scale.png} \\
		\hline
	\end{tabular}

	\caption{Closed cylinder subject to pressure - Evolution of the discretization of the domain and of the solution in terms of von Mises stress through two iterations of adaptive refinement. The results are shown for the stress range 0-37 units for comparison purposes.}
	\label{RefinementStepsClosedCylinder}
\end{figure}

\begin{figure}[!h]
	\centering
	\includegraphics[width=4.5cm]{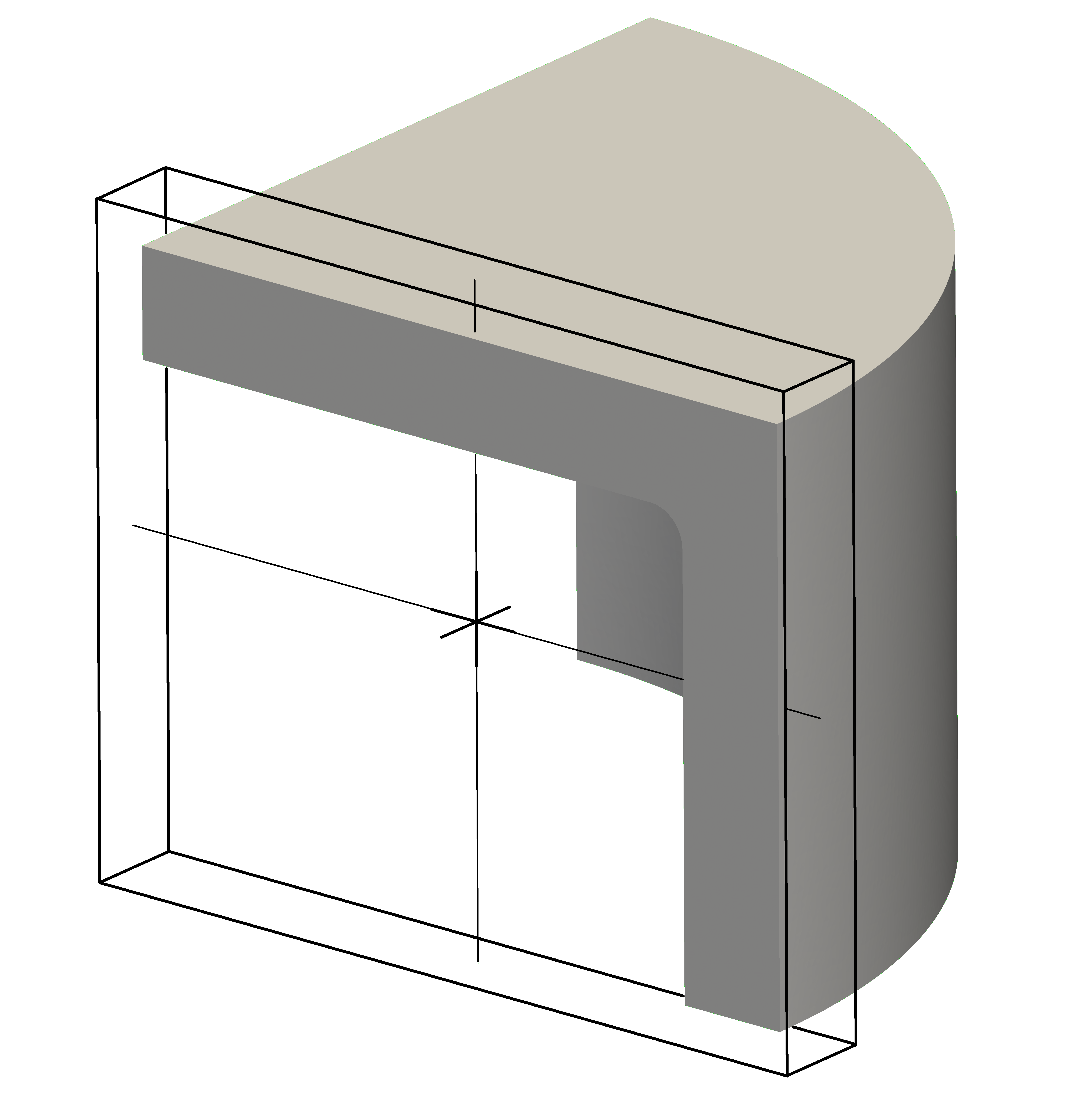}\\
	\caption{Selection of a thick ``slice'' of the domain to allow a closer analysis of the discretization and of the results.}
	\label{ClosedCylinderSlice}
\end{figure}

We show the results obtained for the initial discretization and for two refinement iterations in Figure \ref{RefinementStepsClosedCylinder}. To facilitate the analysis of the results, we also present in this figure the results for a thick ``slice'' of the domain. This allows a closer view of a portion of the domain. The ``slice'' is the intersection of a box and the geometry. The position of the box is shown in Figure \ref{ClosedCylinderSlice}. We focus on the fillet between the closed top and the cylinder since this is the zone where the von Mises stress is the largest. The results are shown for the stress range 0-37 for comparison purposes.

We observe from Figure \ref{RefinementStepsClosedCylinder} that new nodes are placed, at each iteration, in two areas of the domain. The first area is the fillet, where the stress concentration is the highest. The second area is the center of the top surface. It corresponds to the area where the displacement is the greatest. The refinement of the domain is relatively symmetric.

\begin{figure}[h!]
	\centering
	\begin{tabular}{|M{5cm}|c|c|}
		\hline
		                                             & \textbf{Bottom view}                                                                                                  & \textbf{Top view}                                                                                                  \\
		\hline
		\multirow{2}{*}{\begin{minipage}[t]{0.8\columnwidth}  \underline{ Smart cloud collocation solution} \\ $\blacktriangleright$ Refinement iteration \#2 \\ $\blacktriangleright$ 104,314 nodes  \\ $\blacktriangleright$  Approx. inter-node distance: \\ \hspace*{2ex} 0.0077 in dense regions \\ \hspace*{2ex} 0.03 in coarse regions \end{minipage}} &                                                                                                                       &                                                                                                                    \\[-3ex]
		                                             & \includegraphics[width=4.5cm]{ClosedCylinder_Bottom_It2.png} \includegraphics[height=4.5cm]{ClosedCylinder_Scale.png} & \includegraphics[width=4.5cm]{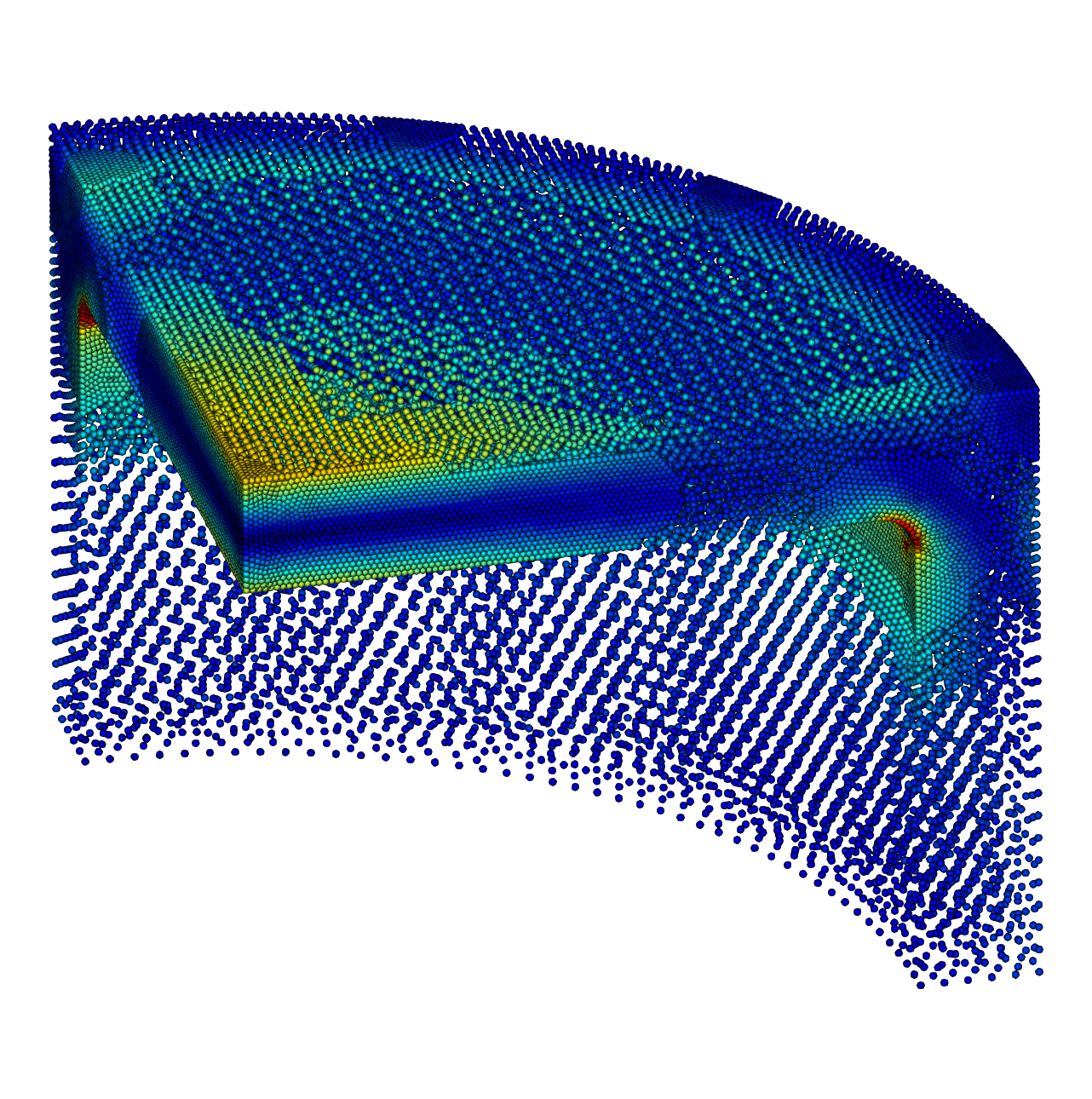} \includegraphics[height=4.5cm]{ClosedCylinder_Scale.png} \\
		\hline
		\multirow{2}{*}{\begin{minipage}[t]{0.8\columnwidth}  \underline{FEA solution} \\ $\blacktriangleright$ 166,193 nodes  \\ $\blacktriangleright$ Approx. inter-node distance:\\ \hspace*{2ex} 0.011 \end{minipage}} &                                                                                                                       &                                                                                                                    \\[-3ex]
		                                             & \includegraphics[width=4.5cm]{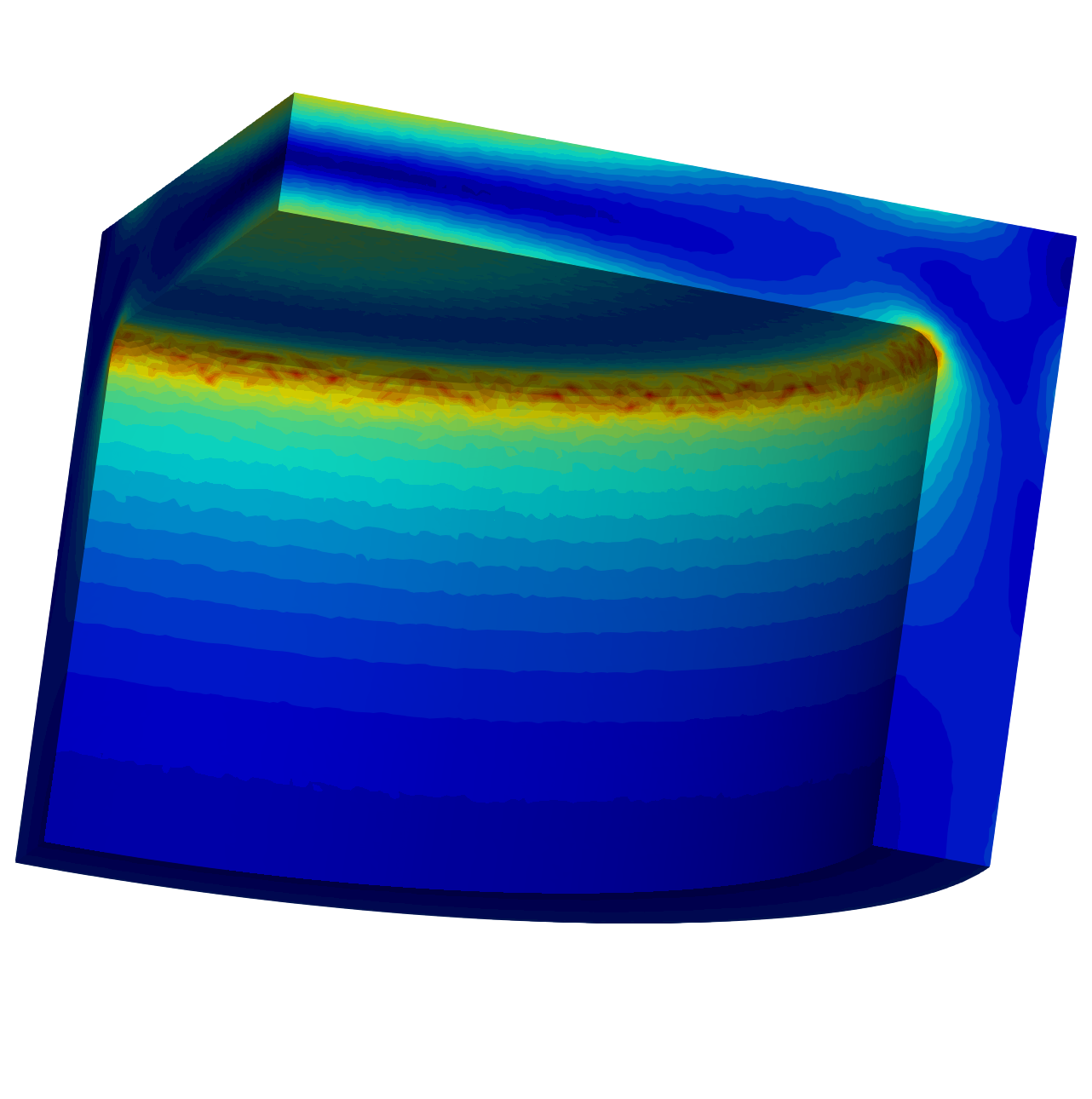} \includegraphics[height=4.5cm]{ClosedCylinder_Scale.png} & \includegraphics[width=4.5cm]{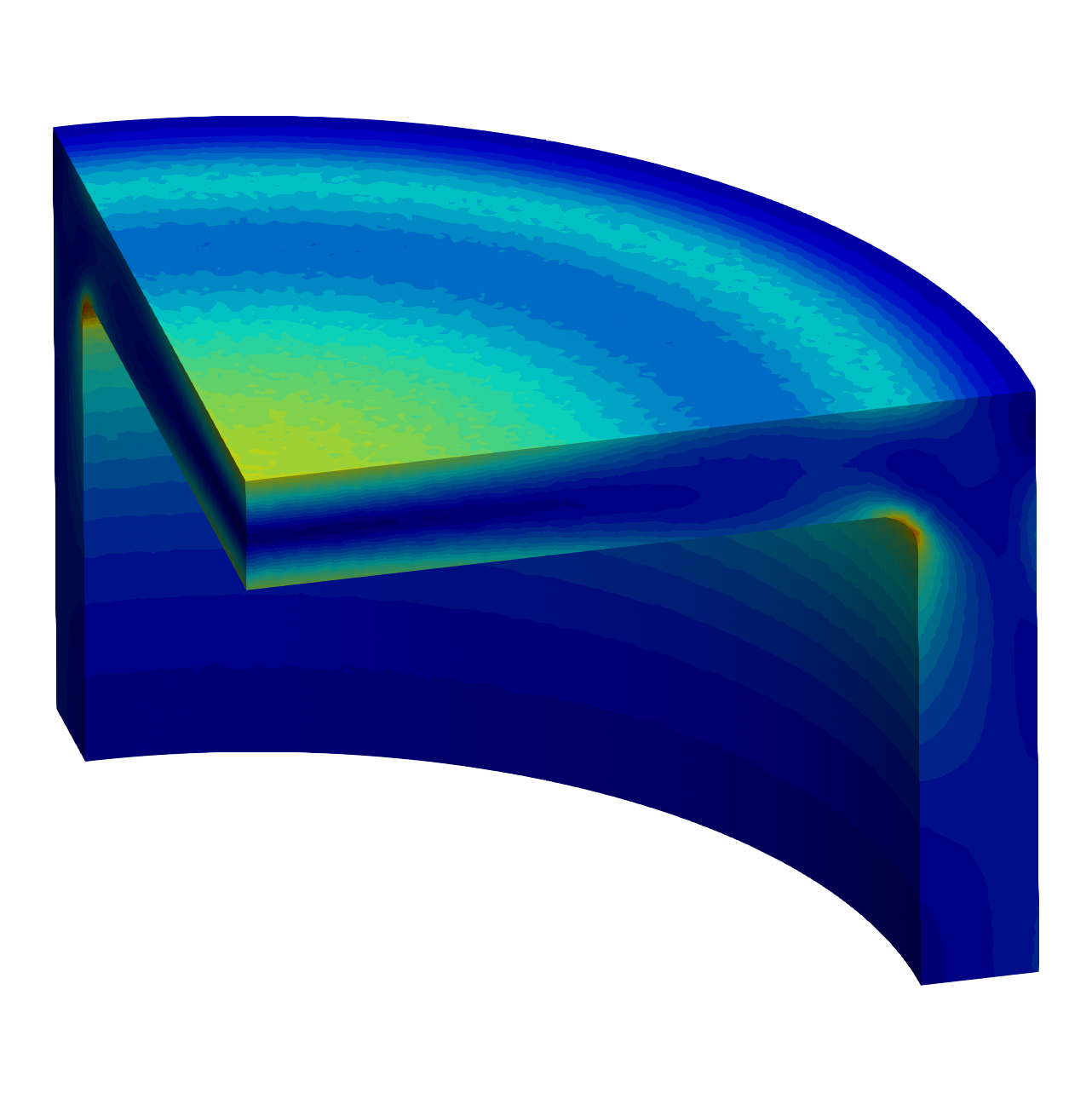} \includegraphics[height=4.5cm]{ClosedCylinder_Scale.png} \\
		\hline
	\end{tabular}

	\caption{Comparison, in terms of von Mises stress, of results obtained with the smart cloud method after two adaptive refinement iterations to results from the reference finite element solution.}
	\label{ClosedCylinderCollocVSFEA}
\end{figure}

We compare in Figure \ref{ClosedCylinderCollocVSFEA} the results obtained from the refined point cloud to a reference finite element solution obtained using code\_aster \cite{CodeASTER}. We used a finite element mesh composed of 166,1933 nodes and 900,039 linear tetrahedra for the purpose of this comparison. The discretization of the domain is uniform. The smart cloud collocation scheme leads to a larger stress concentration in the fillet than the finite element method solution. The von Mises stress in the fillet is approximately 34-37 units for the smart cloud collocation solution and approximately 22-28 units for the finite element solution. The node density in the fillet is larger for the adapted smart cloud than for the finite element discretization, because of the two adaptive refinement iterations. The inter-node  distance is approximately 0.0077 in the fillet for the  smart cloud discretization after two refinement iterations. It is approximately 0.011 for the finite element discretization. Also, for the smart cloud collocation scheme, the strong form of the partial differential equation is solved. These aspects could explain the different stress concentration obtained from both models. We also observe that the stress field in the fillet is smoother for the results obtained with the smart cloud collocation scheme than for the results obtained with the finite element method. Both methods lead to similar results in terms of von Mises stress on the top side of the closed cylinder.

\begin{figure}[!h]
	\centering
	\begin{tikzpicture}[scale=1]
		\begin{axis}[height=7cm,width=11cm, ymode=log, ymin=0.1, ymax=1, xmode=log, xmin=11000, legend entries={ZZ-type indicator}, legend style={ at={(0.5,-0.2)}, anchor=south west,legend columns=1, cells={anchor=west},  font=\footnotesize, rounded corners=2pt,}, legend pos=north west,xlabel=Number of Nodes, ylabel=$L_2$ Relative Error - $\sigma_{VM}$ - Indicator]
			\addplot+[Tblue,mark=triangle*,mark options={fill=Tblue}]   table [x=X, y=L2R-VMS-Est
					, col sep=comma] {ZZIndicator-ClosedCylinder.csv};
			\draw [gray,<-] (axis cs:14827,0.445)+(0,-1mm) -- (axis cs:14827,0.2);
			\node[text=gray](d1) at (axis cs:14827,0.18) {Initial};
			\node[text=gray](d2) [below=-1mm of d1] {discretization};
			\draw [gray,<-] (axis cs:39517,0.413)+(0,-1mm) -- (axis cs:39517,0.2);
			\node[text=gray](d1) at (axis cs:39517,0.18) {It. \#1};
			\draw [gray,<-] (axis cs:104314,0.415)+(0,-1mm) -- (axis cs:104314,0.2);
			\node[text=gray](d1) at (axis cs:104314,0.18) {It. \#2};
		\end{axis}
	\end{tikzpicture}
	\caption{Evolution of the error indicator for the two iterations of adaptive refinement for the problem of the closed cylinder subject to pressure.}
	\label{ClosedCylinderIndicator}
\end{figure}
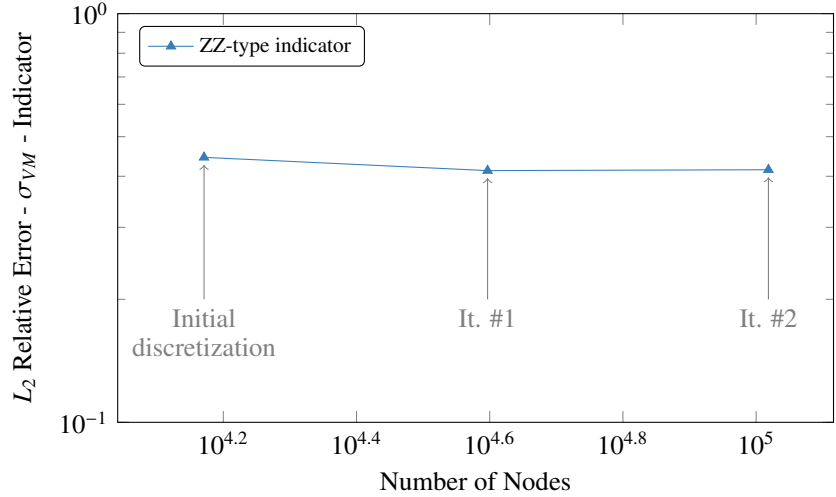

We present in Figure \ref{ClosedCylinderIndicator} the computed error indicator for the initial discretization and the two discretization obtained after two refinement iterations. The error indicator is shown in terms of the $L_2$ relative error norm. We see that the value of the error does not decrease much between the results obtained from the initial discretization and from the discretization after two refinement steps. The error reduction is approximately 6\%. A very slight increase of the error indicator is observed between the first and second refinement steps.

\clearpage

\section{Conclusion}

We presented in this article an integrated smart cloud collocation scheme to solve mechanics problems directly from CAD geometries. We used model adaptivity to speed up convergence of the results, while minimizing the input of the user.

The nodes of the smart cloud contain all the information required for optimum solution of the collocation problem (e.g. reference to the exact geometry, boundary conditions, connections to other boundary nodes). This information is key for adaptive refinement.

We used boundary-type elements to discretize the boundary of the domain. This approach is fast and robust since it relies on the discretization of the simple geometric features that compose the CAD domain. Boundary-type elements are also used to improve the accuracy of the solution using a generalization of the visibility criterion for point collocation methods. A larger number of interior nodes can be generated efficiently using regular node lattices such as triangular lattice for 2D domains or hexagonal close-packed lattice for 3D domains. The CAD file is used to select the nodes of the regular discretization which are in the domain. The proposed method only relies on the selection of a threshold parameter which is used to assess if interior nodes that are close to boundary nodes, are to be included in the collocation model or not. We showed that a threshold of 0.3 leads to good results.

Model adaptivity from CAD using the smart cloud collocation scheme is based on two key steps: error indication and cloud refinement. We presented two error indicators in this article: a ZZ-type error indicator and a residual-type error indicator. Their computation can be parallelized easily. Both allow a successful identification of the areas of the domain where the error is the greatest. We did not select the residual-type error indicator for our adaptive refinement method because of its complexity and computational cost. The adaptive refinement method presented in this article is based on the initial discretization of the CAD geometry. New nodes are added to the initial discretization to generate a refined discretization. Boundary elements are used to speed up the refinement process of the boundary. Once a new node is identified, it is projected onto the actual CAD geometry, thus ensuring that the point cloud discretization remains true to the original CAD geometry, and is not dependent on the element resolution. The required pieces of information about the CAD geometry and/or boundary conditions are added to the new nodes.

Two aspects of the proposed method should be further investigated: the quality of the error indicator and the convergence of the linear system for very large 3D problems. The error indicator considered here allows a proper identification of the refinement areas but cannot be considered as a good estimator. Even if an increase of the error indicator is a sign of an increase of the exact error, we do not consider it as a reliable estimator. The convergence of the linear systems obtained from adapted 3D clouds was difficult to attain. The condition of the matrices is an aspect that should be analyzed in greater detail to speed-up the convergence to the solution.

The use of CAD files at the heart of a collocation method is a great advantage to reduce the steps of the analysis to a minimum and ensure that the key features of the geometry are not lost during the refinement process, in the context of model adaptivity. These aspects are particularly true for domains with significant curvature, re-entrant corners or singularities.

To conclude, the work we have done to date on numerical methods, both within a strong-form framework \cite{Jacquemin2019, Jacquemin2021, Suchde2017, Suchde2018} and within isogeometric boundary element approaches \cite{Simpson2012,Lian2013,Simpson2014,Peng2017,Lian2017,Lian2016,Chen2019,Chen2020,Chen2020b} makes us hopeful that, with some significant work from the community, those approaches may have a strong impact on the way we simulate and optimize engineering systems.

With the advent of imaging approaches, in particular LiDAR and photogrammetry, solving directly from point clouds and adapting them to the simulation that is required becomes an urgent matter. Imaging approaches have been thought out to improve and optimize visual rendering of different scenes. Significant work remains to be done to make those point clouds suitable for physics-based simulations. Based on our experience, we suggest the following research directions in this nascent field:

\begin{itemize}
	\item $\left[ \text{Goal-oriented defeaturing} \right]$ \space Simplify smart point clouds through goal-oriented defeaturing such as done by Rahimi et al. \cite{Rahimi2018} and references therein \cite{Thakur2009,Danglade2013};
	\item $\left[ \text{Goal-oriented error estimation and adaptivity} \right]$ \space Minimize the error on a given quantity of engineering interest through goal-oriented error estimation for smart point clouds \cite{Prudhomme1999,Chamoin2021};
	\item $\left[ \text{Smart clouds to CAD} \right]$ \space Simplify smart point clouds into CAD primitives using classification machine learning approaches \cite{Aghighi2018,Aghighi2020,Li2019,Angles2019,Saporta2021};
	\item $\left[ \text{Properties identification from photogrammetry} \right]$ \space Infer parameter values based on texture and color from photogrammetric images \cite{Zhan2009,Awrangjeb2013,Alshawabkeh2020};
	\item $\left[ \text{Solution process} \right]$ \space Develop preconditioners and parallelization schemes for smart point clouds \cite{petsc-user-ref,Falgout2002,trilinos-website};
	\item $\left[ \text{Integrated simulation pipeline} \right]$ \space Develop open source pipelines based on libraries such as Open CASCADE Technology \cite{opencascade}, {VTK} \cite{Schroeder2006}, {CGAL} \cite{cgal:eb-21b}, Medusa \cite{Medusa}, {Voro++} \cite{Rycroft2009};
	\item $\left[ \text{Multi-scale model reduction approaches} \right]$ \space Defeature complex point clouds in order to accelerate simulations using multi-scale, domain decomposition and model order reduction approaches.
\end{itemize}

We will be pursuing the above directions, both in the context of computational engineering and computational archaeology.

\section*{Acknowledgments}

Thibault Jacquemin and St\'ephane P.A. Bordas would like to acknowledge funding from the Luxembourg National Research Fund (INTER/FWO/15/10318764), and from the European Union’s Horizon 2020 research and innovation program for the DRIVEN project (grant agreement No 811099). Pratik Suchde would like to acknowledge support from the European Union’s Horizon 2020 research and innovation program under the Marie Sk{\l}odowska-Curie Actions grant agreement No. 892761 ``SURFING".

\section*{Data Availability}

The data that support the findings of this study are openly available in \url{https://gitlab.com/tjacquemin/smart-cloud}.

\bibliographystyle{elsarticle-num}
\bibliography{Bibliography}

\end{document}